\documentclass{article}
\usepackage{amssymb}
\usepackage{amsfonts}
\usepackage{amsmath}
\usepackage{hyperref}
\usepackage{geometry}

\newtheorem{theorem}{Theorem}

\newtheorem{corollary}[theorem]{Corollary}

\newtheorem{definition}[theorem]{Definition}

\newtheorem{lemma}[theorem]{Lemma}

\newtheorem{proposition}[theorem]{Proposition}
\newtheorem{remark}[theorem]{Remark}

\numberwithin{theorem}{section}
\input{tcilatex}
\begin{document}

\title{Hecke Eigensystems for $(\func{mod}p)$ Modular Forms of PEL-type and
Algebraic Modular Forms}
\author{Davide A. Reduzzi \\
%EndAName
{\small University of California at Los Angeles, Department of Mathematics}\\
{\scriptsize Math Sciences Building 6363, 520 Portola Plaza, Los Angeles, CA
90095, U.S.A.}\\
{\scriptsize devredu83@math.ucla.edu}}
\maketitle
\tableofcontents

\section{Introduction}

In a letter to J. Tate written in 1987 (\cite{Se96}), J-P. Serre gives an
adelic interpretation of the systems of Hecke eigenvalues occurring in the
space of $(\func{mod}p)$ elliptic modular forms: Serre proves that these
eigensystems coincide with the Hecke eigensystems obtained from locally
constant functions $\mathfrak{B}^{\times }\backslash \mathfrak{\hat{B}}%
^{\times }\rightarrow \overline{\mathbb{F}}_{p}$, where $\mathfrak{B}$ is
the quaternion $%
%TCIMACRO{\U{211a} }%
%BeginExpansion
\mathbb{Q}
%EndExpansion
$-algebra ramified at $p$ and $\infty $, and $\mathfrak{\hat{B}}$ is its
adelization; the algebra $\mathfrak{B}$ naturally appears when evaluating
modular forms at supersingular elliptic curves over $\overline{\mathbb{F}}%
_{p}$.

\noindent This result fits into the more general settings of a $(\func{mod}%
p) $ Langlands philosophy, and it was originally closely related to the
study of Serre conjecture; it can also be regarded as the starting point for
defining the notion of algebraic modular forms, as it was done by B. Gross
in \cite{Gr}.

In \cite{Gh04}, A. Ghitza generalizes Serre's result to $(\func{mod}p)$
Siegel modular forms of genus $g>1$; in this context the algebraic group $%
\mathfrak{B}^{\times }$ is replaced by $GU_{g}(\mathfrak{B})$, an inner form
of $GSp_{2g}$ over $%
%TCIMACRO{\U{211a} }%
%BeginExpansion
\mathbb{Q}
%EndExpansion
$. It is interesting to notice that to obtain the Hecke correspondence for $%
g>1$, it is not enough to restrict modular forms to the supersingular locus
of the Siegel moduli variety over $\overline{\mathbb{F}}_{p}$, since this
locus has positive dimension: one needs to restrict modular forms to the
superspecial locus of the Siegel variety.

It was suggested to us by B. Gross that one could expect results similar to
the ones of Serre and Ghitza in more general settings, the basic idea being
that the systems of $(\func{mod}p)$ Hecke eigenvalues for modular forms
associated to some reductive group $G$ over $%
%TCIMACRO{\U{211a} }%
%BeginExpansion
\mathbb{Q}
%EndExpansion
$\ should coincide with the systems arising from $(\func{mod}p)$ algebraic
modular forms associated to an inner form $I$ of $G$ over $%
%TCIMACRO{\U{211a} }%
%BeginExpansion
\mathbb{Q}
%EndExpansion
$ satisfying some properties: for example $I$ should be compact modulo
center at infinity.

In this paper we investigate the existence of the above correspondence in
the context of modular forms arising from Shimura varieties of PEL-type;
while the procedure we follow generalizes to any PEL-situation of type A or
C, provided the supersingular locus of the corresponding Shimura variety is
non-empty, we present full details only for the case of modular forms
arising from Picard modular varieties (type A).

To state the main result we fix some notation: let $p$ be an odd prime and $%
k=%
%TCIMACRO{\U{211a} }%
%BeginExpansion
\mathbb{Q}
%EndExpansion
(\sqrt{\alpha })$ be a quadratic imaginary field in which $p$ is inert; let $%
g$ be a positive even integer and let $r,s$ be non-negative integers such
that $r+s=g.$ These data define the connected reductive $%
%TCIMACRO{\U{211a} }%
%BeginExpansion
\mathbb{Q}
%EndExpansion
$-group $G=GU_{g}(k;r,s)$ of type A,\ which is the unitary group associated
to the extension $k/%
%TCIMACRO{\U{211a} }%
%BeginExpansion
\mathbb{Q}
%EndExpansion
$ and of signature $(r,s)$; fixing an embedding of $k$ inside the quaternion 
$%
%TCIMACRO{\U{211a} }%
%BeginExpansion
\mathbb{Q}
%EndExpansion
$-algebra $\mathfrak{B}$ of endomorphisms of a fixed superspecial elliptic
curve over $\overline{\mathbb{F}}_{p}$, we can define an inner form $I$ of $%
G $ over $%
%TCIMACRO{\U{211a} }%
%BeginExpansion
\mathbb{Q}
%EndExpansion
$ so that:%
\begin{equation*}
I(%
%TCIMACRO{\U{211a} }%
%BeginExpansion
\mathbb{Q}
%EndExpansion
)=\{X\in GU_{g}(\mathfrak{B}):X\cdot \Phi =\Phi \cdot X\},
\end{equation*}

\noindent where $\Phi $ is the matrix $\left( \QDATOP{-\sqrt{\alpha }\cdot
I_{r}}{0_{s,r}}\QDATOP{0_{r,s}}{\sqrt{\alpha }\cdot I_{s}}\right) $. Let $%
N\geq 3$ be an integer not divisible by $p$ and set:%
\begin{eqnarray*}
U(N) &:&=\ker \left( GU_{g}\left( \mathcal{O}_{k}\otimes _{%
%TCIMACRO{\U{2124} }%
%BeginExpansion
\mathbb{Z}
%EndExpansion
}\widehat{%
%TCIMACRO{\U{2124} }%
%BeginExpansion
\mathbb{Z}
%EndExpansion
}^{p};r,s\right) \rightarrow GU_{g}\left( \mathcal{O}_{k}\otimes _{%
%TCIMACRO{\U{2124} }%
%BeginExpansion
\mathbb{Z}
%EndExpansion
}\frac{\widehat{%
%TCIMACRO{\U{2124} }%
%BeginExpansion
\mathbb{Z}
%EndExpansion
}^{p}}{N\widehat{%
%TCIMACRO{\U{2124} }%
%BeginExpansion
\mathbb{Z}
%EndExpansion
}^{p}};r,s\right) \right) , \\
U_{p} &:&=\ker \left( I(%
%TCIMACRO{\U{2124} }%
%BeginExpansion
\mathbb{Z}
%EndExpansion
_{p})\overset{\pi }{\rightarrow }G(U_{r}\times U_{s})(\mathbb{F}%
_{p^{2}})\right) ,
\end{eqnarray*}

\noindent where $\pi $ denotes the reduction map modulo a uniformizer of $%
\mathfrak{B}$. The main result of the paper, contained in Theorem \ref{main
thm}, is:

\begin{theorem}
\textit{The systems of Hecke eigenvalues arising from }$(r,s)$\textit{%
-unitary PEL-modular forms }$(\func{mod}p)$\textit{\ for the quadratic
imaginary field\ }$k$\textit{, having genus }$g$\textit{, fixed level }$N$%
\textit{\ and any possible weight }$\rho :GL_{g}\rightarrow GL_{m(\rho )}$%
\textit{, are the same as the systems of Hecke eigenvalues arising from }$(%
\func{mod}p)$\textit{\ algebraic modular forms for the group }$I$\textit{\
having level }$U_{p}\times U(N)\subset I(\mathbb{A}_{f})$\textit{\ and any
possible weight }$\rho ^{\prime }:G(U_{r}\times U_{s})\rightarrow
GL_{m^{\prime }(\rho ^{\prime })}$\textit{.}
\end{theorem}

An important tool to prove the above result will be for us the use of a
uniformization result for some isogeny classes on PEL-moduli varieties, due
to M. Rapoport and Th. Zink (Theorem \ref{RZ}). We fix a $p$-divisible group
with additional structure $(\mathbf{X},\overline{\lambda }_{\mathbf{X}},i_{%
\mathbf{X}})$, where $\mathbf{X}=G_{1/2}^{g}$ is superspecial, $\overline{%
\lambda }_{\mathbf{X}}$ is the $%
%TCIMACRO{\U{2124} }%
%BeginExpansion
\mathbb{Z}
%EndExpansion
_{p}$-class of a fixed principal polarization of $\mathbf{X}$, and $i_{%
\mathbf{X}}$ is an action of the order $\mathcal{O}_{k}$ on $\mathbf{X}$. We
can adapt the uniformization result of Rapoport and Zink in order to
uniformize only the finitely many points in the supersingular locus of the
Picard moduli variety having $p$-divisible group \textit{isomorphic} to $(%
\mathbf{X},\overline{\lambda }_{\mathbf{X}},i_{\mathbf{X}})$ (Corollary \ref%
{unif-cor}). This finite set will be our superspecial locus.

\noindent Since the cotangent space of the fixed $p$-divisible group $(%
\mathbf{X},\overline{\lambda }_{\mathbf{X}},i_{\mathbf{X}})$\ has
automorphism group naturally isomorphic to $G(U_{r}\times U_{s})(\mathbb{F}%
_{p^{2}})$, we are able to put differentials into the picture, and realize
all the $(\func{mod}p)$ Hecke eigenvalues of algebraic modular forms for the
group $I$ in the space of superspecial modular forms, and viceversa
(Proposition \ref{alg-superspecial}). Finally, we use an argument that is a
modification of the one appearing in \cite{Gh04} Th. 28,\ to conclude that
the set of superspecial Hecke eigenvalues coincide with the set of all the $(%
\func{mod}p)$ Hecke eigenvalues.

In the more general settings of PEL-varieties over $\overline{\mathbb{F}}%
_{p} $ of type A and C, whose superspecial locus is non-empty, the group $I$
should be taken to be the $%
%TCIMACRO{\U{211a} }%
%BeginExpansion
\mathbb{Q}
%EndExpansion
$-group of automorphisms of a fixed triple $(A_{0},i_{0},\overline{\lambda }%
_{0})$ defining a point in the superspecial locus of the PEL-variety; the
uniformization result of Rapoport and\ Zink can be applied in this general
context, since the $p$-divisible group of a superspecial abelian variety is
basic (cf. \ref{basicity}).

As an application of the main theorem of the paper, we obtain (cf. Theorem %
\ref{counting theorem}) that the number $\mathcal{N}$ of $(\func{mod}p)$
Hecke eigensystems arising from unitary modular forms of signature $(r,s)$,
fixed level $N$ as above and any weight is \textit{finite} and satisfies the
following asymptotic with respect to $p$:%
\begin{equation*}
\mathcal{N}\in O(p^{g^{2}+g+1-rs}).
\end{equation*}

\bigskip

The paper is organized as follows: Section 2 recalls the basic definitions
and results that one needs in order to state a representability theorem for
some moduli functor of $p$-divisible groups, as considered by Rapoport and
Zink in \cite{RZ}; in Section 3 we recall the basic facts on PEL-moduli
varieties and PEL-modular forms; in Section 4 we present the construction of
the uniformization morphism of Rapoport and Zink for basic isogeny classes
in PEL-Shimura varieties, and the modification we need to parametrize our
superspecial locus. The final Section 5 presents the settings and the
computations that allow us to compare $(\func{mod}p)$ Hecke eigensystems for
PEL-modular forms of unitary type, and $(\func{mod}p)$ Hecke eigensystems
for algebraic forms for the group $I$: the main theorem is here stated and
proven, and the number of $(\func{mod}p)$ Hecke eigenvalues is estimated.

The reader might prefer to start reading the paper from section \ref{exp},
and might refer to paragraph \ref{exp1} for the definition of the moduli
problem of $p$-divisible groups that appears in the formulation of Theorem %
\ref{RZ}, and to paragraph \ref{p-adic unif of isogeny classes} for the
description of the Rapoport-Zink uniformization morphism.

\bigskip

\bigskip

\bigskip

\noindent\ I would like to address my gratitude to C. Khare for introducing
me to the topics studied in this paper, and to B.H. Gross for suggesting to
me the problem treated here. I\ would also like to thank Khare for his
valuable supervision, and for sharing with me a note of N. Fakhruddin
solving a problem with the Hecke equivariance of a map considered in \cite%
{Gh04}.

I am grateful to N. Fakhruddin for signaling a problem contained in a
previous version of the paper relatively to the Kodaira-Spencer isomorphism,
and for other valuable remarks. I thank H. Hida and M-H. Nicole for
explaining to me that, in suitable circumstances, the natural map from the
Picard moduli scheme to the Siegel moduli scheme is a closed immersion. I
would like to further thank J. Tilouine for an interesting conversation we
had on the subject of this paper.

\newpage

\section{Moduli of $p$-divisible groups of PEL-type}

\subsection{$p$-divisible groups}

We recall basic facts about $p$-divisible group; the main references are 
\cite{DM}, \cite{F}, \cite{Me} and \cite{RZ}, 2.1-2.8.

\subsubsection{Basic definitions}

Fix a scheme $S$. By an $S$-group we mean a f.p.p.f. sheaf of commutative
groups on the site $(Sch/S)_{f.p.p.f.}$, whose underlying category is the
category of $S$-schemes, endowed with the Grothendieck f.p.p.f. topology. An 
$S$-group that is representable is called an $S$\textit{-}group-scheme. For
any positive integer $n$ and any $S$-group $G$ we denote by $G[n]$ the
kernel of the $S$-morphism $n\cdot 1_{G}:G\rightarrow G$ given my
multiplication by $n$; $G[n]$ is an $S$-group.

Fix a positive prime integer $p.$An $S$-group $G$ is a \textit{%
Barstotti-Tate group} (or a $p$\textit{-divisible group}) if the following
three conditions are satisfied: (1) $G=\lim\limits_{\longrightarrow
}G[p^{n}] $; (2) the morphism $[p]:G\rightarrow G$ is an epimorphism; (3) $%
G[p]$ is a finite locally free group-scheme over $S$.

By the theory of finite group-schemes over a field, the rank of the fiber of 
$G[p]$ at a point $s\in S$ is of the form $p^{h(s)}$, where $h:S\rightarrow 
%TCIMACRO{\U{2124} }%
%BeginExpansion
\mathbb{Z}
%EndExpansion
$ is a locally constant function on $S$; the rank of the fiber of $G[p^{n}]$
at $s$ is $p^{nh(s)}$ for any $n\geq 1$. If $h$ is a constant function
(e.g., when $S=Spec(k)$ for a field $k$), its only value is called the 
\textit{height} $ht(G)$ of $G$.

\bigskip

A morphism $f:G\rightarrow H$ of $p$-divisible groups over $S$ is said to be
an \textit{isogeny} if it is an epimorphism of f.p.p.f. sheaves whose kernel
is representable by a finite locally free $S$-group scheme. \noindent If $S$
is a scheme over $Spec(%
%TCIMACRO{\U{2124} }%
%BeginExpansion
\mathbb{Z}
%EndExpansion
_{p})$ in which $p$ is locally nilpotent, then the kernel of an isogeny $%
f:G\rightarrow H$ is finite of rank $p^{h^{\prime }}$ where $h^{\prime
}:S\rightarrow 
%TCIMACRO{\U{2124} }%
%BeginExpansion
\mathbb{Z}
%EndExpansion
$ is locally constant; if $h^{\prime }$ is constant, its only value is
called the \textit{height }of $f$.

\noindent The $%
%TCIMACRO{\U{2124} }%
%BeginExpansion
\mathbb{Z}
%EndExpansion
$-module $\limfunc{Hom}_{S}(G,H)$ of homomorphisms from $G$ to $H$ is a
torsion-free $%
%TCIMACRO{\U{2124} }%
%BeginExpansion
\mathbb{Z}
%EndExpansion
_{p}$-module. A \textit{quasi-isogeny} $f$ from $G$ to $H$ is global section
of the sheaf $\mathfrak{Hom}_{S}(G,H)\otimes _{%
%TCIMACRO{\U{2124} }%
%BeginExpansion
\mathbb{Z}
%EndExpansion
}%
%TCIMACRO{\U{211a} }%
%BeginExpansion
\mathbb{Q}
%EndExpansion
$ such that any point $s$ of $S$ has a Zariski open neighborhood on which $%
p^{n}f:X\rightarrow Y$ is an isogeny for some positive integer $n=n(s)$. We
denote by $Qisg_{S}(G,H)$ the group of quasi-isogenies from $G$ to $H$.

\noindent We have the following rigidity property (cf. \cite{RZ}, 2.8):

\begin{proposition}
\label{Rigidity of qisogenies}Let $G$ and $H$ be $p$-divisible groups over a
scheme $S$ in which $p$ is locally nilpotent; let $S^{\prime }\subset S$ be
a closed subscheme whose defining sheaf of ideals is locally nilpotent. Then
the canonical homomorphism $Qisg_{S}(G,H)\rightarrow Qisg_{S^{\prime
}}(G_{S^{\prime }},H_{S^{\prime }})$ is bijective.
\end{proposition}

\bigskip

Recall that for any finite flat group scheme $X$ over $S$, the Cartier dual $%
D(X)$ (or $X^{D}$) of $X$ is the finite locally free $S$-group-scheme
defined by $D(X)(T):=\limfunc{Hom}_{T}(G_{T},\mathbb{G}_{m}\times _{S}T)$ ($%
T $ any $S$-scheme). The assignment $D$ induces an additive anti-duality on
the category of finite and locally free $S$-group schemes. Let $G$ be a $p$%
-divisible group over $S$. The \textit{Serre dual} of $G$, denoted by $%
\widehat{G}$ (or $G^{D}$, or $D(G)$), is the $p$-divisible group defined as $%
\widehat{G}:=\lim\limits_{\longrightarrow }D(G\left[ p^{n}\right] )$,
\noindent where the map $D(G[p^{n}])\rightarrow D(G[p^{n+1}])$ in the direct
system is given by $D([p])$ for each $n$.

If $G$ is a $p$-divisible group over $S$, then $\widehat{G}%
[p^{n}]=D(G[p^{n}])$ for any $n$. The assignment $G\rightarrow \widehat{G}$
extends to morphisms in an obvious way and gives rise to an anti-duality in
the category of $p$-divisible groups over $S$ (notice this category is not
abelian) which is compatible with base changes. There is a canonical
isomorphism of $p$-divisible groups $G\overset{\sim }{\rightarrow }\widehat{%
\widehat{G}}.$

\bigskip

An $S$\textit{-polarization} of a $p$-divisible group $G$ over $S$ is an
anti-symmetric $S$-quasi-isogeny $\lambda :G\rightarrow \widehat{G}$. A $%
%TCIMACRO{\U{211a} }%
%BeginExpansion
\mathbb{Q}
%EndExpansion
_{p}$\textit{-homogeneous }$S$\textit{-polarization} of $G$ is the set $%
\overline{\lambda }=%
%TCIMACRO{\U{211a} }%
%BeginExpansion
\mathbb{Q}
%EndExpansion
_{p}^{\times }\lambda $ of $p$-adic non-zero multiples of an $S$%
-polarization of $G$. A \textit{principal }$S$\textit{-polarization} is an $%
S $-polarization that is also an isomorphism of $p$-divisible groups. (A
polarization $\lambda $ of $G$ is anti-symmetric in the sense that $\widehat{%
\lambda }=-\lambda $, where $\widehat{\lambda }$ denotes the Serre dual of $%
\lambda $, viewed in the canonical way as a map from $G$ to $\widehat{G}$;
the reason one requires the polarization to be anti-symmetric can be found
in \cite{O}, Prop. 1.12: a polarization of an abelian variety $A$ over $%
\overline{\mathbb{F}}_{p}$ induces a polarization of the associated $p$%
-divisible group in the sense of the above definition).

Let $O$ be a $\mathbb{Z}_{p}$-algebra with involution $^{\ast }$; an action
of $(O;^{\ast })$\ on $G$ is a homomorphism of $\mathbb{Z}_{p}$-algebras $%
i:O\rightarrow \limfunc{End}_{S}(G)$. If $G$ is endowed with such an action $%
i$, then $\widehat{G}$ is endowed with the dual action $\widehat{i}$ given
by setting $\widehat{i}(a):=i(a^{\ast })^{\widehat{}}$ for any $a\in O$, and
we say that an $S$-polarization $\lambda :G\rightarrow \widehat{G}$ respects
the $O$-action if $\lambda \circ i=\widehat{i}\circ \lambda $.

\subsubsection{Dieudonn\'{e} modules}

Fix a perfect field $k$ of characteristic $p>0$ and denote by $\sigma $ its
absolute Frobenius morphism. Denote by $\mathbb{W}$ the ring scheme over $k$
of Witt vectors, and let $V:\mathbb{W\rightarrow W}$ be the Verschiebung
morphism. The absolute Frobenius morphism $F$ on $\mathbb{W}$ is a
ring-scheme homomorphism and one has $FV=VF=p\cdot id_{\mathbb{W}}$. If $R$
is a $k$-algebra, we denote by $W_{R}$ (or by $W(R)$, or simply by $W$ if no
confusion arises)\ the ring $\mathbb{W(}R\mathbb{)}$.

For the proof of most of the results that follow, cf. \cite{DM} and \cite{F}%
; cf. also \cite{LO}, 1-5 and \cite{CB}.

\bigskip

Let $W$ be the ring $W_{k}$ and denote by $K_{0}$ its quotient field. Also,
denote by $\sigma $ the absolute Frobenius morphism induced by $\sigma
:k\rightarrow k$ on $W$ and on $K_{0}$.

An \textit{isocrystal} over $K_{0}$ is a finite dimensional $K_{0}$-vector
space $D_{0}$ endowed with a Frobenius $K_{0}$-semilinear automorphism $%
F:D_{0}\rightarrow D_{0}$; we denote isocrystals as pairs $(D_{0},F)$ and we
will call the $K_{0}$-automorphism $pF^{-1}$ the Verschiebung of the
isocrystal, and denote it by $V$. The dimension of the vector space $D_{0}$
over $K_{0}$ is called the \textit{height} of the isocrystal; the $p$-adic
valuation of $\det F$ is called the \textit{dimension} of the isocrystal.

Let $K$ be a finite field extension of $K_{0}$. A \textit{filtered isocrystal%
} over $K$ is an isocrystal $(D_{0},F)$ over $K_{0}$ endowed with a
decreasing filtration $Fil(D_{0})^{\bullet }:=\mathcal{F}^{\bullet }$ of the 
$K$-vector space $D_{0}\otimes _{K_{0}}K$ such that $\mathcal{F}^{r}=0$ and $%
\mathcal{F}^{s}=D_{0}\otimes _{K_{0}}K$ for some integers $r,s$ with $r>s$.
Isocrystals over $K_{0}$ form an abelian $%
%TCIMACRO{\U{211a} }%
%BeginExpansion
\mathbb{Q}
%EndExpansion
_{p}$-category, while filtered isocrystal over $K$ do not form an abelian
category.

For $n\in 
%TCIMACRO{\U{2124} }%
%BeginExpansion
\mathbb{Z}
%EndExpansion
$, define the Tate object in the context of isocrystals over $K_{0}$ as $%
\mathbf{1(}n):=(K_{0},p^{n}\sigma )$. As a filtered isocrystal over some $K$%
, we set $Fil(\mathbf{1(}n))^{i}=K$ for $i\leq n$ and $Fil(\mathbf{1(}%
n))^{i}=0$ for $i>n.$ For a given isocrystal $D_{0}$ over $K_{0}$ we define $%
D_{0}(n):=D_{0}\otimes _{K_{0}}\mathbf{1(}n)$.

Let $\lambda =s/r\in 
%TCIMACRO{\U{211a} }%
%BeginExpansion
\mathbb{Q}
%EndExpansion
$ where $r,s\in 
%TCIMACRO{\U{2124} }%
%BeginExpansion
\mathbb{Z}
%EndExpansion
$ are coprime integers and $r>0$. Define $D_{\lambda }$ (also denoted $%
D_{r,s}$) to be the isocrystal $(K_{0}^{r},F:=\Phi \cdot \sigma )$ where $%
\Phi $ is the $r\times r$ matrix having entry $1$ on any slot of the first
upper-diagonal, entry $p^{s}$ in the last row of the first column, and zeros
elsewhere. Notice that $D_{\lambda }$ has height $r$ and dimension $s$.

\noindent Over an algebraically closed field, we understand the category of $%
K_{0}$-isocrystal (cf. \cite{DM}, IV, 4):

\begin{theorem}
If $k=\overline{k}$, the abelian category of $K_{0}$-isocrystals is
semi-simple, and the simple objects of this category are given, up to
isomorphism and without repetition of isomorphism classes, by the $%
D_{\lambda }$ ($\lambda \in 
%TCIMACRO{\U{211a} }%
%BeginExpansion
\mathbb{Q}
%EndExpansion
$).
\end{theorem}

Let $k$ be algebraically closed; if $N$ is an isocrystal over $K_{0}=W(k)[%
\frac{1}{p}]$, the component of the slope $\lambda \in 
%TCIMACRO{\U{211a} }%
%BeginExpansion
\mathbb{Q}
%EndExpansion
$ in $N$ is the sum of the subspaces of $N$ isomorphic to $D_{\lambda }$.
The multiplicity of $\lambda $ is the $K_{0}$-dimension of this component.
The \textit{slope sequence of }$N$ is the non-decreasing sequence $\lambda
_{1}\leq ...\leq \lambda _{h}$ ($h=\dim _{K_{0}}N$) of all slopes of $N$,
each appearing according to its multiplicity. We say that $N$ is \textit{%
isoclinic} if $\lambda _{1}=\lambda _{n}$.

\bigskip

Let $W[F,V]$ be the quotient of the associative free $W$-algebra generated
by the indeterminates $F,V$ with respect to the relations: $Fa=a^{\sigma
}F,Va^{\sigma }=aV,FV=VF=p$ (any $a\in W$). A \textit{Dieudonn\'{e} }$W$%
\textit{-module} is a finitely generated left $W[F,V]$-module. An $F$\textit{%
-lattice \ over }$W$ is a $W$-free module of finite rank endowed with a
Frobenius semilinear injective endomorphism $F$. An $F$-lattice embedded as
a subobject in an isocrystal $(D_{0},F)$ is called a \textit{crystal} if it
is stable under the Verschiebung map.

If $D$ is an $F$-lattice over $W$ with Frobenius $F$ such that $pD\subseteq
FD$, there a unique operator $V:D\rightarrow D$ can be defined on $D$ in
such a way that $D$ becomes a module over $W[F,V]$ (i.e. a Dieudonn\'{e}
module that is finite and free as a $W$-module): if $x\in D$, we set $Vx:=y$
if and only if $px=Fy$. If $D$ is an $F$-lattice over $W$, then $D\otimes
_{W}%
%TCIMACRO{\U{211a} }%
%BeginExpansion
\mathbb{Q}
%EndExpansion
_{p}$ is an isocrystal over $K_{0}$; on the other side, every $F$-lattice
can be realized as an $F$-stable $W$-sublattice of some isocrystal. In
particular, $W$-free finite Dieudonn\'{e} modules are just $W$-sublattices
of isocrystals that are stable under $F$ and $V:=pF^{-1}$.

Let $D$ be a Dieudonn\'{e} module over $W$ that is a finite free $W$-module;
a \textit{polarization} of $D$ is a $W$-bilinear non-degenerate alternating
form $\left\langle ,\right\rangle :D\times D\rightarrow W$ such that $%
\left\langle Fx,y\right\rangle =\left\langle x,Vy\right\rangle ^{\sigma }$
for all $x,y\in D$; a $%
%TCIMACRO{\U{2124} }%
%BeginExpansion
\mathbb{Z}
%EndExpansion
_{p}$\textit{-homogeneous polarization} of $D$ is the equivalence class of $%
%TCIMACRO{\U{2124} }%
%BeginExpansion
\mathbb{Z}
%EndExpansion
_{p}^{\times }$-multiples of a given polarization. A polarization is called
a \textit{principal polarization} if it is a perfect pairing.

If $D_{0}$ is an isocrystal over $K_{0}$, a polarization of $D_{0}$ is a $%
K_{0}$-bilinear non-degenerate alternating pairing of isocrystals $%
\left\langle ,\right\rangle :D_{0}\times D_{0}\rightarrow K_{0}(1)$, where $%
K_{0}(1)=(K_{0},p\sigma )$; this means that $\left\langle ,\right\rangle $
is a $K_{0}$-bilinear non-degenerate alternating pairing such that $%
\left\langle Fx,Fy\right\rangle =p\left\langle x,y\right\rangle ^{\sigma }$
for all $x,y\in D_{0}$. One defines $%
%TCIMACRO{\U{211a} }%
%BeginExpansion
\mathbb{Q}
%EndExpansion
_{p}$-homogeneous polarizations in a similar way as above.

If $D$ is a Dieudonn\'{e} module over $W$ that is finite and free as $W$%
-module, and it is endowed with a polarization $\left\langle ,\right\rangle
:D\times D\rightarrow W$, then $D_{0}:=D[\frac{1}{p}]$ is canonically
endowed with a polarization of isocrystals $\left\langle ,\right\rangle
:D_{0}\times D_{0}\rightarrow K_{0}(1)$. Viceversa, if we start with a
polarization of the isocrystal $D[\frac{1}{p}]$ such that its restriction to 
$D\times D$ takes value into $W$, then a polarization of $D$ remains defined
(it needs not to be principal).

\bigskip

We have the following classification result (cf. \cite{F}, \cite{DM}):

\begin{theorem}
\label{Dieudonne}There is an additive anti-equivalence $M$ of categories
from the category of $p$-divisible groups over $k$ and the category of
Dieudonn\'{e} modules over $W$ that are free and of finite rank as $W$%
-modules. If $G$ is a $p$-divisible group over $k$ we call $M(G)$ its
Dieudonn\'{e} module, and we have:

\begin{enumerate}
\item the height of $G$ is equal to the $W$-rank of $M(G);$

\item if $k^{\prime }/k$ is an extension of perfect fields, there is a
canonical functorial isomorphism of Dieudonn\'{e} modules $M(G_{k^{\prime
}})\simeq M(G)\otimes _{W}W_{k^{\prime }};$

\item there is a canonical functorial isomorphism $M(\widehat{G})\simeq 
\limfunc{Hom}_{W}(M(G),W)$ of Dieudonn\'{e} modules, where the $W[F,V]$%
-structure on the left hand side is given by: $(Ff)(x)=f(Vx)^{\sigma }$ and $%
(Vf)(x)=f(Fx)^{\sigma ^{-1}}$ ( $x\in M(G),$ $f\in \limfunc{Hom}_{W}(M(G),W)$%
);

\item $M(G[n])\simeq M(G)/p^{n}M(G)$ for any positive integer $n$.
\end{enumerate}
\end{theorem}

The functor $M_{\ast }$ defined as $M_{\ast }(G):=M(\widehat{G})$ for any $p$%
-divisible group $G$ over $k$ is called the \textit{covariant Dieudonn\'{e}
functor}; in general we can switch from statements using the contravariant
Dieudonn\'{e} functor into statements with the covariant Dieudonn\'{e}
functor, by using duality of commutative formal groups. In the sequel, we
will sometimes use $M_{\ast }$ instead of $M$, as it is easier to work with
covariant functors than with contravariant.

If $D=M(G)$ for some $p$-divisible group $G$ over $k$, then $D[\frac{1}{p}%
]:=D\otimes _{W}%
%TCIMACRO{\U{211a} }%
%BeginExpansion
\mathbb{Q}
%EndExpansion
_{p}$ is an isocrystals over $K_{0}$ and $D$ is an $F$-lattice over $W$ (and
also a crystal). The functor $M$ of Theorem \ref{Dieudonne} gives rise to an
additive anti-equivalence between the category of $p$-divisible groups over $%
k$ and the category of $F$-lattices $D$ over $W$ such that $pD\subseteq FD$.

\noindent If $G_{1}$ and $G_{2}$ are two $p$-divisible groups over $k$, a
morphism $f:G_{1}\rightarrow G_{2}$ is an isogeny if and only if $%
M(f):M(G_{2})\rightarrow M(G_{1})$ is injective, if and only if $\limfunc{%
coker}M(f)$ is finite, and if and only if $M(f)$ induces an isomorphism $%
M(G_{2})[\frac{1}{p}]\rightarrow M(G_{1})[\frac{1}{p}]$ of isocrystals over $%
K_{0}$. The classification of $p$-divisible groups over $k$ up to isogenies
is equivalent to the classification of isocrystal over $K_{0}$ that contain
a lattice stable under $F$ and $pF^{-1}$. Also, there is a bijection between
the set of quasi-isogenies $G_{1}\rightarrow G_{2}$ and the set of
isomorphisms $M(G_{2})[\frac{1}{p}]\simeq M(G_{1})[\frac{1}{p}]$ of
isocrystals.

\begin{theorem}
\label{Dieudonne-Manin thm}\textbf{(Dieudonn\'{e} - Manin) }Let $k$ be
algebraically closed; then every $p$-divisible group $G$ over $k$ can be
written, up to isogeny, as:%
\begin{equation*}
G\sim \dbigoplus\limits_{0\leq \lambda \leq 1}G_{\lambda }^{\oplus g(\lambda
)},
\end{equation*}

\noindent where $\lambda $ is a rational number between $0$ and $1$
(included), and $g(\lambda )$ is a non-negative integer, equal to zero for
all $\lambda $ except at most a finite number of them. This decomposition
(up to isogeny)\ is uniquely determined by $G$.
\end{theorem}

If $G$ is a $p$-divisible group over $k$ and we write $G\sim
\tbigoplus\nolimits_{0\leq \lambda \leq 1}G_{\lambda }^{\oplus g(\lambda )}$
as in the above theorem, the \textit{slope sequence} of $G$ is the slope
sequence of the associated isocrystal, hence it is given by $0\leq \lambda
_{1}\leq ...\leq \lambda _{h}\leq 1$ (where $h=ht(G)$)\ and $G$ is isoclinic
if and only if $G\sim G_{\lambda }^{g}$ for some slope $\lambda $ and some $%
g\geq 0$. We have that $h=ht(G)$ and $ht(G)=\dim G+\dim \widehat{G}$;
furthermore $\widehat{G}$ has slope sequence $0\leq 1-\lambda _{h}\leq
...\leq 1-\lambda _{1}\leq 1$.

The following result will be crucial later (cf. \cite{DM}, page 92):

\begin{proposition}
\label{isog-isom}Let $k$ be algebraically closed. If the $p$-divisible group 
$G$ over $k$ is isogenous to $G_{1/r}$ (resp. to $G_{(r-1)/r}$), then $G$ is
isomorphic to it.
\end{proposition}

Let $G$ be a $p$-divisible group over $k$ endowed with a left action $i$ of
a $%
%TCIMACRO{\U{2124} }%
%BeginExpansion
\mathbb{Z}
%EndExpansion
_{p}$-algebra with involution $(O,\ast )$; then $M(G)$ is endowed, by
functoriality, with a left action of $O^{opp}$. On the other side, $M_{\ast
}(G)$ is endowed with an action of $O$. The anti-equivalence of Theorem \ref%
{Dieudonne} gives the following:

\begin{proposition}
\label{Rem-polarization} Let $k$ be any perfect field of characteristic $p$.

\begin{enumerate}
\item A polarization $\lambda :G\rightarrow \widehat{G}$ of a $p$-divisible
group over $k$ that it is also an isogeny determines a polarization of $M(%
\widehat{G})=M(G)^{\symbol{94}}$ which is principal if and only if $\lambda $
is principal (in this case, it can then be identified with a principal
polarization of $M(G)$). Viceversa, a polarization (resp. principal
polarization) of $M(\widehat{G})$ determines a polarization (resp. principal
polarization) of $G$. If $\lambda :G\rightarrow \widehat{G}$ is a
polarization (and not necessarily an isogeny), it induces a polarization of $%
M(\widehat{G})[\frac{1}{p}]$, and viceversa.

\item Let $G$ be a $p$-divisible group endowed with the action $i$ of a $%
%TCIMACRO{\U{2124} }%
%BeginExpansion
\mathbb{Z}
%EndExpansion
_{p}$-algebra with involution $(O,\ast )$, and let $\lambda $ be a
polarization of $G$ respecting such an action; then the polarization $%
\left\langle ,\right\rangle $ induced on the left $O^{opp}$-module $M(%
\widehat{G})[\frac{1}{p}]$ is skew-hermitian with respect to the involution
acting on $M(\widehat{G})$, i.e. $\left\langle bf,g\right\rangle
=\left\langle f,b^{\ast }g\right\rangle $ for all $b\in O^{opp},f,g\in M(%
\widehat{G})$.\ The viceversa is also true.
\end{enumerate}
\end{proposition}

\textbf{Proof. }Let $\lambda :G\rightarrow \widehat{G}$ be a polarization
that is also an isogeny of $p$-divisible groups. Then $M(\lambda ):M(%
\widehat{G})\rightarrow M(G)$ is a monomorphism of Dieudonn\'{e} modules.
Using the canonical identification $M(\widehat{G})=\limfunc{Hom}%
\nolimits_{W}(M(G),W)$ coming from Theorem \ref{Dieudonne}, we obtain a
non-degenerate $W$-linear pairing $\left\langle ,\right\rangle :M(\widehat{G}%
)\times M(\widehat{G})\rightarrow W$ defined by $\left\langle
f,g\right\rangle :=f(M(\lambda )(g))$ for any $f,g\in \limfunc{Hom}%
_{W}(M(G),W)$. We have $\left\langle Ff,g\right\rangle =(Ff)(M(\lambda
)(g))=f(VM(\lambda )(g))^{\sigma }$; since $M(\lambda )$ is a map of Dieudonn%
\'{e} modules, it commutes with $V$, so that $\left\langle Ff,g\right\rangle
=f(M(\lambda )(Vg))^{\sigma }=\left\langle f,Vg\right\rangle ^{\sigma }$.
Furthermore, $\left\langle f,g\right\rangle =g(M(\widehat{\lambda })(f))$,
so that $\left\langle ,\right\rangle $ is alternating since $\lambda $ is
anti-symmetric. Notice that $\lambda $ is principal if and only if $%
M(\lambda )$ is an isomorphism, i.e. if and only if $\left\langle
,\right\rangle $ is a perfect $W$-pairing. Using the fact that $M$ is an
anti-equivalence of categories, we obtain that viceversa we can pass from
polarizations of Dieudonn\'{e} modules that are finite and free as $W$%
-modules to polarizations of $p$-divisible groups. The statement for
isocrystals follows after inverting $p.$

For the second statement, let $b\in O$, $f,g\in M(\widehat{G})$ and assume
that $\lambda :G\rightarrow \widehat{G}$ is an isogeny. Notice that the
compatibility of $\lambda $ with $i$ implies that $M(i(b))\circ M(\lambda
)=M(\lambda )\circ M(\widehat{i}(b))$ and similarly for $b^{\ast }$. Then we
have:%
\begin{eqnarray*}
\left\langle bf,g\right\rangle &:&=\left\langle M(\widehat{i}%
(b))(f),g\right\rangle =\left( M(i(b^{\ast })^{\symbol{94}})(f)\right)
(M(\lambda )g)= \\
&=&\left( f\circ M(i(b^{\ast }))\circ M(\lambda )\right) (g)=\left( f\circ
M(\lambda )\circ M(\widehat{i}(b^{\ast }))\right) (g)= \\
&=&f(M(\lambda )(b^{\ast }g))=\left\langle f,b^{\ast }g\right\rangle .
\end{eqnarray*}

In a similar way, one can prove that a skew-hermitian polarization on $M(%
\widehat{G})[\frac{1}{p}]$ induces a $\ast $-compatible polarization of $G$. 
$\blacksquare $

\bigskip

\subsubsection{Isocrystals associated to a $K_{0}$-rational element of a
reductive group}

We still assume that $k$ is a perfect field of positive characteristic $p$,
with $W=\mathbb{W(}k)$ and $K_{0}=W[\frac{1}{p}]$. We also fix a connected
algebraic group $G$ over $%
%TCIMACRO{\U{211a} }%
%BeginExpansion
\mathbb{Q}
%EndExpansion
_{p}$ and an element $b\in G(K_{0})$. If $\rho :G\rightarrow GL(V)$ is a
rational algebraic $%
%TCIMACRO{\U{211a} }%
%BeginExpansion
\mathbb{Q}
%EndExpansion
_{p}$-representation of $G$ of finite dimension we define the isocrystal
associated to $(V;b)$ by the pair:%
\begin{equation*}
\left( V\otimes _{%
%TCIMACRO{\U{211a} }%
%BeginExpansion
\mathbb{Q}
%EndExpansion
_{p}}K_{0},\rho (b)(id_{V}\otimes \sigma )\right) \text{.}
\end{equation*}

\noindent (We denote $\rho (b)(id_{V}\otimes \sigma )$ also by $%
b(id_{V}\otimes \sigma )$ or $b\sigma $ is no confusion arises). An exact
functor remains therefore defined, from the category of rational $%
%TCIMACRO{\U{211a} }%
%BeginExpansion
\mathbb{Q}
%EndExpansion
_{p}$-representation of $G$ of finite dimension, to the category of $K_{0}$%
-isocrystal. Notice that if $g\in G(K_{0})$, and we set $b^{\prime
}:=gb(g^{\sigma })^{-1}$, then the functors associated to $b$ and $b^{\prime
}$ are isomorphic ($b$ and $b^{\prime }$ are in this case said to be $\sigma 
$-conjugate in $G(K_{0})$).

\paragraph{The slope morphism\label{slope morphism}}

Assume that $k$ is algebraically closed. Let $\mathbb{D}$ denote the
diagonalizable pro-algebraic group over $\mathbb{%
%TCIMACRO{\U{211a} }%
%BeginExpansion
\mathbb{Q}
%EndExpansion
}_{p}$ whose character group $X^{\ast }(\mathbb{D})$ (in the sense of
pro-algebraic groups) is isomorphic to $\mathbb{%
%TCIMACRO{\U{211a} }%
%BeginExpansion
\mathbb{Q}
%EndExpansion
}$: we can think of $\mathbb{D}$ as the universal cover of $\mathbb{G}_{m}$
in the sense of quasi-algebraic groups (cf. \cite{Ser}, 7.3). Notice that
the canonical inclusion $%
%TCIMACRO{\U{2124} }%
%BeginExpansion
\mathbb{Z}
%EndExpansion
\subset 
%TCIMACRO{\U{211a} }%
%BeginExpansion
\mathbb{Q}
%EndExpansion
$ induces a morphism $\mathbb{D\rightarrow G}_{m}$ and hence an inclusion $%
\limfunc{Hom}_{%
%TCIMACRO{\U{211a} }%
%BeginExpansion
\mathbb{Q}
%EndExpansion
_{p}}(\mathbb{G}_{m},G)\hookrightarrow \limfunc{Hom}_{%
%TCIMACRO{\U{211a} }%
%BeginExpansion
\mathbb{Q}
%EndExpansion
_{p}}(\mathbb{D},G).$ Furthermore, for any morphism $f\in \limfunc{Hom}_{%
%TCIMACRO{\U{211a} }%
%BeginExpansion
\mathbb{Q}
%EndExpansion
_{p}}(\mathbb{D},G)$ there is a positive integer $n$ such that $nf\in 
\limfunc{Hom}_{%
%TCIMACRO{\U{211a} }%
%BeginExpansion
\mathbb{Q}
%EndExpansion
_{p}}(\mathbb{G}_{m},G)$.

In \cite{Ko} \S 4,\ Kottwitz defines a morphism of over $K_{0}:$%
\begin{equation*}
\nu :\mathbb{D}_{K_{0}}\mathbb{\rightarrow }G_{K_{0}}
\end{equation*}

\noindent characterized by the following property:\ for any $%
%TCIMACRO{\U{211a} }%
%BeginExpansion
\mathbb{Q}
%EndExpansion
_{p}$-rational finite-dimensional representation $\rho :G\rightarrow GL(V)$
of $G$, let $\nu _{\rho }\in \limfunc{Hom}_{K_{0}}(\mathbb{D}%
_{K_{0}},GL(V_{K_{0}}))$ be the only morphism for which the action of $%
\mathbb{D}$ on the isotypical component of the isocrystal $(V\otimes _{%
%TCIMACRO{\U{211a} }%
%BeginExpansion
\mathbb{Q}
%EndExpansion
_{p}}K_{0},\rho (b)(id_{V}\otimes \sigma ))$ of slope $\lambda \in 
%TCIMACRO{\U{211a} }%
%BeginExpansion
\mathbb{Q}
%EndExpansion
$ is given by the character $\lambda \in X^{\ast }(\mathbb{D})$. Then $\nu $
is the only morphism over $K_{0}$ such that $\rho \circ \nu =\nu _{\rho }$
for any representation $\rho $ as above. In other words, $\nu $ is
characterized by the fact that for any representation $V$ of $G$, the
grading induced by $\nu $ on $V\otimes K_{0}$ coincides with the slope
decomposition of $(V\otimes K_{0},b(id_{V}\otimes \sigma ))$. We call such $%
\nu $ the \textit{slope morphism associated to }$b\in G(K_{0})$; that this
definition does not depend upon the choice of $\sigma $-conjugacy class of $%
b $.

\paragraph{Basic elements\label{basicity}}

Let us keep assuming that $k$ is algebraically closed and furthermore that $%
G $ is a connected and reductive algebraic group over $K_{0}$. Following 
\cite{Ko}, 5.1 we say that the element $b\in G(K_{0})$ is \textit{basic} if
the corresponding slope morphism $\nu :\mathbb{D}_{K_{0}}\mathbb{\rightarrow 
}G_{K_{0}}$ factors through the center of $G_{K_{0}}$.

Let $J$ be the functor on $%
%TCIMACRO{\U{211a} }%
%BeginExpansion
\mathbb{Q}
%EndExpansion
_{p}$-algebras defined by:%
\begin{equation*}
J(R):=\{g\in G(R\otimes _{%
%TCIMACRO{\U{211a} }%
%BeginExpansion
\mathbb{Q}
%EndExpansion
_{p}}K_{0}):g(b\sigma )=(b\sigma )g\}\text{,}
\end{equation*}

\noindent for any $%
%TCIMACRO{\U{211a} }%
%BeginExpansion
\mathbb{Q}
%EndExpansion
_{p}$-algebra $R$ (where we see $b\sigma \in G(K_{0})\rtimes \left\langle
\sigma \right\rangle $ and the action by conjugation of $\sigma $ on
elements of $G(K_{0})$ is the one obtained by viewing $\sigma $ as an
automorphism of $K_{0}$). Then $J$ defines a smooth affine group scheme over 
$%
%TCIMACRO{\U{211a} }%
%BeginExpansion
\mathbb{Q}
%EndExpansion
_{p}$ (\cite{RZ}, 1.12) and Kottwitz proves in \cite{Ko} that $b$ is basic
if and only if $J$ is an inner form of $G$. For example, if $V$ is a finite
dimensional $%
%TCIMACRO{\U{211a} }%
%BeginExpansion
\mathbb{Q}
%EndExpansion
_{p}$-vector space and $G=GL(V)$, an element $b\in G(K_{0})$ is basic if and
only if the corresponding isocrystal $(V\otimes _{%
%TCIMACRO{\U{211a} }%
%BeginExpansion
\mathbb{Q}
%EndExpansion
_{p}}K_{0},b\sigma )$ is isoclinic.

\paragraph{Admissibility\label{admissibility}}

We assume here that $k$ is a perfect field of characteristic $p$, with $%
K_{0}=W(k)[\frac{1}{p}]$. Let $K$ be a finite field extension of $K_{0}$. We
say that a $K$-filtered isocrystal $D$ over $K_{0}$ is \textit{admissible}
if the following condition is satisfied: for any $K$-filtered sub-isocrystal 
$D^{\prime }$ of $D$, the rightmost endpoint of the Newton polygon of $%
D^{\prime }$ lies on or above the rightmost endpoint of the Hodge polygon of 
$D^{\prime }$, with equality if $D=D^{\prime }$ (for the definition of
Newton and Hodge polygons see for example \cite{CB}, 8). The category of
admissible $K$-filtered isocrystal over $K_{0}$ is abelian.

Let $B_{cris}$ denote Fontaine's crystalline period ring: it is a $K_{0}$%
-algebra domain, endowed with a continuous action of the absolute Galois
group $G_{K}$ of $K$, such that $B_{cris}^{G_{K}}=K_{0}$; furthermore $%
B_{cris}$ is endowed with a $\sigma $-linear Frobenius endomorphism, and
with a (non Frobenius-stable) $K$-filtration (i.e. a filtration of $K\otimes
_{K_{0}}B_{cris}$) induced by the canonical inclusion:

\begin{equation*}
K\otimes _{K_{0}}B_{cris}\hookrightarrow B_{dR}\text{,}
\end{equation*}

\noindent where $B_{dR}$ is the de Rham period ring, with its canonical
filtration (cf. \cite{Fon}).

Let $\limfunc{Rep}_{%
%TCIMACRO{\U{211a} }%
%BeginExpansion
\mathbb{Q}
%EndExpansion
_{p}}(G_{K})$ be the category of $%
%TCIMACRO{\U{211a} }%
%BeginExpansion
\mathbb{Q}
%EndExpansion
_{p}$-linear continuous representations of $G_{K}$ of finite dimension and,
for any $W\in \limfunc{Rep}_{%
%TCIMACRO{\U{211a} }%
%BeginExpansion
\mathbb{Q}
%EndExpansion
_{p}}(G_{K})$ define $D_{cris}(W):=(W\otimes _{%
%TCIMACRO{\U{211a} }%
%BeginExpansion
\mathbb{Q}
%EndExpansion
_{p}}B_{cris})^{G_{K}}$ and say that $W$ is crystalline if $\dim _{%
%TCIMACRO{\U{211a} }%
%BeginExpansion
\mathbb{Q}
%EndExpansion
_{p}}W=\dim _{K_{0}}D_{cris}(W)$. Let $\limfunc{Rep}_{%
%TCIMACRO{\U{211a} }%
%BeginExpansion
\mathbb{Q}
%EndExpansion
_{p}}^{cris}(G_{K})$ be the subcategory of $\limfunc{Rep}_{%
%TCIMACRO{\U{211a} }%
%BeginExpansion
\mathbb{Q}
%EndExpansion
_{p}}(G_{K})$ consisting of crystalline representations; then $D_{cris}$
defines a fully faithful exact tensor functor from $\limfunc{Rep}_{%
%TCIMACRO{\U{211a} }%
%BeginExpansion
\mathbb{Q}
%EndExpansion
_{p}}^{cris}(G_{K})$ into the category of $K$-filtered isocrystal over $%
K_{0} $. It was proven by Colmez and Fontaine (\cite{FonC}) that this
functor gives rise to an equivalence between the category $\limfunc{Rep}_{%
%TCIMACRO{\U{211a} }%
%BeginExpansion
\mathbb{Q}
%EndExpansion
_{p}}^{cris}(G_{K})$ and the category of $K$-filtered isocrystal over $K_{0}$
that are admissible. Therefore in the sequel an admissible filtered
isocrystal will be the same as an isocrystal coming from a crystalline
representation via the functor $D_{cris}.$

Assume now that $k$ is algebraically closed. Let $G$ be a reductive
connected group defined over $%
%TCIMACRO{\U{211a} }%
%BeginExpansion
\mathbb{Q}
%EndExpansion
_{p}$\ and fix a cocharacter $\mu :\mathbb{G}_{m}\rightarrow G$ \noindent
defined over $K$ (we dropped the subscript $K$). If $V$ is a finite
dimensional $%
%TCIMACRO{\U{211a} }%
%BeginExpansion
\mathbb{Q}
%EndExpansion
_{p}$-rational representation of $G$ and $(V\otimes _{%
%TCIMACRO{\U{211a} }%
%BeginExpansion
\mathbb{Q}
%EndExpansion
_{p}}K_{0},b\sigma )$ is the associated isocrystal over $K_{0}$, we can
construct a $K$-filtered $K_{0}$-isocrystal $\mathcal{I=I}(V;b;\mu )$\
associated to the triple $(V;b;\mu )$ by letting:%
\begin{equation*}
\mathcal{I}:\mathcal{=(}V\otimes _{%
%TCIMACRO{\U{211a} }%
%BeginExpansion
\mathbb{Q}
%EndExpansion
_{p}}K_{0},b\sigma ,V_{K}^{\bullet })
\end{equation*}

\noindent where, for any $i\in 
%TCIMACRO{\U{2124} }%
%BeginExpansion
\mathbb{Z}
%EndExpansion
$, we set:%
\begin{equation*}
V_{K}^{i}=\tbigoplus\nolimits_{j\geq i}(V\otimes _{%
%TCIMACRO{\U{211a} }%
%BeginExpansion
\mathbb{Q}
%EndExpansion
_{p}}K)_{j},
\end{equation*}

\noindent with the subscript $j$ denoting the $j$ weight space for the
action of $\mathbb{G}_{m}$ on $V_{K}:=V\otimes _{%
%TCIMACRO{\U{211a} }%
%BeginExpansion
\mathbb{Q}
%EndExpansion
_{p}}K$ induced by $\mu $.

Keeping the above notation, we say that the pair $(b;\mu )$ is \textit{%
admissible} if for any finite dimensional $%
%TCIMACRO{\U{211a} }%
%BeginExpansion
\mathbb{Q}
%EndExpansion
_{p}$-rational representation $V$ of $G$, the filtered isocrystal $\mathcal{I%
}(V;b;\mu )$ is admissible in the sense explained above. Notice that this is
equivalent to say that $\mathcal{I}(V;b;\mu )$ is admissible for some
faithful finite dimensional $%
%TCIMACRO{\U{211a} }%
%BeginExpansion
\mathbb{Q}
%EndExpansion
_{p}$-rational representation $V$ of $G$ (cf. \cite{RZ}, 1.18).

\subsection{Moduli of $p$-divisible groups of PEL-type\label{exp1}}

We refer for this section to \cite{RZ}, Chapter 3; cf. also \cite{Bou}.

\subsubsection{Definition of local PEL-data\label{def - admiss}}

We define local PEL-data; we will define later on PEL-data of global type,
that will be the starting point for the construction of moduli spaces of
abelian schemes of PEL-type.

Let $B$ be a finite dimensional semi-simple $%
%TCIMACRO{\U{211a} }%
%BeginExpansion
\mathbb{Q}
%EndExpansion
_{p}$-algebra endowed with an involution $^{\ast }$. Let $V\neq \{0\}$ be a
finitely generated left $B$-module and $\left\langle ,\right\rangle :V\times
V\rightarrow 
%TCIMACRO{\U{211a} }%
%BeginExpansion
\mathbb{Q}
%EndExpansion
_{p}$ a non-degenerate, alternating $%
%TCIMACRO{\U{211a} }%
%BeginExpansion
\mathbb{Q}
%EndExpansion
_{p}$-bilinear form which is skew-Hermitian with respect to $^{\ast }$, i.e. 
$\left\langle tv,w\right\rangle =\left\langle v,t^{\ast }w\right\rangle $
for all $v,w$ in $V$ and all $t$ in $B$. These objects define a reductive
algebraic group scheme $G:=G\left( B,^{\ast },V,\left\langle ,\right\rangle
\right) $ over $%
%TCIMACRO{\U{211a} }%
%BeginExpansion
\mathbb{Q}
%EndExpansion
_{p}$ whose $R$-points, for a fixed $%
%TCIMACRO{\U{211a} }%
%BeginExpansion
\mathbb{Q}
%EndExpansion
_{p}$-algebra $R$, are given by:%
\begin{equation*}
G(R)=\left\{ \left( g,s\right) \in GL_{B\otimes _{%
%TCIMACRO{\U{211a} }%
%BeginExpansion
\mathbb{Q}
%EndExpansion
_{p}}R}(V\otimes _{%
%TCIMACRO{\U{211a} }%
%BeginExpansion
\mathbb{Q}
%EndExpansion
_{p}}R)\times \mathbb{G}_{m}\left( R\right) :\left\langle gv,gw\right\rangle
=s\left\langle v,w\right\rangle \text{ }\forall v,w\in V\right\} \text{.}
\end{equation*}

\noindent The map $G(R)\rightarrow R^{\times }$ given by $\left( g,s\right)
\mapsto s$ defines a homomorphism of $%
%TCIMACRO{\U{211a} }%
%BeginExpansion
\mathbb{Q}
%EndExpansion
_{p}$-algebraic groups $c:G\rightarrow \mathbb{G}_{m}$ that is called the 
\textit{similitude character} of $G$; the kernel of $c$ is a reductive $%
%TCIMACRO{\U{211a} }%
%BeginExpansion
\mathbb{Q}
%EndExpansion
_{p}$-subgroup of $G$ denoted by $G_{1}$. By abuse of notation, an element $%
\left( g,s\right) \in G(R)$ will be often denoted by $g$.

Notice that $^{\ast }$ defines on the $%
%TCIMACRO{\U{211a} }%
%BeginExpansion
\mathbb{Q}
%EndExpansion
_{p}$-algebra $C:=\limfunc{End}_{B}V$ an involution $x\longmapsto x^{\ast }$
via the identity $\left\langle xv,w\right\rangle =\left\langle v,x^{\ast
}w\right\rangle \ \ \ (v,w\in V).$

\noindent We therefore have functorial isomorphisms:

\begin{equation*}
G(R)\simeq \left\{ x\in C\otimes _{%
%TCIMACRO{\U{211a} }%
%BeginExpansion
\mathbb{Q}
%EndExpansion
_{p}}R:xx^{\ast }\in R^{\times }\right\} \text{,}
\end{equation*}

\noindent where $R$ is any $%
%TCIMACRO{\U{211a} }%
%BeginExpansion
\mathbb{Q}
%EndExpansion
_{p}$-algebra.

\begin{definition}
\label{PEL p-divisible gps}Fix a tuple $\mathcal{(}B,^{\ast },V,\left\langle
,\right\rangle )$ satisfying the above properties. Let $k$ be an
algebraically closed field of characteristic $p>0$ and set $W=W(k)$, $%
K_{0}=W[\frac{1}{p}]$; denote by $\sigma $ the Frobenius automorphism of $W$
and $K_{0}$.

\begin{enumerate}
\item The datum $\mathcal{D}:\mathcal{=(}B,^{\ast },V,\left\langle
,\right\rangle )$ is called a $%
%TCIMACRO{\U{211a} }%
%BeginExpansion
\mathbb{Q}
%EndExpansion
_{p}$\textbf{-PEL-datum}; the group $G=G\mathcal{(}B,^{\ast },V,\left\langle
,\right\rangle )$ is the algebraic group associated to $\mathcal{D}$. A%
\textbf{\ }$%
%TCIMACRO{\U{211a} }%
%BeginExpansion
\mathbb{Q}
%EndExpansion
_{p}$\textbf{-PEL-datum with integral structure }(or with a PEL-lattice) is
the datum $\mathcal{D}_{\mathcal{O}}:\mathcal{=(}B,^{\ast },V,\left\langle
,\right\rangle ,\mathcal{O}_{B},\Lambda )$, where $\mathcal{(}B,^{\ast
},V,\left\langle ,\right\rangle )$ is a $%
%TCIMACRO{\U{211a} }%
%BeginExpansion
\mathbb{Q}
%EndExpansion
_{p}$-PEL-datum, $\mathcal{O}_{B}$ is a maximal $%
%TCIMACRO{\U{2124} }%
%BeginExpansion
\mathbb{Z}
%EndExpansion
_{p}$-order in $B$ stable under the involution $^{\ast }$, and $\Lambda
\subset V$ is a $%
%TCIMACRO{\U{2124} }%
%BeginExpansion
\mathbb{Z}
%EndExpansion
_{p}$-lattice in $V$ which is also an $\mathcal{O}_{B}$-submodule and which
is self-dual with respect to the pairing $\left\langle ,\right\rangle $
(i.e. the restriction of $\left\langle ,\right\rangle $ to $\Lambda \times
\Lambda $ defines a perfect pairing of $%
%TCIMACRO{\U{2124} }%
%BeginExpansion
\mathbb{Z}
%EndExpansion
_{p}$-modules).

\item A\textbf{\ }$%
%TCIMACRO{\U{211a} }%
%BeginExpansion
\mathbb{Q}
%EndExpansion
_{p}$\textbf{-PEL-datum for moduli of }$p$\textbf{-divisible groups over }$k$
is the datum $\mathcal{D}_{\func{mod}}:\mathcal{=(}B,^{\ast },V,\left\langle
,\right\rangle ,\mathcal{O}_{B},\Lambda ,b,\mu )$, where $\mathcal{(}%
B,^{\ast },V,\left\langle ,\right\rangle ,\mathcal{O}_{B},\Lambda )$ is a $%
%TCIMACRO{\U{211a} }%
%BeginExpansion
\mathbb{Q}
%EndExpansion
_{p}$-PEL-datum with integral structure and associated group $G$, $b$ is a
fixed element of $G(K_{0})$, and $\mu :\mathbb{G}_{m/K}\rightarrow G_{K}$ is
a co-character of $G$ defined over some finite field extension $K$ of $K_{0}$%
. We furthermore require that the following conditions are satisfied:

\begin{enumerate}
\item $(b,\mu )$ is an admissible pair in the sense of \ref{admissibility};

\item the isocrystal $(N,\mathbf{F}):=(V\otimes _{%
%TCIMACRO{\U{211a} }%
%BeginExpansion
\mathbb{Q}
%EndExpansion
_{p}}K_{0},b\sigma )$ associated to $\mathcal{D}_{\func{mod}}$ has slopes in
the interval $[0,1]$;

\item the weight decomposition of $V\otimes _{%
%TCIMACRO{\U{211a} }%
%BeginExpansion
\mathbb{Q}
%EndExpansion
_{p}}K$ with respect to $\mu $ contains only the weights $0$ and $1$: $%
V\otimes _{%
%TCIMACRO{\U{211a} }%
%BeginExpansion
\mathbb{Q}
%EndExpansion
_{p}}K=V_{0}\oplus V_{1}$;

\item if $\nu :\mathbb{D}_{K_{0}}\rightarrow G_{K_{0}}$ denotes the slope
morphism associated to $b$ and $c:G\rightarrow \mathbb{G}_{m}$ is the
similitude character of $G$ over $%
%TCIMACRO{\U{211a} }%
%BeginExpansion
\mathbb{Q}
%EndExpansion
_{p}$, then $c\nu :$ $\mathbb{D}_{K_{0}}\rightarrow \mathbb{G}_{m/K_{0}}$ is
the character of $\mathbb{D}_{K_{0}}$ corresponding to the rational number $%
1 $.
\end{enumerate}

The \textbf{reflex field }(or local Shimura field) of $\mathcal{D}_{\func{mod%
}}$ is the field of definition of the conjugacy class of the co-character $%
\mu $.
\end{enumerate}
\end{definition}

Each of the above $%
%TCIMACRO{\U{211a} }%
%BeginExpansion
\mathbb{Q}
%EndExpansion
_{p}$-PEL-data will be called simple if the algebra $B$ is simple. The $%
%TCIMACRO{\U{211a} }%
%BeginExpansion
\mathbb{Q}
%EndExpansion
_{p}$-PEL-data $\mathcal{D}_{\mathcal{O}}$ or $\mathcal{D}_{\func{mod}}$
have \textit{good reduction} if $B$ is an unramified $%
%TCIMACRO{\U{211a} }%
%BeginExpansion
\mathbb{Q}
%EndExpansion
_{p}$-algebra (i.e. its center is the product of matrix algebras over
unramified field extensions of $%
%TCIMACRO{\U{211a} }%
%BeginExpansion
\mathbb{Q}
%EndExpansion
_{p}$). If the PEL-datum $\mathcal{D}_{\mathcal{O}}$ or $\mathcal{D}_{\func{%
mod}}$ has good reduction at $p$, the associated group $G_{%
%TCIMACRO{\U{211a} }%
%BeginExpansion
\mathbb{Q}
%EndExpansion
_{p}}$ is unramified, i.e. it is quasi-split over $%
%TCIMACRO{\U{211a} }%
%BeginExpansion
\mathbb{Q}
%EndExpansion
_{p}$ and split over an unramified extension of $%
%TCIMACRO{\U{211a} }%
%BeginExpansion
\mathbb{Q}
%EndExpansion
_{p}$ or, equivalently, $G\left( 
%TCIMACRO{\U{211a} }%
%BeginExpansion
\mathbb{Q}
%EndExpansion
_{p}\right) $ has an hyperspecial subgroup. \noindent In particular $G_{%
%TCIMACRO{\U{211a} }%
%BeginExpansion
\mathbb{Q}
%EndExpansion
_{p}}$ has a reductive model $\mathcal{G}$ over $%
%TCIMACRO{\U{2124} }%
%BeginExpansion
\mathbb{Z}
%EndExpansion
_{p}$ such that for any $%
%TCIMACRO{\U{2124} }%
%BeginExpansion
\mathbb{Z}
%EndExpansion
_{p}$-algebra $R$ we have:%
\begin{equation*}
\mathcal{G}\left( R\right) =\left\{ g\in GL_{R\otimes _{%
%TCIMACRO{\U{2124} }%
%BeginExpansion
\mathbb{Z}
%EndExpansion
_{p}}\mathcal{O}_{B}}(\Lambda \otimes _{%
%TCIMACRO{\U{2124} }%
%BeginExpansion
\mathbb{Z}
%EndExpansion
_{p}}R):\left\langle g\lambda _{1},g\lambda _{2}\right\rangle
=c(g)\left\langle \lambda _{1},\lambda _{2}\right\rangle ,c(g)\in R^{\times
}\right\} .
\end{equation*}

The Shimura field $E$ of a fixed datum $\mathcal{D}_{\func{mod}}$ is a
finite extension of $%
%TCIMACRO{\U{211a} }%
%BeginExpansion
\mathbb{Q}
%EndExpansion
_{p}$, contained inside $K$ ($K$ is the field of definition of the
co-character $\mu $). Notice that $E=%
%TCIMACRO{\U{211a} }%
%BeginExpansion
\mathbb{Q}
%EndExpansion
_{p}(tr(\mu (x)):x\in K^{\times }).$ \noindent If $\mathcal{D}_{\func{mod}}$
has good reduction, then $E$ is an unramified extension of $%
%TCIMACRO{\U{211a} }%
%BeginExpansion
\mathbb{Q}
%EndExpansion
_{p}.$

\bigskip

Let us make some comments on the definition of local PEL-datum for moduli of 
$p$-divisible groups (cf. \cite{RZ}, 3.19). We keep the notation from the
above Definition, and fix a local PEL-datum for moduli of $p$-divisible
groups $\mathcal{D}_{\func{mod}}$. We also denote by $\mathbb{D}_{\ast
}(\cdot )$ the covariant functor that assigns to a $p$-divisible group over
a $%
%TCIMACRO{\U{2124} }%
%BeginExpansion
\mathbb{Z}
%EndExpansion
_{p}$-scheme of characteristic $p$ the associated Dieudonn\'{e} crystal (cf. 
\cite{BBM}, \cite{Me}).

$\bullet $ The condition $b)$ on $\mathcal{D}_{\func{mod}}$ is equivalent to
requiring that the $K_{0}$-isocrystal $(N,\mathbf{F})$ comes from a $p$%
-divisible group $\mathbf{X\ }$(defined over $k$) via the \textit{covariant}
Dieudonn\'{e} functor $M_{\ast }$ (cf. Th. \ref{Dieudonne-Manin thm}); such
an $\mathbf{X}$ is uniquely determined up to $k$-quasi-isogeny and we fix
for the sequel such a choice of $\mathbf{X}$.

$\bullet $ Assume that there are a $p$-divisible group $X$ over the ring of
integers $\mathcal{O}_{K}$ of $K$ and a $k$-quasi-isogeny $\varphi :\mathbf{%
X\rightarrow }X_{k}$; we can lift $\varphi $ to a quasi-isogeny of $p$%
-divisible groups over $\limfunc{Spec}(\mathcal{O}_{K}/(p))$:%
\begin{equation*}
\varphi :\mathbf{X\times }_{\limfunc{Spec}k}\limfunc{Spec}\left( \mathcal{O}%
_{K}/(p)\right) \mathbf{\longrightarrow }X\times _{\limfunc{Spf}(\mathcal{O}%
_{K})}\limfunc{Spec}\left( \mathcal{O}_{K}/(p)\right) .
\end{equation*}

\noindent Applying the covariant functor $\mathbb{D}_{\ast }(\cdot )$\ to
the above quasi-isogeny, and then evaluating the corresponding Dieudonn\'{e}
crystals of $\mathcal{O}_{\limfunc{Spec}(\mathcal{O}_{K}/(p))}$-modules at
the PD-thickening $\limfunc{Spec}(\mathcal{O}_{K}/(p))\hookrightarrow 
\limfunc{Spf}(\mathcal{O}_{K})$ belonging to the small crystalline site of $%
\limfunc{Spec}(\mathcal{O}_{K}/(p))$, we obtain an isomorphism of $K$%
-isocrystals:%
\begin{equation*}
\varphi :(N,\mathbf{F})^{(p)}\otimes _{K_{0}}K\overset{\simeq }{%
\longrightarrow }\mathbb{D}_{\ast }\left( X\times _{\limfunc{Spf}(\mathcal{O}%
_{K})}\limfunc{Spec}\left( \mathcal{O}_{K}/(p)\right) \right) _{\limfunc{Spf}%
(\mathcal{O}_{K})}\otimes 
%TCIMACRO{\U{211a} }%
%BeginExpansion
\mathbb{Q}
%EndExpansion
,
\end{equation*}

\noindent where we used the fact that, by assumption, $M_{\ast }(\mathbf{X}%
)[p^{-1}]=(N,\mathbf{F)}$.

\noindent If we assume that, under this identification, the Hodge filtration
of the crystal $\mathbb{D}_{\ast }\left( X\times _{\limfunc{Spf}(\mathcal{O}%
_{K})}\limfunc{Spec}\left( \mathcal{O}_{K}/(p)\right) \right) _{\limfunc{Spf}%
(\mathcal{O}_{K})}\otimes 
%TCIMACRO{\U{211a} }%
%BeginExpansion
\mathbb{Q}
%EndExpansion
$\ corresponds to the filtration induced by $\mu $ on $V\otimes _{%
%TCIMACRO{\U{211a} }%
%BeginExpansion
\mathbb{Q}
%EndExpansion
_{p}}K$: $0\rightarrow V_{1}\rightarrow V\otimes _{%
%TCIMACRO{\U{211a} }%
%BeginExpansion
\mathbb{Q}
%EndExpansion
_{p}}K\rightarrow V_{0}\rightarrow 0,$\noindent then the conditions $%
a),b),c) $ of the above definition of $\mathcal{D}_{\func{mod}}$ are
automatically satisfied: for the pair $(b,\mu )$ is, in this case of "good
reduction", admissible in the sense of \autoref{admissibility}. \noindent In
the later application of these constructions, we will be in the "good
reduction" case, so we will not need to check conditions $a),b),c)$ on the
PEL-data we will have in hands.

$\bullet $ Let $\chi :G\rightarrow \mathbb{G}_{m}$ be a $%
%TCIMACRO{\U{211a} }%
%BeginExpansion
\mathbb{Q}
%EndExpansion
_{p}$-rational character of $G$. Then $(K_{0}(\chi ),b\sigma )$ is a one
dimensional $K_{0}$-isocrystal whose only slope is $\limfunc{ord}_{p}\chi
(b) $, since $(b\sigma )(W)=\chi (b)W=p^{\limfunc{ord}_{p}\chi (b)}W$. If
the pair $(b,\mu )$ is admissible, the $K$-filtered $K_{0}$-isocrystal
associated to $(K_{0}(\chi ),b\sigma )$ and $\mu $ need to satisfy the
admissibility equation $\tsum\nolimits_{i\in 
%TCIMACRO{\U{2124} }%
%BeginExpansion
\mathbb{Z}
%EndExpansion
}i\cdot \dim _{K}K\left( \chi \right) ^{i}=\limfunc{ord}\nolimits_{p}\chi
(b),$\noindent\ where $K\left( \chi \right) ^{i}$ is the $i^{th}$ term in
the filtration that is induced by $\mu $ on $K\left( \chi \right) $. We
conclude that $\left\langle \mu ,\chi \right\rangle =\limfunc{ord}%
\nolimits_{p}\chi (b).$

\noindent In particular:

$\bullet $ taking $\chi =\det_{%
%TCIMACRO{\U{211a} }%
%BeginExpansion
\mathbb{Q}
%EndExpansion
_{p}}$, defined by viewing $G\subset GL(V)$, we obtain $\dim _{K}V_{1}=%
\limfunc{ord}\nolimits_{p}(\det_{K_{0}}(b;V\otimes _{%
%TCIMACRO{\U{211a} }%
%BeginExpansion
\mathbb{Q}
%EndExpansion
_{p}}K_{0}))$;

$\bullet $ taking $\chi =c$ to be the similitude factor of $G$, condition $%
d) $ from the definition of $\mathcal{D}_{\func{mod}}$ implies that $%
\left\langle \mu ,c\right\rangle =\limfunc{ord}\nolimits_{p}c(b)=1.$ Hence
the subspaces $V_{0}$ and $V_{1}$ ov $V\otimes _{%
%TCIMACRO{\U{211a} }%
%BeginExpansion
\mathbb{Q}
%EndExpansion
_{p}}K$ are isotropic for the pairing induced by $\left\langle
,\right\rangle $ on $V\otimes _{%
%TCIMACRO{\U{211a} }%
%BeginExpansion
\mathbb{Q}
%EndExpansion
_{p}}K$. (In fact if $v\in V_{0}$, we have $\left\langle v,v\right\rangle
=\left\langle \mu (x)v,\mu (x)v\right\rangle =x\left\langle v,v\right\rangle 
$ for all $x\in K^{\times }$; similarly for $V_{1}$).

\bigskip

Fix a local PEL-datum for moduli of $p$-divisible groups $\mathcal{D}_{\func{%
mod}}$; denote by $(N,\mathbf{F})$ the associated isocrystal. The action of $%
B$ on $V$ by left multiplication induces an action of $B$ on $N$.
Furthermore, since $k$ is algebraically closed and $\limfunc{ord}%
\nolimits_{p}c(b)=1$, we can write $c(b)=p\cdot u\sigma (u)^{-1}$ for some $%
u\in W^{\times }$. Define a map $\Psi :N\times N\longrightarrow K_{0}$
\noindent by setting $\Psi (v,w):=u^{-1}\cdot \left\langle v,w\right\rangle $
for all $v,w\in N$; $\Psi $ is a non-degenerate $K_{0}$-bilinear pairing,
which is alternating and skew-Hermitian with respect to $^{\ast }$.
Furthermore, $\Psi (\mathbf{F}v,\mathbf{F}w)=p\Psi (v,w)^{\sigma }$, so that 
$\Psi $ defines a polarization of isocrystals $N\times N\longrightarrow 
\mathbf{1}(1)$.

\noindent Notice that any other choice of $u$ gives rise to a multiple of $%
\Psi $ by an element of $%
%TCIMACRO{\U{211a} }%
%BeginExpansion
\mathbb{Q}
%EndExpansion
_{p}^{\times }$, so that our local PEL-datum gives rise to a well defined
triple $(N,\mathbf{F},%
%TCIMACRO{\U{211a} }%
%BeginExpansion
\mathbb{Q}
%EndExpansion
_{p}^{\times }\Psi )$ \noindent that is to a $%
%TCIMACRO{\U{211a} }%
%BeginExpansion
\mathbb{Q}
%EndExpansion
_{p}$-homogeneously polarized $K_{0}$-isocrystal endowed with an action of $%
B $ (the polarization form is skew-Hermitian with respect to the involution
of $B$).

\subsubsection{The moduli functor for $p$-divisible groups\label%
{moduli-functor-p div def.}}

Fix a $%
%TCIMACRO{\U{211a} }%
%BeginExpansion
\mathbb{Q}
%EndExpansion
_{p}$-PEL-datum for moduli of $p$-divisible groups $\mathcal{D}_{\func{mod}}:%
\mathcal{=(}B,^{\ast },V,\left\langle ,\right\rangle ,\mathcal{O}%
_{B},\Lambda ,b,\mu )$ over the algebraically closed field $k$ of
characteristic $p$. Denote by $G$ the associated algebraic $%
%TCIMACRO{\U{211a} }%
%BeginExpansion
\mathbb{Q}
%EndExpansion
_{p}$-group, and let $(N,\mathbf{F},%
%TCIMACRO{\U{211a} }%
%BeginExpansion
\mathbb{Q}
%EndExpansion
_{p}^{\times }\Psi )$ be the $%
%TCIMACRO{\U{211a} }%
%BeginExpansion
\mathbb{Q}
%EndExpansion
_{p}$-homogeneously polarized $K_{0}$-isocrystal with $B$-action associated
to $\mathcal{D}_{\func{mod}}$, as in the previous section.

Pick a $p$-divisible group $\mathbf{X}$ over $k$ whose $K_{0}$-isocrystal
constructed via $M_{\ast }$ is isomorphic to $(N,\mathbf{F})$; by Prop. \ref%
{Rem-polarization}, $B$ acts on the group of quasi-isogenies of $\mathbf{X}$%
, so that we have a $%
%TCIMACRO{\U{211a} }%
%BeginExpansion
\mathbb{Q}
%EndExpansion
_{p}$-algebra homomorphism:%
\begin{equation*}
i_{\mathbf{X}}:B\rightarrow \limfunc{End}(\mathbf{X)\otimes }_{%
%TCIMACRO{\U{2124} }%
%BeginExpansion
\mathbb{Z}
%EndExpansion
_{p}}\mathbf{%
%TCIMACRO{\U{211a} }%
%BeginExpansion
\mathbb{Q}
%EndExpansion
;}
\end{equation*}

\noindent furthermore $%
%TCIMACRO{\U{211a} }%
%BeginExpansion
\mathbb{Q}
%EndExpansion
_{p}^{\times }\Psi $ induces a $%
%TCIMACRO{\U{211a} }%
%BeginExpansion
\mathbb{Q}
%EndExpansion
_{p}$-homogeneous polarization $\lambda _{\mathbf{X}}:(\mathbf{X},i_{\mathbf{%
X}})\mathbf{\rightarrow (}\widehat{\mathbf{X}},i_{\widehat{\mathbf{X}}}$) of 
$\mathbf{X}$ that respect the action of $B$ (this polarization needs not to
be principal). We have therefore associated to $\mathcal{D}_{\func{mod}}$
the triple $(\mathbf{X},i_{\mathbf{X}},\overline{\lambda }_{\mathbf{X}})$
consisting on a $%
%TCIMACRO{\U{211a} }%
%BeginExpansion
\mathbb{Q}
%EndExpansion
_{p}$-homogeneously polarized $p$-divisible group over $k$ endowed with an
action of $B$. Such a triple is unique only up to quasi-isogenies, and it is
assumed fixed in the rest of this paragraph.

Denote by $\breve{E}$ the complete unramified extension of the local Shimura
field $E$ of $\mathcal{D}_{\func{mod}}$, which has residue field $k$ and is
contained inside $K$, the field of definition of the quasi-character $\mu $;
we have $\breve{E}=EK_{0}.$ Let $\mathcal{O}_{\breve{E}}$ be the ring of
integers of $\breve{E}$ and denote by $NILP_{\mathcal{O}_{\breve{E}}}$ the
category of locally noetherian $\limfunc{Spec}\mathcal{O}_{\breve{E}}$%
-schemes $(S,\mathcal{O}_{S})$ such that the ideal sheaf $p\mathcal{O}_{S}$
is locally nilpotent; for such an $S$, we denote by $\overline{S}$ the
closed subscheme of $S$ defined by $p\mathcal{O}_{S}$: it is a scheme over $%
\limfunc{Spec}k$. For example, we could take $S=\limfunc{Spec}k=\overline{S}$%
.

The moduli problem for $p$-divisible groups associated to the above data, as
defined in \cite{RZ}, 3.21, is the following:

\begin{definition}
\label{Shi p-div}Fix a local PEL-datum $\mathcal{D}_{\func{mod}}$ for moduli
of $p$-divisible groups over $k$. Denote by $E$ its local Shimura field and
let $(\mathbf{X},i_{\mathbf{X}},\overline{\lambda }_{\mathbf{X}})$ be a
choice of $%
%TCIMACRO{\U{211a} }%
%BeginExpansion
\mathbb{Q}
%EndExpansion
_{p}$-homogeneously polarized $p$-divisible group over $k$ with $B$-action
associated via $M_{\ast }$ to the triple $(N,\mathbf{F},%
%TCIMACRO{\U{211a} }%
%BeginExpansion
\mathbb{Q}
%EndExpansion
_{p}^{\times }\Psi )$ defined by $\mathcal{D}_{\func{mod}}$.

We let $\mathcal{\breve{M}}$ be the contravariant functor from $NILP_{%
\mathcal{O}_{\breve{E}}}$ to $SETS$ defined as follows: if $S$ is a scheme
in $NILP_{\mathcal{O}_{\breve{E}}}$, $\mathcal{\breve{M}}(S)$ consists of
the equivalence classes of tuples $(X,i,\overline{\lambda };\rho )$ where:

\begin{enumerate}
\item $X$ is a $p$-divisible group over $S$;

\item $i:\mathcal{O}_{B}\rightarrow \limfunc{End}X$ is a $\mathcal{%
%TCIMACRO{\U{2124} }%
%BeginExpansion
\mathbb{Z}
%EndExpansion
}_{p}$-algebra homomorphism satisfying the determinant condition, i.e. we
require an equality of polynomial functions: $\det_{\mathcal{O}_{S}}(a,%
\limfunc{Lie}X_{S})=\det_{K}(a,V_{0})$ for all $a\in \mathcal{O}_{B}$;

\item $\lambda :(X,i)\rightarrow (\widehat{X},\widehat{i}$) is a principal
polarization of $(X,i)$, and $\overline{\lambda }$ is the corresponding $%
%TCIMACRO{\U{211a} }%
%BeginExpansion
\mathbb{Q}
%EndExpansion
_{p}$-homogeneous (principal)\ polarization;

\item $\rho :(\mathbf{X},i_{\mathbf{X}})_{\overline{S}}\rightarrow (X,i)_{%
\overline{S}}$ is a quasi isogeny of $p$-divisible groups over $\overline{S}$
that respects the $\mathcal{O}_{B}$-structure and such that $\widehat{\rho }%
\circ \lambda _{\overline{S}}\circ \rho \in 
%TCIMACRO{\U{211a} }%
%BeginExpansion
\mathbb{Q}
%EndExpansion
_{p}^{\times }(\lambda _{\mathbf{X}})_{\overline{S}}$.
\end{enumerate}

\noindent \noindent Two tuples $(X,i,\overline{\lambda };\rho )$,$(X^{\prime
},i^{\prime },\overline{\lambda }^{\prime };\rho ^{\prime })\in \mathcal{%
\breve{M}}(S)$ are said to be equivalent if the $\overline{S}$-quasi-isogeny 
$\rho ^{\prime }\circ \rho ^{-1}$ lifts to an isomorphism $%
f:(X,i)\rightarrow (X^{\prime },i^{\prime })$ of $p$-divisible groups over $%
S $ with $\mathcal{O}_{B}$-action, such that $\widehat{f}\circ \lambda
^{\prime }\circ f\in 
%TCIMACRO{\U{2124} }%
%BeginExpansion
\mathbb{Z}
%EndExpansion
_{p}^{\times }\lambda $.
\end{definition}

Notice that the definition of the functor $\mathcal{\breve{M}}$ depends - up
to isomorphism - only upon the isocrystal $(N,\mathbf{F},%
%TCIMACRO{\U{211a} }%
%BeginExpansion
\mathbb{Q}
%EndExpansion
_{p}^{\times }\Psi )$, and not upon the choice of $p$-divisible group $(%
\mathbf{X},i_{\mathbf{X}},\overline{\lambda }_{\mathbf{X}})$ that is
associated to $(N,\mathbf{F},%
%TCIMACRO{\U{211a} }%
%BeginExpansion
\mathbb{Q}
%EndExpansion
_{p}^{\times }\Psi )$ via covariant Dieudonn\'{e} theory.

It was proven by Rapoport ant Zink that the above functor is representable
by a formal scheme: this is the content of Theorem 3.25 of \cite{RZ}:

\begin{theorem}
The functor $\mathcal{\breve{M}}$ defined above is representable by a formal
scheme (still denoted by $\mathcal{\breve{M}}$) which is formally locally of
finite type over $\limfunc{Spf}\mathcal{O}_{_{\breve{E}}}$. Furthermore, if
the original local PEL-datum $\mathcal{D}_{\func{mod}}$ has good reduction,
we have $\breve{E}=K_{0}$ and the formal scheme representing $\mathcal{%
\breve{M}}$ is formally smooth over $\limfunc{Spf}W.$
\end{theorem}

\paragraph{Determinant condition \label{det condition abs}}

The determinant condition required in the above definition is due to
Kottwitz. A precise formulation, a bit hidden behind the one we gave, is
contained in \cite{Kot92}, 5, and also \cite{RZ}, 3.23: let $\mathcal{V}$ be
the $%
%TCIMACRO{\U{2124} }%
%BeginExpansion
\mathbb{Z}
%EndExpansion
_{p}$-scheme whose values in any $%
%TCIMACRO{\U{2124} }%
%BeginExpansion
\mathbb{Z}
%EndExpansion
_{p}$-algebra $R$ are given by $\mathcal{O}_{B}\otimes _{%
%TCIMACRO{\U{2124} }%
%BeginExpansion
\mathbb{Z}
%EndExpansion
_{p}}R$; fix an $\mathcal{O}_{B}$-invariant $\mathcal{O}_{K}$-lattice $%
\Gamma \subset V_{0}$ and define a map of $\mathcal{O}_{K}$-schemes $%
\det\nolimits_{K}(\cdot ,V_{0}):\mathcal{V}_{\mathcal{O}_{K}}\longrightarrow 
\mathbb{A}_{\mathcal{O}_{K}}^{1}$ \noindent by setting, for any $\mathcal{O}%
_{K}$-algebra $R$ and any $a\in \mathcal{O}_{B}\otimes _{%
%TCIMACRO{\U{2124} }%
%BeginExpansion
\mathbb{Z}
%EndExpansion
_{p}}R$, $\det_{K}(a,V_{0}):=\det (a;\Gamma \otimes _{\mathcal{O}_{K}}R).$
It is clear that this map does not depend upon the choice of lattice $\Gamma 
$ and is defined over $\mathcal{O}_{E}$ ($E$ is the local Shimura field of
the datum), so that it defines a morphism of $S$-schemes $%
\det\nolimits_{K}(\cdot ,V_{0}):\mathcal{V}_{S}\longrightarrow \mathbb{A}%
_{S}^{1}.$\noindent \noindent\ In a similar way, $\det_{\mathcal{O}%
_{S}}(\cdot ,\limfunc{Lie}X)$ can be seen as a morphism of $S$-schemes:%
\begin{equation*}
\det\nolimits_{\mathcal{O}_{S}}(\cdot ,\limfunc{Lie}X):\mathcal{V}%
_{S}\longrightarrow \mathbb{A}_{S}^{1}.
\end{equation*}

\noindent The determinant condition is the requirement that these two
morphisms of schemes over $S$ coincide. \noindent By the definition of $%
\mathcal{V}$, we could rephrase this condition by requiring that for any $S$%
-scheme $S^{\prime }$ and any $a\in \mathcal{O}_{B}\otimes \mathcal{O}%
_{S^{\prime }}$ we have an identity of polynomial functions $\det_{\mathcal{O%
}_{S^{\prime }}}(a,\limfunc{Lie}X_{S^{\prime }})=\det_{K}(a,V_{0})$; notice
that by fixing a basis for the $%
%TCIMACRO{\U{2124} }%
%BeginExpansion
\mathbb{Z}
%EndExpansion
_{p}$-free module $\mathcal{O}_{B}$ - say, of rank $t$ - then $%
\det_{K}(a,V_{0})$ for a variable $a\in \mathcal{O}_{B}$\ can be written as
a polynomial of the algebra $\mathcal{O}_{E}[X_{1},...,X_{t}]$.

\newpage

\section{Moduli of abelian schemes of PEL-type}

\subsection{Global PEL-data}

We recall the definition of a PEL-data over $%
%TCIMACRO{\U{211a} }%
%BeginExpansion
\mathbb{Q}
%EndExpansion
$, following \cite{Kot92}.

\subsubsection{Definition of a PEL-data over $%
%TCIMACRO{\U{211a} }%
%BeginExpansion
\mathbb{Q}
%EndExpansion
\label{def of global PEL}$}

Let $B$ be a finite dimensional semi-simple $%
%TCIMACRO{\U{211a} }%
%BeginExpansion
\mathbb{Q}
%EndExpansion
$-algebra endowed with a positive involution $^{\ast }$ (positivity means
that $Tr_{B/%
%TCIMACRO{\U{211a} }%
%BeginExpansion
\mathbb{Q}
%EndExpansion
}\left( xx^{\ast }\right) >0$ for any $x$ in $B-\{0\}$). Let $V\neq \{0\}$
be a finitely generated left $B$-module and $\left\langle ,\right\rangle
:V\times V\rightarrow 
%TCIMACRO{\U{211a} }%
%BeginExpansion
\mathbb{Q}
%EndExpansion
$ a non-degenerate, alternating $%
%TCIMACRO{\U{211a} }%
%BeginExpansion
\mathbb{Q}
%EndExpansion
$-bilinear form which is skew-hermitian with respect to $^{\ast }$, i.e. $%
\left\langle bv,w\right\rangle =\left\langle v,b^{\ast }w\right\rangle $ for
all $v,w$ in $V$ and all $b$ in $B$. These objects define a reductive\
algebraic group $G$ over $%
%TCIMACRO{\U{211a} }%
%BeginExpansion
\mathbb{Q}
%EndExpansion
$ whose $R$-points, for a fixed $%
%TCIMACRO{\U{211a} }%
%BeginExpansion
\mathbb{Q}
%EndExpansion
$-algebra $R$, are given by:%
\begin{equation*}
G(R)=\left\{ g\in GL_{B\otimes _{%
%TCIMACRO{\U{211a} }%
%BeginExpansion
\mathbb{Q}
%EndExpansion
}R}(V\otimes _{%
%TCIMACRO{\U{211a} }%
%BeginExpansion
\mathbb{Q}
%EndExpansion
}R):\left\langle gv,gw\right\rangle =c(g)\left\langle v,w\right\rangle \text{
}\forall v,w\in V;c(g)\in R^{\times }\right\} \text{.}
\end{equation*}

\noindent The map $G(R)\rightarrow R^{\times }$ given by $g\mapsto c(g)$
defines a homomorphism of $%
%TCIMACRO{\U{211a} }%
%BeginExpansion
\mathbb{Q}
%EndExpansion
$-algebraic groups $c:G\rightarrow \mathbb{G}_{m}$ (the similitude character
of $G$); the kernel of $c$ is a reductive $%
%TCIMACRO{\U{211a} }%
%BeginExpansion
\mathbb{Q}
%EndExpansion
$-subgroup of $G$ denoted by $G_{1}$.

\noindent Notice that $^{\ast }$ defines on the $%
%TCIMACRO{\U{211a} }%
%BeginExpansion
\mathbb{Q}
%EndExpansion
$-algebra $C:=\limfunc{End}_{B}V$ an involution $x\longmapsto x^{\ast }$ via
the identity $\left\langle xv,w\right\rangle =\left\langle v,x^{\ast
}w\right\rangle \ \ $($v,w\in V$)$.$

\noindent We therefore have functorial isomorphisms:

\begin{equation*}
G(R)\simeq \left\{ x\in C\otimes _{%
%TCIMACRO{\U{211a} }%
%BeginExpansion
\mathbb{Q}
%EndExpansion
}R:xx^{\ast }\in R^{\times }\right\} \text{,}
\end{equation*}

\noindent where $R$ is any $%
%TCIMACRO{\U{211a} }%
%BeginExpansion
\mathbb{Q}
%EndExpansion
$-algebra.

\begin{definition}
Fix a tuple $\mathcal{(}B,^{\ast },V,\left\langle ,\right\rangle )$
satisfying the above properties, and let $p$ be a fixed prime number.

\begin{enumerate}
\item The datum $\mathcal{D=(}B,^{\ast },V,\left\langle ,\right\rangle )$ is
called a $%
%TCIMACRO{\U{211a} }%
%BeginExpansion
\mathbb{Q}
%EndExpansion
$\textbf{-PEL datum}; the group $G$ is the algebraic group associated to $%
\mathcal{D}$. A $%
%TCIMACRO{\U{211a} }%
%BeginExpansion
\mathbb{Q}
%EndExpansion
$\textbf{-PEL datum with integral structure at }$p$ (or with a $p$-adic
PEL-lattice) is the datum $\mathcal{D}_{\mathcal{O}}:=(B,^{\ast
},V,\left\langle ,\right\rangle ,\mathcal{O}_{B},\Lambda )$ where $\mathcal{(%
}B,^{\ast },V,\left\langle ,\right\rangle )$ is a $%
%TCIMACRO{\U{211a} }%
%BeginExpansion
\mathbb{Q}
%EndExpansion
$-PEL datum; $\mathcal{O}_{B}$ is a $%
%TCIMACRO{\U{2124} }%
%BeginExpansion
\mathbb{Z}
%EndExpansion
_{(p)}$-order of $B$ stable under the involution $^{\ast }$ and such that $%
\mathcal{O}_{B}\mathcal{\otimes }_{%
%TCIMACRO{\U{2124} }%
%BeginExpansion
\mathbb{Z}
%EndExpansion
}%
%TCIMACRO{\U{2124} }%
%BeginExpansion
\mathbb{Z}
%EndExpansion
_{p}$ is a maximal order in $B\mathcal{\otimes }_{%
%TCIMACRO{\U{211a} }%
%BeginExpansion
\mathbb{Q}
%EndExpansion
}%
%TCIMACRO{\U{211a} }%
%BeginExpansion
\mathbb{Q}
%EndExpansion
_{p}$; $\Lambda \subset V\mathcal{\otimes }_{%
%TCIMACRO{\U{211a} }%
%BeginExpansion
\mathbb{Q}
%EndExpansion
}%
%TCIMACRO{\U{211a} }%
%BeginExpansion
\mathbb{Q}
%EndExpansion
_{p}$ is a $%
%TCIMACRO{\U{2124} }%
%BeginExpansion
\mathbb{Z}
%EndExpansion
_{p}$-lattice and an $\mathcal{O}_{B}$-submodule such that the restriction
of $\left\langle ,\right\rangle \otimes _{%
%TCIMACRO{\U{211a} }%
%BeginExpansion
\mathbb{Q}
%EndExpansion
}%
%TCIMACRO{\U{211a} }%
%BeginExpansion
\mathbb{Q}
%EndExpansion
_{p}$ to $\Lambda \times \Lambda $ gives a perfect pairing of $%
%TCIMACRO{\U{2124} }%
%BeginExpansion
\mathbb{Z}
%EndExpansion
_{p}$-modules.

\item A $%
%TCIMACRO{\U{211a} }%
%BeginExpansion
\mathbb{Q}
%EndExpansion
$\textbf{-PEL datum for moduli of abelian schemes (at }$p$\textbf{)} is the
datum $\mathcal{D}_{\func{mod}}:=(B,^{\ast },V,\left\langle ,\right\rangle ,%
\mathcal{O}_{B},\Lambda ,h,K^{p},\nu )$ where $(B,^{\ast },V,\left\langle
,\right\rangle ,\mathcal{O}_{B},\Lambda )$ is a PEL-datum with integral
structure at $p$; $K^{p}\subset G(\mathbb{A}_{f}^{p})$ is an open compact
subgroup of $G(\mathbb{A}_{f}^{p})$; $\nu :\overline{%
%TCIMACRO{\U{211a} }%
%BeginExpansion
\mathbb{Q}
%EndExpansion
}\hookrightarrow \overline{%
%TCIMACRO{\U{211a} }%
%BeginExpansion
\mathbb{Q}
%EndExpansion
}_{p}$ is an embedding of fields; $h:%
%TCIMACRO{\U{2102} }%
%BeginExpansion
\mathbb{C}
%EndExpansion
\rightarrow \limfunc{End}_{B}V\otimes _{%
%TCIMACRO{\U{211a} }%
%BeginExpansion
\mathbb{Q}
%EndExpansion
}%
%TCIMACRO{\U{211d} }%
%BeginExpansion
\mathbb{R}
%EndExpansion
$ is an $%
%TCIMACRO{\U{211d} }%
%BeginExpansion
\mathbb{R}
%EndExpansion
$-algebra homomorphism such that:

\begin{enumerate}
\item $h\left( \overline{z}\right) =h\left( z\right) ^{\ast }$ for all $z$
in $%
%TCIMACRO{\U{2102} }%
%BeginExpansion
\mathbb{C}
%EndExpansion
$ (i.e. $h$ is a $^{\ast }$-homomorphism);

\item the symmetric $%
%TCIMACRO{\U{211d} }%
%BeginExpansion
\mathbb{R}
%EndExpansion
$-bilinear form $\left( ,\right) :V_{%
%TCIMACRO{\U{211d} }%
%BeginExpansion
\mathbb{R}
%EndExpansion
}\times V_{%
%TCIMACRO{\U{211d} }%
%BeginExpansion
\mathbb{R}
%EndExpansion
}\rightarrow 
%TCIMACRO{\U{211d} }%
%BeginExpansion
\mathbb{R}
%EndExpansion
$ defined by $\left( v,w\right) :=\left\langle v,h(\sqrt{-1})w\right\rangle $
is positive definite.

\noindent \noindent Such a map $h$ is called a \textbf{polarization} for the
PEL-datum $\mathcal{(}B,^{\ast },V,\left\langle ,\right\rangle )$.
\end{enumerate}
\end{enumerate}
\end{definition}

\noindent Each of the above PEL-data over $%
%TCIMACRO{\U{211a} }%
%BeginExpansion
\mathbb{Q}
%EndExpansion
$ will be called simple if $B$ is a simple $%
%TCIMACRO{\U{211a} }%
%BeginExpansion
\mathbb{Q}
%EndExpansion
$-algebra. The PEL-data $\mathcal{D}_{\mathcal{O}}$ or $\mathcal{D}_{\func{%
mod}}$ will be said to have \textit{good reduction at }$p$ if the algebra $%
B\otimes _{%
%TCIMACRO{\U{211a} }%
%BeginExpansion
\mathbb{Q}
%EndExpansion
}%
%TCIMACRO{\U{211a} }%
%BeginExpansion
\mathbb{Q}
%EndExpansion
_{p}$ is unramified and, in case $\limfunc{End}_{B}V\otimes _{%
%TCIMACRO{\U{211a} }%
%BeginExpansion
\mathbb{Q}
%EndExpansion
}%
%TCIMACRO{\U{211d} }%
%BeginExpansion
\mathbb{R}
%EndExpansion
$ has a factor isomorphic to $M_{n}(\mathbb{H})$ for some $n>0$, then $p$ is
odd (here $\mathbb{H}$ denotes the division algebra of real quaternions).
(Cf. \cite{Wed99}, 1.4).

Let $\mathcal{D}$ $\mathcal{=(}B,^{\ast },V,\left\langle ,\right\rangle )$
be a PEL-datum endowed with a polarization $h:%
%TCIMACRO{\U{2102} }%
%BeginExpansion
\mathbb{C}
%EndExpansion
\rightarrow \limfunc{End}_{B}V\otimes _{%
%TCIMACRO{\U{211a} }%
%BeginExpansion
\mathbb{Q}
%EndExpansion
}%
%TCIMACRO{\U{211d} }%
%BeginExpansion
\mathbb{R}
%EndExpansion
$; denote by $G$ the associated algebraic group. Denote tensoring (or
extension of scalars)\ with a subscript. \noindent Since $h$ is a $^{\ast }$%
-homomorphism, $h(z)\in G(%
%TCIMACRO{\U{211d} }%
%BeginExpansion
\mathbb{R}
%EndExpansion
)$ for any $z\in 
%TCIMACRO{\U{2102} }%
%BeginExpansion
\mathbb{C}
%EndExpansion
^{\times }$. Define the map:%
\begin{eqnarray*}
\text{ \ \ \ \ }h_{%
%TCIMACRO{\U{2102} }%
%BeginExpansion
\mathbb{C}
%EndExpansion
} &:&%
%TCIMACRO{\U{2102} }%
%BeginExpansion
\mathbb{C}
%EndExpansion
\times 
%TCIMACRO{\U{2102} }%
%BeginExpansion
\mathbb{C}
%EndExpansion
\longrightarrow \limfunc{End}\nolimits_{B_{%
%TCIMACRO{\U{211d} }%
%BeginExpansion
\mathbb{R}
%EndExpansion
}}V_{%
%TCIMACRO{\U{211d} }%
%BeginExpansion
\mathbb{R}
%EndExpansion
}\otimes _{%
%TCIMACRO{\U{211d} }%
%BeginExpansion
\mathbb{R}
%EndExpansion
}%
%TCIMACRO{\U{2102} }%
%BeginExpansion
\mathbb{C}
%EndExpansion
\\
(z_{1},z_{2}) &\mapsto &\mathbf{1}\otimes \frac{z_{1}+z_{2}}{2}+h(\sqrt{-1}%
)\otimes \frac{z_{1}-z_{2}}{2\sqrt{-1}}.
\end{eqnarray*}

\noindent If we view $%
%TCIMACRO{\U{2102} }%
%BeginExpansion
\mathbb{C}
%EndExpansion
\times 
%TCIMACRO{\U{2102} }%
%BeginExpansion
\mathbb{C}
%EndExpansion
$ as a $%
%TCIMACRO{\U{2102} }%
%BeginExpansion
\mathbb{C}
%EndExpansion
$-algebra via the diagonal embedding, $h_{%
%TCIMACRO{\U{2102} }%
%BeginExpansion
\mathbb{C}
%EndExpansion
}$ is a $%
%TCIMACRO{\U{2102} }%
%BeginExpansion
\mathbb{C}
%EndExpansion
$-algebra map such that $h_{%
%TCIMACRO{\U{2102} }%
%BeginExpansion
\mathbb{C}
%EndExpansion
}(z_{1},z_{2})\in G(%
%TCIMACRO{\U{2102} }%
%BeginExpansion
\mathbb{C}
%EndExpansion
)$ if $z_{1},z_{2}\in 
%TCIMACRO{\U{2102} }%
%BeginExpansion
\mathbb{C}
%EndExpansion
^{\times }$. Let $\mathbb{S}:=\limfunc{Res}\nolimits_{%
%TCIMACRO{\U{2102} }%
%BeginExpansion
\mathbb{C}
%EndExpansion
/%
%TCIMACRO{\U{211d} }%
%BeginExpansion
\mathbb{R}
%EndExpansion
}(\mathbb{G}_{m/%
%TCIMACRO{\U{2102} }%
%BeginExpansion
\mathbb{C}
%EndExpansion
})$ with the usual identifications $\mathbb{S}(%
%TCIMACRO{\U{211d} }%
%BeginExpansion
\mathbb{R}
%EndExpansion
)=\mathbb{%
%TCIMACRO{\U{2102} }%
%BeginExpansion
\mathbb{C}
%EndExpansion
}^{\times }$ and $\mathbb{S}(%
%TCIMACRO{\U{2102} }%
%BeginExpansion
\mathbb{C}
%EndExpansion
)=(\mathbb{%
%TCIMACRO{\U{2102} }%
%BeginExpansion
\mathbb{C}
%EndExpansion
\otimes }_{%
%TCIMACRO{\U{211d} }%
%BeginExpansion
\mathbb{R}
%EndExpansion
}\mathbb{%
%TCIMACRO{\U{2102} }%
%BeginExpansion
\mathbb{C}
%EndExpansion
)}^{\times }\simeq $ $%
%TCIMACRO{\U{2102} }%
%BeginExpansion
\mathbb{C}
%EndExpansion
^{\times }\times 
%TCIMACRO{\U{2102} }%
%BeginExpansion
\mathbb{C}
%EndExpansion
^{\times }$, where $\mathbb{%
%TCIMACRO{\U{2102} }%
%BeginExpansion
\mathbb{C}
%EndExpansion
}$ is embedded into $\mathbb{%
%TCIMACRO{\U{2102} }%
%BeginExpansion
\mathbb{C}
%EndExpansion
\otimes }_{%
%TCIMACRO{\U{211d} }%
%BeginExpansion
\mathbb{R}
%EndExpansion
}\mathbb{%
%TCIMACRO{\U{2102} }%
%BeginExpansion
\mathbb{C}
%EndExpansion
}$ by $z\mapsto z\otimes \overline{z}$; by the above considerations, a
morphism of real algebraic groups $\mathfrak{h:}$ $\mathbb{S}\rightarrow G_{%
%TCIMACRO{\U{211d} }%
%BeginExpansion
\mathbb{R}
%EndExpansion
}$ \noindent remains defined by setting:%
\begin{equation*}
\mathfrak{h}(%
%TCIMACRO{\U{2102} }%
%BeginExpansion
\mathbb{C}
%EndExpansion
):=h_{%
%TCIMACRO{\U{2102} }%
%BeginExpansion
\mathbb{C}
%EndExpansion
|%
%TCIMACRO{\U{2102} }%
%BeginExpansion
\mathbb{C}
%EndExpansion
^{\times }\times 
%TCIMACRO{\U{2102} }%
%BeginExpansion
\mathbb{C}
%EndExpansion
^{\times }};\text{ \ }\mathfrak{h}(%
%TCIMACRO{\U{211d} }%
%BeginExpansion
\mathbb{R}
%EndExpansion
):=h_{|%
%TCIMACRO{\U{2102} }%
%BeginExpansion
\mathbb{C}
%EndExpansion
^{\times }},
\end{equation*}

\noindent for if $z\in 
%TCIMACRO{\U{2102} }%
%BeginExpansion
\mathbb{C}
%EndExpansion
^{\times }$ we have $\mathfrak{h}(%
%TCIMACRO{\U{2102} }%
%BeginExpansion
\mathbb{C}
%EndExpansion
)(z,\overline{z})=(\func{Re}z\cdot \mathbf{1}+\func{Im}z\cdot h(\sqrt{-1}%
))\otimes 1$ that we identify with $\mathfrak{h}(%
%TCIMACRO{\U{211d} }%
%BeginExpansion
\mathbb{R}
%EndExpansion
)(z)$. \noindent We denote the above morphism of algebraic groups by $h$ and
no confusion should arise.

\noindent By definition, our original map $h:%
%TCIMACRO{\U{2102} }%
%BeginExpansion
\mathbb{C}
%EndExpansion
\rightarrow \limfunc{End}_{B}V\otimes _{%
%TCIMACRO{\U{211a} }%
%BeginExpansion
\mathbb{Q}
%EndExpansion
}%
%TCIMACRO{\U{211d} }%
%BeginExpansion
\mathbb{R}
%EndExpansion
$ endows $V_{%
%TCIMACRO{\U{211d} }%
%BeginExpansion
\mathbb{R}
%EndExpansion
}$ with a complex structure, hence it defined a Hodge structure of type $%
(0,-1),(-1,0)$ on the vector space $V_{%
%TCIMACRO{\U{2102} }%
%BeginExpansion
\mathbb{C}
%EndExpansion
}=V^{(0,-1)}\oplus V^{(-1,0)}$, where:%
\begin{equation*}
V^{(0,-1)}=\{v\in V_{%
%TCIMACRO{\U{2102} }%
%BeginExpansion
\mathbb{C}
%EndExpansion
}:h(%
%TCIMACRO{\U{2102} }%
%BeginExpansion
\mathbb{C}
%EndExpansion
)(z_{1},z_{2})v=z_{2}v\text{, }z_{1},z_{2}\in 
%TCIMACRO{\U{2102} }%
%BeginExpansion
\mathbb{C}
%EndExpansion
^{\times }\},\text{ \ }V^{(-1,0)}=\overline{V^{(0,-1)}}.\text{\ }
\end{equation*}

\noindent Observe that $%
%TCIMACRO{\U{2102} }%
%BeginExpansion
\mathbb{C}
%EndExpansion
\subset \mathbb{%
%TCIMACRO{\U{2102} }%
%BeginExpansion
\mathbb{C}
%EndExpansion
\otimes }_{%
%TCIMACRO{\U{211d} }%
%BeginExpansion
\mathbb{R}
%EndExpansion
}\mathbb{%
%TCIMACRO{\U{2102} }%
%BeginExpansion
\mathbb{C}
%EndExpansion
}$ acts on $V^{(0,-1)}$ via the character $z\mapsto \overline{z}$, while it
acts on $V^{(-1,0)}$ via the character $z\mapsto z$.

Now let $\mu :\mathbb{G}_{m/%
%TCIMACRO{\U{2102} }%
%BeginExpansion
\mathbb{C}
%EndExpansion
}\rightarrow G_{%
%TCIMACRO{\U{2102} }%
%BeginExpansion
\mathbb{C}
%EndExpansion
}$ be the map induced by the assignment $z\mapsto h(%
%TCIMACRO{\U{2102} }%
%BeginExpansion
\mathbb{C}
%EndExpansion
)(z,1)$ for $z\in 
%TCIMACRO{\U{2102} }%
%BeginExpansion
\mathbb{C}
%EndExpansion
^{\times }$; $V^{(0,-1)}$ is the subspace of $V_{%
%TCIMACRO{\U{2102} }%
%BeginExpansion
\mathbb{C}
%EndExpansion
}$ on which $\mathbb{G}_{m/%
%TCIMACRO{\U{2102} }%
%BeginExpansion
\mathbb{C}
%EndExpansion
}$ acts - via $\mu $ - through the trivial character (weight zero), and $%
V^{(-1,0)}$ is the subspace of $V_{%
%TCIMACRO{\U{2102} }%
%BeginExpansion
\mathbb{C}
%EndExpansion
}$ on which $\mathbb{G}_{m/%
%TCIMACRO{\U{2102} }%
%BeginExpansion
\mathbb{C}
%EndExpansion
}$ acts through the identity character (weight one), so that we will\ write
the decomposition $V_{%
%TCIMACRO{\U{2102} }%
%BeginExpansion
\mathbb{C}
%EndExpansion
}=V^{(0,-1)}\oplus V^{(-1,0)}$ as $V_{%
%TCIMACRO{\U{2102} }%
%BeginExpansion
\mathbb{C}
%EndExpansion
}=V_{%
%TCIMACRO{\U{2102} }%
%BeginExpansion
\mathbb{C}
%EndExpansion
,0}\oplus V_{%
%TCIMACRO{\U{2102} }%
%BeginExpansion
\mathbb{C}
%EndExpansion
,1}.$

\noindent Both $V_{%
%TCIMACRO{\U{2102} }%
%BeginExpansion
\mathbb{C}
%EndExpansion
,0}$ and $V_{%
%TCIMACRO{\U{2102} }%
%BeginExpansion
\mathbb{C}
%EndExpansion
,1}$ have an action of $B_{%
%TCIMACRO{\U{2102} }%
%BeginExpansion
\mathbb{C}
%EndExpansion
}$, so that we obtain a semisimple complex representation of the $%
%TCIMACRO{\U{211a} }%
%BeginExpansion
\mathbb{Q}
%EndExpansion
$-algebras $B$:%
\begin{equation*}
\rho :B\rightarrow \limfunc{End}\nolimits_{%
%TCIMACRO{\U{2102} }%
%BeginExpansion
\mathbb{C}
%EndExpansion
}V_{%
%TCIMACRO{\U{2102} }%
%BeginExpansion
\mathbb{C}
%EndExpansion
,0}.
\end{equation*}

The \textit{reflex field} (or \textit{global Shimura field}) of the
PEL-datum $\mathcal{D=(}B,^{\ast },V,\left\langle ,\right\rangle )$ endowed
with a polarization $h$ is the field of definition $E=E(\mathcal{D},h)$ of
the isomorphism class of the complex representation $\rho :B\rightarrow 
\limfunc{End}_{%
%TCIMACRO{\U{2102} }%
%BeginExpansion
\mathbb{C}
%EndExpansion
}V_{%
%TCIMACRO{\U{2102} }%
%BeginExpansion
\mathbb{C}
%EndExpansion
,0}$. If $\nu :\overline{%
%TCIMACRO{\U{211a} }%
%BeginExpansion
\mathbb{Q}
%EndExpansion
}\hookrightarrow \overline{%
%TCIMACRO{\U{211a} }%
%BeginExpansion
\mathbb{Q}
%EndExpansion
}_{p}$ is an embedding of fields, the $\nu $-adic completion of $E$ will be
denote by $E_{\nu }$ and called the $\nu $\textit{-adic Shimura field}
associated to our set of data.

\noindent Fix embeddings $\overline{%
%TCIMACRO{\U{211a} }%
%BeginExpansion
\mathbb{Q}
%EndExpansion
}\hookrightarrow 
%TCIMACRO{\U{2102} }%
%BeginExpansion
\mathbb{C}
%EndExpansion
$ and $\nu :\overline{%
%TCIMACRO{\U{211a} }%
%BeginExpansion
\mathbb{Q}
%EndExpansion
}\hookrightarrow \overline{%
%TCIMACRO{\U{211a} }%
%BeginExpansion
\mathbb{Q}
%EndExpansion
}_{p}$; keeping the above assumptions, we have that $E=%
%TCIMACRO{\U{211a} }%
%BeginExpansion
\mathbb{Q}
%EndExpansion
(Tr(\rho (b)):b\in B)$ is a number field; we can characterize $E$ as the
field of definition of the $G^{0}(%
%TCIMACRO{\U{2102} }%
%BeginExpansion
\mathbb{C}
%EndExpansion
)$-conjugacy class of the co-character $\mu :\mathbb{G}_{m/%
%TCIMACRO{\U{2102} }%
%BeginExpansion
\mathbb{C}
%EndExpansion
}\rightarrow G_{%
%TCIMACRO{\U{2102} }%
%BeginExpansion
\mathbb{C}
%EndExpansion
}$. On the other side, the field of definition of the $G^{0}(%
%TCIMACRO{\U{211a} }%
%BeginExpansion
\mathbb{Q}
%EndExpansion
_{p})$-conjugacy class of $\mu $ coincide with $E_{\nu }$ and is a finite
extension of $%
%TCIMACRO{\U{211a} }%
%BeginExpansion
\mathbb{Q}
%EndExpansion
_{p}$.

Notice that the decomposition $V_{%
%TCIMACRO{\U{2102} }%
%BeginExpansion
\mathbb{C}
%EndExpansion
}=V_{%
%TCIMACRO{\U{2102} }%
%BeginExpansion
\mathbb{C}
%EndExpansion
,0}\oplus V_{%
%TCIMACRO{\U{2102} }%
%BeginExpansion
\mathbb{C}
%EndExpansion
,1}$ induced by the polarization $h$ is defined over a finite extension $%
K^{\prime }$of $%
%TCIMACRO{\U{211a} }%
%BeginExpansion
\mathbb{Q}
%EndExpansion
_{p}$, so that we can write $V_{K^{\prime }}=V_{K^{\prime },0}\oplus
V_{K^{\prime },1}$, \noindent \noindent where $V_{K^{\prime },0}$ is the $%
K^{\prime }$-subspace of $V_{K^{\prime }}$ on which $\mathbb{G}_{m/K^{\prime
}}$ acts trivially, and similarly for $V_{K^{\prime },1}$.

We have in conclusion:

\begin{lemma}
\label{p-div PEL datum}Let $\mathcal{D}:\mathcal{=(}B,^{\ast
},V,\left\langle ,\right\rangle )$ be a $%
%TCIMACRO{\U{211a} }%
%BeginExpansion
\mathbb{Q}
%EndExpansion
$\textit{-PEL datum} with associated group $G$. Define $B_{p}:=B\otimes _{%
%TCIMACRO{\U{211a} }%
%BeginExpansion
\mathbb{Q}
%EndExpansion
}%
%TCIMACRO{\U{211a} }%
%BeginExpansion
\mathbb{Q}
%EndExpansion
_{p}$, $V_{p}:=V\otimes _{%
%TCIMACRO{\U{211a} }%
%BeginExpansion
\mathbb{Q}
%EndExpansion
}%
%TCIMACRO{\U{211a} }%
%BeginExpansion
\mathbb{Q}
%EndExpansion
_{p}$, $\left\langle ,\right\rangle _{p}:=\left\langle ,\right\rangle
\otimes _{%
%TCIMACRO{\U{211a} }%
%BeginExpansion
\mathbb{Q}
%EndExpansion
}%
%TCIMACRO{\U{211a} }%
%BeginExpansion
\mathbb{Q}
%EndExpansion
_{p}$, $G_{p}=G_{%
%TCIMACRO{\U{211a} }%
%BeginExpansion
\mathbb{Q}
%EndExpansion
_{p}}$.

\begin{enumerate}
\item The tuple $\mathcal{D}_{p}:\mathcal{=(}B_{p},^{\ast
},V_{p},\left\langle ,\right\rangle _{p})$ is a $%
%TCIMACRO{\U{211a} }%
%BeginExpansion
\mathbb{Q}
%EndExpansion
_{p}$-PEL datum with associated group $G_{p}$. If $\mathcal{D}_{\mathcal{O}%
}:=(B,^{\ast },V,\left\langle ,\right\rangle ,\mathcal{O}_{B},\Lambda )$ is
a $%
%TCIMACRO{\U{211a} }%
%BeginExpansion
\mathbb{Q}
%EndExpansion
$-PEL datum with integral structure at $p$, and $\mathcal{O}_{B_{p}}:=%
\mathcal{O}_{B}\otimes _{%
%TCIMACRO{\U{2124} }%
%BeginExpansion
\mathbb{Z}
%EndExpansion
}%
%TCIMACRO{\U{2124} }%
%BeginExpansion
\mathbb{Z}
%EndExpansion
_{p}$, then $\mathcal{D}_{\mathcal{O}_{p}}:=(B_{p},^{\ast
},V_{p},\left\langle ,\right\rangle _{p},\mathcal{O}_{B_{p}},\Lambda )$ is a 
$%
%TCIMACRO{\U{211a} }%
%BeginExpansion
\mathbb{Q}
%EndExpansion
_{p}$-PEL datum with integral structure.

\item If $\mathcal{D}_{\func{mod}}:=(B,^{\ast },V,\left\langle
,\right\rangle ,\mathcal{O}_{B},\Lambda ,h,K^{p},\nu )$ is a $%
%TCIMACRO{\U{211a} }%
%BeginExpansion
\mathbb{Q}
%EndExpansion
$-PEL datum for moduli of abelian schemes at $p$, then the $p$-adic field of
definition of the $G^{0}(%
%TCIMACRO{\U{211a} }%
%BeginExpansion
\mathbb{Q}
%EndExpansion
_{p})$-conjugacy class of the associated cocharacter $\mu :\mathbb{G}_{m/%
%TCIMACRO{\U{2102} }%
%BeginExpansion
\mathbb{C}
%EndExpansion
}\rightarrow G_{%
%TCIMACRO{\U{2102} }%
%BeginExpansion
\mathbb{C}
%EndExpansion
}$ coincide with the $\nu $-adic completion $E_{\nu }$ of the Shimura field $%
E$ of $\mathcal{D}_{\func{mod}}$; furthermore there is a finite field
extension $K^{\prime }$ of $%
%TCIMACRO{\U{211a} }%
%BeginExpansion
\mathbb{Q}
%EndExpansion
_{p}$ such that the weight decomposition $V_{%
%TCIMACRO{\U{2102} }%
%BeginExpansion
\mathbb{C}
%EndExpansion
}=V_{%
%TCIMACRO{\U{2102} }%
%BeginExpansion
\mathbb{C}
%EndExpansion
,0}\oplus V_{%
%TCIMACRO{\U{2102} }%
%BeginExpansion
\mathbb{C}
%EndExpansion
,1}$ under the action of $\mu $ is defined over $K^{\prime }$.

If any of the above global PEL-data is simple, then the corresponding local
PEL-datum is simple; if $\mathcal{D}_{\mathcal{O}}$ has good reduction at $p$%
, then $\mathcal{D}_{\mathcal{O}_{p}}$ has good reduction.
\end{enumerate}
\end{lemma}

Let $\mathcal{D=(}B,^{\ast },V,\left\langle ,\right\rangle )$ be a global
PEL-datum with polarization $h$ e cocharacter $\mu :\mathbb{G}_{m/%
%TCIMACRO{\U{2102} }%
%BeginExpansion
\mathbb{C}
%EndExpansion
}\rightarrow G_{%
%TCIMACRO{\U{2102} }%
%BeginExpansion
\mathbb{C}
%EndExpansion
}$ as before. If $K^{\prime }$ is a finite field extension of $%
%TCIMACRO{\U{211a} }%
%BeginExpansion
\mathbb{Q}
%EndExpansion
_{p}$ such that $V_{K^{\prime }}=V_{K^{\prime },0}\oplus V_{K^{\prime },1}$
is the weight decomposition of $V_{K^{\prime }}$ under $\mu $, then we can
define a morphism of schemes over $\mathcal{O}_{K^{\prime }}$ as in \ref{det
condition abs}:%
\begin{equation*}
\det\nolimits_{K^{\prime }}(\cdot ,V_{K^{\prime },0}):\mathcal{V}_{\mathcal{O%
}_{K^{\prime }}}\longrightarrow \mathbb{A}_{\mathcal{O}_{K^{\prime }}}^{1}.
\end{equation*}

\noindent This morphism is actually defined over the ring $\mathcal{O}%
_{E}\otimes _{%
%TCIMACRO{\U{2124} }%
%BeginExpansion
\mathbb{Z}
%EndExpansion
}%
%TCIMACRO{\U{2124} }%
%BeginExpansion
\mathbb{Z}
%EndExpansion
_{(p)}$, embedded in $\mathcal{O}_{E_{\nu }}$ via $\nu $.

\subsubsection{Data of type $A$ and $C$\label{Examples}}

Fix a $%
%TCIMACRO{\U{211a} }%
%BeginExpansion
\mathbb{Q}
%EndExpansion
$-PEL datum $\mathcal{D}_{\func{mod}}=(B,^{\ast },V,\left\langle
,\right\rangle ,\mathcal{O}_{B},\Lambda ,h,K^{p},\nu )$ for moduli of
abelian schemes (at $p$); let $G$ be the corresponding algebraic group and $%
E $ the Shimura field. Assume that $\mathcal{D}_{\func{mod}}$ is simple, so
that $B\simeq M_{N}(D)$ for some skew field $D$ of characteristic zero. Let $%
F$ be the center of $B$, and $F_{0}$ be the subfield of $F$ fixed by the
involution $^{\ast }$; $F$ is a number field endowed with a positive
involution, so that $F_{0}$ is totally real and either $F=F_{0}$ (and $%
^{\ast }$ is called an involution of the \textit{first kind}), or $F/F_{0}$
is a quadratic totally complex extension (and $^{\ast }$ is said to be of
the \textit{second kind}). Then there is a reductive algebraic group $G_{0}$%
\ over $F_{0}$ such that $G_{1}=\limfunc{Res}\nolimits_{F_{0}/%
%TCIMACRO{\U{211a} }%
%BeginExpansion
\mathbb{Q}
%EndExpansion
}G_{0}$, namely $G_{0}(R)=\{x\in \limfunc{End}\nolimits_{B}V\otimes
_{F_{0}}R:xx^{\ast }=1\}$ for any $F_{0}$-algebra $R$. Set $%
m:=[F:F_{0}]\cdot \left( \dim _{F}(\limfunc{End}\nolimits_{B}V)\right)
^{1/2}.$

\noindent Since $\limfunc{End}\nolimits_{B}V\otimes _{%
%TCIMACRO{\U{211a} }%
%BeginExpansion
\mathbb{Q}
%EndExpansion
}%
%TCIMACRO{\U{211d} }%
%BeginExpansion
\mathbb{R}
%EndExpansion
$ is endowed with a complex structure given by $h$, we have $m=2n$ for some
integer $n$. We have three possible situations:

\begin{description}
\item[(A)] if $^{\ast }$ is of the second kind, then $G_{0}$ is an inner
form of the quasi-split unitary group over $F_{0}$ associated to the
quadratic imaginary extension $F/F_{0}$. Over an algebraically closure of $%
F_{0}$, the group $G_{0}$ is of type $A_{n-1}$; in this case $\limfunc{End}%
\nolimits_{B}V\otimes _{%
%TCIMACRO{\U{211a} }%
%BeginExpansion
\mathbb{Q}
%EndExpansion
}%
%TCIMACRO{\U{211d} }%
%BeginExpansion
\mathbb{R}
%EndExpansion
\simeq M_{n}(%
%TCIMACRO{\U{2102} }%
%BeginExpansion
\mathbb{C}
%EndExpansion
)^{[F_{0}:%
%TCIMACRO{\U{211a} }%
%BeginExpansion
\mathbb{Q}
%EndExpansion
]}$;

\item[(C)] if $^{\ast }$ is of the first kind and $\limfunc{End}%
\nolimits_{B}V\otimes _{%
%TCIMACRO{\U{211a} }%
%BeginExpansion
\mathbb{Q}
%EndExpansion
}%
%TCIMACRO{\U{211d} }%
%BeginExpansion
\mathbb{R}
%EndExpansion
\simeq M_{2n}(%
%TCIMACRO{\U{211d} }%
%BeginExpansion
\mathbb{R}
%EndExpansion
)^{[F_{0}:%
%TCIMACRO{\U{211a} }%
%BeginExpansion
\mathbb{Q}
%EndExpansion
]}$, then over an algebraically closure of $F_{0}$, the group $G_{0}$ is a
symplectic group in $2n$ variables: it is of type $C_{n}$;

\item[(D)] if $^{\ast }$ is of the first kind and $\limfunc{End}%
\nolimits_{B}V\otimes _{%
%TCIMACRO{\U{211a} }%
%BeginExpansion
\mathbb{Q}
%EndExpansion
}%
%TCIMACRO{\U{211d} }%
%BeginExpansion
\mathbb{R}
%EndExpansion
\simeq M_{n}(\mathbb{H})^{[F_{0}:%
%TCIMACRO{\U{211a} }%
%BeginExpansion
\mathbb{Q}
%EndExpansion
]}$, then over an algebraically closure of $F_{0}$, the group $G_{0}$ is an
orthogonal group in $2n$ variables: it is of type $D_{n}$.
\end{description}

\ \noindent As remarked in \cite{Kot92}, 7, if the PEL-datum falls into
cases $A$ or $C$, the associated reductive group $G$ is connected, while in
case $D$ it has $2^{[F_{0}:%
%TCIMACRO{\U{211a} }%
%BeginExpansion
\mathbb{Q}
%EndExpansion
]}>1$ connected components. Furthermore, in case $A$ and $C$, the derived
subgroup $G^{\prime }$ is simply-connected. If $G$ is of type $A$ with $n$
even (notation as above) or of type $C$, then $G$ satisfies the Hasse
principle, while in case $D$, the group $G$ does not satisfy the Hasse
principle. Notice finally that, in case $D$, in order to guarantee good
reduction of our PEL-datum, we need to exclude the prime $p=2$. For these
reasons we won't consider PEL-data of type $D$.

\bigskip

It is enough for us to consider simple PEL-data as above, in which $B$ is a
division $%
%TCIMACRO{\U{211a} }%
%BeginExpansion
\mathbb{Q}
%EndExpansion
$-algebra, endowed with the positive involution $^{\ast }$. Define $%
d^{2}:=[B:F]$, $e:=[F:%
%TCIMACRO{\U{211a} }%
%BeginExpansion
\mathbb{Q}
%EndExpansion
]$, $e_{0}:=[F_{0}:%
%TCIMACRO{\U{211a} }%
%BeginExpansion
\mathbb{Q}
%EndExpansion
]$, where as above $F=Z(B)$ and $F_{0}=F^{\ast =id}.$

By Albert's classification of division algebras with positive involutions,
we only have four possibilities for $(B,^{\ast })$ (cf. \cite{Mum}, 21, Th.
2), that reduce, if we exclude the case of PEL-data of type $D$, to the
following three:

\begin{description}
\item[(C-I)] $B=F=F_{0}$ is a totally real number field, with the trivial
involution $^{\ast }=id_{B}$;

\item[(C-II)] $F=F_{0}$ is a totally real number field, and $B$ is a
quaternion division algebra over $F$ such that $B\otimes _{F,\iota }%
%TCIMACRO{\U{211d} }%
%BeginExpansion
\mathbb{R}
%EndExpansion
\simeq M_{2}(%
%TCIMACRO{\U{211d} }%
%BeginExpansion
\mathbb{R}
%EndExpansion
)$ for any real embedding $\iota :F\hookrightarrow 
%TCIMACRO{\U{211d} }%
%BeginExpansion
\mathbb{R}
%EndExpansion
$; the involution $^{\ast }$ on such a $B$ is given by conjugating the
natural involution $x\mapsto Tr_{B/F}(x)-x$ of $B$ by some element $a\in
B^{\times }$ such that $a^{2}\in F$ is totally negative in $F$; any such map
is a positive involution on $B$. In this case we can choose an isomorphism $%
B\otimes _{%
%TCIMACRO{\U{211a} }%
%BeginExpansion
\mathbb{Q}
%EndExpansion
}%
%TCIMACRO{\U{211d} }%
%BeginExpansion
\mathbb{R}
%EndExpansion
\simeq M_{2}(%
%TCIMACRO{\U{211d} }%
%BeginExpansion
\mathbb{R}
%EndExpansion
)^{e_{0}}$ carrying the involution $^{\ast }$ into the involution $%
(X_{1},...,X_{e_{0}})\mapsto (X_{1}^{t},...,X_{e_{0}}^{t})$;

\item[(A)] $F$ is a totally imaginary quadratic extension of the totally
real field $F_{0}$, with complex conjugation $c$ and, for any finite place $%
v $ of $F$ we have $inv_{v}B=0$ if $v=cv$, and $inv_{v}B+inv_{cv}B=0$
otherwise; there is also a positive involution $^{\prime }$ on $B$ and an
isomorphism $B\otimes _{%
%TCIMACRO{\U{211a} }%
%BeginExpansion
\mathbb{Q}
%EndExpansion
}%
%TCIMACRO{\U{211d} }%
%BeginExpansion
\mathbb{R}
%EndExpansion
\simeq M_{d}(%
%TCIMACRO{\U{2102} }%
%BeginExpansion
\mathbb{C}
%EndExpansion
)^{e_{0}}$ which carries $^{\prime }$ into the involution $%
(X_{1},...,X_{e_{0}})\mapsto (\overline{X}_{1}^{t},...,\overline{X}%
_{e_{0}}^{t})$. Given such an involution $^{\prime }$, any other positive
involution is obtained by conjugating it through an element $a\in D$ such
that $a=a^{\prime }$ and such that the image of $a$ via the above
isomorphism is of the form $(A_{1},...,A_{e_{0}})$ where each matrix $A_{i}$
is Hermitian and positive definite.
\end{description}

\bigskip

\noindent Let $2g=\dim _{%
%TCIMACRO{\U{211a} }%
%BeginExpansion
\mathbb{Q}
%EndExpansion
}V$ and define, for each of the above cases, the integer $r$ as follows (cf. 
\cite{Yu}, 2):%
\begin{equation*}
r:=\left\{ 
\begin{array}{cc}
\frac{g}{e_{0}} & \text{(Type C-I)} \\ 
\frac{g}{2e_{0}} & \text{(Type C-II)} \\ 
\frac{g}{d^{2}e_{0}} & \text{(Type A).}%
\end{array}%
\right.
\end{equation*}

\noindent Being $B$ a simple algebra over $%
%TCIMACRO{\U{211a} }%
%BeginExpansion
\mathbb{Q}
%EndExpansion
$, $B_{%
%TCIMACRO{\U{2102} }%
%BeginExpansion
\mathbb{C}
%EndExpansion
}:=B\otimes _{%
%TCIMACRO{\U{211a} }%
%BeginExpansion
\mathbb{Q}
%EndExpansion
}%
%TCIMACRO{\U{2102} }%
%BeginExpansion
\mathbb{C}
%EndExpansion
$ is a semisimple algebra over $%
%TCIMACRO{\U{2102} }%
%BeginExpansion
\mathbb{C}
%EndExpansion
$ and $V_{%
%TCIMACRO{\U{2102} }%
%BeginExpansion
\mathbb{C}
%EndExpansion
,0}$ becomes a finite dimensional $B_{%
%TCIMACRO{\U{2102} }%
%BeginExpansion
\mathbb{C}
%EndExpansion
}$-module, so that it decomposes as the direct sum of\ irreducible $B_{%
%TCIMACRO{\U{2102} }%
%BeginExpansion
\mathbb{C}
%EndExpansion
}$-modules (the irreducible $B_{%
%TCIMACRO{\U{2102} }%
%BeginExpansion
\mathbb{C}
%EndExpansion
}$-modules are in one to one correspondence with the irreducible modules of
each simple constituent of $B_{%
%TCIMACRO{\U{2102} }%
%BeginExpansion
\mathbb{C}
%EndExpansion
}$).

Using this fact, one sees (cf. \cite{Shi}, 5) that in cases C-I and C-II
(and also in case D), the representation $\rho :B\rightarrow \limfunc{End}%
\nolimits_{%
%TCIMACRO{\U{2102} }%
%BeginExpansion
\mathbb{C}
%EndExpansion
}V_{%
%TCIMACRO{\U{2102} }%
%BeginExpansion
\mathbb{C}
%EndExpansion
,0}$ is uniquely determined and has to be a multiple of a reduced $%
%TCIMACRO{\U{2102} }%
%BeginExpansion
\mathbb{C}
%EndExpansion
$-representation of $B$ over $%
%TCIMACRO{\U{211a} }%
%BeginExpansion
\mathbb{Q}
%EndExpansion
$, so that the reflex field is $E=%
%TCIMACRO{\U{211a} }%
%BeginExpansion
\mathbb{Q}
%EndExpansion
$. In case A, the situation is less rigid: we know that $F$ is a quadratic
imaginary extension of the totally real field $F_{0}$; let $\left\{ \tau
_{i}\right\} _{i=1}^{e_{0}}$ be a CM-type of $F$, and let $n_{i}$ (resp. $%
\overline{n}_{i}$) be the multiplicity of the standard representation of $%
B\otimes _{F,\tau _{i}}%
%TCIMACRO{\U{2102} }%
%BeginExpansion
\mathbb{C}
%EndExpansion
$ (resp. $B\otimes _{F,c\tau _{i}}%
%TCIMACRO{\U{2102} }%
%BeginExpansion
\mathbb{C}
%EndExpansion
$) in $V_{%
%TCIMACRO{\U{2102} }%
%BeginExpansion
\mathbb{C}
%EndExpansion
,0}$, for $i=1,...,e_{0}$ (notice that the standard representations of $%
B\otimes _{F,\tau _{i}}%
%TCIMACRO{\U{2102} }%
%BeginExpansion
\mathbb{C}
%EndExpansion
$ and of $B\otimes _{F,c\tau _{i}}%
%TCIMACRO{\U{2102} }%
%BeginExpansion
\mathbb{C}
%EndExpansion
$ have dimension $q$ over $%
%TCIMACRO{\U{2102} }%
%BeginExpansion
\mathbb{C}
%EndExpansion
$). Any fixed collection of pairs $(n_{i},\overline{n}_{i})$ such that $%
n_{i}+\overline{n}_{i}=dr$ ($1\leq i\leq e_{0}$)\ gives rise in the obvious
way to a unique finite dimensional complex representation of $B$; $\rho $
has to be isomorphic to a representation of this form.

\bigskip

We now present two cases that will be interesting for us: the first is of
type C-I, the third is of type A. We will always assume fixed a choice of
square root of $-1$ in $%
%TCIMACRO{\U{2102} }%
%BeginExpansion
\mathbb{C}
%EndExpansion
$, denoted by $\sqrt{-1}$.

\paragraph{PEL-datum of type C-I}

Let $B=F$ be a totally real, finite, Galois extension of $%
%TCIMACRO{\U{211a} }%
%BeginExpansion
\mathbb{Q}
%EndExpansion
$ of degree $f$. Denote by $\left\{ \tau _{1},...,\tau _{f}\right\} $ the
distinct embeddings of $F$ into $%
%TCIMACRO{\U{211d} }%
%BeginExpansion
\mathbb{R}
%EndExpansion
$ and endow $F$ with the identity involution. Assume that the fixed prime $p$
is inert in $F/%
%TCIMACRO{\U{211a} }%
%BeginExpansion
\mathbb{Q}
%EndExpansion
$, so that $F\otimes _{%
%TCIMACRO{\U{211a} }%
%BeginExpansion
\mathbb{Q}
%EndExpansion
}%
%TCIMACRO{\U{211a} }%
%BeginExpansion
\mathbb{Q}
%EndExpansion
_{p}$ is the unramified extension of $%
%TCIMACRO{\U{211a} }%
%BeginExpansion
\mathbb{Q}
%EndExpansion
_{p}$ of degree $f$ in a fixed algebraic closure of $%
%TCIMACRO{\U{211a} }%
%BeginExpansion
\mathbb{Q}
%EndExpansion
_{p}$. Let $V=F^{2g}$ for a fixed integer $g>0$ and denote by $\left\langle
,\right\rangle $ the map $V\times V\rightarrow 
%TCIMACRO{\U{211a} }%
%BeginExpansion
\mathbb{Q}
%EndExpansion
$ defined by setting $\left\langle v,w\right\rangle :=Tr_{F/%
%TCIMACRO{\U{211a} }%
%BeginExpansion
\mathbb{Q}
%EndExpansion
}\left( v^{t}J_{2g}w\right) $ \noindent for all $v,w\in V$.

Identifying $C=\limfunc{End}_{F}V$ with the matrix algebra $M_{2g}(F)$, the
induced positive involution on $C$ is given by $A\mapsto A^{\ast
}:=J_{2g}^{-1}A^{t}J_{2g}$. We identify the $%
%TCIMACRO{\U{211d} }%
%BeginExpansion
\mathbb{R}
%EndExpansion
$-algebras $F\otimes _{%
%TCIMACRO{\U{211a} }%
%BeginExpansion
\mathbb{Q}
%EndExpansion
}%
%TCIMACRO{\U{211d} }%
%BeginExpansion
\mathbb{R}
%EndExpansion
$ and $%
%TCIMACRO{\U{211d} }%
%BeginExpansion
\mathbb{R}
%EndExpansion
^{f}$ by the map $x\otimes r\mapsto (\tau _{1}(x)r,...,\tau _{f}(x)r)$; an
identification $C_{%
%TCIMACRO{\U{211d} }%
%BeginExpansion
\mathbb{R}
%EndExpansion
}=M_{2g}(F)\otimes _{%
%TCIMACRO{\U{211a} }%
%BeginExpansion
\mathbb{Q}
%EndExpansion
}%
%TCIMACRO{\U{211d} }%
%BeginExpansion
\mathbb{R}
%EndExpansion
\simeq M_{2g}(%
%TCIMACRO{\U{211d} }%
%BeginExpansion
\mathbb{R}
%EndExpansion
)^{\oplus f}$ remains therefore defined via the embeddings $\tau _{i}.$
Notice that the involution $^{\ast }$ acts componentwise on $M_{2g}(%
%TCIMACRO{\U{211d} }%
%BeginExpansion
\mathbb{R}
%EndExpansion
)^{\oplus f}$, i.e. $(X_{1},...,X_{f})^{\ast
}=(J_{2g}^{-1}X_{1}^{t}J_{2g},...,J_{2g}^{-1}X_{f}^{t}J_{2g}).$

\noindent Let $h:%
%TCIMACRO{\U{2102} }%
%BeginExpansion
\mathbb{C}
%EndExpansion
\rightarrow M_{2g}(%
%TCIMACRO{\U{211d} }%
%BeginExpansion
\mathbb{R}
%EndExpansion
)^{\oplus f}$ be the $%
%TCIMACRO{\U{211d} }%
%BeginExpansion
\mathbb{R}
%EndExpansion
$-algebra map defined by the assignment $a+b\sqrt{-1}\longmapsto
(a+J_{2g}^{t}b)^{\oplus f}.$The algebraic $%
%TCIMACRO{\U{211a} }%
%BeginExpansion
\mathbb{Q}
%EndExpansion
$-group $G$ associated to the above data is isomorphic to the reductive
connected group $GSp_{2g}(F)_{/%
%TCIMACRO{\U{211a} }%
%BeginExpansion
\mathbb{Q}
%EndExpansion
}$ of symplectic similitudes of $F$: if $R$ is a $%
%TCIMACRO{\U{211a} }%
%BeginExpansion
\mathbb{Q}
%EndExpansion
$-algebra 
\begin{equation*}
GSp_{2g}\left( F\right) (R)=\left\{ A\in GL_{2g}(F\otimes _{%
%TCIMACRO{\U{211a} }%
%BeginExpansion
\mathbb{Q}
%EndExpansion
}R):A^{t}J_{2g}A=c(A)J_{2g},c(A)\in R^{\times }\right\} \text{.}
\end{equation*}

\noindent Furthermore, $G_{1}=Sp_{2g}(F)$ and $G_{0}=Sp_{2g}(F)_{/F}$ is a
group of type $C_{g}$ when viewed over $\overline{F}$;\ $G$ satisfies the
Hasse principle.

By making the identification $V_{%
%TCIMACRO{\U{2102} }%
%BeginExpansion
\mathbb{C}
%EndExpansion
}\simeq (%
%TCIMACRO{\U{211d} }%
%BeginExpansion
\mathbb{R}
%EndExpansion
^{2g})^{\oplus f}\otimes _{%
%TCIMACRO{\U{211d} }%
%BeginExpansion
\mathbb{R}
%EndExpansion
}%
%TCIMACRO{\U{2102} }%
%BeginExpansion
\mathbb{C}
%EndExpansion
$, we write, for any $z_{1},z_{2}\in 
%TCIMACRO{\U{2102} }%
%BeginExpansion
\mathbb{C}
%EndExpansion
^{\times }$:%
\begin{equation*}
h(%
%TCIMACRO{\U{2102} }%
%BeginExpansion
\mathbb{C}
%EndExpansion
)(z_{1},z_{2})=(1,...,1)\otimes \frac{z_{1}+z_{2}}{2}%
+(J_{2g}^{t},...,J_{2g}^{t})\otimes \frac{z_{1}-z_{2}}{2\sqrt{-1}}\in M_{2g}(%
%TCIMACRO{\U{211d} }%
%BeginExpansion
\mathbb{R}
%EndExpansion
)^{\oplus f}\otimes _{%
%TCIMACRO{\U{211d} }%
%BeginExpansion
\mathbb{R}
%EndExpansion
}%
%TCIMACRO{\U{2102} }%
%BeginExpansion
\mathbb{C}
%EndExpansion
.
\end{equation*}

\noindent If $\{e_{1},...,e_{g},f_{1},...,f_{g}\}$ is the standard ordered
basis of $%
%TCIMACRO{\U{211d} }%
%BeginExpansion
\mathbb{R}
%EndExpansion
^{2g}$, denote by $\left\{ e_{i}B_{h},f_{j}B_{h}:1\leq i,j\leq g,1\leq h\leq
f\right\} $ the corresponding canonical ordered basis of $(%
%TCIMACRO{\U{211d} }%
%BeginExpansion
\mathbb{R}
%EndExpansion
^{2g})^{\oplus f}$, so that for example $e_{i}B_{h}$ is the vector $%
(0,...,e_{i},...,0)\in $ $(%
%TCIMACRO{\U{211d} }%
%BeginExpansion
\mathbb{R}
%EndExpansion
^{2g})^{\oplus f}$ where $e_{i}$ appears in position $h$. We have:%
\begin{eqnarray*}
V_{%
%TCIMACRO{\U{2102} }%
%BeginExpansion
\mathbb{C}
%EndExpansion
,0} &=&\left\langle e_{i}B_{h}\otimes 1+f_{i}B_{h}\otimes \sqrt{-1}:1\leq
i\leq g,1\leq h\leq f\right\rangle , \\
V_{%
%TCIMACRO{\U{2102} }%
%BeginExpansion
\mathbb{C}
%EndExpansion
,1} &=&\left\langle e_{i}B_{h}\otimes 1-f_{i}B_{h}\otimes \sqrt{-1}:1\leq
i\leq g,1\leq h\leq f\right\rangle .
\end{eqnarray*}

If $x\in F$, then $\rho (x)$ acts on the vector $e_{i}B_{h}\otimes
1+f_{i}B_{h}\otimes \sqrt{-1}\in V_{%
%TCIMACRO{\U{2102} }%
%BeginExpansion
\mathbb{C}
%EndExpansion
,0}$ as multiplication by $(\tau _{1}(x),...,\tau _{f}(x))\otimes 1\in F_{%
%TCIMACRO{\U{211d} }%
%BeginExpansion
\mathbb{R}
%EndExpansion
}\otimes _{%
%TCIMACRO{\U{211d} }%
%BeginExpansion
\mathbb{R}
%EndExpansion
}%
%TCIMACRO{\U{2102} }%
%BeginExpansion
\mathbb{C}
%EndExpansion
=%
%TCIMACRO{\U{211d} }%
%BeginExpansion
\mathbb{R}
%EndExpansion
^{f}\otimes _{%
%TCIMACRO{\U{211d} }%
%BeginExpansion
\mathbb{R}
%EndExpansion
}%
%TCIMACRO{\U{2102} }%
%BeginExpansion
\mathbb{C}
%EndExpansion
$, so that:%
\begin{equation*}
\rho (x)\left( e_{i}B_{h}\otimes 1+f_{i}B_{h}\otimes \sqrt{-1}\right) =\tau
_{h}(x)\left( e_{i}B_{h}\otimes 1+f_{i}B_{h}\otimes \sqrt{-1}\right) ,
\end{equation*}

\noindent for all $1\leq i\leq g$ and $1\leq h\leq f$. Referring this linear
transformation the ordered basis of $V_{%
%TCIMACRO{\U{2102} }%
%BeginExpansion
\mathbb{C}
%EndExpansion
,0}$ given by:

\begin{eqnarray*}
&&e_{1}B_{1}\otimes 1+f_{1}B_{1}\otimes \sqrt{-1},...,e_{1}B_{f}\otimes
1+f_{1}B_{f}\otimes \sqrt{-1}, \\
&&e_{2}B_{1}\otimes 1+f_{2}B_{1}\otimes \sqrt{-1},...,e_{2}B_{f}\otimes
1+f_{2}B_{f}\otimes \sqrt{-1},...,
\end{eqnarray*}

\noindent we obtain the matrix form of the representation $\rho
:F\rightarrow M_{gf}(%
%TCIMACRO{\U{2102} }%
%BeginExpansion
\mathbb{C}
%EndExpansion
)$:%
\begin{equation*}
\rho :x\longmapsto \limfunc{diag}(\tau _{1}(x),\tau _{2}(x),...,\tau
_{f}(x))^{\oplus g}\text{ \ \ \ \ \ (}x\in F\text{);}
\end{equation*}

\noindent $\rho $ is the $g$-fold multiple of the reduced representation of $%
F$ over $%
%TCIMACRO{\U{211a} }%
%BeginExpansion
\mathbb{Q}
%EndExpansion
$, and the reflex field for our PEL-datum is $E=%
%TCIMACRO{\U{211a} }%
%BeginExpansion
\mathbb{Q}
%EndExpansion
$. Observe that $\det \rho (x)=\left( Nm_{F/%
%TCIMACRO{\U{211a} }%
%BeginExpansion
\mathbb{Q}
%EndExpansion
}(x)\right) ^{g}$. If $\left\{ b_{1},...,b_{f}\right\} $ is a $%
%TCIMACRO{\U{2124} }%
%BeginExpansion
\mathbb{Z}
%EndExpansion
_{(p)}$-basis of the free $%
%TCIMACRO{\U{2124} }%
%BeginExpansion
\mathbb{Z}
%EndExpansion
_{(p)}$-module $\mathcal{O}_{B}:=\mathcal{O}_{F}\otimes _{%
%TCIMACRO{\U{2124} }%
%BeginExpansion
\mathbb{Z}
%EndExpansion
}%
%TCIMACRO{\U{2124} }%
%BeginExpansion
\mathbb{Z}
%EndExpansion
_{(p)}\simeq \mathcal{O}_{F,(p)}$, then the determinant polynomial function
is: 
\begin{equation*}
f(X_{1},...,X_{f}):=\prod\nolimits_{i=1}^{f}\left( \tau
_{1}(b_{1})X_{1}+...+\tau _{f}(b_{f})X_{f}\right) ^{g}=\left( Nm_{F/%
%TCIMACRO{\U{211a} }%
%BeginExpansion
\mathbb{Q}
%EndExpansion
}(b_{1}X_{1}+...+b_{f}X_{f})\right) ^{g}.
\end{equation*}

\noindent Notice that this polynomial belongs to $%
%TCIMACRO{\U{2124} }%
%BeginExpansion
\mathbb{Z}
%EndExpansion
_{(p)}[X_{1},...,X_{f}]$.\ If we set $\Lambda :=\left( \mathcal{O}%
_{B}\otimes _{%
%TCIMACRO{\U{2124} }%
%BeginExpansion
\mathbb{Z}
%EndExpansion
_{(p)}}%
%TCIMACRO{\U{2124} }%
%BeginExpansion
\mathbb{Z}
%EndExpansion
_{p}\right) ^{2g}$ and pick any compact open subgroup $K^{p}$ of $G(\widehat{%
%TCIMACRO{\U{2124} }%
%BeginExpansion
\mathbb{Z}
%EndExpansion
}^{p})$ and any choice of embedding $\nu :\overline{%
%TCIMACRO{\U{211a} }%
%BeginExpansion
\mathbb{Q}
%EndExpansion
}\hookrightarrow \overline{%
%TCIMACRO{\U{211a} }%
%BeginExpansion
\mathbb{Q}
%EndExpansion
}_{p},$\ we have all the information necessary to define a simple $%
%TCIMACRO{\U{211a} }%
%BeginExpansion
\mathbb{Q}
%EndExpansion
$-PEL datum with good reduction at $p$. We denote this datum by $\mathcal{D}%
_{2g,p}^{Sp(F)}$ and notice that $G$ is defined over $%
%TCIMACRO{\U{2124} }%
%BeginExpansion
\mathbb{Z}
%EndExpansion
$, by setting for any $%
%TCIMACRO{\U{2124} }%
%BeginExpansion
\mathbb{Z}
%EndExpansion
$-algebra $R$: $G(R)=\left\{ A\in GL_{2g}(\mathcal{O}_{F}\otimes _{%
%TCIMACRO{\U{2124} }%
%BeginExpansion
\mathbb{Z}
%EndExpansion
}R):A^{t}J_{2g}A=c(A)J_{2g},c(A)\in R^{\times }\right\} .$

\paragraph{PEL-datum of type $A$}

Let $B=k$ be a quadratic imaginary field, say $k=%
%TCIMACRO{\U{211a} }%
%BeginExpansion
\mathbb{Q}
%EndExpansion
(\sqrt{\alpha })$, where $\alpha $ is a negative square-free integer; fix an
embedding $\tau :k\hookrightarrow 
%TCIMACRO{\U{2102} }%
%BeginExpansion
\mathbb{C}
%EndExpansion
$ such that $\tau (\sqrt{\alpha })=\sqrt{-1}\sqrt{-\alpha }$; via this
embedding we make the identification $k\otimes _{%
%TCIMACRO{\U{211a} }%
%BeginExpansion
\mathbb{Q}
%EndExpansion
}%
%TCIMACRO{\U{211d} }%
%BeginExpansion
\mathbb{R}
%EndExpansion
\simeq 
%TCIMACRO{\U{2102} }%
%BeginExpansion
\mathbb{C}
%EndExpansion
$. Assume that the fixed prime $p$ is inert in the extension $k/%
%TCIMACRO{\U{211a} }%
%BeginExpansion
\mathbb{Q}
%EndExpansion
$, so that $k\otimes _{%
%TCIMACRO{\U{211a} }%
%BeginExpansion
\mathbb{Q}
%EndExpansion
}%
%TCIMACRO{\U{211a} }%
%BeginExpansion
\mathbb{Q}
%EndExpansion
_{p}=%
%TCIMACRO{\U{211a} }%
%BeginExpansion
\mathbb{Q}
%EndExpansion
_{p}(\sqrt{\alpha })$ is the quadratic unramified extension of $%
%TCIMACRO{\U{211a} }%
%BeginExpansion
\mathbb{Q}
%EndExpansion
_{p}$ in a fixed algebraic closure of $%
%TCIMACRO{\U{211a} }%
%BeginExpansion
\mathbb{Q}
%EndExpansion
_{p}$.

Let $x\mapsto \overline{x}$ ($x\in k$) denote the non-trivial field
automorphism of $k$: it is a positive involution on $k$. Set $V=k^{g}$ for a
positive even integer $g=2n$; fix two non-negative integers $r$ and $s$
whose sum is $g$ and let: 
\begin{equation*}
H:=H_{r,s}:=%
\begin{pmatrix}
-\sqrt{\alpha }I_{r} & 0_{r,s} \\ 
0_{s,r} & \sqrt{\alpha }I_{s}%
\end{pmatrix}%
.
\end{equation*}

\noindent Let us denote by $\left\langle ,\right\rangle :V\times
V\rightarrow 
%TCIMACRO{\U{211a} }%
%BeginExpansion
\mathbb{Q}
%EndExpansion
$ the map defined by setting $\left\langle v,w\right\rangle :=Tr_{k/%
%TCIMACRO{\U{211a} }%
%BeginExpansion
\mathbb{Q}
%EndExpansion
}\left( \overline{v}^{t}Hw\right) $ \noindent for every $v,w\in V$. Notice
that $\left\langle ,\right\rangle $ is a $%
%TCIMACRO{\U{211a} }%
%BeginExpansion
\mathbb{Q}
%EndExpansion
$-bilinear non-degenerate skew-Hermitian pairing. Letting $C=\limfunc{End}%
_{k}V=M_{g}(k)$, we find that the involution induced by $^{\ast }$ on $C$ is 
$A\mapsto A^{\ast }:=H^{-1}\overline{A}^{t}H,$ where $\overline{A}$ is the
matrix obtained from $A$ by applying the field automorphism $^{\ast }$ to
the entries.

\noindent Let $h:%
%TCIMACRO{\U{2102} }%
%BeginExpansion
\mathbb{C}
%EndExpansion
\rightarrow C\otimes _{%
%TCIMACRO{\U{211a} }%
%BeginExpansion
\mathbb{Q}
%EndExpansion
}%
%TCIMACRO{\U{211d} }%
%BeginExpansion
\mathbb{R}
%EndExpansion
=M_{g}(k)\otimes _{%
%TCIMACRO{\U{211a} }%
%BeginExpansion
\mathbb{Q}
%EndExpansion
}%
%TCIMACRO{\U{211d} }%
%BeginExpansion
\mathbb{R}
%EndExpansion
$ be the $%
%TCIMACRO{\U{211d} }%
%BeginExpansion
\mathbb{R}
%EndExpansion
$-algebra homomorphism defined by:%
\begin{equation*}
a+b\sqrt{-1}\longmapsto 1\otimes a-H\otimes \frac{b}{\sqrt{-\alpha }}.
\end{equation*}

The algebraic group $G$ associated to the above data is identified with the $%
%TCIMACRO{\U{211a} }%
%BeginExpansion
\mathbb{Q}
%EndExpansion
$-group $GU_{g}(k;r,s)$: for any $%
%TCIMACRO{\U{211a} }%
%BeginExpansion
\mathbb{Q}
%EndExpansion
$-algebra $R$ we have%
\begin{equation*}
GU_{g}\left( k\otimes _{%
%TCIMACRO{\U{211a} }%
%BeginExpansion
\mathbb{Q}
%EndExpansion
}R;r,s\right) =\left\{ A\in GL_{g}(k\otimes _{%
%TCIMACRO{\U{211a} }%
%BeginExpansion
\mathbb{Q}
%EndExpansion
}R):\overline{A}^{t}HA=c(A)H,c(A)\in R^{\times }\right\} .
\end{equation*}

\noindent Furthermore, $G_{1}=U_{g}(k;H)$ is an inner form of the
quasi-split unitary group over $%
%TCIMACRO{\U{211a} }%
%BeginExpansion
\mathbb{Q}
%EndExpansion
$ associated to the extension $k/%
%TCIMACRO{\U{211a} }%
%BeginExpansion
\mathbb{Q}
%EndExpansion
$, hence it is a group of type $A_{g-1}$ when viewed over $\overline{%
%TCIMACRO{\U{211a} }%
%BeginExpansion
\mathbb{Q}
%EndExpansion
}$;\ being $g$ even, $G$ satisfies the weak Hasse principle. We also notice $%
G$ is connected.

The $%
%TCIMACRO{\U{211d} }%
%BeginExpansion
\mathbb{R}
%EndExpansion
$-algebra $C_{%
%TCIMACRO{\U{211d} }%
%BeginExpansion
\mathbb{R}
%EndExpansion
}$ is isomorphic to $M_{g}(%
%TCIMACRO{\U{2102} }%
%BeginExpansion
\mathbb{C}
%EndExpansion
)$ via the fixed embedding $k\hookrightarrow 
%TCIMACRO{\U{2102} }%
%BeginExpansion
\mathbb{C}
%EndExpansion
$, and by making the identification $V_{%
%TCIMACRO{\U{2102} }%
%BeginExpansion
\mathbb{C}
%EndExpansion
}=%
%TCIMACRO{\U{2102} }%
%BeginExpansion
\mathbb{C}
%EndExpansion
^{g}\otimes _{%
%TCIMACRO{\U{211d} }%
%BeginExpansion
\mathbb{R}
%EndExpansion
}%
%TCIMACRO{\U{2102} }%
%BeginExpansion
\mathbb{C}
%EndExpansion
$ we write, for any $z_{1},z_{2}\in 
%TCIMACRO{\U{2102} }%
%BeginExpansion
\mathbb{C}
%EndExpansion
^{\times }$:

\begin{equation*}
h(%
%TCIMACRO{\U{2102} }%
%BeginExpansion
\mathbb{C}
%EndExpansion
)(z_{1},z_{2})=1\otimes \frac{z_{1}+z_{2}}{2}-\frac{H}{\sqrt{-\alpha }}%
\otimes \frac{z_{1}-z_{2}}{2\sqrt{-1}}\in G(%
%TCIMACRO{\U{2102} }%
%BeginExpansion
\mathbb{C}
%EndExpansion
).
\end{equation*}

\noindent If $\{e_{1},...,e_{g},f_{1},...,f_{g}\}$ denotes the standard
ordered basis of $%
%TCIMACRO{\U{2102} }%
%BeginExpansion
\mathbb{C}
%EndExpansion
^{2g}$, it is easy to see that:%
\begin{eqnarray*}
V_{%
%TCIMACRO{\U{2102} }%
%BeginExpansion
\mathbb{C}
%EndExpansion
,0} &=&\left\langle \sqrt{-1}e_{i}\otimes 1-e_{i}\otimes \sqrt{-1},\sqrt{-1}%
f_{j}\otimes 1+f_{j}\otimes \sqrt{-1}:1\leq i\leq r,1\leq j\leq
s\right\rangle , \\
V_{%
%TCIMACRO{\U{2102} }%
%BeginExpansion
\mathbb{C}
%EndExpansion
,1} &=&\left\langle \sqrt{-1}e_{i}\otimes 1+e_{i}\otimes \sqrt{-1},\sqrt{-1}%
f_{j}\otimes 1-f_{j}\otimes \sqrt{-1}:1\leq i\leq r,1\leq j\leq
s\right\rangle ,
\end{eqnarray*}

\noindent and the representation $\rho :B\rightarrow \limfunc{End}\nolimits_{%
%TCIMACRO{\U{2102} }%
%BeginExpansion
\mathbb{C}
%EndExpansion
}(V_{%
%TCIMACRO{\U{2102} }%
%BeginExpansion
\mathbb{C}
%EndExpansion
,0})\simeq M_{g}(%
%TCIMACRO{\U{2102} }%
%BeginExpansion
\mathbb{C}
%EndExpansion
)$ is obviously induced by the assignment:

\begin{equation*}
\sqrt{\alpha }\longmapsto 
\begin{pmatrix}
-\sqrt{\alpha }I_{r} & 0_{r,s} \\ 
0_{s,r} & \sqrt{\alpha }I_{s}%
\end{pmatrix}%
.
\end{equation*}

The reflex field for our PEL-datum is $E=%
%TCIMACRO{\U{211a} }%
%BeginExpansion
\mathbb{Q}
%EndExpansion
$ if $r=s(=n)$, and it is $E=k$ otherwise. The multiplicity of $k\otimes
_{k,\tau }%
%TCIMACRO{\U{2102} }%
%BeginExpansion
\mathbb{C}
%EndExpansion
\simeq 
%TCIMACRO{\U{2102} }%
%BeginExpansion
\mathbb{C}
%EndExpansion
$ in $\rho $ is equal to $s$ and the multiplicity of $k\otimes _{k,\overline{%
\tau }}%
%TCIMACRO{\U{2102} }%
%BeginExpansion
\mathbb{C}
%EndExpansion
$ is $r$. \noindent Moreover:%
\begin{equation*}
f(X_{1},X_{2}):=\det (X_{1}+\sqrt{\alpha }X_{2};V_{%
%TCIMACRO{\U{2102} }%
%BeginExpansion
\mathbb{C}
%EndExpansion
,0})=(X_{1}-\sqrt{\alpha }X_{2})^{r}(X_{1}+\sqrt{\alpha }X_{2})^{s}\in 
\mathcal{O}_{E}[X_{1},X_{2}]\text{.}
\end{equation*}

Finally, set $\mathcal{O}_{B}:=\mathcal{O}_{k}\otimes _{%
%TCIMACRO{\U{2124} }%
%BeginExpansion
\mathbb{Z}
%EndExpansion
}%
%TCIMACRO{\U{2124} }%
%BeginExpansion
\mathbb{Z}
%EndExpansion
_{(p)}\simeq \mathcal{O}_{k,(p)}$, $\Lambda :=\left( \mathcal{O}_{B}\otimes
_{%
%TCIMACRO{\U{2124} }%
%BeginExpansion
\mathbb{Z}
%EndExpansion
_{(p)}}%
%TCIMACRO{\U{2124} }%
%BeginExpansion
\mathbb{Z}
%EndExpansion
_{p}\right) ^{2g}$ (notice that $p$ does not divide $\alpha $ in $%
%TCIMACRO{\U{2124} }%
%BeginExpansion
\mathbb{Z}
%EndExpansion
$, so that our pairing $\left\langle ,\right\rangle $ restricts to a perfect
pairing on $\Lambda ^{2}$). Then for any compact open subgroup $K^{p}$ of $G(%
\widehat{%
%TCIMACRO{\U{2124} }%
%BeginExpansion
\mathbb{Z}
%EndExpansion
}^{p})$ and for any choice of embedding $\nu :\overline{%
%TCIMACRO{\U{211a} }%
%BeginExpansion
\mathbb{Q}
%EndExpansion
}\hookrightarrow \overline{%
%TCIMACRO{\U{211a} }%
%BeginExpansion
\mathbb{Q}
%EndExpansion
}_{p}$\ we have all the information necessary to define a simple $%
%TCIMACRO{\U{211a} }%
%BeginExpansion
\mathbb{Q}
%EndExpansion
$-PEL datum with good reduction at $p$. Denote this datum by $\mathcal{D}%
_{(r,s);p}^{U}$ and notice that $G$ is defined over $%
%TCIMACRO{\U{2124} }%
%BeginExpansion
\mathbb{Z}
%EndExpansion
$, by setting for any $%
%TCIMACRO{\U{2124} }%
%BeginExpansion
\mathbb{Z}
%EndExpansion
$-algebra $R$: $G(R)=\left\{ A\in GL_{g}(\mathcal{O}_{k}\otimes _{%
%TCIMACRO{\U{2124} }%
%BeginExpansion
\mathbb{Z}
%EndExpansion
}R):\overline{A}^{t}HA=c(A)H,c(A)\in R^{\times }\right\} $.

\subsection{The moduli functor for abelian schemes}

\subsubsection{Abelian schemes up to prime-to-$p$ isogenies}

We assume fixed in this paragraph a global PEL-datum for moduli of abelian
schemes (at $p$) $\mathcal{D}_{\func{mod}}=(B,^{\ast },V,\left\langle
,\right\rangle ,\mathcal{O}_{B},\Lambda ,h,K^{p},\nu )$; we write $G$ for
the associated algebraic group, and $E$ for the reflex field; we furthermore
assume that our datum has good reduction at $p$. Let us fix a locally
noetherian base scheme $S$; we recall some definitions following \cite{RZ},
6.3 and \cite{Lan}, 1.3.1 and 1.3.2.

The category of \textit{abelian }$O_{B}$\textit{-schemes over }$S$\textit{\
up to isogeny of order prime to }$p$, denoted by $AV_{\mathcal{O}_{B}/S}$ or
simply by $AV$, is defined as follows: its objects are pairs $(A,i)$ where $%
A $ is an abelian scheme over $S$, and $i$ is a homomorphism of $%
%TCIMACRO{\U{2124} }%
%BeginExpansion
\mathbb{Z}
%EndExpansion
_{(p)}$-algebras $i:\mathcal{O}_{B}\mathcal{\rightarrow }\limfunc{End}%
A\otimes _{%
%TCIMACRO{\U{2124} }%
%BeginExpansion
\mathbb{Z}
%EndExpansion
}%
%TCIMACRO{\U{2124} }%
%BeginExpansion
\mathbb{Z}
%EndExpansion
_{(p)}$. A morphism $f:\left( A_{1},i_{1}\right) \rightarrow (A_{2},i_{2})$
in $AV$ is an element of the group $\limfunc{Hom}_{\mathcal{O}%
_{B}}(A_{1},A_{2})\otimes _{%
%TCIMACRO{\U{2124} }%
%BeginExpansion
\mathbb{Z}
%EndExpansion
}%
%TCIMACRO{\U{2124} }%
%BeginExpansion
\mathbb{Z}
%EndExpansion
_{(p)}$, where $\limfunc{Hom}_{\mathcal{O}_{B}}(A_{1},A_{2})$ is the module
of morphisms of abelian $S$-schemes that respect the action of $\mathcal{O}%
_{B}$.

An \textit{isogeny }$\varphi :\left( A_{1},i_{1}\right) \rightarrow
(A_{2},i_{2})$\textbf{\ }in $AV$ is a quasi-isogeny of abelian $S$-schemes $%
A_{1}\rightarrow A_{2}$ which is also a morphism of $AV$; its kernel\textbf{%
\ }is the kernel of the corresponding isogeny of $p$-divisible groups $%
A_{1}(p)\rightarrow A_{2}(p)$, so that $\ker \varphi $ is a finite locally
free group scheme whose order is locally a power of $p$, and all isogenies
have degree a power of $p$ (locally). A quasi-isogeny in $AV$ is a
quasi-isogeny of abelian schemes that respects the action of $\mathcal{O}%
_{B} $.

If $\left( A,i\right) $ is an object of $AV$, then we define an object of $%
AV $ by $\left( A,i\right) \symbol{94}:=(\widehat{A},\widehat{i})$, where $%
\widehat{A}$ is the dual abelian scheme of $A$, and $\widehat{i}:\mathcal{O}%
_{B}\mathcal{\rightarrow }\limfunc{End}\widehat{A}\otimes _{%
%TCIMACRO{\U{2124} }%
%BeginExpansion
\mathbb{Z}
%EndExpansion
}%
%TCIMACRO{\U{2124} }%
%BeginExpansion
\mathbb{Z}
%EndExpansion
_{(p)}$ is given by $\widehat{i}(b)=i\left( b^{\ast }\right) \symbol{94}$.
If $\varphi :\left( A_{1},i_{1}\right) \rightarrow (A_{2},i_{2})$ is an
isogeny in $AV$, then the dual quasi-isogeny $\widehat{\varphi }:\widehat{A}%
_{2}\rightarrow \widehat{A}_{1}$ is an isogeny in $AV$, called the dual
isogeny of $\varphi $ in $AV$.

A \textit{polarization}\textbf{\ }of $\left( A,i\right) $ in $AV_{\mathcal{O}%
_{B}/S}$ is a quasi-isogeny $\lambda :\left( A,i\right) \rightarrow (%
\widehat{A},\widehat{i})$ in $AV$ such that there exists a positive integer $%
n$ for which $n\lambda $ is induced by an ample line bundle on $A$. Such a $%
\lambda $ is called a principal polarization if furthermore $\lambda $ is an
isomorphism in $AV$. A $%
%TCIMACRO{\U{211a} }%
%BeginExpansion
\mathbb{Q}
%EndExpansion
$-homogeneous\textbf{\ }(resp. $%
%TCIMACRO{\U{2124} }%
%BeginExpansion
\mathbb{Z}
%EndExpansion
_{(p)}$-homogeneous)\ polarization $\overline{\lambda }:\left( A,i\right)
\rightarrow (\widehat{A},\widehat{i})$ is the set of (locally on $S$) $%
%TCIMACRO{\U{211a} }%
%BeginExpansion
\mathbb{Q}
%EndExpansion
^{\times }$-multiples (resp. $%
%TCIMACRO{\U{2124} }%
%BeginExpansion
\mathbb{Z}
%EndExpansion
_{(p)}^{\times }$-multiples) of a polarization $\lambda $ of $(A,i)$ in $AV$%
; such a set is called a principal $%
%TCIMACRO{\U{211a} }%
%BeginExpansion
\mathbb{Q}
%EndExpansion
$-homogeneous\textbf{\ }(resp. $%
%TCIMACRO{\U{2124} }%
%BeginExpansion
\mathbb{Z}
%EndExpansion
_{(p)}$-homogeneous)\ polarization if there is an element $\lambda \in 
\overline{\lambda }$ that is a principal polarization in $AV$.

An isogeny of polarized (resp. $%
%TCIMACRO{\U{211a} }%
%BeginExpansion
\mathbb{Q}
%EndExpansion
$-homogeneously polarized; $%
%TCIMACRO{\U{2124} }%
%BeginExpansion
\mathbb{Z}
%EndExpansion
_{(p)}$-homogeneously polarized)\ abelian varieties $\varphi
:(A_{1},i_{1};\lambda _{1})\rightarrow (A_{2},i_{2};\lambda _{2})$ in $AV$
is an isogeny $\varphi :\left( A_{1},i_{1}\right) \rightarrow (A_{2},i_{2})$
in $AV$ such that $\widehat{\varphi }\circ \lambda _{2}\circ \varphi
=\lambda _{1}$ (resp. $\widehat{\varphi }\circ \lambda _{2}\circ \varphi \in 
%TCIMACRO{\U{211a} }%
%BeginExpansion
\mathbb{Q}
%EndExpansion
^{\times }\lambda _{1}$; resp. $\widehat{\varphi }\circ \lambda _{2}\circ
\varphi \in 
%TCIMACRO{\U{2124} }%
%BeginExpansion
\mathbb{Z}
%EndExpansion
_{(p)}^{\times }\lambda _{1}$).

\bigskip

Notice that if $\lambda :\left( A,i\right) \rightarrow (\widehat{A},\widehat{%
i})$ is a polarization in $AV$, then $i\left( b^{\ast }\right) =i\left(
b\right) ^{\dag }$ for any $b$ in $\mathcal{O}_{B}$, where $^{\dag }$
denotes the Rosati involution induced by $\lambda $.

An isomorphism $\varphi :\left( A_{1},i_{1}\right) \rightarrow (A_{2},i_{2})$
in $AV$ is of the form $f\otimes r$ where $f$ is an isogeny of abelian
schemes endowed with $\mathcal{O}_{B}$-action whose degree is prime to $p$,
and $r\in 
%TCIMACRO{\U{2124} }%
%BeginExpansion
\mathbb{Z}
%EndExpansion
_{(p)}^{\times }$. We can also say that $\varphi $ is an isogeny of $AV$ of
prime-to-$p$ degree. Furthermore, if $\varphi :(A_{1},i_{1};\lambda
_{1})\rightarrow (A_{2},i_{2};\lambda _{2})$ is an isomorphism of polarized
abelian varieties in $AV$, then $i_{2}$ and $\lambda _{2}$ are determined
uniquely by $\varphi $, $i_{1}$ and $\lambda _{1}$. If $\varphi
:(A_{1},i_{1})\rightarrow (A_{2},i_{2})$ is an isomorphism\ in $AV$ and $%
\lambda _{2}:\left( A_{2},i_{2}\right) \rightarrow (\widehat{A_{2}},\widehat{%
i_{2}})$ is a principal polarization in $AV$, then $\widehat{\varphi }\circ
\lambda _{2}\circ \varphi $ is a principal polarization of $\left(
A_{1},i_{1}\right) $ and $\lambda _{2}^{-1}$ is a principal polarization of $%
(\widehat{A_{2}},\widehat{i_{2}})$.

\paragraph{Level structure}

Recall that our PEL-datum comes with an open compact subgroup $K^{p}$ of $G(%
\mathbb{A}_{f}^{p})$ ($\mathbb{A}_{f}^{p}$ denotes the ring of finite ad\`{e}%
les over $%
%TCIMACRO{\U{211a} }%
%BeginExpansion
\mathbb{Q}
%EndExpansion
$ with trivial $p$-component). Let $(A,i;\lambda )$ be a principally
polarized abelian scheme in $AV_{\mathcal{O}_{B}/S}$, where the base scheme $%
S$ is assumed to be a connected locally noetherian scheme over $\mathcal{O}%
_{E}\otimes _{%
%TCIMACRO{\U{2124} }%
%BeginExpansion
\mathbb{Z}
%EndExpansion
}%
%TCIMACRO{\U{2124} }%
%BeginExpansion
\mathbb{Z}
%EndExpansion
_{(p)}$. Let $s$ be a geometric point of $S$ and consider:%
\begin{equation*}
H_{1}(A_{s},\mathbb{A}_{f}^{p})=\left( \tprod\limits_{l\neq
p}T_{l}(A_{s})\right) \otimes _{%
%TCIMACRO{\U{2124} }%
%BeginExpansion
\mathbb{Z}
%EndExpansion
}%
%TCIMACRO{\U{211a} }%
%BeginExpansion
\mathbb{Q}
%EndExpansion
\text{,}
\end{equation*}

\noindent the Tate $\mathbb{A}_{f}^{p}$-module of the abelian variety $A_{s}$%
. It is endowed with a continuous action of $\pi _{1}(S,s)$. The action of $%
\mathcal{O}_{B}$ on $A$ endows $H_{1}(A_{s},\mathbb{A}_{f}^{p})$ with a
structure of $B$-module, and the principal polarization $\lambda $ of $(A,i)$
induces a canonical skew-symmetric $\mathbb{A}_{f}^{p}$-pairing (the Weil
pairing):%
\begin{equation*}
H_{1}(A_{s},\mathbb{A}_{f}^{p})\times H_{1}(A_{s},\mathbb{A}%
_{f}^{p})\rightarrow \mathbb{A}_{f}^{p}(1)
\end{equation*}

\noindent which is non-degenerate and skew-Hermitian with respect to $^{\ast
}$. On the other side, by definition of PEL-datum, $V_{\mathbb{A}%
_{f}^{p}}:=V\otimes _{%
%TCIMACRO{\U{211a} }%
%BeginExpansion
\mathbb{Q}
%EndExpansion
}\mathbb{A}_{f}^{p}$ is endowed with an action of $B$ and a skew-Hermitian
(with respect to $^{\ast }$) non-degenerate $\mathbb{A}_{f}^{p}$-pairing
with values in $\mathbb{A}_{f}^{p}$.

A \textit{level structure of type}\textbf{\ }$K^{p}$ on $(A,i;\lambda )$ is
the left $K^{p}$-orbit $\overline{\alpha }$ of an isomorphism $\alpha
:H_{1}(A_{s},\mathbb{A}_{f}^{p})\rightarrow V_{\mathbb{A}_{f}^{p}}$ of
skew-Hermitian $B$-modules such that $\overline{\alpha }$ is fixed by $\pi
_{1}(S,s)$. Here by isomorphisms of skew-Hermitian $B$-modules we mean an
isomorphism of $B$-modules carrying one alternating form into a $(\mathbb{A}%
_{f}^{p})^{\times }$-multiple of the other.

Assume \ now that the group $G_{/%
%TCIMACRO{\U{211a} }%
%BeginExpansion
\mathbb{Q}
%EndExpansion
}$ associated to the PEL-datum $\mathcal{D}_{\func{mod}}$ has a model $G_{/%
%TCIMACRO{\U{2124} }%
%BeginExpansion
\mathbb{Z}
%EndExpansion
}$ over $%
%TCIMACRO{\U{2124} }%
%BeginExpansion
\mathbb{Z}
%EndExpansion
$ (but it does not need to be smooth over $%
%TCIMACRO{\U{2124} }%
%BeginExpansion
\mathbb{Z}
%EndExpansion
$). Let $N\geq 1$ be an integer not divisible by the prime $p$ and let $%
(A,i;\lambda )$ be a principally polarized abelian scheme in $AV$. A \textit{%
principal level-}$N$\textit{\ structure} on $(A,i;\lambda )$ is a level
structure of type: 
\begin{equation*}
U\left( N\right) =\limfunc{Ker}\left( G(\widehat{%
%TCIMACRO{\U{2124} }%
%BeginExpansion
\mathbb{Z}
%EndExpansion
}^{p})\rightarrow G(\widehat{%
%TCIMACRO{\U{2124} }%
%BeginExpansion
\mathbb{Z}
%EndExpansion
}^{p}/N\widehat{%
%TCIMACRO{\U{2124} }%
%BeginExpansion
\mathbb{Z}
%EndExpansion
}^{p})=G(%
%TCIMACRO{\U{2124} }%
%BeginExpansion
\mathbb{Z}
%EndExpansion
/N%
%TCIMACRO{\U{2124} }%
%BeginExpansion
\mathbb{Z}
%EndExpansion
)\right) \text{.}
\end{equation*}

\noindent Notice that $U(N)=\tprod\nolimits_{l\neq p}U_{l}(N),$\noindent
where $U_{l}(N)=G\left( 
%TCIMACRO{\U{2124} }%
%BeginExpansion
\mathbb{Z}
%EndExpansion
_{l}\right) $ if $l\neq p$ and $l\nmid N$, and $U_{l}(N)=\limfunc{Ker}\left(
G\left( 
%TCIMACRO{\U{2124} }%
%BeginExpansion
\mathbb{Z}
%EndExpansion
_{l}\right) \rightarrow G(%
%TCIMACRO{\U{2124} }%
%BeginExpansion
\mathbb{Z}
%EndExpansion
_{l}/l^{n_{l}}%
%TCIMACRO{\U{2124} }%
%BeginExpansion
\mathbb{Z}
%EndExpansion
_{l}\right) )$, where $n_{l}=\limfunc{ord}_{l}N$. If $K^{p}$ is a compact
open subgroup of $G(\mathbb{A}_{f}^{p})$ contained inside $U\left( N\right) $
for some $N\geq 3$ not divisible by $p$, then "Serre's lemma" implies that $%
K^{p}$ is neat in the sense of Pink (cf. \cite{Lan} 1.4.1.9-10).

\paragraph{The determinant condition}

Let $V_{%
%TCIMACRO{\U{2102} }%
%BeginExpansion
\mathbb{C}
%EndExpansion
}=V_{%
%TCIMACRO{\U{2102} }%
%BeginExpansion
\mathbb{C}
%EndExpansion
,0}\oplus V_{%
%TCIMACRO{\U{2102} }%
%BeginExpansion
\mathbb{C}
%EndExpansion
,1}$ be the Hodge decomposition of $V_{%
%TCIMACRO{\U{2102} }%
%BeginExpansion
\mathbb{C}
%EndExpansion
}$ as in \ref{def of global PEL}, where $%
%TCIMACRO{\U{2102} }%
%BeginExpansion
\mathbb{C}
%EndExpansion
^{\times }$ act on $V_{%
%TCIMACRO{\U{2102} }%
%BeginExpansion
\mathbb{C}
%EndExpansion
,0}$ via $\mu $\ through the trivial character. We have recalled above the
morphism of schemes $\det\nolimits_{K^{\prime }}(\cdot ,V_{K^{\prime },0}):%
\mathcal{V}_{\mathcal{O}_{K^{\prime }}}\longrightarrow \mathbb{A}_{\mathcal{O%
}_{K^{\prime }}}^{1}$, that is defined over $\mathcal{O}_{E}\otimes _{%
%TCIMACRO{\U{2124} }%
%BeginExpansion
\mathbb{Z}
%EndExpansion
}%
%TCIMACRO{\U{2124} }%
%BeginExpansion
\mathbb{Z}
%EndExpansion
_{(p)}\overset{\nu }{\hookrightarrow }\mathcal{O}_{E_{\nu }}$; let us denote
this morphism by $\det\nolimits_{E}(\cdot ,V_{0})$. If $S$ is a locally
noetherian scheme over $\mathcal{O}_{E}\otimes _{%
%TCIMACRO{\U{2124} }%
%BeginExpansion
\mathbb{Z}
%EndExpansion
}%
%TCIMACRO{\U{2124} }%
%BeginExpansion
\mathbb{Z}
%EndExpansion
_{(p)}$, we can therefore define a morphism of $S$-schemes $%
\det\nolimits_{E}(\cdot ,V_{0}):\mathcal{V}_{S}\longrightarrow \mathbb{A}%
_{S}^{1}.$

\noindent Similarly (cf. Remark \ref{det condition abs}) we have, for any
object $(A,i)$ in $AV_{\mathcal{O}_{B}/S}$, a well defined morphism of $S$%
-schemes $\det\nolimits_{\mathcal{O}_{S}}(\cdot ,\limfunc{Lie}A):\mathcal{V}%
_{S}\longrightarrow \mathbb{A}_{S}^{1}.$

We say that $(A,i)$ satisfies the \textit{Kottwitz\ determinant condition}
if for any locally noetherian $S$-scheme $S^{\prime }$ we have:%
\begin{equation*}
\det\nolimits_{\mathcal{O}_{S^{\prime }}}(a,\limfunc{Lie}A_{S^{\prime
}})=\det\nolimits_{E}(a,V_{0})\text{ \ \ far all }a\in \mathcal{O}%
_{B}\otimes \mathcal{O}_{S^{\prime }}.
\end{equation*}

Let us fix a basis $\left\{ b_{1},...,b_{t}\right\} $ of the $%
%TCIMACRO{\U{2124} }%
%BeginExpansion
\mathbb{Z}
%EndExpansion
_{(p)}$-free module $\mathcal{O}_{B}$ and let $\{X_{1},...,X_{t}\}$ be
indeterminates. Set $f(X_{1},...,X_{t}):=\det (b_{1}X_{1}+...+b_{t}X_{t};V_{%
%TCIMACRO{\U{2102} }%
%BeginExpansion
\mathbb{C}
%EndExpansion
,0}).$ \noindent This is a homogeneous polynomial of degree $\dim _{%
%TCIMACRO{\U{2102} }%
%BeginExpansion
\mathbb{C}
%EndExpansion
}V_{%
%TCIMACRO{\U{2102} }%
%BeginExpansion
\mathbb{C}
%EndExpansion
,0}$ in the indeterminates $X_{1},...,X_{t}$ with coefficients in $\mathcal{O%
}_{E}\otimes _{%
%TCIMACRO{\U{2124} }%
%BeginExpansion
\mathbb{Z}
%EndExpansion
}%
%TCIMACRO{\U{2124} }%
%BeginExpansion
\mathbb{Z}
%EndExpansion
_{(p)}$. \noindent On the other side, $\limfunc{Lie}A$ is as a locally free $%
\mathcal{O}_{S}$-module with an action of $\mathcal{O}_{B}$; hence it makes
sense to consider the polynomial $g(X_{1},...,X_{t}):=\det
(b_{1}X_{1}+...+b_{t}X_{t};\limfunc{Lie}A)$, \noindent which is homogeneous
of degree $\dim _{S}A$ with coefficients in the ring of global sections of $%
\mathcal{O}_{S}$. Since we are assuming that $S$ is a scheme over $\mathcal{O%
}_{E}\otimes _{%
%TCIMACRO{\U{2124} }%
%BeginExpansion
\mathbb{Z}
%EndExpansion
}%
%TCIMACRO{\U{2124} }%
%BeginExpansion
\mathbb{Z}
%EndExpansion
_{(p)}$, the condition $f=g$ makes sense: this is equivalent to the above
defined determinant condition.

\subsubsection{The moduli problem}

Following \cite{Kot92}, one defines the following moduli problem:

\begin{definition}
Let $\mathcal{D}\mathcal{=}(B,^{\ast },V,\left\langle ,\right\rangle ,%
\mathcal{O}_{B},\Lambda ,h,K^{p},\nu )$ be a $%
%TCIMACRO{\U{211a} }%
%BeginExpansion
\mathbb{Q}
%EndExpansion
$-PEL-datum with good reduction at $p$, having Shimura field $E$ and
associated group $G$. The moduli problem $\mathbf{M:=M}(\mathcal{D})$
associated to the above data is the contravariant functor from the category $%
SCH_{\mathcal{O}_{E}\otimes 
%TCIMACRO{\U{2124} }%
%BeginExpansion
\mathbb{Z}
%EndExpansion
_{(p)}}$ of locally noetherian schemes over $\mathcal{O}_{E}\otimes _{%
%TCIMACRO{\U{2124} }%
%BeginExpansion
\mathbb{Z}
%EndExpansion
}%
%TCIMACRO{\U{2124} }%
%BeginExpansion
\mathbb{Z}
%EndExpansion
_{(p)}$ to the category of sets defined as follows: if $S$ is an object of $%
SCH_{\mathcal{O}_{E}\otimes 
%TCIMACRO{\U{2124} }%
%BeginExpansion
\mathbb{Z}
%EndExpansion
_{p}}$, then $\mathbf{M}(S)$ is the set of isomorphism classes of tuples $%
(A,i,\overline{\lambda },\overline{\alpha })$ where:

\begin{enumerate}
\item $(A,i)$ is an object in $AV$ satisfying the determinant condition;

\item $\overline{\lambda }:(A,i)\rightarrow (\widehat{A},\widehat{i})$ is a $%
%TCIMACRO{\U{211a} }%
%BeginExpansion
\mathbb{Q}
%EndExpansion
$-homogeneous principal polarization in $AV$;

\item $\overline{\alpha }$ is a level structure of type $K^{p}$ on $(A,i,%
\overline{\lambda }).$
\end{enumerate}

\noindent Here we consider two tuples $(A_{1},i_{1},\overline{\lambda }_{1},%
\overline{\alpha }_{1})$ and $(A_{2},i_{2},\overline{\lambda }_{2},\overline{%
\alpha }_{2})$ as above isomorphic if there is an isomorphism $%
f:(A_{1},i_{1};\overline{\lambda }_{1})\rightarrow (A_{2},i_{2};\overline{%
\lambda }_{2})$ of $%
%TCIMACRO{\U{211a} }%
%BeginExpansion
\mathbb{Q}
%EndExpansion
$-homogeneously principally polarized abelian schemes in $AV$ carrying $%
\overline{\alpha }_{1}$ into $\overline{\alpha }_{2}$ (in the sense that $%
\alpha _{2}\circ H_{1}(f,\mathbb{A}_{f}^{p})\circ \alpha _{1}^{-1}\in K^{p}$
and $c(\alpha _{2})\cdot c(\alpha _{1})^{-1}\in r\cdot c(K^{p})$, where $c$
denotes the similitude factor homomorphism, and $r\in 
%TCIMACRO{\U{2124} }%
%BeginExpansion
\mathbb{Z}
%EndExpansion
_{(p)}^{\times }$ is such that $r\cdot \hat{f}\circ \lambda _{2}\circ
f=\lambda _{1}$).
\end{definition}

Notice that any abelian scheme over some scheme $S$\ that satisfies the
conditions given above has relative dimension over $S$ equal to $\dim _{%
%TCIMACRO{\U{2102} }%
%BeginExpansion
\mathbb{C}
%EndExpansion
}V_{%
%TCIMACRO{\U{2102} }%
%BeginExpansion
\mathbb{C}
%EndExpansion
,0}=\frac{1}{2}\dim _{%
%TCIMACRO{\U{211a} }%
%BeginExpansion
\mathbb{Q}
%EndExpansion
}V$, in virtue of the determinant condition. Furthermore, the possible
determinant conditions that we can impose are subject to rigid constraints;
in particular they are only finitely many.

We have the following result (cf. \cite{Kot92}; \cite{Lan}, Ch. 2, 1.4.1.14,
7.2.3.10):

\begin{theorem}
\label{Representability}Let $\mathcal{D}$ be a global PEL-datum with Shimura
field $E$; assume $\mathcal{D}$ had good reduction at $p$ and assume $K^{p}$
is neat. Then the associated functor $\mathbf{M:=M}(\mathcal{D})$ is
represented by a quasi-projective smooth separated scheme $\mathcal{S}_{%
\mathcal{D},K^{p}}$ over $\mathcal{O}_{E}\otimes _{%
%TCIMACRO{\U{2124} }%
%BeginExpansion
\mathbb{Z}
%EndExpansion
}%
%TCIMACRO{\U{2124} }%
%BeginExpansion
\mathbb{Z}
%EndExpansion
_{(p)}$ which is of finite type.
\end{theorem}

We will later need the following:

\begin{theorem}
\label{lifting with separability}Let $\mathcal{D}$ be a global PEL-datum
with Shimura field $E$ and let $k$ be the algebraic closure of the residue
field of $E$ at $p$; assume $\mathcal{D}$ had good reduction at $p$ and $%
K^{p}$ is neat; let us denote by $\mathcal{S}_{\mathcal{D},K^{p}}$ the
scheme representing $\mathbf{M}(\mathcal{D})$. The canonical map $\mathcal{S}%
_{\mathcal{D},K^{p}}(W(k))\rightarrow \mathcal{S}_{\mathcal{D},K^{p}}(k)$ is
surjective.
\end{theorem}

The proof of the last result can be found in \cite{Lan}, 2.2.4.16, 2.3.2.1,
where it is shown that the deformation functor $ART(W(k))\rightarrow SETS$
associated to a fixed point $\left[ \xi \right] \in \mathcal{S}_{\mathcal{D}%
,K^{p}}(k)$ is pro-representable and formally smooth, and that one can apply
Grothendieck's Formal Existence Theorem to guarantee that we can algebraize
the formal scheme defined on $W(k)$ by a projective system of deformations
of $\left[ \xi \right] $ over the rings $W_{n}(k)$'s. The main hypothesis
necessary to prove the above theorem is that the polarization of $\xi $ is
separable.

\paragraph{Hecke action\label{Hecke}}

Let us only consider open compact subgroups $K^{p}$ of $G(\mathbb{A}%
_{f}^{p}) $ that are small enough, so that each element in $\mathbf{M}(S)$
has no non-trivial automorphisms, for any $S$ in $SCH_{\mathcal{O}%
_{E}\otimes 
%TCIMACRO{\U{2124} }%
%BeginExpansion
\mathbb{Z}
%EndExpansion
_{p}}$: for example, if $\mathcal{D}$ has a $%
%TCIMACRO{\U{2124} }%
%BeginExpansion
\mathbb{Z}
%EndExpansion
$-model, we consider only open compact subgroups contained inside $U(N)$ for
some integer $N\geq 3$. If $K_{1}^{p}\subseteq K_{2}^{p}$\ are two such open
compact subgroups of $G(\mathbb{A}_{f}^{p})$, then the transition map $%
\mathcal{S}_{\mathcal{D},K_{1}^{p}}\rightarrow \mathcal{S}_{\mathcal{D}%
,K_{2}^{p}}$ induced by $(A,i,\overline{\lambda },K_{1}^{p}\alpha )\mapsto
(A,i,\overline{\lambda },K_{2}^{p}\alpha )$\ is a finite \'{e}tale covering
which is Galois, with Galois group $K_{2}^{p}/K_{1}^{p}$, if $K_{1}^{p}$ is
normal in $K_{2}^{p}$. Denote by $\mathcal{S}_{\mathcal{D}}$ the projective
systems of the family of schemes $\{\mathcal{S}_{\mathcal{D}%
,K^{p}}\}_{K^{p}} $ where the $K^{p}$'s are small enough; we define the
following natural \textit{Hecke action} of $G(\mathbb{A}_{f}^{p})$ on $%
\mathcal{S}_{\mathcal{D}}$: if $g\in G(\mathbb{A}_{f}^{p})$, then $g$ acts
on the right on $\mathcal{S}_{\mathcal{D}}$ via the isomorphism:%
\begin{equation*}
g:\mathcal{S}_{\mathcal{D},K^{p}}\rightarrow \mathcal{S}_{\mathcal{D}%
,g^{-1}K^{p}g}
\end{equation*}

\noindent defined by $[(A,i,\overline{\lambda },\overline{\alpha })]\cdot
g:=[(A,i,\overline{\lambda },\overline{g^{-1}\circ \alpha })]$.

\subsubsection{Modular forms of PEL-type}

We define modular forms of PEL-type, as a generalization of Siegel modular
forms (cf. \cite{Gor}, 5.1). We will keep the notation of the previous
sections.

Let $\mathcal{D=(}B,^{\ast },V,\left\langle ,\right\rangle ,\mathcal{O}%
_{B},\Lambda ,h,K^{p},\nu \mathcal{)}$ be a simple PEL-datum for moduli of
abelian schemes, having good reduction at the fixed prime $p$; let $G$ be
the associated algebraic group and assume it has a model over $%
%TCIMACRO{\U{2124} }%
%BeginExpansion
\mathbb{Z}
%EndExpansion
$. Let $E$ be the reflex field of $\mathcal{D}$ and let $g=\dim _{%
%TCIMACRO{\U{2102} }%
%BeginExpansion
\mathbb{C}
%EndExpansion
}V_{%
%TCIMACRO{\U{2102} }%
%BeginExpansion
\mathbb{C}
%EndExpansion
,0}$; fix an integer $N\geq 3$ not divisible by $p$ and assume $K^{p}=U(N)$.
Let us denote by $\mathcal{S}_{\mathcal{D},N}:=\mathcal{S}_{\mathcal{D}%
,U(N)} $ the quasi-projective smooth scheme over $\mathcal{O}_{E}\otimes _{%
%TCIMACRO{\U{2124} }%
%BeginExpansion
\mathbb{Z}
%EndExpansion
}%
%TCIMACRO{\U{2124} }%
%BeginExpansion
\mathbb{Z}
%EndExpansion
_{(p)}$ representing the functor $\mathbf{M=M}(\mathcal{D})$. Let $\pi :%
\mathcal{X}_{\mathcal{D},N}\rightarrow \mathcal{S}_{\mathcal{D},N}$
\noindent be the corresponding universal abelian scheme over $\mathcal{S}_{%
\mathcal{D},N}$, and let $0:\mathcal{S}_{\mathcal{D},N}\rightarrow \mathcal{X%
}_{\mathcal{D},N}$ be its zero section. Denote by $\Omega _{\mathcal{X}_{%
\mathcal{D},N}/\mathcal{S}_{\mathcal{D},N}}^{1}$ the sheaf of relative
invariant differentials of $\mathcal{X}_{\mathcal{D},N}$ over $\mathcal{S}_{%
\mathcal{D},N}$: it is a locally free sheaf of $\mathcal{O}_{\mathcal{X}_{%
\mathcal{D},N}}$-modules over $\mathcal{X}_{\mathcal{D},N}$, having rank $g$%
. Its pull-back via the zero section $0$, i.e. the sheaf of relative
cotangent vectors at the origin of $\mathcal{X}_{\mathcal{D},N}:$ 
\begin{equation*}
\mathbb{E}:=0^{\ast }\left( \Omega _{\mathcal{X}_{\mathcal{D},N}/\mathcal{S}%
_{\mathcal{D},N}}^{1}\right) ,
\end{equation*}

\noindent is the \textit{Hodge bundle} of the PEL-scheme $\mathcal{S}_{%
\mathcal{D},N}$; it is a locally free sheaf of $\mathcal{O}_{\mathcal{S}_{%
\mathcal{D},N}}$-modules over $\mathcal{S}_{\mathcal{D},N}$ and its rank
equals $g$. If $\rho :GL_{g}\rightarrow GL_{m}$ is an $\mathcal{O}%
_{E}\otimes _{%
%TCIMACRO{\U{2124} }%
%BeginExpansion
\mathbb{Z}
%EndExpansion
}%
%TCIMACRO{\U{2124} }%
%BeginExpansion
\mathbb{Z}
%EndExpansion
_{(p)}$-representation of the algebraic group $GL_{g}$, we denote by $%
\mathbb{E}_{\rho }$ the locally free sheaf of rank $m$ on $\mathcal{S}_{%
\mathcal{D},N}$ obtained by twisting $\mathbb{E}$ via $\rho $ (cf. \cite%
{Gh04} 2.2.1).

\begin{definition}
Let $\mathcal{S}_{\mathcal{D},N}$ and $\rho \ $be as above. For any $%
\mathcal{O}_{E}\otimes _{%
%TCIMACRO{\U{2124} }%
%BeginExpansion
\mathbb{Z}
%EndExpansion
}%
%TCIMACRO{\U{2124} }%
%BeginExpansion
\mathbb{Z}
%EndExpansion
_{(p)}$-algebra $\mathfrak{R}$ we define the $\mathfrak{R}$-module:%
\begin{equation*}
M_{\rho }(\mathcal{D};\mathfrak{R}):=H^{0}(\mathcal{S}_{\mathcal{D}%
,N}\otimes _{\mathcal{O}_{E}\otimes _{%
%TCIMACRO{\U{2124} }%
%BeginExpansion
\mathbb{Z}
%EndExpansion
}%
%TCIMACRO{\U{2124} }%
%BeginExpansion
\mathbb{Z}
%EndExpansion
_{(p)}}\mathfrak{R;}\mathbb{E}_{\rho }\otimes \mathfrak{R}),
\end{equation*}

\noindent and we call it the space of \textbf{PEL-modular forms over }$%
\mathfrak{R}$\textbf{\ of} \textbf{weight }$\rho $\textbf{\ relative to the
moduli problem }$M(\mathcal{D})$. If no confusion arises we also denote this
space as $M_{\rho }^{g}(N;\mathfrak{R})$ and we say that the modular forms
in $M_{\rho }^{g}(N;\mathfrak{R})$ have genus $g$ and (full) level $N$.
\end{definition}

(For brevity, if we are given a morphism of schemes $X\rightarrow Y$ with a
section $e:Y\rightarrow X$, we set $\mathfrak{t}_{X/Y}^{\ast }:=e^{\ast
}\Omega _{X/Y}^{1}$; $\mathfrak{t}_{X/Y}^{\ast }$ is called the sheaf of
relative cotangent vectors along the section $e$. Notice that $\mathbb{E=}%
\mathfrak{t}_{\mathcal{X}_{\mathcal{D},N}/\mathcal{S}_{\mathcal{D},N}}^{\ast
}$).

With the above assumptions, we have: 
\begin{equation*}
\mathbb{E}_{\rho }\otimes \mathfrak{R=}0^{\ast }(\Omega _{\mathcal{X}_{%
\mathcal{D},N}\otimes \mathfrak{R}/\mathcal{S}_{\mathcal{D},N}\otimes 
\mathfrak{R}}^{1})_{\rho }=(\mathfrak{t}_{\mathcal{X}_{\mathcal{D},N}\otimes 
\mathfrak{R}/\mathcal{S}_{\mathcal{D},N}\otimes \mathfrak{R}}^{\ast })_{\rho
}.
\end{equation*}

\noindent Fix $f\in M_{\rho }(\mathcal{D};\mathfrak{R})$ and let $\xi :=(A,i,%
\overline{\lambda },\overline{\alpha })/S/\mathfrak{R}$ be a tuple such that 
$[\xi ]$ is an element of $\mathfrak{(}\mathcal{S}_{\mathcal{D},N}\otimes _{%
\mathcal{O}_{E}\otimes _{%
%TCIMACRO{\U{2124} }%
%BeginExpansion
\mathbb{Z}
%EndExpansion
}%
%TCIMACRO{\U{2124} }%
%BeginExpansion
\mathbb{Z}
%EndExpansion
_{(p)}}\mathfrak{R)(}S)$, for some $\mathfrak{R}$-scheme $S$; let $\phi
:S\rightarrow \mathcal{S}_{\mathcal{D},N}\otimes _{\mathcal{O}_{E}\otimes _{%
%TCIMACRO{\U{2124} }%
%BeginExpansion
\mathbb{Z}
%EndExpansion
}%
%TCIMACRO{\U{2124} }%
%BeginExpansion
\mathbb{Z}
%EndExpansion
_{(p)}}\mathfrak{R}$ be the morphism of $\mathfrak{R}$-schemes parametrizing
the element $[\xi ]$. We have:%
\begin{eqnarray*}
f([\xi ]) &\in &\phi ^{\ast }(\mathbb{E\otimes \mathfrak{R})}_{\rho }\mathbb{%
(}S\mathbb{)}= \\
&=&\phi ^{\ast }(0^{\ast }\Omega _{\mathcal{X}_{\mathcal{D},N}\otimes 
\mathfrak{R}/\mathcal{S}_{\mathcal{D},N}\otimes \mathfrak{R}}^{1})_{\rho
}(S)\simeq (\mathfrak{t}_{A/S}^{\ast })_{\rho }(S),
\end{eqnarray*}%
where the last isomorphism depends upon the choice of representative for $%
[\xi ]$.\ We conclude that for $\mathcal{D},p,N,\rho ,\mathfrak{R}$ as
above, a PEL-modular forms over $\mathfrak{R}$\ of weight $\rho $\ relative
to the moduli problem $M(\mathcal{D})$ is a rule $f$ that assigns, to any
tuple $\xi :=(A,i,\overline{\lambda },\overline{\alpha })/S/\mathfrak{R}$
such that $[\xi ]$ is an element of $\mathfrak{(}\mathcal{S}_{\mathcal{D}%
,N}\otimes _{\mathcal{O}_{E}\otimes _{%
%TCIMACRO{\U{2124} }%
%BeginExpansion
\mathbb{Z}
%EndExpansion
}%
%TCIMACRO{\U{2124} }%
%BeginExpansion
\mathbb{Z}
%EndExpansion
_{(p)}}\mathfrak{R)(}S)$ an element $f(\xi )$ of $(\mathfrak{t}_{A/S}^{\ast
})_{\rho }(S)$ in such a way that the rule $f$ is compatible isomorphisms
and commute with base change.

Pick a modular form $f$ described as in the above proposition and let $\xi
:=(A,i,\overline{\lambda },\overline{\alpha })/S/\mathfrak{R}$ and $\xi
^{\prime }:=(A^{\prime },i^{\prime },\overline{\lambda }^{\prime },\overline{%
\alpha }^{\prime })/S/\mathfrak{R}$ be two isomorphic tuples, say $\varphi
:\xi \rightarrow \xi ^{\prime }$ is an $S$-isomorphism; let $\varphi _{\rho
}^{\ast }:(\mathfrak{t}_{A^{\prime }/S}^{\ast })_{\rho }(S)\rightarrow (%
\mathfrak{t}_{A/S}^{\ast })_{\rho }(S)$ be the induced isomorphism on the
cotangent spaces. The compatibility of $f$ with respect to isomorphism means
that $\varphi _{\rho }^{\ast }(f(\xi ^{\prime }))=f(\xi )$. Let $\xi $ be as
above and fix an $S$-scheme $S^{\prime }$; denote by $c:\xi \otimes
_{S}S^{\prime }\rightarrow S^{\prime }$ the canonical map. The commutativity
of $f$ with base extension means that $c_{\rho }^{\ast }(f(\xi ))=f(\xi
\otimes _{S}S^{\prime })$ as elements of $(\mathfrak{t}_{A\otimes
_{S}S^{\prime }/S^{\prime }}^{\ast })_{\rho }(S^{\prime })$.

\paragraph{Modular forms mod $p$}

Let us assume that $\mathfrak{R}=\overline{\mathbb{F}}_{p}$ is a fixed
algebraic closure of the field with $p$ elements. If $A/\overline{\mathbb{F}}%
_{p}$ is an abelian variety of dimension $g$, we denote for brevity $%
\mathfrak{t}_{A/\overline{\mathbb{F}}_{p}}^{\ast }(\limfunc{Spec}\overline{%
\mathbb{F}}_{p})$ by $\mathfrak{t}_{A/\overline{\mathbb{F}}_{p}}^{\ast }$:
it is a vector space over $\overline{\mathbb{F}}_{p}$ of dimension $g$.

Let $A$ and $A^{\prime }$ be two abelian varieties over $\overline{\mathbb{F}%
}_{p}$ of dimension $g$, and let $\rho :GL_{g}\rightarrow GL_{m}$ as above;
let us fix isomorphisms of vector spaces $\gamma :\mathfrak{t}_{A/\overline{%
\mathbb{F}}_{p}}^{\ast }\rightarrow \overline{\mathbb{F}}_{p}^{g}$ and $%
\gamma ^{\prime }:\mathfrak{t}_{A^{\prime }/\overline{\mathbb{F}}_{p}}^{\ast
}\rightarrow \overline{\mathbb{F}}_{p}^{g}$. If $\varphi :A\rightarrow
A^{\prime }$ is a morphism of abelian varieties over $\overline{\mathbb{F}}%
_{p}$, then $\varphi ^{\ast }:\mathfrak{t}_{A^{\prime }/\overline{\mathbb{F}}%
_{p}}^{\ast }\rightarrow \mathfrak{t}_{A/\overline{\mathbb{F}}_{p}}^{\ast }$
is defined. By definition of the functor $\left( \cdot \right) _{\rho }$, we
can make the identifications $(\mathfrak{t}_{A/\overline{\mathbb{F}}%
_{p}}^{\ast })_{\rho }=(\mathfrak{t}_{A^{\prime }/\overline{\mathbb{F}}%
_{p}}^{\ast })_{\rho }=\overline{\mathbb{F}}_{p}^{m}$, so that $\varphi
_{\rho }^{\ast }:=\rho (\gamma \circ \varphi ^{\ast }\circ \gamma ^{\prime
-1})$ is a morphism $(\mathfrak{t}_{A^{\prime }/\overline{\mathbb{F}}%
_{p}}^{\ast })_{\rho }\rightarrow (\mathfrak{t}_{A/\overline{\mathbb{F}}%
_{p}}^{\ast })_{\rho }$.

Let us compute $H^{0}((\mathcal{S}_{\mathcal{D},N}\otimes \overline{\mathbb{F%
}}_{p})(\overline{\mathbb{F}}_{p})\mathfrak{;}\mathbb{E}_{\rho }\otimes 
\overline{\mathbb{F}}_{p})$. Fix a modular form $f$ and a tuple $\xi =(A,i,%
\overline{\lambda },\overline{\alpha })/\overline{\mathbb{F}}_{p}$ such that 
$[\xi ]$ is an element of $\mathfrak{(}\mathcal{S}_{\mathcal{D},N}\otimes 
\overline{\mathbb{F}}_{p}\mathfrak{)(}\overline{\mathbb{F}}_{p})$; fix also
an isomorphism of vector spaces $\gamma :\mathfrak{t}_{A/\overline{\mathbb{F}%
}_{p}}^{\ast }\rightarrow \overline{\mathbb{F}}_{p}^{g}$ that will allow us
to identify these two spaces in the sequel (the choice of $\gamma $ will not
be influent). If $\eta =\left( \eta _{1},...,\eta _{g}\right) $ is an
ordered basis for $\mathfrak{t}_{A/\overline{\mathbb{F}}_{p}}^{\ast }$ over $%
\overline{\mathbb{F}}_{p}$, denote by the same symbol the matrix $\eta
=[\eta _{1}|...|\eta _{g}]\in GL_{g}(\overline{\mathbb{F}}_{p})$ obtained by
placing the vectors $\eta _{i}\in \overline{\mathbb{F}}_{p}^{g}$ as columns.
Write $\underline{\nu }=[\nu _{1}|...|\nu _{m}]:=\rho (\eta )$, so that we
can find a unique column vector $x\in \overline{\mathbb{F}}_{p}^{m}$ such
that $f(\xi )=\sum\nolimits_{j=1}^{m}x_{j}\nu _{j}=\rho (\eta )\cdot x.$

We can define an assignment $\widetilde{f}$ on tuples $(A,i,\overline{%
\lambda },\overline{\alpha },\eta )$, where $(A,i,\overline{\lambda },%
\overline{\alpha })/\overline{\mathbb{F}}_{p}$ is as above and $\eta $ is an
ordered basis for $\mathfrak{t}_{A/\overline{\mathbb{F}}_{p}}^{\ast }$ over $%
\overline{\mathbb{F}}_{p}$, by setting:%
\begin{equation*}
\widetilde{f}:(A,i,\overline{\lambda },\overline{\alpha },\eta )\longmapsto x%
\text{ \ }\Longleftrightarrow \text{ \ }f(A,i,\overline{\lambda },\overline{%
\alpha })=\rho (\eta )\cdot x.\text{\ }
\end{equation*}

\noindent Notice that if $M\in GL_{g}(\overline{\mathbb{F}}_{p})$, then:%
\begin{equation*}
\widetilde{f}(A,i,\overline{\lambda },\overline{\alpha },\eta M)=\rho
(M)^{-1}\cdot \widetilde{f}(A,i,\overline{\lambda },\overline{\alpha },\eta
).
\end{equation*}

Assume that we are given another tuple $\xi ^{\prime }=(A^{\prime
},i^{\prime },\overline{\lambda }^{\prime },\overline{\alpha }^{\prime })/%
\overline{\mathbb{F}}_{p}$ such that $[\xi ^{\prime }]$ is an element of $%
\mathfrak{(}\mathcal{S}_{\mathcal{D},N}\otimes \overline{\mathbb{F}}_{p}%
\mathfrak{)(}\overline{\mathbb{F}}_{p})$; fix an isomorphism of vector
spaces $\gamma ^{\prime }:\mathfrak{t}_{A^{\prime }/\overline{\mathbb{F}}%
_{p}}^{\ast }\rightarrow \overline{\mathbb{F}}_{p}^{g}$ and pick an ordered
basis $\eta ^{\prime }$ of invariant differentials for $A^{\prime }$. Assume
there is an isomorphism $\varphi :(\xi ,\eta )\rightarrow (\xi ^{\prime
},\eta ^{\prime })$, i.e. an isomorphism of abelian varieties (up to
primo-to-$p$ isogeny) with additional structure $\varphi :\xi \rightarrow
\xi ^{\prime }$ such that $\varphi ^{\ast }\eta ^{\prime }=\eta $. Notice
that by functoriality of the $\rho $-twisting we have $\rho (\eta )=\varphi
_{\rho }^{\ast }\cdot \rho (\eta ^{\prime })$, since $\varphi ^{\ast }\eta
^{\prime }=\eta $. We compute, using the fact that $f$ is compatible with
isomorphisms:%
\begin{eqnarray*}
\widetilde{f}(\xi ,\eta ) &=&\rho (\eta )^{-1}\cdot f(\xi )= \\
&=&(\varphi _{\rho }^{\ast }\cdot \rho (\eta ^{\prime }))^{-1}\cdot \varphi
_{\rho }^{\ast }\cdot f(\xi ^{\prime })= \\
&=&\rho (\eta ^{\prime })^{-1}\cdot f(\xi ^{\prime })=\widetilde{f}(\xi
^{\prime },\eta ^{\prime }).
\end{eqnarray*}

We have shown:

\begin{proposition}
\label{down to earth mod forms}Let $\mathcal{D},p,N,\rho $ be as above; then
a PEL-modular forms over $\overline{\mathbb{F}}_{p}$\ of weight $\rho $\
relative to the moduli problem $M(\mathcal{D})$ is a rule $f$ that assigns,
to any tuple $(A,i,\overline{\lambda },\overline{\alpha },\eta )/\overline{%
\mathbb{F}}_{p}$ such that $[(A,i,\overline{\lambda },\overline{\alpha })]$
is an element of $\mathfrak{(}\mathcal{S}_{\mathcal{D},N}\otimes \overline{%
\mathbb{F}}_{p}\mathfrak{)(}\overline{\mathbb{F}}_{p})$ and $\eta $ is an
ordered basis for $\mathfrak{t}_{A/\overline{\mathbb{F}}_{p}}^{\ast }$ over $%
\overline{\mathbb{F}}_{p}$, an element $f(A,i,\overline{\lambda },\overline{%
\alpha },\eta )\in \overline{\mathbb{F}}_{p}^{m}$ of in such a way that:

\begin{description}
\item[(a)] $f(A,i,\overline{\lambda },\overline{\alpha },\eta M)=\rho
(M)^{-1}\cdot f(A,i,\overline{\lambda },\overline{\alpha },\eta )$ for all $%
M\in GL_{g}(\overline{\mathbb{F}}_{p});$

\item[(b)] if $(A,i,\overline{\lambda },\overline{\alpha },\eta )\simeq
(A^{\prime },i^{\prime },\overline{\lambda }^{\prime },\overline{\alpha }%
^{\prime },\eta ^{\prime })$ then $f(A,i,\overline{\lambda },\overline{%
\alpha },\eta )=f(A^{\prime },i^{\prime },\overline{\lambda }^{\prime },%
\overline{\alpha }^{\prime },\eta ^{\prime })$.
\end{description}
\end{proposition}

\paragraph{Hilbert-Siegel modular forms\label{Siegel}}

Let $N\geq 3$ be an integer prime to $p$, $F$ a totally real Galois
extension of $%
%TCIMACRO{\U{211a} }%
%BeginExpansion
\mathbb{Q}
%EndExpansion
$ having degree $f$ in which $p$ is unramified. We consider the functor $%
\mathbf{M}(\mathcal{D}_{2g,p}^{Sp(F)})$ associated to the PEL-datum $%
\mathcal{D}_{2g,p}^{Sp(F)}$ with good reduction at $p$\ (cf. \ref{Examples}%
), with $K^{p}=U(N)$; in this case the Shimura field is $E=%
%TCIMACRO{\U{211a} }%
%BeginExpansion
\mathbb{Q}
%EndExpansion
$. Fix a rational $%
%TCIMACRO{\U{2124} }%
%BeginExpansion
\mathbb{Z}
%EndExpansion
_{(p)}$-representation $\rho :GL_{fg}\rightarrow GL_{m}$ and let $\mathfrak{R%
}$ be any $%
%TCIMACRO{\U{2124} }%
%BeginExpansion
\mathbb{Z}
%EndExpansion
_{(p)}$-algebra. The $\mathfrak{R}$-module:%
\begin{equation*}
M_{\rho }^{fg}(F;N;\mathfrak{R})_{HS}:=M_{\rho }(\mathcal{D}_{2g,p}^{Sp(F)};%
\mathfrak{R})
\end{equation*}

\noindent is the space of \textit{Hilbert-Siegel }$\mathfrak{R}$\textit{%
-modular forms} of genus $fg$, level $N$, weight $\rho $ and relative to the
field $F$. If $\mathfrak{R=%
%TCIMACRO{\U{2102} }%
%BeginExpansion
\mathbb{C}
%EndExpansion
}$, the Hermitian symmetric domain associated to the corresponding Shimura
variety $\mathcal{S}_{\mathcal{D}_{2g,p}^{Sp(F)},N}$ is the product $%
\mathfrak{h}_{g}^{f}$ of $f$ copies of the genus-$g$ Siegel upper half plane 
$\mathfrak{h}_{g}$. If furthermore $f=1$, we obtain the classical space of
Siegel modular forms of genus $g$. (Cf. \cite{FC}, Ch.V and \cite{Lan}
1.4.1-3 for the comparison between $\mathbf{M}(\mathcal{D}_{2g}^{Sp(%
%TCIMACRO{\U{211a} }%
%BeginExpansion
\mathbb{Q}
%EndExpansion
)})$\ and the classical Mumford's functor $\mathcal{A}_{g,1,N}$).

\paragraph{Unitary modular forms\label{quasi-split}}

We consider the functor $\mathbf{M}(\mathcal{D}_{(r,s),p}^{U})$ associated
to the PEL-datum $\mathcal{D}_{(r,s),p}^{U}$ with good reduction at $p$ as
defined in \ref{Examples}. Let $N\geq 3$ denotes an integer prime to $p$.
Recall that the PEL-datum comes with a quadratic imaginary field $k$ in
which $p$ is inert. For any $%
%TCIMACRO{\U{2124} }%
%BeginExpansion
\mathbb{Z}
%EndExpansion
_{(p)}$-algebra $\mathfrak{R}$, we call 
\begin{equation*}
M_{\rho }^{(r,s)}(k;N;\mathfrak{R})_{U}:=M_{\rho }(\mathcal{D}_{(r,s),p}^{U};%
\mathfrak{R})
\end{equation*}%
the space of\textbf{\ }\textit{unitary }$\mathfrak{R}$\textit{-modular forms
of signature }$(r,s)$\textit{\ for the field }$k$, having genus $g:=r+s$,
level $N$ and weight $\rho .$ Over the complex numbers, these forms can be
constructed analytically starting from the Hermitian domain $\mathfrak{h}%
_{s,r}:=\{z\in M_{r\times s}(%
%TCIMACRO{\U{2102} }%
%BeginExpansion
\mathbb{C}
%EndExpansion
):1-\overline{z}^{t}z$ positive Hermitian$\}$ (the Picard space).

\newpage

\section{Uniformization results for the supersingular and the superspecial
loci}

We recall a uniformization result for isogeny classes in a PEL-moduli space
due to Rapoport and Zink (\cite{RZ}, Ch. 6); we then present a modification
of this result that allows us to parametrize the superspecial locus.

\bigskip

\subsection{The result of Rapoport and Zink\label{p-adic unif of isogeny
classes}}

We fix some notation. Let $\mathcal{D=(}B,^{\ast },V,\left\langle
,\right\rangle ,\mathcal{O}_{B},\Lambda ,h,K^{p},\nu \mathcal{)}$ be a
simple $%
%TCIMACRO{\U{211a} }%
%BeginExpansion
\mathbb{Q}
%EndExpansion
$-PEL-datum for moduli of abelian schemes with good reduction at $p$, and
neat level $K^{p}$. Let $G$ the associated reductive group over $%
%TCIMACRO{\U{211a} }%
%BeginExpansion
\mathbb{Q}
%EndExpansion
$, and $E$ the Shimura field. The completion $E_{\nu }$ of $E$ at $\nu $
coincide with the field of definition of the $G^{0}(%
%TCIMACRO{\U{211a} }%
%BeginExpansion
\mathbb{Q}
%EndExpansion
_{p})$-conjugacy class of $\mu $. Let $k=\overline{\mathbb{F}}_{p}$ be a
fixed algebraic closure of the residue field of $E_{\nu }$, and let $W=W(%
\overline{\mathbb{F}}_{p})$, $K_{0}=W[\frac{1}{p}]$ and $\sigma $ the
Frobenius morphism of $W$; fix a finite extension $K$ of $K_{0}$ such that $%
\mu $ is defined over $K$, so that the corresponding weight decomposition $%
V_{K}=V_{K,0}\oplus V_{K,1}$ is also defined over $K$. Set $\breve{E}=E_{\nu
}K_{0}$ as in \ref{moduli-functor-p div def.}: since we are in the good
reduction case, we have $\breve{E}=K_{0}$. Set $B_{p}:=B\otimes _{%
%TCIMACRO{\U{211a} }%
%BeginExpansion
\mathbb{Q}
%EndExpansion
}%
%TCIMACRO{\U{211a} }%
%BeginExpansion
\mathbb{Q}
%EndExpansion
_{p}$, $V_{p}:=V\otimes _{%
%TCIMACRO{\U{211a} }%
%BeginExpansion
\mathbb{Q}
%EndExpansion
}%
%TCIMACRO{\U{211a} }%
%BeginExpansion
\mathbb{Q}
%EndExpansion
_{p}$, $\left\langle ,\right\rangle _{p}:=\left\langle ,\right\rangle
\otimes _{%
%TCIMACRO{\U{211a} }%
%BeginExpansion
\mathbb{Q}
%EndExpansion
}%
%TCIMACRO{\U{211a} }%
%BeginExpansion
\mathbb{Q}
%EndExpansion
_{p}$, $G_{p}=G_{%
%TCIMACRO{\U{211a} }%
%BeginExpansion
\mathbb{Q}
%EndExpansion
_{p}}$, $\mathcal{O}_{B_{p}}:=\mathcal{O}_{B}\otimes _{%
%TCIMACRO{\U{2124} }%
%BeginExpansion
\mathbb{Z}
%EndExpansion
}%
%TCIMACRO{\U{2124} }%
%BeginExpansion
\mathbb{Z}
%EndExpansion
_{p}$.

\subsubsection{From global PEL-data to local PEL-data}

Let $\mathcal{S}_{\mathcal{D},K^{p}}$ be the quasi-projective smooth scheme
over $\mathcal{O}_{E}\otimes _{%
%TCIMACRO{\U{2124} }%
%BeginExpansion
\mathbb{Z}
%EndExpansion
}%
%TCIMACRO{\U{2124} }%
%BeginExpansion
\mathbb{Z}
%EndExpansion
_{(p)}$ representing the functor $\mathbf{M=M}(\mathcal{D})$ of Theorem \ref%
{Representability}, and let us fix a point $[(A_{0},i_{0},\overline{\lambda }%
_{0},\overline{\alpha }_{0})]\in \mathcal{S}_{\mathcal{D},K^{p}}\left( 
\overline{\mathbb{F}}_{p}\right) $, where $\lambda _{0}$ denotes a principal
polarization of $(A_{0},i_{0})$. Correspondingly we have a $p$-divisible
group $\mathbf{X:=}A_{0}(p)$ over $\overline{\mathbb{F}}_{p}$, endowed with
the action $i_{\mathbf{X}}:$ $\mathcal{O}_{B_{p}}\rightarrow \limfunc{End}%
\mathbf{X}$ of $\mathcal{O}_{B_{p}}$\ induced by $i_{0}$; the principal
polarization $\lambda _{0}$ induces a principal polarization $\lambda _{%
\mathbf{X}}:A_{0}(p)\rightarrow \widehat{A_{0}}(p)=\widehat{A_{0}(p)}$ of $p$%
-divisible groups (cf. \cite{O}, 1.8, 1.12) respecting the $\mathcal{O}%
_{B_{p}}$-action; $\lambda _{\mathbf{X}}$ is well defined up to a constant
in $%
%TCIMACRO{\U{211a} }%
%BeginExpansion
\mathbb{Q}
%EndExpansion
_{p}^{\times }$. The triple $(\mathbf{X},i_{\mathbf{X}},\overline{\lambda }_{%
\mathbf{X}})$ is well defined modulo isomorphisms by the given point $%
[(A_{0},i_{0},\overline{\lambda }_{0},\overline{\alpha }_{0})]\in \mathcal{S}%
_{\mathcal{D},K^{p}}\left( \overline{\mathbb{F}}_{p}\right) .$

By covariant Dieudonn\'{e} theory, we can associate to $(\mathbf{X},i_{%
\mathbf{X}},\overline{\lambda }_{\mathbf{X}})$ the isocrystal $(N:=M_{\ast }(%
\mathbf{X})[\frac{1}{p}],\mathbf{F})$ over $K_{0}$ endowed with an action of 
$B_{p}$ and with a non-degenerate bilinear form of isocrystals $\Psi
:N\times N\rightarrow \mathbf{1}(1)$ (well defined up to an element of $%
%TCIMACRO{\U{211a} }%
%BeginExpansion
\mathbb{Q}
%EndExpansion
_{p}^{\times }$) that is skew-Hermitian with respect to $^{\ast }$. The
quasi-isogeny class of the principally polarized $B_{p}$-isocrystal $(N,%
\mathbf{F},%
%TCIMACRO{\U{211a} }%
%BeginExpansion
\mathbb{Q}
%EndExpansion
_{p}^{\times }\Psi )$ depends upon the isomorphism class of $(\mathbf{X},i_{%
\mathbf{X}},\overline{\lambda }_{\mathbf{X}})$.

We fix an isomorphism of $B\otimes _{%
%TCIMACRO{\U{211a} }%
%BeginExpansion
\mathbb{Q}
%EndExpansion
}K_{0}$-modules $N\simeq V\otimes _{%
%TCIMACRO{\U{211a} }%
%BeginExpansion
\mathbb{Q}
%EndExpansion
}K_{0}$ that respects the skew-symmetric forms on both sides. We then write
the action of Frobenius on the right hand side as $\mathbf{F}=b\otimes
\sigma $ for a unique element $b\in G_{p}(K_{0})$. By construction we have,
for any $x,y\in V\otimes _{%
%TCIMACRO{\U{211a} }%
%BeginExpansion
\mathbb{Q}
%EndExpansion
}K_{0}$:%
\begin{equation*}
c(b)\left\langle x,y\right\rangle =\left\langle b\otimes \sigma \cdot
x,b\otimes \sigma \cdot y\right\rangle ^{\sigma ^{-1}}=p\left\langle
x,y\right\rangle ,
\end{equation*}

\noindent so that $c(b)=p$. The isocrystal $V\otimes _{%
%TCIMACRO{\U{211a} }%
%BeginExpansion
\mathbb{Q}
%EndExpansion
}K_{0}$ has slopes in the interval $[0,1]$, and in the decomposition of the $%
K$-vector space $V\otimes _{%
%TCIMACRO{\U{211a} }%
%BeginExpansion
\mathbb{Q}
%EndExpansion
}K$ under the co-character $\mu $ only the weights $0$ and $1$ appear.
\noindent Finally, the pair $(b,\mu )$ is admissible in the sense of
Definition \ref{admissibility}. For, since $\lambda _{0}$ is a separable
polarization, the point $\left[ (A_{0},i_{0},\lambda _{0},\alpha _{0})\right]
\in \mathcal{S}_{\mathcal{D},K^{p}}\left( \overline{\mathbb{F}}_{p}\right) $
is liftable to a point $[(\widetilde{A}_{0},\widetilde{i}_{0},\overline{%
\widetilde{\lambda }}_{0},\overline{\widetilde{\alpha }}_{0})]\in \mathcal{S}%
_{\mathcal{D},K^{p}}\left( W\right) $, by Theorem \ref{lifting with
separability}; this implies that the $p$-divisible group $(\mathbf{X},i_{%
\mathbf{X}},\overline{\lambda }_{\mathbf{X}})$ over $\overline{\mathbb{F}}%
_{p}$ can be correspondingly lifted to a $p$-divisible group $(\widetilde{%
\mathbf{X}},i_{\widetilde{\mathbf{X}}},\overline{\lambda }_{\widetilde{%
\mathbf{X}}})$ over $W$, so that the $K$-filtered isocrystal over $K_{0}$
given by $(V\otimes _{%
%TCIMACRO{\U{211a} }%
%BeginExpansion
\mathbb{Q}
%EndExpansion
}K_{0},b\otimes \sigma ,\left\{ V_{K,1}\subset V_{K}\right\} )$ is
associated to a $p$-divisible group over $W\subset \mathcal{O}_{K}$\bigskip
. By the considerations we made in \ref{def - admiss}, this implies that $%
(b,\mu )$ is admissible, since we are in the good reduction case.

We conclude that the choice of a global PEL-datum $\mathcal{D}$ plus a fixed
point $[(A_{0},i_{0},\overline{\lambda }_{0},\overline{\alpha }_{0})]\in 
\mathcal{S}_{\mathcal{D},K^{p}}\left( \overline{\mathbb{F}}_{p}\right) $
determines a simple $%
%TCIMACRO{\U{211a} }%
%BeginExpansion
\mathbb{Q}
%EndExpansion
_{p}$-PEL datum for moduli of $p$-divisible groups: 
\begin{equation*}
\mathcal{D}_{p}:\mathcal{=(}B_{p},^{\ast },V_{p},\left\langle ,\right\rangle
_{p},\mathcal{O}_{B_{p}},\Lambda ,b,\mu ),
\end{equation*}%
in the sense of Definition \ref{PEL p-divisible gps}, having good reduction
at $p$ and Shimura field equal to $E_{\nu }$. Denote by $\mathcal{\breve{M}}$
the contravariant functor from $NILP_{\mathcal{O}_{\overset{\smile }{E_{\nu }%
}}}$ to $SETS$ defined in \ref{Shi p-div} starting from the $%
%TCIMACRO{\U{211a} }%
%BeginExpansion
\mathbb{Q}
%EndExpansion
_{p}$-homogeneously principally polarized $p$-divisible group $(\mathbf{X}%
,i_{\mathbf{X}},\overline{\lambda }_{\mathbf{X}})$ endowed with the action
of $\mathcal{O}_{B_{p}}$: $\mathcal{\breve{M}}$ is representable by a formal
scheme $\mathcal{\breve{M}}$ which is formally locally of finite type over $%
\limfunc{Spf}\mathcal{O}_{_{\overset{\smile }{E_{\nu }}}}$. Furthermore $%
\overset{\smile }{E_{\nu }}=K_{0}$ and the $\mathcal{\breve{M}}$ is formally
smooth over $\limfunc{Spf}W.$

\paragraph{The groups $I$ and $J\label{groups I and J}$}

Denote by $J(%
%TCIMACRO{\U{211a} }%
%BeginExpansion
\mathbb{Q}
%EndExpansion
_{p})$ the group of automorphisms of the polarized $B\otimes _{%
%TCIMACRO{\U{211a} }%
%BeginExpansion
\mathbb{Q}
%EndExpansion
}K_{0}$-isocrystal $(V\otimes _{%
%TCIMACRO{\U{211a} }%
%BeginExpansion
\mathbb{Q}
%EndExpansion
}K_{0};b\otimes \sigma )$: this means that $J(%
%TCIMACRO{\U{211a} }%
%BeginExpansion
\mathbb{Q}
%EndExpansion
_{p})$ is the group of $K_{0}$-automorphisms of the isocrystal $(V\otimes _{%
%TCIMACRO{\U{211a} }%
%BeginExpansion
\mathbb{Q}
%EndExpansion
}K_{0};b\otimes \sigma )$ equivariant for the action of $B_{p}$ and
preserving the polarization form induced by $\left\langle ,\right\rangle $
up to a non-zero scalar in $%
%TCIMACRO{\U{211a} }%
%BeginExpansion
\mathbb{Q}
%EndExpansion
_{p}$. $J(%
%TCIMACRO{\U{211a} }%
%BeginExpansion
\mathbb{Q}
%EndExpansion
_{p})$ is the group of $%
%TCIMACRO{\U{211a} }%
%BeginExpansion
\mathbb{Q}
%EndExpansion
_{p}$-rational points of an algebraic group $J$ defined over $%
%TCIMACRO{\U{211a} }%
%BeginExpansion
\mathbb{Q}
%EndExpansion
_{p}$. By covariant Dieudonn\'{e} theory, $J(%
%TCIMACRO{\U{211a} }%
%BeginExpansion
\mathbb{Q}
%EndExpansion
_{p})$ is isomorphic to the group $J^{\prime }(%
%TCIMACRO{\U{211a} }%
%BeginExpansion
\mathbb{Q}
%EndExpansion
_{p})$\ of quasi-isogenies $f:(\mathbf{X},i_{\mathbf{X}})\mathbf{\rightarrow 
}(\mathbf{X},i_{\mathbf{X}})$ of $p$-divisible groups over $\overline{%
\mathbb{F}}_{p}$\ such that $\widehat{f}\circ \lambda _{\mathbf{X}}\circ
f\in 
%TCIMACRO{\U{211a} }%
%BeginExpansion
\mathbb{Q}
%EndExpansion
_{p}^{\times }\lambda _{\mathbf{X}}$. The isomorphism $J(%
%TCIMACRO{\U{211a} }%
%BeginExpansion
\mathbb{Q}
%EndExpansion
_{p})\rightarrow J^{\prime }(%
%TCIMACRO{\U{211a} }%
%BeginExpansion
\mathbb{Q}
%EndExpansion
_{p})$ depends upon the choice of isomorphism of $B\otimes _{%
%TCIMACRO{\U{211a} }%
%BeginExpansion
\mathbb{Q}
%EndExpansion
}K_{0}$-modules $N\simeq V\otimes _{%
%TCIMACRO{\U{211a} }%
%BeginExpansion
\mathbb{Q}
%EndExpansion
}K_{0}$ that respects the skew-symmetric forms on both sides (recall that $N$
is the isocrystal associated to $\mathbf{X}=A_{0}(p)$), hence it depends
upon the choice of $b\in G_{p}(K_{0})$. We have fixed such an isomorphism,
so that we think $J(%
%TCIMACRO{\U{211a} }%
%BeginExpansion
\mathbb{Q}
%EndExpansion
_{p})=J^{\prime }(%
%TCIMACRO{\U{211a} }%
%BeginExpansion
\mathbb{Q}
%EndExpansion
_{p}).$ The isomorphism class of $J(%
%TCIMACRO{\U{211a} }%
%BeginExpansion
\mathbb{Q}
%EndExpansion
_{p})$ only depends upon the quasi-isogeny class of $(\mathbf{X},i_{\mathbf{X%
}},\overline{\lambda }_{\mathbf{X}})$.

The group $J(%
%TCIMACRO{\U{211a} }%
%BeginExpansion
\mathbb{Q}
%EndExpansion
_{p})$ acts on $\mathcal{\breve{M}}$ from the left by the rule:%
\begin{equation*}
g\cdot \lbrack (X,i,\overline{\lambda };\rho )]:=[(X,i,\overline{\lambda }%
;\rho \circ g^{-1})],
\end{equation*}%
where $[(X,i,\overline{\lambda };\rho )]\in \mathcal{\breve{M}}(S)$ for some
scheme $S$ in $NILP_{\mathcal{O}_{\overset{\smile }{E_{\nu }}}}$.

Denote by $I(%
%TCIMACRO{\U{211a} }%
%BeginExpansion
\mathbb{Q}
%EndExpansion
)$ the group of quasi-isogenies of $(A_{0},i_{0},\overline{\lambda }_{0})$,
i.e. the group of quasi-isogenies of $\overline{\mathbb{F}}_{p}$-abelian
variety $(A_{0},i_{0})\rightarrow (A_{0},i_{0})$ that send $\overline{%
\lambda }_{0}$ into itself; $I\left( 
%TCIMACRO{\U{211a} }%
%BeginExpansion
\mathbb{Q}
%EndExpansion
\right) $ is the group of rational points of an algebraic group $I$ defined
over $%
%TCIMACRO{\U{211a} }%
%BeginExpansion
\mathbb{Q}
%EndExpansion
$ and the isomorphism class of $I$ only depends upon the quasi-isogeny class
of $(A_{0},i_{0},\overline{\lambda }_{0})$.

The group $I(%
%TCIMACRO{\U{211a} }%
%BeginExpansion
\mathbb{Q}
%EndExpansion
)$ acts by quasi-isogenies on the tuple $(\mathbf{X},i_{\mathbf{X}},%
\overline{\lambda }_{\mathbf{X}})$, hence on $((N,\mathbf{F}),i,%
%TCIMACRO{\U{211a} }%
%BeginExpansion
\mathbb{Q}
%EndExpansion
_{p}^{\times }\Psi )$,\ defining a morphism $I(%
%TCIMACRO{\U{211a} }%
%BeginExpansion
\mathbb{Q}
%EndExpansion
)\longrightarrow J(%
%TCIMACRO{\U{211a} }%
%BeginExpansion
\mathbb{Q}
%EndExpansion
_{p})$ factoring through $I(%
%TCIMACRO{\U{211a} }%
%BeginExpansion
\mathbb{Q}
%EndExpansion
)\hookrightarrow I(%
%TCIMACRO{\U{211a} }%
%BeginExpansion
\mathbb{Q}
%EndExpansion
_{p})$. Such a monomorphism depends upon the choice of an isomorphism $%
J^{\prime }(%
%TCIMACRO{\U{211a} }%
%BeginExpansion
\mathbb{Q}
%EndExpansion
_{p})\simeq J(%
%TCIMACRO{\U{211a} }%
%BeginExpansion
\mathbb{Q}
%EndExpansion
_{p})$. Since we fixed an isomorphism of isocrystal (with additional
structure) $N\simeq V\otimes _{%
%TCIMACRO{\U{211a} }%
%BeginExpansion
\mathbb{Q}
%EndExpansion
}K_{0}$ at the beginning of this section, we denote by $\alpha _{0,p}:I(%
%TCIMACRO{\U{211a} }%
%BeginExpansion
\mathbb{Q}
%EndExpansion
)\longrightarrow J(%
%TCIMACRO{\U{211a} }%
%BeginExpansion
\mathbb{Q}
%EndExpansion
_{p})$ the corresponding map.

Notice that $I(%
%TCIMACRO{\U{211a} }%
%BeginExpansion
\mathbb{Q}
%EndExpansion
)$ acts by skew-Hermitian symplectic $B$-equivariant similitudes on the
space $H_{1}(A_{0},\mathbb{A}_{f}^{p})$, that we identified with $V\otimes _{%
%TCIMACRO{\U{211a} }%
%BeginExpansion
\mathbb{Q}
%EndExpansion
}\mathbb{A}_{f}^{p}$ through the isomorphism $\alpha _{0}$, and hence we
also have a homomorphism $\alpha _{0}^{p}:I(%
%TCIMACRO{\U{211a} }%
%BeginExpansion
\mathbb{Q}
%EndExpansion
)\rightarrow G(\mathbb{A}_{f}^{p})$ depending of the choice of a
representative $\alpha _{0}$ for the class $\overline{\alpha }_{0}$ (hence
well defined only modulo $K^{p}$).

In conclusion, once $b$ and $\alpha _{0}$ are fixed, we have well defined
group homomorphisms $\alpha _{0,p}:I(%
%TCIMACRO{\U{211a} }%
%BeginExpansion
\mathbb{Q}
%EndExpansion
)\longrightarrow J(%
%TCIMACRO{\U{211a} }%
%BeginExpansion
\mathbb{Q}
%EndExpansion
_{p})$ and $\alpha _{0}^{p}:I(%
%TCIMACRO{\U{211a} }%
%BeginExpansion
\mathbb{Q}
%EndExpansion
)\longrightarrow G(\mathbb{A}_{f}^{p})$. \noindent We can make the
identification:%
\begin{equation*}
J(%
%TCIMACRO{\U{211a} }%
%BeginExpansion
\mathbb{Q}
%EndExpansion
_{p})=\{g\in G(K_{0}):g\cdot b\sigma =b\sigma \cdot g,c(g)\in 
%TCIMACRO{\U{211a} }%
%BeginExpansion
\mathbb{Q}
%EndExpansion
_{p}^{\times }\}.
\end{equation*}

\subsubsection{The uniformization morphism}

The formal scheme $\mathcal{\breve{M}}$ is not defined, in general, over $%
\mathcal{O}_{E_{\nu }}$, i.e. over the (ring of integer of the)\ local
Shimura field. There is a suitable completion of $\mathcal{\breve{M}}$ that
can be written as $\mathcal{M\otimes }_{\limfunc{Spf}(\mathcal{O}_{E_{\nu
}})}\limfunc{Spf}(\mathcal{O}_{\breve{E}_{\nu }})$ for a pro-formal scheme $%
\mathcal{M}$ over $Spf(\mathcal{O}_{E_{\nu }})$. The uniformization result
of Rapoport and Zink that we recall here can be stated in terms of $\mathcal{%
M}.$

To this purpose, let $s\in 
%TCIMACRO{\U{2124} }%
%BeginExpansion
\mathbb{Z}
%EndExpansion
$ such that $\mathbb{D}_{K_{0}}\overset{s}{\mathbb{\rightarrow }}\mathbb{D}%
_{K_{0}}\overset{\nu }{\mathbb{\rightarrow }}G_{K_{0}}$ factors, via the
canonical projection $\mathbb{D}_{K_{0}}\rightarrow \mathbb{G}_{m/K_{0}}$
induced by $%
%TCIMACRO{\U{2124} }%
%BeginExpansion
\mathbb{Z}
%EndExpansion
\subset 
%TCIMACRO{\U{211a} }%
%BeginExpansion
\mathbb{Q}
%EndExpansion
$, through a map $s\nu :\mathbb{G}_{m/K_{0}}\rightarrow G_{K_{0}}$. Since $%
\overline{\mathbb{F}}_{p}$ is algebraically closed, and if we assume $G$ to
be connected, we can assume that $b$ satisfies the decency equation $%
(b\sigma )^{s}=((s\nu )(p))\sigma ^{s}$ \noindent in $G(K_{0})\rtimes
\left\langle \sigma \right\rangle $, and that $b\in G(%
%TCIMACRO{\U{211a} }%
%BeginExpansion
\mathbb{Q}
%EndExpansion
_{p^{s}})$ (\cite{RZ}, 1.8-9). Let $\gamma _{s}:=p^{s}\cdot (s\nu
(p))^{-1}\in G(K_{0})$; since $\gamma _{s}\cdot b\sigma =b\sigma \cdot
\gamma _{s}$ and $c(\gamma _{s})\in 
%TCIMACRO{\U{211a} }%
%BeginExpansion
\mathbb{Q}
%EndExpansion
_{p}^{\times }$, we have $\gamma _{s}\in J(%
%TCIMACRO{\U{211a} }%
%BeginExpansion
\mathbb{Q}
%EndExpansion
_{p}).$ By what we saw above, $\gamma _{s}$ acts upon the functor $\mathcal{%
\breve{M}}$, and we define $\mathcal{M}_{s}$ to be the sheaf associated to
the functor:%
\begin{equation*}
S\longmapsto \mathcal{\breve{M}}(S)/\gamma _{s}^{%
%TCIMACRO{\U{2124} }%
%BeginExpansion
\mathbb{Z}
%EndExpansion
}
\end{equation*}

\noindent ($S\in Nilp_{\mathcal{O}_{\breve{E}_{\nu }}}$). By \cite{RZ},
3.42, $\mathcal{M}_{s}$ is represented by a formal scheme locally of finite
type over $\limfunc{Spf}(\mathcal{O}_{\breve{E}_{\nu }})$; notice that $J(%
%TCIMACRO{\U{211a} }%
%BeginExpansion
\mathbb{Q}
%EndExpansion
_{p})$ continues to act on $\mathcal{M}_{s}$.

The formal schemes $\mathcal{\breve{M}}$ and $\mathcal{M}_{s}$ have a
natural Weyl descent datum relative to the extension $\mathcal{O}_{\breve{E}%
_{\nu }}/\mathcal{O}_{E_{\nu }}$ (cf. \cite{RZ}, 3.45-46); this descent
datum is effective on each $\mathcal{M}_{s}$ for $s$ as above. Then the
projective limit of the $\mathcal{M}_{s}$ is a completion of $\mathcal{%
\breve{M}}$ that comes as base-change from a pro-formal scheme $\mathcal{M}$
defined over $\limfunc{Spf}(\mathcal{O}_{E_{\nu }}):$%
\begin{equation*}
\mathcal{M=}\underset{\longleftarrow }{\lim }\mathcal{M}_{s}\text{ }/\text{ }%
\limfunc{Spf}(\mathcal{O}_{E_{\nu }}).
\end{equation*}

\noindent The action of $J(%
%TCIMACRO{\U{211a} }%
%BeginExpansion
\mathbb{Q}
%EndExpansion
_{p})$ on $\mathcal{\breve{M}}$ commutes with the descent datum, and then
gives an action of $J(%
%TCIMACRO{\U{211a} }%
%BeginExpansion
\mathbb{Q}
%EndExpansion
_{p})$ on $\mathcal{M}.$

\bigskip

From now on, we shall assume that the group $G$ is connected and satisfies
the Hasse principle, unless otherwise stated; this is because under this
assumption, the uniformization result has an easier formulation. Recall that
we fixed $[(A_{0},i_{0},\overline{\lambda }_{0},\overline{\alpha }_{0})]\in 
\mathcal{S}_{\mathcal{D},K^{p}}(\overline{\mathbb{F}}_{p})$; we also keep
the notation introduced at the beginning of the section.

\begin{definition}
We let $Z([(A_{0},i_{0},\overline{\lambda }_{0},\overline{\alpha }%
_{0})])\subseteq \mathcal{S}_{\mathcal{D},K^{p}}(\overline{\mathbb{F}}_{p})$
to be the set of points $[(A,i,\overline{\lambda },\overline{\alpha })]\in 
\mathcal{S}_{\mathcal{D},K^{p}}\left( \overline{\mathbb{F}}_{p}\right) $
such that there exists an isogeny $(A_{0},i_{0})\rightarrow (A,i)$ in $AV$
sending the $%
%TCIMACRO{\U{211a} }%
%BeginExpansion
\mathbb{Q}
%EndExpansion
$-homogeneous polarization $\overline{\lambda }_{0}$ into $\overline{\lambda 
}$ (no condition on the level structure).
\end{definition}

Notice that $Z([(A_{0},i_{0},\overline{\lambda }_{0},\overline{\alpha }%
_{0})])$ is also the set of points $[(A,i,\overline{\lambda },\overline{%
\alpha })]\in \mathcal{S}_{\mathcal{D},K^{p}}\left( \overline{\mathbb{F}}%
_{p}\right) $ whose underlying abelian scheme $(A,i,\overline{\lambda })$ is
quasi-isogenous to $(A_{0},i_{0},\overline{\lambda }_{0})$. By a result of
Rapoport and Richartz (cf. \cite{RR}, \cite{RZ} 6.26-27), one knows:

\begin{theorem}
\textbf{(Rapoport, Richartz) }If the $p$-divisible group $\mathbf{X}$ of $%
(A_{0},i_{0},\overline{\lambda }_{0})$ is basic, then $Z([(A_{0},i_{0},%
\overline{\lambda }_{0},\overline{\alpha }_{0})])$ is the set of $\overline{%
\mathbb{F}}_{p}$-valued points of a closed subset $Z$ of $\mathcal{S}_{%
\mathcal{D},K^{p}}\otimes \overline{\mathbb{F}}_{p}$.
\end{theorem}

We can now state:

\begin{theorem}
\label{RZ}\textbf{(Rapoport, Zink) }Let us fix $[(A_{0},i_{0},\overline{%
\lambda }_{0},\overline{\alpha }_{0})]\in \mathcal{S}_{\mathcal{D}%
,K^{p}}\left( \overline{\mathbb{F}}_{p}\right) $ such that the $p$-divisible
group $\mathbf{X}$ of $(A_{0},i_{0},\overline{\lambda }_{0})$ is basic;\
denote by $Z\subseteq \mathcal{S}_{\mathcal{D},K^{p}}$ the closed subspace
whose $\overline{\mathbb{F}}_{p}$-points are the points of $Z([(A_{0},i_{0},%
\overline{\lambda }_{0},\overline{\alpha }_{0})])\subseteq \mathcal{S}_{%
\mathcal{D},K^{p}}(\overline{\mathbb{F}}_{p}).$\ Let $\widehat{\mathcal{S}}_{%
\mathcal{D},C^{p}/Z}$ be the formal completion of the scheme $\mathcal{S}_{%
\mathcal{D},K^{p}}$ along $Z$. If the algebraic group $G$ associated to the
PEL-datum $\mathcal{D}$ is connected and satisfies the Hasse principle, then
there is a canonical isomorphism of formal schemes over $\limfunc{Spf}%
\mathcal{O}_{E_{\nu }}$:%
\begin{equation*}
\vartheta _{K^{p}}:I(%
%TCIMACRO{\U{211a} }%
%BeginExpansion
\mathbb{Q}
%EndExpansion
)\backslash \mathcal{M\times }G(\mathbb{A}_{f}^{p})/K^{p}\rightarrow 
\widehat{\mathcal{S}}_{\mathcal{D},K^{p}/Z},
\end{equation*}

where $I(%
%TCIMACRO{\U{211a} }%
%BeginExpansion
\mathbb{Q}
%EndExpansion
)$ acts on $\mathcal{M}$ via $\alpha _{0,p}$,\ and on $G(\mathbb{A}_{f}^{p})$
via $\alpha _{0}^{p}$. The system of morphisms $\{\vartheta
_{K^{p}}\}_{K^{p}}$ is equivariant with respect to the right Hecke $G(%
\mathbb{A}_{f}^{p})$-action on the projective systems of both sides above.
Furthermore, $I$ is an inner form of $G$, and we have canonical
identifications $I(\mathbb{A}_{f}^{p})=G(\mathbb{A}_{f}^{p})$, $J(%
%TCIMACRO{\U{211a} }%
%BeginExpansion
\mathbb{Q}
%EndExpansion
_{p})=I(%
%TCIMACRO{\U{211a} }%
%BeginExpansion
\mathbb{Q}
%EndExpansion
_{p})$; $I(%
%TCIMACRO{\U{211d} }%
%BeginExpansion
\mathbb{R}
%EndExpansion
)$ is compact modulo its center.
\end{theorem}

Recall that the right action of $G(\mathbb{A}_{f}^{p})$ on the projective
system $\{\widehat{\mathcal{S}}_{\mathcal{D},K^{p}/Z}\}_{K^{p}}$\ is defined
at \ref{Hecke}. In a similar way, notice that if $K_{1}^{p}\subseteq
K_{2}^{p}$ are open compact subgroups of $G(\mathbb{A}_{f}^{p})$, we have a
transition map $I(%
%TCIMACRO{\U{211a} }%
%BeginExpansion
\mathbb{Q}
%EndExpansion
)\backslash \mathcal{M\times }G(\mathbb{A}_{f}^{p})/K_{1}^{p}\rightarrow I(%
%TCIMACRO{\U{211a} }%
%BeginExpansion
\mathbb{Q}
%EndExpansion
)\backslash \mathcal{M\times }G(\mathbb{A}_{f}^{p})/K_{2}^{p}$; we obtain
correspondingly a projective system of formal schemes. If $g\in G(\mathbb{A}%
_{f}^{p})$, we define a morphism:%
\begin{equation*}
g:I(%
%TCIMACRO{\U{211a} }%
%BeginExpansion
\mathbb{Q}
%EndExpansion
)\backslash \mathcal{M\times }G(\mathbb{A}_{f}^{p})/K^{p}\rightarrow I(%
%TCIMACRO{\U{211a} }%
%BeginExpansion
\mathbb{Q}
%EndExpansion
)\backslash \mathcal{M\times }G(\mathbb{A}_{f}^{p})/g^{-1}K^{p}g
\end{equation*}%
induced by the map $x\cdot K^{p}\mapsto g^{-1}xg\cdot g^{-1}K^{p}g$ ($x\in G(%
\mathbb{A}_{f}^{p})$). The equivariance of the isomorphisms $\vartheta
_{K^{p}}$'s stated in the above theorem is with respect to these two right
actions of $G(\mathbb{A}_{f}^{p})$ (that we will call the action of the 
\textit{Hecke operators} of $G(\mathbb{A}_{f}^{p})$).

\noindent The statement about $I(\mathbb{A}_{f}^{p})$ and $I(%
%TCIMACRO{\U{211a} }%
%BeginExpansion
\mathbb{Q}
%EndExpansion
_{p})$ in the above theorem is a consequence of the basicity of $\mathbf{X}$
(cf. \cite{RZ}, 6.29).

To prove Theorem \ref{RZ}, Rapoport and Zink first define a morphism of
functors over $NILP_{\mathcal{O}_{\overset{\smile }{E_{\nu }}}}:$%
\begin{equation*}
\Theta _{K^{p}}:I(%
%TCIMACRO{\U{211a} }%
%BeginExpansion
\mathbb{Q}
%EndExpansion
)\backslash \mathcal{\breve{M}}\times G(\mathbb{A}_{f}^{p})/K^{p}\rightarrow 
\mathcal{S}_{\mathcal{D},K^{p}}\otimes _{\mathcal{O}_{E_{\nu }}}\mathcal{O}_{%
\overset{\smile }{E_{\nu }}},
\end{equation*}

\noindent and then they show the crucial:

\begin{proposition}
\label{key}Assume that $G$ is connected and satisfies the Hasse principle;
the morphism of functors $\Theta _{K^{p}}$ defines a canonical isomorphism:%
\begin{equation*}
\Theta _{K^{p}}(\mathcal{\overline{\mathbb{F}}}_{p}):I(%
%TCIMACRO{\U{211a} }%
%BeginExpansion
\mathbb{Q}
%EndExpansion
)\backslash \mathcal{\breve{M}(\overline{\mathbb{F}}}_{p}\mathcal{)\times }G(%
\mathbb{A}_{f}^{p})/K^{p}\rightarrow Z(\mathcal{\overline{\mathbb{F}}}_{p}).
\end{equation*}
\end{proposition}

Both the domain and the codomain of $\Theta _{K^{p}}$ have a Weyl descent
datum to $\mathcal{O}_{E_{\nu }}$; one can prove that $\Theta _{K^{p}}$ is
compatible with such a descent, and then deduce Theorem \ref{RZ} from the
above Proposition. The morphism $\Theta _{K^{p}}$ is called the
uniformization morphism in \cite{RZ}, 6.15.

(In some special situations, $Z$ coincides with the whole special fiber of
the scheme $\mathcal{S}_{\mathcal{D},K^{p}}$ over $%
%TCIMACRO{\U{2124} }%
%BeginExpansion
\mathbb{Z}
%EndExpansion
_{p}$; in this cases, the result of Rapoport and Zink gives a $p$-adic
uniformization of the whole Shimura variety at $p$. This happens, for
example, in the case considered by Cherednik in \cite{Ch}, where $B$ is a
rational quaternion algebra unramified at infinity, $G=B^{\times }$ and $p$
a prime at which $B$ is ramified).

\subsubsection{Description of the uniformization morphism over $\overline{%
\mathbb{F}}_{p}$\label{uniformization}}

We keep the assumptions of the previous sections. In particular, we have
fixed an element $[(A_{0},i_{0},\overline{\lambda }_{0},\overline{\alpha }%
_{0})]\in \mathcal{S}_{\mathcal{D},K^{p}}\left( \overline{\mathbb{F}}%
_{p}\right) $ that gives rise to the $p$-divisible group $(\mathbf{X},i_{%
\mathbf{X}},\overline{\lambda }_{\mathbf{X}})$ over $\overline{\mathbb{F}}%
_{p}$\ and then to the isocrystal $((N,\mathbf{F}),i,%
%TCIMACRO{\U{211a} }%
%BeginExpansion
\mathbb{Q}
%EndExpansion
_{p}^{\times }\Psi )$ over $K_{0}$; we also have fixed an isomorphism $%
N\simeq V\otimes _{%
%TCIMACRO{\U{211a} }%
%BeginExpansion
\mathbb{Q}
%EndExpansion
}K_{0}$ of $B\otimes _{%
%TCIMACRO{\U{211a} }%
%BeginExpansion
\mathbb{Q}
%EndExpansion
}K_{0}$-modules that respects the skew-symmetric forms on both sides, and
such that $\mathbf{F}=b\otimes \sigma $. We still assume that $\mathbf{X}$
is basic and that $G$ is connected and satisfies the Hasse principle.

We recall the definition of the uniformization morphism of functors over $%
NILP_{\mathcal{O}_{\overset{\smile }{E_{\nu }}}}$:%
\begin{equation*}
\Theta _{K^{p}}:I(%
%TCIMACRO{\U{211a} }%
%BeginExpansion
\mathbb{Q}
%EndExpansion
)\backslash \mathcal{\breve{M}}(S)\times G(\mathbb{A}_{f}^{p})/K^{p}%
\rightarrow \mathcal{S}_{\mathcal{D},K^{p}}\otimes _{\mathcal{O}_{E_{\nu }}}%
\mathcal{O}_{\overset{\smile }{E_{\nu }}}
\end{equation*}

\noindent that appears in Proposition \ref{key}. We need:

\begin{lemma}
\label{push-forward}Let $S\in NILP_{\mathcal{%
%TCIMACRO{\U{2124} }%
%BeginExpansion
\mathbb{Z}
%EndExpansion
}_{p}}$, and let $A^{\prime }$ be an object in $AV_{\mathcal{O}_{B}/S}$ ;
denote by $X^{\prime }=A^{\prime }(p)$ the $p$-divisible group of $A^{\prime
}$. For any quasi-isogeny $\xi :X^{\prime }\rightarrow X^{\prime \prime }$
of $p$-divisible groups over $S$ that respects the $\mathcal{O}_{B_{p}}$%
-action, there exists an element $A^{\prime \prime }$ of $AV_{\mathcal{O}%
_{B}/S}$ whose $p$-divisible group is $X^{\prime \prime }$ and a
quasi-isogeny $\xi :A^{\prime }\rightarrow A^{\prime \prime }$ of $AV_{%
\mathcal{O}_{B}/S}$ inducing $\xi :X^{\prime }\rightarrow X^{\prime \prime }$%
. Furthermore the arrow $\xi :A^{\prime }\rightarrow A^{\prime \prime }$ in $%
AV_{\mathcal{O}_{B}/S}$ is uniquely determined; we denote $A^{\prime \prime
} $ by $\xi _{\ast }A^{\prime }$. This construction is functorial, i.e. $%
\left( \xi _{2}\xi _{1}\right) _{\ast }A^{\prime }=\xi _{2\ast }\left( \xi
_{1\ast }A^{\prime }\right) $.
\end{lemma}

Under the hypothesis of the above Lemma, if $A^{\prime }$ comes with a
polarization $\lambda :A^{\prime }\rightarrow \widehat{A}^{\prime }$, then $%
\xi $ defines a polarization $\xi _{\ast }\lambda :=\left( \xi ^{-1}\right)
^{\wedge }\lambda \xi ^{-1}$ on $A^{\prime \prime }$.\ If furthermore $%
A^{\prime }$ comes with a rigidification $\alpha :H_{1}(A^{\prime },\mathbb{A%
}_{f}^{p})\rightarrow V\otimes _{%
%TCIMACRO{\U{211a} }%
%BeginExpansion
\mathbb{Q}
%EndExpansion
}\mathbb{A}_{f}^{p}$ (i.e. a symplectic $\mathcal{O}$-equivariant
isomorphism), then $A^{\prime \prime }$ comes with the rigidification $\xi
_{\ast }\alpha :=\alpha \circ H_{1}(\xi ^{-1},\mathbb{A}_{f}^{p}).$

Let $S$ be a fixed scheme in $NILP_{\mathcal{O}_{\overset{\smile }{E_{\nu }}%
}}$; we denote by $(\widetilde{\mathbf{X}},\widetilde{i}_{\mathbf{X}},%
\widetilde{\overline{\lambda }}_{\mathbf{X}})$ a \textit{fixed} lifting of $(%
\mathbf{X},i_{\mathbf{X}},\overline{\lambda }_{\mathbf{X}})$\ to $\mathcal{O}%
_{\overset{\smile }{E_{\nu }}}$, and we let $(\widetilde{A}_{0},\widetilde{i}%
_{0},\widetilde{\overline{\lambda }}_{0},\widetilde{\overline{\alpha }}_{0})$
be the corresponding lifting of $(A_{0},i_{0},\overline{\lambda }_{0},%
\overline{\alpha }_{0})$ over $\mathcal{O}_{\overset{\smile }{E_{\nu }}}$%
(cf. Prop. \ref{lifting with separability}). We can consider base changes of
these objects to $S$ and to $\overline{S}$ (which is the closed subscheme of 
$S$ defined by the ideal sheaf $p\mathcal{O}_{S}$): we denote it by the
subscripts $\cdot _{S}$ and $\cdot _{\overline{S}}$ respectively.

Consider a $p$-divisible group with additional structure $[(X,i,\overline{%
\lambda };\rho )]\in \mathcal{\breve{M}}(S)$; $\rho :(\mathbf{X},i_{\mathbf{X%
}})_{\overline{S}}\rightarrow (X,i)_{\overline{S}}$ is an $\overline{S}$%
-quasi-isogeny of $p$-divisible groups with $\mathcal{O}_{B_{p}}$-action
such that $\widehat{\rho }\circ \lambda _{X}\circ \rho \in 
%TCIMACRO{\U{211a} }%
%BeginExpansion
\mathbb{Q}
%EndExpansion
_{p}^{\times }\lambda _{\mathbf{X}}$; by Prop. \ref{Rigidity of qisogenies},$%
\ \rho $ lifts uniquely to a quasi-isogeny $\widetilde{\rho }:(\widetilde{%
\mathbf{X}}_{S},\widetilde{i}_{\mathbf{X}})\rightarrow (X,i)$\ of $p$%
-divisible groups over $S$ with an action of $\mathcal{O}_{B_{p}}$. By the
above lemma, we obtain therefore an abelian scheme $\widetilde{\rho }_{\ast
}(\widetilde{A}_{0/S})$ over $S$ endowed with an action $\widetilde{\rho }%
_{\ast }(\widetilde{i}_{0})$ of $\mathcal{O}_{B}$, a polarization $%
\widetilde{\rho }_{\ast }(\widetilde{\overline{\lambda }}_{0})$ and a level
structure $\widetilde{\rho }_{\ast }(\widetilde{\overline{\alpha }}_{0})$,
such that:%
\begin{equation*}
\lbrack (\widetilde{\rho }_{\ast }(\widetilde{A}_{0/S}),\widetilde{\rho }%
_{\ast }(\widetilde{i}_{0}),\widetilde{\rho }_{\ast }(\widetilde{\overline{%
\lambda }}_{0}),\widetilde{\rho }_{\ast }(\widetilde{\overline{\alpha }}%
_{0}))]\in \mathcal{S}_{\mathcal{D},K^{p}}\otimes _{\mathcal{O}_{E_{\nu }}}%
\mathcal{O}_{\overset{\smile }{E_{\nu }}}\left( S\right) .
\end{equation*}%
We define a morphism of functors $\Theta _{K^{p}}$ over $NILP_{\mathcal{O}_{%
\overset{\smile }{E_{\nu }}}}$ by letting, for any $S\in \limfunc{Obj}(NILP_{%
\mathcal{O}_{\overset{\smile }{E_{\nu }}}})$: 
\begin{eqnarray*}
\Theta _{K^{p}}(S) &:&\mathcal{\breve{M}}(S)\times G(\mathbb{A}%
_{f}^{p})/K^{p}\rightarrow \mathcal{S}_{\mathcal{D},K^{p}}(S), \\
\lbrack (X,i,\overline{\lambda };\rho )]\times gK^{p} &\longmapsto &[(%
\widetilde{\rho }_{\ast }(\widetilde{A}_{0/S}),\widetilde{\rho }_{\ast }(%
\widetilde{i}_{0}),\widetilde{\rho }_{\ast }(\widetilde{\overline{\lambda }}%
_{0}),g^{-1}\cdot \widetilde{\rho }_{\ast }(\widetilde{\overline{\alpha }}%
_{0}))].
\end{eqnarray*}

\bigskip

(We need to explain better the meaning of the notation $g^{-1}\cdot 
\widetilde{\rho }_{\ast }(\widetilde{\overline{\alpha }}_{0})$. More
properly, we assumed \underline{fixed} a representative $\alpha _{0}$ for $%
\overline{\alpha }_{0}$, and $g^{-1}\cdot \widetilde{\rho }_{\ast }(%
\widetilde{\overline{\alpha }}_{0})$ denotes the $K^{p}$-class of the
isomorphism $g^{-1}\circ \widetilde{\alpha }_{0}\circ H_{1}(\widetilde{\rho }%
^{-1})$; this class does not depend on the choice of representative for the
coset $gK^{p}$. In the sequel, we will use without any further explanation
the notation just introduced).

We shall use the shorter notation:%
\begin{equation*}
\Theta _{K^{p}}(S):[(X,i,\overline{\lambda };\rho )]\times gK^{p}\longmapsto
\lbrack (\widetilde{\rho }_{\ast }(\widetilde{A}_{0/S},\widetilde{i}_{0},%
\widetilde{\overline{\lambda }}_{0}),g^{-1}\cdot \widetilde{\rho }_{\ast }(%
\widetilde{\overline{\alpha }}_{0}))]
\end{equation*}

\noindent Notice that $\Theta _{K^{p}}(S)$ is well defined (all the choices
we made above give the same tuples, up to isomorphism). Furthermore, $%
\{\Theta _{K^{p}}\}_{K^{p}}$ is equivariant with respect to the right $G(%
\mathbb{A}_{f}^{p})$-action on the projective systems of both sides above.

\noindent We now define a left action of $I(%
%TCIMACRO{\U{211a} }%
%BeginExpansion
\mathbb{Q}
%EndExpansion
)$ on the left hand side. Fix $S\in \limfunc{Obj}(NILP_{\mathcal{O}_{\overset%
{\smile }{E_{\nu }}}})$, $[(X,i,\overline{\lambda };\rho )]\in $ $\mathcal{%
\breve{M}}(S)$, $g\in G(\mathbb{A}_{f}^{p})$ and $\xi \in I(%
%TCIMACRO{\U{211a} }%
%BeginExpansion
\mathbb{Q}
%EndExpansion
)$; we set:%
\begin{equation*}
\xi \cdot \left( \lbrack (X,i,\overline{\lambda };\rho )]\times
gK^{p}\right) :=[(X,i,\overline{\lambda };\rho \circ \alpha _{0,p}(\xi
^{-1})_{\overline{S}})]\times \alpha _{0}^{p}(\xi )gK^{p}\text{,}
\end{equation*}

\noindent where $\alpha _{0,p}:I(%
%TCIMACRO{\U{211a} }%
%BeginExpansion
\mathbb{Q}
%EndExpansion
)\rightarrow J(%
%TCIMACRO{\U{211a} }%
%BeginExpansion
\mathbb{Q}
%EndExpansion
_{p})$ and $\alpha _{0}^{p}:I(%
%TCIMACRO{\U{211a} }%
%BeginExpansion
\mathbb{Q}
%EndExpansion
)\rightarrow G(\mathbb{A}_{f}^{p})$ are the homomorphisms defined in \ref%
{p-adic unif of isogeny classes}. It is easy to check that $\Theta _{K^{p}}$
is invariant under the left action of $I(%
%TCIMACRO{\U{211a} }%
%BeginExpansion
\mathbb{Q}
%EndExpansion
)$ just defined: in fact since $\alpha _{0}^{p}(\xi ^{-1}):(\mathbf{X},i_{%
\mathbf{X}},\overline{\lambda }_{\mathbf{X}})\rightarrow (\mathbf{X},i_{%
\mathbf{X}},\overline{\lambda }_{\mathbf{X}})$ is a quasi-isogeny coming
from $\xi ^{-1}:(A_{0},i_{0},\overline{\lambda }_{0})\rightarrow
(A_{0},i_{0},\overline{\lambda }_{0})$, we have canonical identifications 
\begin{eqnarray*}
\alpha _{0}^{p}(\xi ^{-1})_{\ast }(A_{0},i_{0},\overline{\lambda }_{0})
&=&(A_{0},i_{0},\overline{\lambda }_{0}), \\
\left( \rho \circ \alpha _{0,p}(\xi ^{-1})_{\overline{S}}\right) _{\ast
}^{\sim }(\widetilde{A}_{0},\widetilde{i}_{0},\widetilde{\overline{\lambda }}%
_{0}) &=&\widetilde{\rho }_{\ast }(\widetilde{A}_{0},\widetilde{i}_{0},%
\widetilde{\overline{\lambda }}_{0}).
\end{eqnarray*}
Under these identifications, the level structures associated to $\xi \cdot
\left( \lbrack (X,i,\overline{\lambda };\rho )]\times gK^{p}\right) $ and $%
[(X,i,\overline{\lambda };\rho )]\times gK^{p}$ coincide.

We call the \textbf{uniformization morphism} 
\begin{equation*}
\Theta _{K^{p}}:I(%
%TCIMACRO{\U{211a} }%
%BeginExpansion
\mathbb{Q}
%EndExpansion
)\backslash \mathcal{\breve{M}}\times G(\mathbb{A}_{f}^{p})/C^{p}\rightarrow 
\mathcal{S}_{\mathcal{D},K^{p}}\otimes _{\mathcal{O}_{E_{\nu }}}\mathcal{O}_{%
\overset{\smile }{E_{\nu }}}
\end{equation*}%
the map induced by $\Theta _{K^{p}}$ defined above modulo the left action of 
$I(%
%TCIMACRO{\U{211a} }%
%BeginExpansion
\mathbb{Q}
%EndExpansion
)$ on $\mathcal{\breve{M}}\times G(\mathbb{A}_{f}^{p})/K^{p}$.

Let us now consider the map on geometric points $\Theta _{K^{p}}(\mathcal{%
\overline{\mathbb{F}}}_{p})$: it is clear that its image coincide with $Z(%
\mathcal{\overline{\mathbb{F}}}_{p})$. Furthermore $\Theta _{K^{p}}(\mathcal{%
\overline{\mathbb{F}}}_{p})$ is injective: assume that we start with two
points $[(X_{j},i_{j},\overline{\lambda }_{j};\rho _{j})]\in $ $\mathcal{%
\breve{M}(\overline{\mathbb{F}}}_{p}\mathcal{)}$ ($j=1,2$) in ; we obtain
correspondingly quasi-isogenies of abelian varieties $\rho
_{j,ab}:A_{0}\rightarrow A_{j}:=\rho _{j,\ast }(A_{0})$ ($j=1,2$). If $%
f:A_{1}\rightarrow A_{2}$ is an isomorphism of abelian varieties with
additional structure, then $\xi :=\rho _{2,ab}^{-1}\circ f\circ \rho
_{1,ab}\in I(%
%TCIMACRO{\U{211a} }%
%BeginExpansion
\mathbb{Q}
%EndExpansion
)$ induces an element $\alpha _{0,p}(\xi )\in J(%
%TCIMACRO{\U{211a} }%
%BeginExpansion
\mathbb{Q}
%EndExpansion
_{p})$ such that $f\circ \rho _{1}=\rho _{2}\circ \alpha _{0,p}(\xi )$, so
that $(X_{1},i_{1},\overline{\lambda }_{1};\rho _{1})\simeq (X_{2},i_{2},%
\overline{\lambda }_{2};\rho _{2})$. Keeping track of the level structures,
we conclude the injectivity of $\Theta _{K^{p}}(\mathcal{\overline{\mathbb{F}%
}}_{p})$. This proves Proposition \ref{key}.

\subsection{Restriction of the moduli problems to the superspecial locus}

We recall the notion of superspecial and supersingular abelian varieties
over $\overline{\mathbb{F}}_{p}$; then we apply the results of the previous
paragraph to the superspecial locus of our moduli scheme.

\subsubsection{Supersingular and superspecial abelian varieties\label%
{superspecial-supersingular}}

Let $A$ be an abelian variety of dimension $g\geq 1$ over the field $%
\overline{\mathbb{F}}_{p}$; denote as usual by $A(p)$ its Barsotti-Tate
group: it has dimension $g$, height $2g$ and $\#A[p](\overline{\mathbb{F}}%
_{p})\leq $ $p^{g}$. We have $\widehat{A(p)}\simeq \widehat{A}(p)$, where $%
\widehat{A}$ denotes the dual abelian variety, hence $A(p)$ is isogenous to
its Serre dual. This implies that if $0\leq \lambda _{1}\leq ...\leq \lambda
_{2g}\leq 1$ is the slope sequence of $A(p)$, then $\lambda _{i}+\lambda
_{2g-i}=1$ for all $1\leq i\leq 2g$. In particular, if $g=1$ the only two
possibilities for the isomorphism class of $A(p)$ are $G_{0}\oplus G_{1}$
and $G_{1/2}$; we say that $A$ is ordinary in the first case, otherwise $A$
is said to be supersingular. In higher dimension the situation is richer.

An abelian variety $A$ over $\overline{\mathbb{F}}_{p}$ is said to be 
\textit{supersingular} if $A(p)$ is isogenous to $G_{1/2}^{\oplus g}$ (over $%
\overline{\mathbb{F}}_{p}$); it is said to be \textit{superspecial} if $A(p)$
is isomorphic to $G_{1/2}^{\oplus g}$ (over $\overline{\mathbb{F}}_{p}$). An
abelian variety $A^{\prime }$ over a finite extension of $\mathbb{F}_{p}$ is
said to be supersingular (resp. superspecial) if its base change to $%
\overline{\mathbb{F}}_{p}$ is supersingular (resp. superspecial).

\bigskip

If $g=1$, there is no difference between supersingular and superspecial
abelian varieties. We recall some properties of supersingular elliptic
curves (cf.:\ \cite{Sil}, Ch.\ V; \cite{T}; \cite{G}, 2.2.).

An elliptic curve $E$ defined over $\overline{\mathbb{F}}_{p}$ is
supersingular if and only if its endomorphism ring $\limfunc{End}E:=\limfunc{%
End}_{\overline{\mathbb{F}}_{p}}E$ is isomorphic to a maximal order in the
quaternion $%
%TCIMACRO{\U{211a} }%
%BeginExpansion
\mathbb{Q}
%EndExpansion
$-algebra ramified exactly at the places $\{p,\infty \}$; if $E/\overline{%
\mathbb{F}}_{p}$ is supersingular, then there exists a unique (up to $%
\mathbb{F}_{p^{2}}$-isomorphism) elliptic curve $E^{\prime }$ defined over $%
\mathbb{F}_{p^{2}}$, such that $E=E^{\prime }\otimes _{\mathbb{F}_{p^{2}}}%
\overline{\mathbb{F}}_{p}$ and such that the geometric Frobenius $E^{\prime
}\rightarrow E^{\prime (p^{2})}=E^{\prime }$ equals $[-p]$. Furthermore, $%
\limfunc{End}_{\mathbb{F}_{p^{2}}}E^{\prime }=\limfunc{End}_{\overline{%
\mathbb{F}}_{p}}E$ and the association $E\mapsto E^{\prime }$ is \textit{%
functorial}. The cotangent space $\omega (E)$ of $E$ has a canonical $%
\mathbb{F}_{p^{2}}$-structure. If $E$ is a supersingular elliptic curve over 
$\overline{\mathbb{F}}_{p}$ and $E^{\prime }$ is its canonical model over $%
\mathbb{F}_{p^{2}}$, then $E^{\prime }(p)$ is a canonical model of $E(p)$
whose covariant Dieduonn\'{e} module $A_{1/2}^{\prime }:=M_{\ast }(E^{\prime
}(p))$ (over $W(\mathbb{F}_{p^{2}})$)\ is isomorphic to:%
\begin{equation*}
\left( W(\mathbb{F}_{p^{2}})^{2};F=%
\begin{pmatrix}
0 & 1 \\ 
-p & 0%
\end{pmatrix}%
\sigma ,V=%
\begin{pmatrix}
0 & -1 \\ 
p & 0%
\end{pmatrix}%
\sigma ^{-1}\right) ,
\end{equation*}%
where $\sigma $ denotes the Frobenius morphism of $W(\mathbb{F}_{p^{2}}).$
\noindent Furthermore, $A_{1/2}:=M_{\ast }(E(p))\simeq \frac{W(\overline{%
\mathbb{F}}_{p})[F,V]}{W(\overline{\mathbb{F}}_{p})[F,V](F-V)}$ (cf. Th. \ref%
{Dieudonne}), and the ring of Dieudonn\'{e} module endomorphisms $\limfunc{%
End}A_{1/2}=\limfunc{End}A_{1/2}^{\prime }$ is isomorphic to the maximal
order in the quaternion division algebra over $%
%TCIMACRO{\U{211a} }%
%BeginExpansion
\mathbb{Q}
%EndExpansion
_{p}$ (cf. \cite{Gh04}, Corollary 7).

There are finitely many isomorphism classes of supersingular elliptic curves
over $\overline{\mathbb{F}}_{p}$, and they all are isogenous. It is a result
of Deligne and Ogus that if $g\geq 2$ and $E_{1},...,E_{2g}$ are
supersingular elliptic curves over $\overline{\mathbb{F}}_{p}$, then we have
an isomorphism over $\overline{\mathbb{F}}_{p}$: $E\times ...\times
E_{g}\simeq E_{g+1}\times ...\times E_{2g}.$For references on the proof of
this and the following results, cf. \cite{LO}, 1.6.:

\begin{proposition}
\label{Superspecial}Let $A$ be an abelian variety over $\overline{\mathbb{F}}%
_{p}$ of dimension $g\geq 2$. Then the following are equivalent:

\begin{enumerate}
\item $A$ is supersingular;

\item all $2g$ slopes of the Newton polygon of $A$ are equal to $1/2$;

\item $A$ is isogenous over $\overline{\mathbb{F}}_{p}$ to $E^{g}$ for some
(any) supersingular elliptic curve $E$ over $\overline{\mathbb{F}}_{p}$.
\end{enumerate}

\noindent Furthermore, the following are equivalent:

\begin{enumerate}
\item $A$ is superspecial;

\item $M_{\ast }(A(p))\simeq A_{1/2}^{\oplus g}$ as Dieduonn\'{e} $W(%
\overline{\mathbb{F}}_{p})$-modules;

\item $A$ is isomorphic over $\overline{\mathbb{F}}_{p}$ to $E^{g}$ for some
(any) supersingular elliptic curve $E$ over $\overline{\mathbb{F}}_{p}.$
\end{enumerate}
\end{proposition}

If $A$ is a superspecial abelian variety over $\overline{\mathbb{F}}_{p}$ of
dimension $g\geq 2$, then we see that $A$ has a canonical model $A^{\prime }$
over $\mathbb{F}_{p^{2}}$, in which the geometric Frobenius equals $\left[ -p%
\right] $. Furthermore the association $A\mapsto A^{\prime }$ is functorial.
If $A=E^{g}$ is a superspecial abelian variety over $\overline{\mathbb{F}}%
_{p}$ of dimension $g\geq 1$, then $A$ comes with a canonical principal
polarization $A\rightarrow \widehat{A}$ induced from the canonical
polarization of the elliptic curve $E$. For this reason, in this case we
will identify $A$ and $\widehat{A}$. As a consequence, $M(A(p))$ and $%
M_{\ast }(A(p))$ are canonically isomorphic as Dieudonn\'{e} modules.

\subsubsection{A variant of the moduli problem for $p$-divisible groups\label%
{exp}}

Let us fix integers $g\geq 1$, $N\geq 3$ and a prime number $p$ not dividing 
$N$; let us denote by $\mathcal{A}_{g,1,N}$ the Siegel moduli scheme
associated to the PEL-datum $\mathcal{D}_{2g,p}^{Sp(%
%TCIMACRO{\U{211a} }%
%BeginExpansion
\mathbb{Q}
%EndExpansion
)}$ having good reduction at $p$, and to the choice of principal level $U(N)$
(cf. \ref{Siegel}). If $g=1$, $\mathcal{A}_{1,1,N}(\overline{\mathbb{F}}%
_{p}) $ contains a finite number of supersingular elliptic curves, which
form an isogeny class; if $g>1$, the supersingular abelian varieties living
in $\mathcal{A}_{g,1,N}(\overline{\mathbb{F}}_{p})$ define a closed subset
of positive dimension (cf. \cite{LO}, 4.9), hence "too big" for our
purposes. For this reason, in \cite{Gh04}, the author needs to consider
superspecial abelian varieties, instead of supersingular varieties, in order
to construct a map from geometric to algebraic eigenforms of Siegel type.
The above situation also occurs in other cases of PEL-type.

\bigskip

Fix a prime number $p$ and denote by $W=W\left( \overline{\mathbb{F}}%
_{p}\right) $ the ring of Witt vectors of $\overline{\mathbb{F}}_{p}$ and by 
$K_{0}$ its fraction field, endowed with the Frobenius automorphism $\sigma $%
. We assume fixed a simple $%
%TCIMACRO{\U{211a} }%
%BeginExpansion
\mathbb{Q}
%EndExpansion
_{p}$-PEL-datum with good reduction at $p$ for moduli of $p$-divisible
groups over $\overline{\mathbb{F}}_{p}:$ $\mathcal{D}_{\func{mod}}\mathcal{=}%
(B_{p},^{\ast },V_{p},\left\langle ,\right\rangle ,\mathcal{O}%
_{B_{p}},\Lambda ,b,\mu )$. We denote by $K$ the finite extension of $K_{0}$
over which the co-character $\mu $ is defined, by $E_{p}$ the Shimura field
of the datum; we let $\breve{E}_{p}:=E_{p}K_{0}$.

The above data define an isocrystal $(N,\mathbf{F}):=\left( V_{p}\otimes _{%
%TCIMACRO{\U{211a} }%
%BeginExpansion
\mathbb{Q}
%EndExpansion
_{p}}K_{0},b\sigma \right) $ endowed with an action $i$ of $B_{p}$ and a
skew-hermitian non-degenerate form of isocrystals $\Psi :N\times
N\rightarrow \mathbf{1}(1)$. By our assumptions, this isocrystal comes - via
the covariant Dieudonn\'{e}\ functor - from some $p$-divisible group over $%
\overline{\mathbb{F}}_{p}$ - endowed with action of $B_{p}$ and polarization
- that is uniquely determined only up to isogeny. We fix a \textit{choice of 
\textit{isomorphism} class} $(\mathbf{X},i_{\mathbf{X}},\overline{\lambda }_{%
\mathbf{X}})$ of polarized $p$-divisible group over $\overline{\mathbb{F}}%
_{p}$\ (with Shimura field $E_{p}$) associated to $(N,i,%
%TCIMACRO{\U{211a} }%
%BeginExpansion
\mathbb{Q}
%EndExpansion
_{p}^{\times }\Psi )$. We furthermore \textit{assume that }$\lambda _{%
\mathbf{X}}$\textit{\ is a principal polarization} and that $i_{\mathbf{X}}$
comes from an action of $\mathcal{O}_{B_{p}}$ on $\mathbf{X}$\ (this will
always be automatically true in our applications).

\begin{definition}
\label{modified functor}Let us fix $(\mathbf{X},i_{\mathbf{X}},\overline{%
\lambda }_{\mathbf{X}})$ as above.\ The set: 
\begin{equation*}
\mathcal{M}^{\prime }(\overline{\mathbb{F}}_{p}):=\mathcal{M}_{(\mathbf{X}%
,i_{\mathbf{X}},\overline{\lambda }_{\mathbf{X}})}^{\prime }(\overline{%
\mathbb{F}}_{p})
\end{equation*}%
is the collection of equivalence classes of quasi-isogenies $\rho :(\mathbf{X%
},i_{\mathbf{X}},\overline{\lambda }_{\mathbf{X}})\rightarrow (\mathbf{X},i_{%
\mathbf{X}},\overline{\lambda }_{\mathbf{X}})$ of the $p$-divisible group $%
\mathbf{X}$ over $\overline{\mathbb{F}}_{p}$ that respect the $\mathcal{O}%
_{B_{p}}$-structure and such that $\widehat{\rho }\circ \lambda _{\mathbf{X}%
}\circ \rho \in 
%TCIMACRO{\U{211a} }%
%BeginExpansion
\mathbb{Q}
%EndExpansion
_{p}^{\times }\lambda _{\mathbf{X}}$. Two quasi-isogenies $\rho $ and $\rho
^{\prime }$ are said to be equivalent if the $\overline{\mathbb{F}}_{p}$%
-quasi-isogeny $f:=\rho ^{\prime }\circ \rho ^{-1}$ is an isomorphism $(%
\mathbf{X},i_{\mathbf{X}},\overline{\lambda }_{\mathbf{X}})\rightarrow (%
\mathbf{X},i_{\mathbf{X}},\overline{\lambda }_{\mathbf{X}})$ of $p$%
-divisible groups over $\overline{\mathbb{F}}_{p}$ with $\mathcal{O}_{B_{p}}$%
-action, such that $\widehat{f}\circ \lambda _{\mathbf{X}}\circ f\in 
%TCIMACRO{\U{2124} }%
%BeginExpansion
\mathbb{Z}
%EndExpansion
_{p}^{\times }\lambda _{\mathbf{X}}$.
\end{definition}

The set $\mathcal{M}^{\prime }(\overline{\mathbb{F}}_{p})$ is non-empty
since $[(\mathbf{X},i_{\mathbf{X}},\overline{\lambda }_{\mathbf{X}})]\in 
\mathcal{M}^{\prime }(\overline{\mathbb{F}}_{p})$. Furthermore $\mathcal{M}%
^{\prime }(\overline{\mathbb{F}}_{p})\subseteq \mathcal{\breve{M}}(\overline{%
\mathbb{F}}_{p})$ is closed. Notice also that in the definition of $\mathcal{%
M}^{\prime }(\overline{\mathbb{F}}_{p})$ we can forget about the determinant
condition that appears in the definition of $\mathcal{\breve{M}}$, since it
is automatically satisfied.

Let $J(%
%TCIMACRO{\U{211a} }%
%BeginExpansion
\mathbb{Q}
%EndExpansion
_{p})$ denote the group of quasi-isogenies $\rho :(\mathbf{X},i_{\mathbf{X}})%
\mathbf{\rightarrow }(\mathbf{X},i_{\mathbf{X}})$ over $\overline{\mathbb{F}}%
_{p}$\ such that $\widehat{\rho }\circ \lambda _{\mathbf{X}}\circ \rho \in 
%TCIMACRO{\U{211a} }%
%BeginExpansion
\mathbb{Q}
%EndExpansion
_{p}^{\times }\lambda _{\mathbf{X}}$, as in \ref{groups I and J}, and let $J(%
%TCIMACRO{\U{2124} }%
%BeginExpansion
\mathbb{Z}
%EndExpansion
_{p})$ be the subgroup of isomorphisms $(\mathbf{X},i_{\mathbf{X}})\mathbf{%
\rightarrow }(\mathbf{X},i_{\mathbf{X}})$ preserving the polarization form
up to a factor in $%
%TCIMACRO{\U{2124} }%
%BeginExpansion
\mathbb{Z}
%EndExpansion
_{p}^{\times }$. Although the space $\mathcal{\breve{M}}(\overline{\mathbb{F}%
}_{p})$ is somehow mysterious, we have a better understanding of $\mathcal{M}%
^{\prime }(\overline{\mathbb{F}}_{p})$:

\begin{proposition}
\label{obvious}There is a natural bijection $\mathcal{M}^{\prime }(\overline{%
\mathbb{F}}_{p})\simeq J(%
%TCIMACRO{\U{211a} }%
%BeginExpansion
\mathbb{Q}
%EndExpansion
_{p})/J(%
%TCIMACRO{\U{2124} }%
%BeginExpansion
\mathbb{Z}
%EndExpansion
_{p})$.
\end{proposition}

\textbf{Proof. }Clear from the definition of $\mathcal{M}^{\prime }(%
\overline{\mathbb{F}}_{p})$:\ the bijection associates to $\left[ \rho %
\right] \in \mathcal{M}^{\prime }(\overline{\mathbb{F}}_{p})$ the left coset 
$\rho ^{-1}J(%
%TCIMACRO{\U{2124} }%
%BeginExpansion
\mathbb{Z}
%EndExpansion
_{p})$, where we view $\rho ^{-1}$ as an element of $J(%
%TCIMACRO{\U{211a} }%
%BeginExpansion
\mathbb{Q}
%EndExpansion
_{p})$.\ $\blacksquare $

\subsubsection{Uniformization of the superspecial locus}

We now apply the above results to the situation we are interested in. For
convenience we recall some notation already introduced. Fix $\mathcal{D=(}%
B,^{\ast },V,\left\langle ,\right\rangle ,\mathcal{O}_{B},\Lambda
,h,K^{p},\nu \mathcal{)}$ a simple $%
%TCIMACRO{\U{211a} }%
%BeginExpansion
\mathbb{Q}
%EndExpansion
$-PEL-datum for moduli of abelian schemes with good reduction at $p$, and
neat level $K^{p}$. Denote by $G$ the associated algebraic group, and assume
it is connected and satisfies the Hasse principle. Let $E$ be the Shimura
field of $\mathcal{D}$ and $E_{\nu }$ the completion of $E$ at $\nu $. Let $%
\overline{\mathbb{F}}_{p}$ be a fixed algebraic closure of the residue field
of $E_{\nu }$, and $W=W(\overline{\mathbb{F}}_{p})$, $K_{0}=W[\frac{1}{p}]$
and $\sigma $ the Frobenius morphism of $W$; fix a finite extension $K$ of $%
K_{0}$ such that $\mu $ is defined over $K$; set $\breve{E}_{\nu }=E_{\nu
}K_{0}$. The objects $B_{p}$, $V_{p}$, $\left\langle ,\right\rangle _{p}$, $%
G_{p}$, $\mathcal{O}_{B_{p}}$ are defined as before.

Let $\mathcal{S}_{\mathcal{D},K^{p}}$ be the quasi-projective smooth scheme
over $\mathcal{O}_{E}\otimes _{%
%TCIMACRO{\U{2124} }%
%BeginExpansion
\mathbb{Z}
%EndExpansion
}%
%TCIMACRO{\U{2124} }%
%BeginExpansion
\mathbb{Z}
%EndExpansion
_{(p)}$ representing $\mathbf{M}(\mathcal{D})$; we see $\mathcal{S}_{%
\mathcal{D},K^{p}}$ as a scheme over $\mathcal{O}_{\breve{E}_{\nu }}$.
Suppose that the common dimension of the abelian schemes parametrized by $%
\mathcal{S}_{\mathcal{D},K^{p}}$ is $g:=\dim _{%
%TCIMACRO{\U{2102} }%
%BeginExpansion
\mathbb{C}
%EndExpansion
}V_{%
%TCIMACRO{\U{2102} }%
%BeginExpansion
\mathbb{C}
%EndExpansion
,0}\geq 2$; fix a supersingular elliptic curve $E_{0}$ over $\overline{%
\mathbb{F}}_{p}$, and denote its canonical model over $\mathbb{F}_{p^{2}}$
by $E_{0}^{\prime }$. Let $A_{0}=E_{0}^{g}$ be the corresponding
superspecial abelian variety over $\overline{\mathbb{F}}_{p}$, endowed with
the \textit{identity} principal polarization $\lambda _{0}:=id_{E_{0}}^{g}$
(we identify canonically $E_{0}$ and $\widehat{E}_{0}$).

\textit{Assume} that the moduli scheme $\mathcal{S}_{\mathcal{D},K^{p}}$
contains a point of the form $[(A_{0},i_{0},\overline{\lambda }_{0},%
\overline{\alpha }_{0})]\in \mathcal{S}_{\mathcal{D},K^{p}}\left( \overline{%
\mathbb{F}}_{p}\right) $ that we fix; the $p$-divisible group $\mathbf{X=}%
A_{0}(p)$ over $\overline{\mathbb{F}}_{p}$ is isomorphic to $G_{1/2}^{\oplus
g}$\ and is endowed with the action $i_{\mathbf{X}}$ of $\mathcal{O}_{B_{p}}$
induced by $i_{0}$ and the principal polarization $\lambda _{\mathbf{X}%
}:A_{0}(p)\rightarrow \widehat{A_{0}}(p)\simeq G_{1/2}^{\oplus g}$. The
triple $(\mathbf{X},i_{\mathbf{X}},\overline{\lambda }_{\mathbf{X}})$ is
well defined modulo isomorphisms by the given point $[(A_{0},i_{0},\overline{%
\lambda }_{0},\overline{\alpha }_{0})]\in \mathcal{S}_{\mathcal{D}%
,K^{p}}\left( \overline{\mathbb{F}}_{p}\right) .$ As usual, we associate to $%
(\mathbf{X},i_{\mathbf{X}},\overline{\lambda }_{\mathbf{X}})$ the isocrystal 
$(N:=M_{\ast }(\mathbf{X})[\frac{1}{p}],\mathbf{F})$ over $K_{0}$ endowed
with an action of $B_{p}$ and with a non-degenerate bilinear form of
isocrystals $\Psi :N\times N\rightarrow \mathbf{1}(1)$. We fix an
isomorphism of $B\otimes _{%
%TCIMACRO{\U{211a} }%
%BeginExpansion
\mathbb{Q}
%EndExpansion
}K_{0}$-modules $N\simeq V\otimes _{%
%TCIMACRO{\U{211a} }%
%BeginExpansion
\mathbb{Q}
%EndExpansion
}K_{0}$ that respects the skew-symmetric forms on both sides and we then
write the action of Frobenius on the right hand side as $\mathbf{F}=b\otimes
\sigma $ for some $b\in G_{p}(K_{0})$. Since $N$ is isoclinic, the slope
morphism associated to $G$ and $b$ over $K_{0}$ has image contained inside
the center of $G$, so that $b$ is basic in the sense of \ref{basicity}.

We have in hands also a simple $%
%TCIMACRO{\U{211a} }%
%BeginExpansion
\mathbb{Q}
%EndExpansion
_{p}$-PEL datum for moduli of $p$-divisible groups $\mathcal{D}_{p}:\mathcal{%
=(}B_{p},^{\ast },V_{p},\left\langle ,\right\rangle _{p},\mathcal{O}%
_{B_{p}},\Lambda ,b,\mu ),$ having good reduction at $p$ and Shimura field
equal to $E_{\nu }$. The closed subscheme $\mathcal{M}^{\prime }(\overline{%
\mathbb{F}}_{p}):=\mathcal{M}_{(\mathbf{X},i_{\mathbf{X}},\overline{\lambda }%
_{\mathbf{X}})}^{\prime }(\overline{\mathbb{F}}_{p})$ of $\mathcal{\breve{M}}%
(\overline{\mathbb{F}}_{p})$ is then defined and identified with $J(%
%TCIMACRO{\U{211a} }%
%BeginExpansion
\mathbb{Q}
%EndExpansion
_{p})/J(%
%TCIMACRO{\U{2124} }%
%BeginExpansion
\mathbb{Z}
%EndExpansion
_{p})$.

\begin{definition}
We let $Z^{\prime }(\overline{\mathbb{F}}_{p}):=Z^{\prime }([(A_{0},i_{0},%
\overline{\lambda }_{0},\overline{\alpha }_{0})])(\overline{\mathbb{F}}%
_{p})\subseteq \mathcal{S}_{\mathcal{D},K^{p}}(\overline{\mathbb{F}}_{p})$
be the set of points $[(A,i,\overline{\lambda },\overline{\alpha })]\in
Z\left( \overline{\mathbb{F}}_{p}\right) $ such that the principally
polarized $p$-divisible group $(A(p),i,\overline{\lambda })$ of $(A,i,%
\overline{\lambda })$ is isomorphic to $(\mathbf{X},i_{\mathbf{X}},\overline{%
\lambda }_{\mathbf{X}})$. We call $Z^{\prime }(\overline{\mathbb{F}}_{p})$
the \textbf{superspecial locus} associated to $(\mathbf{X},i_{\mathbf{X}},%
\overline{\lambda }_{\mathbf{X}})$.
\end{definition}

The set $Z^{\prime }(\overline{\mathbb{F}}_{p})$ is a closed subset of $Z(%
\overline{\mathbb{F}}_{p})$; furthermore if the class $[(A,i,\overline{%
\lambda },\overline{\alpha })]$ belongs to $Z^{\prime }\left( \overline{%
\mathbb{F}}_{p}\right) $, the $p$-divisible group of $A$ is isomorphic to $%
G_{1/2}^{\oplus g}$, so that $A\simeq A_{0}$ is superspecial.

\begin{remark}
\label{gh is ok}Assume that $\mathcal{D}$ is the PEL-datum of type $C$
defined in \ref{Examples}, with $e=1$ (so that $G(%
%TCIMACRO{\U{211a} }%
%BeginExpansion
\mathbb{Q}
%EndExpansion
)=GSp_{2g}(%
%TCIMACRO{\U{211a} }%
%BeginExpansion
\mathbb{Q}
%EndExpansion
)$); let $A_{0}$ be as above a fixed superspecial abelian variety of
dimension $g>1$ over $\overline{\mathbb{F}}_{p}$. In \cite{Ek}, it is shown
that the isomorphism classes of principal polarizations on $A_{0}$ form a
single genus class:\ this means in particular that if $\lambda $ and $%
\lambda ^{\prime }$ are two principal polarizations on $A_{0}$, then the $p$%
-adic polarizations associated to $\lambda $ and $\lambda ^{\prime }$
respectively on the Dieudonn\'{e} module of $A_{0}$ are isomorphic. Hence $%
(A_{0}(p),\overline{\lambda })\simeq (A_{0}(p),\overline{\lambda }^{\prime
}) $ as principally polarized $p$-divisible groups over $\overline{\mathbb{F}%
}_{p}$.

As a consequence, in the case $\mathcal{D}$ is of type $C$ with $e=1$, we
have: 
\begin{eqnarray*}
Z^{\prime }(\overline{\mathbb{F}}_{p}) &=&\{[(A_{0},\overline{\lambda },%
\overline{\alpha })]:\lambda \text{ a principal polarization on }A_{0}\text{,%
} \\
&&\text{ \ \ \ \ \ \ \ \ \ \ \ \ \ \ \ \ \ \ \ \ \ \ }\alpha \text{ a }K^{p}%
\text{-level structure on }A_{0}\}
\end{eqnarray*}
\end{remark}

\begin{proposition}
\label{unif - prop}The uniformization morphism $\Theta _{K^{p}}(\mathcal{%
\overline{\mathbb{F}}}_{p})$ of Proposition \ref{key} induces a canonical
Hecke-equivariant isomorphism:%
\begin{equation*}
\Theta _{K^{p}}^{\prime }(\mathcal{\overline{\mathbb{F}}}_{p}):I(%
%TCIMACRO{\U{211a} }%
%BeginExpansion
\mathbb{Q}
%EndExpansion
)\backslash \mathcal{M}^{\prime }\mathcal{(\overline{\mathbb{F}}}_{p}%
\mathcal{)\times }G(\mathbb{A}_{f}^{p})/K^{p}\longrightarrow Z^{\prime }(%
\mathcal{\overline{\mathbb{F}}}_{p}).
\end{equation*}

\noindent We call $\Theta _{K^{p}}^{\prime }(\mathcal{\overline{\mathbb{F}}}%
_{p})$ the \textbf{uniformization morphism for the superspecial locus}%
.\noindent
\end{proposition}

\textbf{Proof. }First recall that under our assumptions on $G$, and by the
basicity of $b$, the map $\Theta _{K^{p}}(\mathcal{\overline{\mathbb{F}}}%
_{p})$ is a well defined Hecke-equivariant isomorphism. The action of $I(%
%TCIMACRO{\U{211a} }%
%BeginExpansion
\mathbb{Q}
%EndExpansion
)$ on $\mathcal{\breve{M}(\overline{\mathbb{F}}}_{p}\mathcal{)}$ determines
an action on $\mathcal{M}^{\prime }\mathcal{(\overline{\mathbb{F}}}_{p}%
\mathcal{)\subseteq \breve{M}(\overline{\mathbb{F}}}_{p}\mathcal{)}$, so
that we obtain a natural injective Hecke-equivariant map:%
\begin{equation*}
I(%
%TCIMACRO{\U{211a} }%
%BeginExpansion
\mathbb{Q}
%EndExpansion
)\backslash \mathcal{M}^{\prime }\mathcal{(\overline{\mathbb{F}}}_{p}%
\mathcal{)\times }G(\mathbb{A}_{f}^{p})/K^{p}\hookrightarrow I(%
%TCIMACRO{\U{211a} }%
%BeginExpansion
\mathbb{Q}
%EndExpansion
)\backslash \mathcal{\breve{M}(\overline{\mathbb{F}}}_{p}\mathcal{)\times }G(%
\mathbb{A}_{f}^{p})/K^{p}.
\end{equation*}%
\ Define $\Theta _{K^{p}}^{\prime }(\mathcal{\overline{\mathbb{F}}}_{p})$\
by precomposing this last map with $\Theta _{K^{p}}(\mathcal{\overline{%
\mathbb{F}}}_{p})$. In order to determine the image of $\Theta
_{K^{p}}^{\prime }(\mathcal{\overline{\mathbb{F}}}_{p})$, we follow the
construction of the uniformization morphism over the field $\mathcal{%
\overline{\mathbb{F}}}_{p}$ (cf. \ref{uniformization} and \cite{RZ},
6.13-14).

Pick an element $\left[ \rho \right] \in \mathcal{M}^{\prime }\mathcal{(%
\overline{\mathbb{F}}}_{p}\mathcal{)}$; the quasi-isogeny $\rho :(\mathbf{X}%
,i_{\mathbf{X}},\overline{\lambda }_{\mathbf{X}})\rightarrow (\mathbf{X},i_{%
\mathbf{X}},\overline{\lambda }_{\mathbf{X}})$ determines, by Lemma \ref%
{push-forward}, a principally polarized abelian variety:%
\begin{equation*}
(\rho _{\ast }A_{0},\rho _{\ast }i_{0},\rho _{\ast }\overline{\lambda }_{0})
\end{equation*}
in $AV_{\mathcal{O}_{B}}$, whose $p$-divisible group is isomorphic to $(%
\mathbf{X},i_{\mathbf{X}},\overline{\lambda }_{\mathbf{X}})$,\ so that the
image of $\overline{\Theta }_{K^{p}}^{\prime }(\mathcal{\overline{\mathbb{F}}%
}_{p})$ is contained inside the superspecial locus.

Viceversa, let $[(A,i,\overline{\lambda },\overline{\alpha })]\in Z^{\prime
}(\mathcal{\overline{\mathbb{F}}}_{p})$ and choose a quasi-isogeny of
principally polarized abelian varieties $\rho :(A_{0},i_{0},\overline{%
\lambda }_{0})\rightarrow (A,i,\overline{\lambda })$. Then $\rho $ defines a
quasi-isogeny of the corresponding $p$-divisible groups $\rho :(\mathbf{X}%
,i_{\mathbf{X}},\overline{\lambda }_{\mathbf{X}})\rightarrow (A(p),i,%
\overline{\lambda })$. Precomposing $\rho $ with an isomorphism $\mu
:(A(p),i,\overline{\lambda })\rightarrow (\mathbf{X},i_{\mathbf{X}},%
\overline{\lambda }_{\mathbf{X}})$ we obtain an element $[\mu \circ \rho
]\in \mathcal{M}^{\prime }\mathcal{(\overline{\mathbb{F}}}_{p}\mathcal{)}$
such that $(\mu \circ \rho )_{\ast }(A_{0},i_{0},\overline{\lambda }%
_{0})=\mu _{\ast }(A,i,\overline{\lambda })\simeq (A,i,\overline{\lambda })$%
. Let now $g\in G(\mathbb{A}_{f}^{p})$ defined by $g:=\left( \mu \circ \rho
\right) _{\ast }\overline{\alpha }_{0}\circ \overline{\alpha }^{-1}$; the
pre-image of $[(A,i,\overline{\lambda },\overline{\alpha }]\in Z^{\prime }(%
\mathcal{\overline{\mathbb{F}}}_{p})$ under $\Theta _{K^{p}}^{\prime }(%
\mathcal{\overline{\mathbb{F}}}_{p})$ is the $I(%
%TCIMACRO{\U{211a} }%
%BeginExpansion
\mathbb{Q}
%EndExpansion
)$-class represented by $\left[ \mu \circ \rho \right] \times gK^{p}$.\ $%
\blacksquare $

\begin{remark}
\label{Important}We have seen that the group $I(%
%TCIMACRO{\U{211a} }%
%BeginExpansion
\mathbb{Q}
%EndExpansion
):=(\limfunc{End}_{\mathcal{O}_{B}}(A_{0},\overline{\lambda }_{0})\otimes _{%
%TCIMACRO{\U{2124} }%
%BeginExpansion
\mathbb{Z}
%EndExpansion
}%
%TCIMACRO{\U{211a} }%
%BeginExpansion
\mathbb{Q}
%EndExpansion
)^{\times }$ acts on the left upon $\mathcal{M}^{\prime }(\overline{\mathbb{F%
}}_{p})$ through the map $\alpha _{0,p}:I(%
%TCIMACRO{\U{211a} }%
%BeginExpansion
\mathbb{Q}
%EndExpansion
)\rightarrow J(%
%TCIMACRO{\U{211a} }%
%BeginExpansion
\mathbb{Q}
%EndExpansion
_{p})$. We have therefore the canonical identification:%
\begin{equation*}
I(%
%TCIMACRO{\U{211a} }%
%BeginExpansion
\mathbb{Q}
%EndExpansion
)\backslash \mathcal{M}^{\prime }(\overline{\mathbb{F}}_{p})\overset{\simeq }%
{\rightarrow }I(%
%TCIMACRO{\U{211a} }%
%BeginExpansion
\mathbb{Q}
%EndExpansion
)\backslash J(%
%TCIMACRO{\U{211a} }%
%BeginExpansion
\mathbb{Q}
%EndExpansion
_{p})/J(%
%TCIMACRO{\U{2124} }%
%BeginExpansion
\mathbb{Z}
%EndExpansion
_{p}),
\end{equation*}

\noindent where the action of $I(%
%TCIMACRO{\U{211a} }%
%BeginExpansion
\mathbb{Q}
%EndExpansion
)$ on $\mathcal{M}^{\prime }(\overline{\mathbb{F}}_{p})$ is the one
described in \ref{groups I and J}, so that (cf. proof of Proposition \ref%
{obvious}) the action of $I(%
%TCIMACRO{\U{211a} }%
%BeginExpansion
\mathbb{Q}
%EndExpansion
)$ on the coset space $J(%
%TCIMACRO{\U{211a} }%
%BeginExpansion
\mathbb{Q}
%EndExpansion
_{p})/J(%
%TCIMACRO{\U{2124} }%
%BeginExpansion
\mathbb{Z}
%EndExpansion
_{p})$ is given by $x\cdot gJ(%
%TCIMACRO{\U{2124} }%
%BeginExpansion
\mathbb{Z}
%EndExpansion
_{p})=(M_{\ast }(x)\cdot g)J(%
%TCIMACRO{\U{2124} }%
%BeginExpansion
\mathbb{Z}
%EndExpansion
_{p})$, for all $x\in I(%
%TCIMACRO{\U{211a} }%
%BeginExpansion
\mathbb{Q}
%EndExpansion
)$ and all $g\in J(%
%TCIMACRO{\U{211a} }%
%BeginExpansion
\mathbb{Q}
%EndExpansion
_{p})$. More easily, we will write $x\cdot gJ(%
%TCIMACRO{\U{2124} }%
%BeginExpansion
\mathbb{Z}
%EndExpansion
_{p})=xgJ(%
%TCIMACRO{\U{2124} }%
%BeginExpansion
\mathbb{Z}
%EndExpansion
_{p}).$
\end{remark}

\begin{corollary}
\label{unif-cor}There is a canonical Hecke equivariant isomorphism:%
\begin{equation*}
\Theta _{K^{p}}^{\prime }(\mathcal{\overline{\mathbb{F}}}_{p}):I(%
%TCIMACRO{\U{211a} }%
%BeginExpansion
\mathbb{Q}
%EndExpansion
)\backslash \left( J(%
%TCIMACRO{\U{211a} }%
%BeginExpansion
\mathbb{Q}
%EndExpansion
_{p})/J(%
%TCIMACRO{\U{2124} }%
%BeginExpansion
\mathbb{Z}
%EndExpansion
_{p})\mathcal{\times }G(\mathbb{A}_{f}^{p})/K^{p}\right) \rightarrow
Z^{\prime }(\mathcal{\overline{\mathbb{F}}}_{p}),
\end{equation*}

\noindent where the action of $I(%
%TCIMACRO{\U{211a} }%
%BeginExpansion
\mathbb{Q}
%EndExpansion
)$ on $J(%
%TCIMACRO{\U{211a} }%
%BeginExpansion
\mathbb{Q}
%EndExpansion
_{p})/J(%
%TCIMACRO{\U{2124} }%
%BeginExpansion
\mathbb{Z}
%EndExpansion
_{p})$ is the natural one, described in Remark \ref{Important}. Furthermore, 
$Z^{\prime }(\mathcal{\overline{\mathbb{F}}}_{p})$ is a finite set.
\end{corollary}

\textbf{Proof. }We just need to show the finiteness of $Z^{\prime }(\mathcal{%
\overline{\mathbb{F}}}_{p})$. By Theorem \ref{RZ} we have canonical
identifications $I(\mathbb{A}_{f}^{p})=G(\mathbb{A}_{f}^{p})$ and $J(%
%TCIMACRO{\U{211a} }%
%BeginExpansion
\mathbb{Q}
%EndExpansion
_{p})=I(%
%TCIMACRO{\U{211a} }%
%BeginExpansion
\mathbb{Q}
%EndExpansion
_{p})$, so that we can rewrite the domain of the morphism $\Theta
_{K^{p}}^{\prime }(\mathcal{\overline{\mathbb{F}}}_{p})$ as: 
\begin{equation*}
I(%
%TCIMACRO{\U{211a} }%
%BeginExpansion
\mathbb{Q}
%EndExpansion
)\backslash \left( I(%
%TCIMACRO{\U{211a} }%
%BeginExpansion
\mathbb{Q}
%EndExpansion
_{p})/I(%
%TCIMACRO{\U{2124} }%
%BeginExpansion
\mathbb{Z}
%EndExpansion
_{p})\mathcal{\times }I(\mathbb{A}_{f}^{p})/C^{p}\right) =I(%
%TCIMACRO{\U{211a} }%
%BeginExpansion
\mathbb{Q}
%EndExpansion
)\backslash I(\mathbb{A}_{f})/C,
\end{equation*}%
where $C^{p}$ is the image of $K^{p}$ in $I(\mathbb{A}_{f}^{p})$ and $C=I(%
%TCIMACRO{\U{2124} }%
%BeginExpansion
\mathbb{Z}
%EndExpansion
_{p})\times C^{p}$. By the proof of 6.23 in \cite{RZ}, $I(%
%TCIMACRO{\U{211a} }%
%BeginExpansion
\mathbb{Q}
%EndExpansion
)$ is a discrete subgroup of $I(%
%TCIMACRO{\U{211a} }%
%BeginExpansion
\mathbb{Q}
%EndExpansion
_{p})\mathcal{\times }I(\mathbb{A}_{f}^{p})=I(\mathbb{A}_{f})$; by
Proposition 1.4 of \cite{Gr}, the quotient space $I(%
%TCIMACRO{\U{211a} }%
%BeginExpansion
\mathbb{Q}
%EndExpansion
)\backslash I(\mathbb{A}_{f})$ is therefore compact, so that $I(%
%TCIMACRO{\U{211a} }%
%BeginExpansion
\mathbb{Q}
%EndExpansion
)\backslash I(\mathbb{A}_{f})/C$ is finite. $\blacksquare $

\newpage

\section{Comparing Hecke eigensystems}

We apply the results of the last section to study systems of Hecke
eigenvalues coming from unitary modular forms. We begin by recalling some
notation: although this notation will be assumed fixed in the rest of the
chapter, we will specify in the following paragraphs additional requirements
on the objects considered below.

Fix an algebraic closure $\overline{%
%TCIMACRO{\U{211a} }%
%BeginExpansion
\mathbb{Q}
%EndExpansion
}$ (resp. $\overline{%
%TCIMACRO{\U{211a} }%
%BeginExpansion
\mathbb{Q}
%EndExpansion
}_{p}$) of $%
%TCIMACRO{\U{211a} }%
%BeginExpansion
\mathbb{Q}
%EndExpansion
$ (resp. $%
%TCIMACRO{\U{211a} }%
%BeginExpansion
\mathbb{Q}
%EndExpansion
_{p}$) and an embedding $\overline{%
%TCIMACRO{\U{211a} }%
%BeginExpansion
\mathbb{Q}
%EndExpansion
}\subset 
%TCIMACRO{\U{2102} }%
%BeginExpansion
\mathbb{C}
%EndExpansion
$. Fix a simple $%
%TCIMACRO{\U{211a} }%
%BeginExpansion
\mathbb{Q}
%EndExpansion
$-PEL-datum $\mathcal{D=(}B,^{\ast },V,\left\langle ,\right\rangle ,\mathcal{%
O}_{B},\Lambda ,h,K^{p},\nu \mathcal{)}$ for moduli of abelian schemes with
good reduction at $p$, and neat level $K^{p}$; denote by $G$ the associated
algebraic group, and assume it is connected and satisfies the Hasse
principle. Let $\mu $ be the co-character associated to $h$, and let $%
E\subset \overline{%
%TCIMACRO{\U{211a} }%
%BeginExpansion
\mathbb{Q}
%EndExpansion
}$ be the Shimura field of $\mathcal{D}$. Let $\overline{\mathbb{F}}_{p}$ be
a fixed algebraic closure of the residue field of $E_{\nu }\subset \overline{%
%TCIMACRO{\U{211a} }%
%BeginExpansion
\mathbb{Q}
%EndExpansion
}_{p}$; set $W=W(\overline{\mathbb{F}}_{p})$, $K_{0}=W[\frac{1}{p}]\subset 
\overline{%
%TCIMACRO{\U{211a} }%
%BeginExpansion
\mathbb{Q}
%EndExpansion
}_{p}$ and denote by $\sigma $ the Frobenius morphism of $W$; recall that $%
E_{\nu }\subset K_{0}$.\ Fix a finite extension $K\subset \overline{%
%TCIMACRO{\U{211a} }%
%BeginExpansion
\mathbb{Q}
%EndExpansion
}_{p}$ of $K_{0}$ such that $\mu $ is defined over $K$; set $\breve{E}_{\nu
}=E_{\nu }K_{0}=K_{0}$. Define $B_{p}$, $V_{p}$, $\left\langle
,\right\rangle _{p}$, $G_{p}$, $\mathcal{O}_{B_{p}}$ as usual.

Let $\mathcal{S}_{\mathcal{D},K^{p}}$ be the quasi-projective smooth scheme
over $\mathcal{O}_{\breve{E}_{\nu }}$ defined in Th. \ref{Representability};
suppose that the common relative dimension of the abelian schemes
parametrized by $\mathcal{S}_{\mathcal{D},K^{p}}$ is $g:=\dim _{%
%TCIMACRO{\U{2102} }%
%BeginExpansion
\mathbb{C}
%EndExpansion
}V_{%
%TCIMACRO{\U{2102} }%
%BeginExpansion
\mathbb{C}
%EndExpansion
,0}\geq 2$; fix a supersingular elliptic curve $E_{0}$ over $\overline{%
\mathbb{F}}_{p}$, and denote its canonical model over $\mathbb{F}_{p^{2}}$
by $E_{0}^{\prime }$. Let $A_{0}:=E_{0}^{g}$ be the corresponding
superspecial abelian variety over $\overline{\mathbb{F}}_{p}$, endowed with
the identity principal polarization $\lambda _{0}:=id_{E_{0}}^{g}$ (we will
always identify canonically $E_{0}$ and $\widehat{E}_{0}$). Denote by $%
\mathfrak{R}$ the $%
%TCIMACRO{\U{2124} }%
%BeginExpansion
\mathbb{Z}
%EndExpansion
$-algebra $\limfunc{End}E_{0}=\limfunc{End}_{\mathbb{F}_{p^{2}}}E_{0}^{%
\prime }$ of endomorphisms of the supersingular elliptic curve $E_{0}$ over $%
\mathcal{\overline{\mathbb{F}}}_{p}$; $\mathfrak{R}$ is a maximal order in $%
\mathfrak{B}:=\limfunc{End}^{0}E_{0}=\limfunc{End}E_{0}\otimes _{%
%TCIMACRO{\U{2124} }%
%BeginExpansion
\mathbb{Z}
%EndExpansion
}%
%TCIMACRO{\U{211a} }%
%BeginExpansion
\mathbb{Q}
%EndExpansion
$, where $\mathfrak{B}$ is a quaternion algebra over $%
%TCIMACRO{\U{211a} }%
%BeginExpansion
\mathbb{Q}
%EndExpansion
$ whose ramification set is $\{p,\infty \}$. If $v$ is a place of $%
%TCIMACRO{\U{211a} }%
%BeginExpansion
\mathbb{Q}
%EndExpansion
$, we denote by $\mathfrak{B}_{v}$ the $%
%TCIMACRO{\U{211a} }%
%BeginExpansion
\mathbb{Q}
%EndExpansion
_{v}$-algebra $\mathfrak{B}\otimes _{%
%TCIMACRO{\U{211a} }%
%BeginExpansion
\mathbb{Q}
%EndExpansion
}%
%TCIMACRO{\U{211a} }%
%BeginExpansion
\mathbb{Q}
%EndExpansion
_{v}$; we also denote by $\overline{\cdot }$ the canonical involution of $%
\mathfrak{B}$ (cf. \cite{Vig}).

Assume that $\mathcal{S}_{\mathcal{D},K^{p}}$ contains a point of the form $%
[(A_{0},i_{0},\overline{\lambda }_{0},\overline{\alpha }_{0})]\in \mathcal{S}%
_{\mathcal{D},K^{p}}\left( \overline{\mathbb{F}}_{p}\right) $ for some $%
i_{0},\overline{\lambda }_{0}$ and $\overline{\alpha }_{0}$. Let us fix such
a point; the $p$-divisible group $\mathbf{X:=}A_{0}(p)$ over $\overline{%
\mathbb{F}}_{p}$ is isomorphic to $G_{1/2}^{\oplus g}$\ and is endowed with
the action $i_{\mathbf{X}}$ of $\mathcal{O}_{B_{p}}$ induced by $i_{0}$ and
the principal polarization $\lambda _{\mathbf{X}}:A_{0}(p)\rightarrow 
\widehat{A_{0}}(p)\simeq G_{1/2}^{\oplus g}$. By covariant Dieudonn\'{e}
theory, we associate to $(\mathbf{X},i_{\mathbf{X}},\overline{\lambda }_{%
\mathbf{X}})$ a Dieudonn\'{e} module $\mathbf{M}:=M_{\ast }(\mathbf{X})$
over $W$, endowed with an action $i_{\mathbf{M}}$ of $\mathcal{O}_{B_{p}}$
and a principal polarization $e_{\mathbf{M}}:\mathbf{M\times M\rightarrow }W$
of Dieudonn\'{e} modules, which is skew-Hermitian with respect to $^{\ast }$%
, and well defined only up to a scalar factor in $%
%TCIMACRO{\U{2124} }%
%BeginExpansion
\mathbb{Z}
%EndExpansion
_{p}^{\times }$. By inverting $p$, we obtain an isocrystal $(N:=\mathbf{M}[%
\frac{1}{p}],\mathbf{F})$ over $K_{0}$ endowed with an action of $B_{p}$ and
with a non-degenerate bilinear form of isocrystals $\Psi :N\times
N\rightarrow \mathbf{1}(1)$. We fix an isomorphism of $B\otimes _{%
%TCIMACRO{\U{211a} }%
%BeginExpansion
\mathbb{Q}
%EndExpansion
}K_{0}$-modules $N\simeq V\otimes _{%
%TCIMACRO{\U{211a} }%
%BeginExpansion
\mathbb{Q}
%EndExpansion
}K_{0}$ that respects the skew-symmetric forms on both sides and we then
write the action of Frobenius on the right hand side as $\mathbf{F}=b\otimes
\sigma $ for some $b\in G_{p}(K_{0})$. Since $N$ is isoclinic, $b$ is basic
in the sense of \ref{basicity}.

We have a simple $%
%TCIMACRO{\U{211a} }%
%BeginExpansion
\mathbb{Q}
%EndExpansion
_{p}$-PEL datum for moduli of $p$-divisible groups $\mathcal{D}_{p}\mathcal{%
=(}B_{p},^{\ast },V_{p},\left\langle ,\right\rangle _{p},\mathcal{O}%
_{B_{p}},\Lambda ,b,\mu ),$having good reduction at $p$ and Shimura field
equal to $E_{\nu }$. Associated to $\mathcal{D}_{p}$ we have the moduli
functor $\mathcal{\breve{M}}$. The fixed choice of $(\mathbf{X},i_{\mathbf{X}%
},\overline{\lambda }_{\mathbf{X}})$ gives rise to a closed subset $\mathcal{%
M}^{\prime }(\overline{\mathbb{F}}_{p})\subseteq \mathcal{\breve{M}}(%
\overline{\mathbb{F}}_{p})$ (cf. Def. \ref{modified functor}), and a
uniformization Hecke-equivariant isomorphism of finite sets:%
\begin{equation*}
\Theta _{K^{p}}^{\prime }(\mathcal{\overline{\mathbb{F}}}_{p}):I(%
%TCIMACRO{\U{211a} }%
%BeginExpansion
\mathbb{Q}
%EndExpansion
)\backslash \mathcal{M}^{\prime }\mathcal{(\overline{\mathbb{F}}}_{p}%
\mathcal{)\times }G(\mathbb{A}_{f}^{p})/K^{p}\longrightarrow Z^{\prime }(%
\mathcal{\overline{\mathbb{F}}}_{p}),
\end{equation*}

\noindent or:%
\begin{equation*}
\Theta _{K^{p}}^{\prime }(\mathcal{\overline{\mathbb{F}}}_{p}):I(%
%TCIMACRO{\U{211a} }%
%BeginExpansion
\mathbb{Q}
%EndExpansion
)\backslash \left( J(%
%TCIMACRO{\U{211a} }%
%BeginExpansion
\mathbb{Q}
%EndExpansion
_{p})/J(%
%TCIMACRO{\U{2124} }%
%BeginExpansion
\mathbb{Z}
%EndExpansion
_{p})\mathcal{\times }G(\mathbb{A}_{f}^{p})/K^{p}\right) \longrightarrow
Z^{\prime }(\mathcal{\overline{\mathbb{F}}}_{p}).
\end{equation*}

Since $A_{0}=E_{0}^{g}$, we obtain canonical isomorphisms $\limfunc{End}%
A_{0}=M_{g}(\mathfrak{R})$ and $\limfunc{End}^{0}A_{0}=M_{g}(\mathfrak{B})$.
Under these identifications, and the canonical isomorphism $A_{0}=\widehat{A}%
_{0}$, the principal polarization $\lambda _{0}=id_{E_{0}}^{g}$ coincides
with the identity matrix $I_{g}\in M_{g}(\mathfrak{R})$, so that the
auto-quasi-isogenies of the principally (homogeneously) polarized abelian
variety $(A_{0},\overline{\lambda }_{0})$ are identified with the elements
of the unitary quaternion similitude group:%
\begin{equation*}
GU_{g}(\mathfrak{B};I_{g}):=\{X\in GL_{g}(\mathfrak{B}):X^{\ast }X=c(X)\cdot
I_{g},\text{ }c(X)\in 
%TCIMACRO{\U{211a} }%
%BeginExpansion
\mathbb{Q}
%EndExpansion
^{\times }\},
\end{equation*}%
\noindent where $X^{\ast }:=\overline{X}^{t}$. Similarly, the automorphisms
of the pair $(A_{0},\overline{\lambda }_{0})$ are given by $GU_{g}(\mathfrak{%
R};I_{g})$, and the automorphisms of $(A_{0},\overline{\lambda }_{0})$
viewed as a polarized abelian variety up to prime-to-$p$ isogeny are given
by $GU_{g}(\mathfrak{R}\otimes _{%
%TCIMACRO{\U{2124} }%
%BeginExpansion
\mathbb{Z}
%EndExpansion
}%
%TCIMACRO{\U{2124} }%
%BeginExpansion
\mathbb{Z}
%EndExpansion
_{(p)};I_{g})$.

Notice that $GU_{g}(\mathfrak{B};I_{g})$ defines an algebraic group over $%
%TCIMACRO{\U{211a} }%
%BeginExpansion
\mathbb{Q}
%EndExpansion
$, which is reductive, since it is a form of the reductive group $GSp_{2g/%
%TCIMACRO{\U{211a} }%
%BeginExpansion
\mathbb{Q}
%EndExpansion
}$,\ and such that $U_{g}(\mathfrak{B};I_{g})$ is compact at infinity (since 
$U_{g}(\mathfrak{B}_{\infty };I_{g})\subset O(4g)$).

\subsection{The settings}

Let us assume that $\mathcal{D=D}_{(r,s),p}^{U}$ is the PEL-datum of type A
constructed in \ref{Examples}: in particular, $\mathcal{D}$ is associated to
a quadratic imaginary field $k=%
%TCIMACRO{\U{211a} }%
%BeginExpansion
\mathbb{Q}
%EndExpansion
\left( \sqrt{\alpha }\right) \subset \overline{%
%TCIMACRO{\U{211a} }%
%BeginExpansion
\mathbb{Q}
%EndExpansion
}$ ($\alpha \in 
%TCIMACRO{\U{2124} }%
%BeginExpansion
\mathbb{Z}
%EndExpansion
_{<0}$ square free) in which $p$ is inert. By our previous assumptions, we
are given an embedding $\tau :k\hookrightarrow 
%TCIMACRO{\U{2102} }%
%BeginExpansion
\mathbb{C}
%EndExpansion
$ of fields such that $\tau (\sqrt{\alpha })=\sqrt{-1}\sqrt{-\alpha }$ for
our fixed choice of square root of $-1$; we denote by $\overline{\cdot }$
the non-trivial field automorphism of $k$. Recall that we set $V=k^{g}$, $%
g=2n$; $r,s$ are two positive integers such that $r+s=g$. The algebraic $%
%TCIMACRO{\U{211a} }%
%BeginExpansion
\mathbb{Q}
%EndExpansion
$-group associated to $\mathcal{D}_{(r,s),p}^{U}$ is the unitary group
relative to $k$ and to the form:%
\begin{equation*}
H=%
\begin{pmatrix}
-\sqrt{\alpha }I_{r} & 0_{r,s} \\ 
0_{s,r} & \sqrt{\alpha }I_{s}%
\end{pmatrix}%
,
\end{equation*}%
i.e. it is $G=GU_{g}(k;r,s)$; $G$ is connected, satisfies the Hasse
principle and can be defined over $%
%TCIMACRO{\U{2124} }%
%BeginExpansion
\mathbb{Z}
%EndExpansion
$. We have shown in \ref{Examples} that the reflex field of our datum is $E=%
%TCIMACRO{\U{211a} }%
%BeginExpansion
\mathbb{Q}
%EndExpansion
$ if $r=s=n$, and $E=k$ otherwise. The determinant polynomial is:%
\begin{equation*}
f(X_{1},X_{2})=(X_{1}-\sqrt{\alpha }X_{2})^{r}(X_{1}+\sqrt{\alpha }%
X_{2})^{s}\in \mathcal{O}_{E}[X_{1},X_{2}]\text{.}
\end{equation*}

We take $K^{p}:=U(N)$ for a fixed integer $N\geq 3$ not divisible by $p$.
Assume for simplicity that $p\neq 2$.

\bigskip

The embedding $\nu :\overline{%
%TCIMACRO{\U{211a} }%
%BeginExpansion
\mathbb{Q}
%EndExpansion
}\hookrightarrow \overline{%
%TCIMACRO{\U{211a} }%
%BeginExpansion
\mathbb{Q}
%EndExpansion
}_{p}$ that we fixed identifies $k_{\nu }$ with the only degree two
unramified extension of $%
%TCIMACRO{\U{211a} }%
%BeginExpansion
\mathbb{Q}
%EndExpansion
_{p}$ inside $K_{0}$; if $\mathbb{F}_{p^{2}}$ denotes the residue field of $%
k_{\nu }$,\ we obtain an embedding $\mathbb{F}_{p^{2}}\subset $ $\mathcal{%
\overline{\mathbb{F}}}_{p}$.

\begin{proposition}
There is a supersingular elliptic curve $E_{0}$ over $\mathcal{\overline{%
\mathbb{F}}}_{p}$ whose endomorphism ring $\limfunc{End}E_{0}$ contains an
element $\varphi _{\alpha }$ such that:

\begin{enumerate}
\item $\varphi _{\alpha }^{2}=\alpha ;$

\item the tangent map $\limfunc{Lie}\varphi _{\alpha }:\limfunc{Lie}%
E_{0}\rightarrow \limfunc{Lie}E_{0}$ is multiplication by the scalar $\sqrt{%
\alpha }(\func{mod}p)\in \mathbb{F}_{p^{2}}\subset \mathcal{\overline{%
\mathbb{F}}}_{p}$, where here $\sqrt{\alpha }\in k$ is viewed as an element
of $k_{\nu }$ via $\nu .$
\end{enumerate}
\end{proposition}

\textbf{Proof. }Let $E$ be any fixed supersingular elliptic curve over $%
\mathcal{\overline{\mathbb{F}}}_{p}$. By \cite{Vig}, Th\'{e}or\`{e}me 3.8,
page 78, there is an embedding of $%
%TCIMACRO{\U{211a} }%
%BeginExpansion
\mathbb{Q}
%EndExpansion
$-algebras $j:k\hookrightarrow \limfunc{End}^{0}E$; there is a maximal order 
$R$ of $\limfunc{End}^{0}E$ containing $j(\sqrt{\alpha })$: in fact we can
write $\limfunc{End}^{0}E=j(k)\oplus j(k)u$ for some $u\in \limfunc{End}%
^{0}E $ (cf. \cite{Vig}, Corollaire 2.2, page 6), and the left order of the
ideal $\mathfrak{%
%TCIMACRO{\U{2124} }%
%BeginExpansion
\mathbb{Z}
%EndExpansion
+%
%TCIMACRO{\U{2124} }%
%BeginExpansion
\mathbb{Z}
%EndExpansion
}j\mathfrak{(}\sqrt{\alpha })+%
%TCIMACRO{\U{2124} }%
%BeginExpansion
\mathbb{Z}
%EndExpansion
u+%
%TCIMACRO{\U{2124} }%
%BeginExpansion
\mathbb{Z}
%EndExpansion
j(\sqrt{\alpha })u$ of $\limfunc{End}^{0}E$ clearly contains $j\mathfrak{(}%
\sqrt{\alpha })$. By work of Deuring, it is known that there is a an
elliptic curve $E_{0}$ over $\mathcal{\overline{\mathbb{F}}}_{p}$ and a
quasi-isogeny $f:E_{0}\rightarrow E$ such that $R=f\circ \limfunc{End}%
E_{0}\circ f^{-1}$, so that $\limfunc{End}E_{0}$ contains an element $%
\varphi _{\alpha }^{\prime }$ whose square equals $\alpha $ (cf. \cite{Wed}
2.15).

The tangent morphism $\limfunc{Lie}\varphi _{\alpha }^{\prime }$ can be
canonically identified with an element of $\mathcal{\overline{\mathbb{F}}}%
_{p}$; since $\left( \limfunc{Lie}\varphi _{\alpha }^{\prime }\right)
^{2}=\alpha \func{mod}p$, we have $\limfunc{Lie}\varphi _{\alpha }^{\prime
}=\pm \sqrt{\alpha }(\func{mod}p)\in \mathbb{F}_{p^{2}}$ (recall we are
given a fixed embedding $\mathbb{F}_{p^{2}}\subset $ $\mathcal{\overline{%
\mathbb{F}}}_{p}$). We define $\varphi _{\alpha }:=\pm \varphi _{\alpha
}^{\prime }$ depending on $\limfunc{Lie}\varphi _{\alpha }^{\prime }$ being
equal to $\pm \sqrt{\alpha }(\func{mod}p)$ respectively. The pair $%
(E_{0},\varphi _{\alpha })$ we just constructed satisfies the requirement of
the proposition. $\blacksquare $

\bigskip

Fix a pair $(E_{0},\varphi _{\alpha })$ over $\mathcal{\overline{\mathbb{F}}}%
_{p}$ as in the above proposition (recall that the choice of isomorphism
class of $E_{0}$ will not be relevant later on, since $g\geq 2$) and set $%
A_{0}=E_{0}^{g}$, $\lambda _{0}=id_{E_{0}}^{g}$; recall that $\mathfrak{R}$
is the $%
%TCIMACRO{\U{2124} }%
%BeginExpansion
\mathbb{Z}
%EndExpansion
$-algebra $\limfunc{End}E_{0}$ (containing $\varphi _{\alpha }$), and $%
\mathfrak{B}:=\limfunc{End}^{0}E_{0}.$ Let $\iota :k\hookrightarrow 
\mathfrak{B}$\ be the $%
%TCIMACRO{\U{211a} }%
%BeginExpansion
\mathbb{Q}
%EndExpansion
$-algebra homomorphism such that $\iota (\sqrt{\alpha })=\varphi _{\alpha }$%
. Since $p$ is odd, we have $\mathcal{O}_{B}:=\mathcal{O}_{k,(p)}=%
%TCIMACRO{\U{2124} }%
%BeginExpansion
\mathbb{Z}
%EndExpansion
_{(p)}[\sqrt{\alpha }]$. Define a $%
%TCIMACRO{\U{2124} }%
%BeginExpansion
\mathbb{Z}
%EndExpansion
_{(p)}$-algebra monomorphism $i_{0}:\mathcal{O}_{k,(p)}\mathcal{%
\hookrightarrow }\limfunc{End}(A_{0})\otimes _{%
%TCIMACRO{\U{2124} }%
%BeginExpansion
\mathbb{Z}
%EndExpansion
}%
%TCIMACRO{\U{2124} }%
%BeginExpansion
\mathbb{Z}
%EndExpansion
_{(p)}=M_{g}(\mathfrak{R}\otimes _{%
%TCIMACRO{\U{2124} }%
%BeginExpansion
\mathbb{Z}
%EndExpansion
}%
%TCIMACRO{\U{2124} }%
%BeginExpansion
\mathbb{Z}
%EndExpansion
_{(p)})$ by requiring that: 
\begin{equation*}
\sqrt{\alpha }\mapsto 
\begin{pmatrix}
-\varphi _{\alpha }I_{r} & 0_{r,s} \\ 
0_{s,r} & \varphi _{\alpha }I_{s}%
\end{pmatrix}%
.
\end{equation*}

\noindent On the dual variety $\hat{A}_{0}=A_{0}$ we have the dual action $%
\hat{\imath}_{0}:\mathcal{O}_{k,(p)}\mathcal{\hookrightarrow }M_{g}(%
\mathfrak{R}\otimes _{%
%TCIMACRO{\U{2124} }%
%BeginExpansion
\mathbb{Z}
%EndExpansion
}%
%TCIMACRO{\U{2124} }%
%BeginExpansion
\mathbb{Z}
%EndExpansion
_{(p)})$\ defined by $\hat{\imath}_{0}(b)=i_{0}(\overline{b})^{\ast }$ for
any $b\in \mathcal{O}_{k,(p)}$ (recall that for $X\in M_{g}(\mathfrak{B})$
we have set $X^{\ast }:=\overline{X}^{t}$). Since $k$ is embedded into $%
\mathfrak{B}$ by $\iota $, the conjugation on $\mathfrak{B}$ induces on $%
\iota (k)$ the only non-trivial automorphism (cf. \cite{Vig}, Corollaire
2.2, page 6);\ we therefore have:%
\begin{equation*}
\iota _{0}(\overline{\sqrt{\alpha }})=-i_{0}(\sqrt{\alpha })=i_{0}(\sqrt{%
\alpha })^{\ast },
\end{equation*}%
so that $\lambda _{0}\circ i_{0}(b)=\hat{\imath}_{0}(b)\circ \lambda _{0}$
for any $b\in \mathcal{O}_{k,(p)}$, and $\lambda _{0}$ is a principal
polarization for $(A_{0},i_{0})$ (equivalently, $i_{0}(b^{\ast
})=i_{0}(b)^{\dag }$ for any $b\in \mathcal{O}_{k,(p)}$, where $\dag $
denotes the Rosati involution associated to $\lambda _{0}$).

Fix an ordered basis $\{t_{1},...,t_{g}\}$\ for the $\mathcal{\overline{%
\mathbb{F}}}_{p}$-vector space $\limfunc{Lie}A_{0}=\left( \limfunc{Lie}%
E_{0}\right) ^{g}$ such that $\left\{ t_{i}\right\} $ is (the natural image
of) a basis for the Lie algebra of the simple $i^{th}$ factor of $%
A_{0}=E_{0}^{g}$ ($1\leq i\leq g$). With respect to $\{t_{1},...,t_{g}\}$,
by construction $\limfunc{Lie}\iota _{0}(\sqrt{\alpha })$ acts on $\limfunc{%
Lie}A_{0}$ via the matrix:%
\begin{equation*}
\begin{pmatrix}
-\sqrt{\alpha }(\func{mod}p)\cdot I_{r} & 0_{r,s} \\ 
0_{s,r} & \sqrt{\alpha }(\func{mod}p)\cdot I_{s}%
\end{pmatrix}%
\in GL_{g}(\mathbb{F}_{p^{2}})\subset GL_{g}(\mathcal{\overline{\mathbb{F}}}%
_{p}).
\end{equation*}

\noindent We conclude that the fixed pair $(A_{0},i_{0})$ satisfies
Kottwitz' determinant condition, since we have the following equalities in $%
\mathbb{F}_{p^{2}}$: 
\begin{eqnarray*}
\det (X_{1}+\sqrt{\alpha }X_{2};\limfunc{Lie}A_{0}) &=&\det 
\begin{pmatrix}
(X_{1}-\sqrt{\alpha }X_{2})I_{r} & 0_{r,s} \\ 
0_{s,r} & (X_{1}+\sqrt{\alpha }X_{2})I_{s}%
\end{pmatrix}%
(\func{mod}p) \\
&=&(X_{1}-\sqrt{\alpha }X_{2})^{r}(X_{1}+\sqrt{\alpha }X_{2})^{s}(\func{mod}%
p)= \\
&=&f(X_{1},X_{2})\text{ }(\func{mod}p).
\end{eqnarray*}

\noindent We fix a $U(N)$-orbit of an isomorphism $\alpha _{0}:H_{1}(A_{0},%
\mathbb{A}_{f}^{p})\rightarrow V\otimes _{%
%TCIMACRO{\U{211a} }%
%BeginExpansion
\mathbb{Q}
%EndExpansion
}\mathbb{A}_{f}^{p}$ of skew-hermitian modules with $k$-action. By
definition of our moduli variety, we have determined a point:%
\begin{equation*}
\lbrack (A_{0},i_{0},\overline{\lambda }_{0},\overline{\alpha }_{0})]\in 
\mathcal{S}_{\mathcal{D},U(N)}(\mathcal{\overline{\mathbb{F}}}_{p})
\end{equation*}

\noindent that we consider \textit{fixed} for the remaining of this section.
In particular, we have a canonical Hecke-equivariant isomorphism associated
to $(\mathbf{X},i_{\mathbf{X}},\overline{\lambda }_{\mathbf{X}})$ from
Corollary \ref{unif-cor}:%
\begin{equation*}
\Theta _{U(N)}^{\prime }(\mathcal{\overline{\mathbb{F}}}_{p}):I(%
%TCIMACRO{\U{211a} }%
%BeginExpansion
\mathbb{Q}
%EndExpansion
)\backslash \left( J(%
%TCIMACRO{\U{211a} }%
%BeginExpansion
\mathbb{Q}
%EndExpansion
_{p})/J(%
%TCIMACRO{\U{2124} }%
%BeginExpansion
\mathbb{Z}
%EndExpansion
_{p})\mathcal{\times }G(\mathbb{A}_{f}^{p})/U(N)\right) \rightarrow
Z^{\prime }(\mathcal{\overline{\mathbb{F}}}_{p}).
\end{equation*}

\noindent Recall that both sides above are finite sets.

\subsubsection{The groups $I$, $J$ and $G\label{I, J, G and phi}$}

We recall the nature of the algebraic groups $I,J,G$ appearing in the domain
of the uniformization morphism. Recall that $G=GU_{g}(k;r,s)$ is the
reductive group over $%
%TCIMACRO{\U{211a} }%
%BeginExpansion
\mathbb{Q}
%EndExpansion
$ associated to the quadratic imaginary field extension $k/%
%TCIMACRO{\U{211a} }%
%BeginExpansion
\mathbb{Q}
%EndExpansion
$ and the Hermitian form determine by the matrix $H$; furthermore: 
\begin{equation*}
U(N):=\ker \left( GU_{g}\left( \mathcal{O}_{k}\otimes _{%
%TCIMACRO{\U{2124} }%
%BeginExpansion
\mathbb{Z}
%EndExpansion
}\widehat{%
%TCIMACRO{\U{2124} }%
%BeginExpansion
\mathbb{Z}
%EndExpansion
}^{p};r,s\right) \rightarrow GU_{g}\left( \mathcal{O}_{k}\otimes _{%
%TCIMACRO{\U{2124} }%
%BeginExpansion
\mathbb{Z}
%EndExpansion
}\frac{\widehat{%
%TCIMACRO{\U{2124} }%
%BeginExpansion
\mathbb{Z}
%EndExpansion
}^{p}}{N\widehat{%
%TCIMACRO{\U{2124} }%
%BeginExpansion
\mathbb{Z}
%EndExpansion
}^{p}};r,s\right) \right)
\end{equation*}%
is a compact open subgroup of $G(\mathbb{A}_{f}^{p})$.

By definition, $I(%
%TCIMACRO{\U{211a} }%
%BeginExpansion
\mathbb{Q}
%EndExpansion
)$ is the group of auto-$\mathcal{O}_{k,(p)}$-quasi-isogenies of the
homogeneously principally polarized abelian variety $(A_{0},\overline{%
\lambda }_{0})$; let us denote by $\Phi $ the matrix:%
\begin{equation*}
\Phi :=\Phi _{\alpha }=%
\begin{pmatrix}
-\varphi _{\alpha }I_{r} & 0_{r,s} \\ 
0_{s,r} & \varphi _{\alpha }I_{s}%
\end{pmatrix}%
\in M_{g}(\mathfrak{R}).
\end{equation*}%
Then we have $I(%
%TCIMACRO{\U{211a} }%
%BeginExpansion
\mathbb{Q}
%EndExpansion
)=\{X\in GU_{g}(\mathfrak{B};I_{g}):X\Phi =\Phi X\},$ \noindent and for any $%
%TCIMACRO{\U{211a} }%
%BeginExpansion
\mathbb{Q}
%EndExpansion
$-algebra $S$ we have $I(S)=\{X\in GU_{g}(\mathfrak{B}\otimes _{%
%TCIMACRO{\U{211a} }%
%BeginExpansion
\mathbb{Q}
%EndExpansion
}S;I_{g}):X\Phi =\Phi X\}.${}

On the other side, $J(%
%TCIMACRO{\U{211a} }%
%BeginExpansion
\mathbb{Q}
%EndExpansion
_{p})$ is the group of $K_{0}$-automorphisms of the homogeneously
principally polarized isocrystal with $k$-action $(N,\mathbf{F;}i,%
%TCIMACRO{\U{211a} }%
%BeginExpansion
\mathbb{Q}
%EndExpansion
_{p}^{\times }\Psi )$; it has a compact subgroup $J(%
%TCIMACRO{\U{2124} }%
%BeginExpansion
\mathbb{Z}
%EndExpansion
_{p})$ that is given by the $W$-automorphisms of the homogeneously
principally polarized Dieudonn\'{e} module $(\mathbf{M,}i_{\mathbf{M}},%
\overline{e}_{\mathbf{M}})$ endowed with $\mathcal{O}_{k,(p)}$-action. Since 
$\mathbf{M}=M_{\ast }(A_{0}(p))\simeq A_{1/2}^{\oplus g}$ as principally
polarized Dieudonn\'{e} modules, where $A_{1/2}^{\oplus g}$ is endowed with
the product polarization coming from the polarization $\left( \QTATOP{0}{-1}%
\QTATOP{1}{0}\right) $ on $A_{1/2}$, one deduces (\cite{Gh04}, Cor. 10),
keeping track of the action of $\mathcal{O}_{k,(p)}$:%
\begin{eqnarray*}
J(%
%TCIMACRO{\U{211a} }%
%BeginExpansion
\mathbb{Q}
%EndExpansion
_{p}) &\simeq &\{X\in GU_{g}(\mathfrak{B}_{p};I_{g}):X\Phi =\Phi X\}=I(%
%TCIMACRO{\U{211a} }%
%BeginExpansion
\mathbb{Q}
%EndExpansion
_{p}), \\
J(%
%TCIMACRO{\U{2124} }%
%BeginExpansion
\mathbb{Z}
%EndExpansion
_{p}) &\simeq &\{X\in GU_{g}(\mathfrak{R}_{p};I_{g}):X\Phi =\Phi X\}=:I(%
%TCIMACRO{\U{2124} }%
%BeginExpansion
\mathbb{Z}
%EndExpansion
_{p}),
\end{eqnarray*}

\noindent where $\mathfrak{R}_{p}$ denotes the unique maximal order of the
skew-field $\mathfrak{B}_{p}$. This result can also be deduced from \cite{RZ}
6.29; the above isomorphisms are canonical.

We deduce from Corollary \ref{unif-cor}:

\begin{proposition}
There is a canonical Hecke-equivariant isomorphism of finite sets:%
\begin{equation*}
\Theta _{U(N)}^{\prime }(\mathcal{\overline{\mathbb{F}}}_{p}):I(%
%TCIMACRO{\U{211a} }%
%BeginExpansion
\mathbb{Q}
%EndExpansion
)\backslash I(\mathbb{A}_{f})/U_{N}\rightarrow Z^{\prime }(\mathcal{%
\overline{\mathbb{F}}}_{p}),
\end{equation*}

\noindent where $U_{N}:=J(%
%TCIMACRO{\U{2124} }%
%BeginExpansion
\mathbb{Z}
%EndExpansion
_{p})\times U(N)$ is viewed as an open compact subgroup of $I(\mathbb{A}%
_{f})=J(%
%TCIMACRO{\U{211a} }%
%BeginExpansion
\mathbb{Q}
%EndExpansion
_{p})\times G(\mathbb{A}_{f}^{p})$.
\end{proposition}

If $\mathbb{L}_{q}$ is a finite field of cardinality $q$, there is, up to
isomorphism, a unique Hermitian space of dimension $m$ associated to the
quadratic extension $\mathbb{L}_{q^{2}}/\mathbb{L}_{q}$ (cf. \cite{L}). We
denote the associated unitary group - defined over $\mathbb{L}_{q}$ - by $%
U_{m}$; in particular: 
\begin{equation*}
GU_{m}(\mathbb{L}_{q^{2}})=\{X\in GL_{m}(\mathbb{L}_{q^{2}}):X^{\ast
}X=c(X)\cdot I_{m},c(X)\in \mathbb{L}_{q}^{\times }\}.
\end{equation*}%
We will also need to consider the algebraic group $G(U_{m_{1}}\times
U_{m_{2}})\subset GU_{m_{1}}\times GU_{m_{2}}$ defined over $\mathbb{L}_{q}$%
, whose $\mathbb{L}_{q}$-points are: 
\begin{eqnarray*}
G(U_{m_{1}}\times U_{m_{2}})(\mathbb{L}_{q^{2}}) &=&\left\{ g=\left( 
\begin{array}{cc}
X & 0_{m_{1},m_{2}} \\ 
0_{m_{2},m_{1}} & Y%
\end{array}%
\right) :g\in GU_{m_{1}+m_{2}}(\mathbb{L}_{q^{2}})\right\} = \\
&=&\left\{ \left( 
\begin{array}{cc}
X & 0_{m_{1},m_{2}} \\ 
0_{m_{2},m_{1}} & Y%
\end{array}%
\right) :X^{\ast }X=cI_{r},Y^{\ast }Y=cI_{s},c\in \mathbb{L}_{q}^{\times
}\right\} .
\end{eqnarray*}

(Notice the slight abuse of notation). Via the embedding $\iota
:k\hookrightarrow \mathfrak{B}$ that we fixed, we obtain a natural
epimorphism $\mathfrak{R}_{p}\twoheadrightarrow \mathbb{F}_{p^{2}}$, indeed
we can write $\mathfrak{R}_{p}=%
%TCIMACRO{\U{2124} }%
%BeginExpansion
\mathbb{Z}
%EndExpansion
_{p}[\varphi _{\alpha }]\oplus 
%TCIMACRO{\U{2124} }%
%BeginExpansion
\mathbb{Z}
%EndExpansion
_{p}[\varphi _{\alpha }]\Pi $ for a choice of uniformizer $\Pi $ such that $%
\Pi ^{2}=p$, so that: $\frac{\mathfrak{R}_{p}}{(\Pi )}\cong \frac{%
%TCIMACRO{\U{2124} }%
%BeginExpansion
\mathbb{Z}
%EndExpansion
_{p}[\varphi _{\alpha }]}{(p)}\overset{\iota }{\underset{\simeq }{%
\longrightarrow }}\frac{%
%TCIMACRO{\U{2124} }%
%BeginExpansion
\mathbb{Z}
%EndExpansion
_{p}[\sqrt{\alpha }]}{(p)}=\mathbb{F}_{p^{2}}.$

\begin{lemma}
\label{G(p)}Let $G(p):=G(U_{r}\times U_{s})(\mathbb{F}_{p^{2}})$; there is a
short exact sequence of groups (defining $U_{p}$):%
\begin{equation*}
1\rightarrow U_{p}\rightarrow J(%
%TCIMACRO{\U{2124} }%
%BeginExpansion
\mathbb{Z}
%EndExpansion
_{p})\overset{\pi }{\rightarrow }G(p)\rightarrow 1\text{,}
\end{equation*}

\noindent where the map $\pi :J(%
%TCIMACRO{\U{2124} }%
%BeginExpansion
\mathbb{Z}
%EndExpansion
_{p})\rightarrow G(p)$ is induced by the canonical epimorphism $\mathfrak{R}%
_{p}\twoheadrightarrow \mathbb{F}_{p^{2}}$ arising from our fixed embedding $%
\iota :k\hookrightarrow \mathfrak{B}$.
\end{lemma}

\textbf{Proof}. By previous considerations, we have a natural identification:%
\begin{equation*}
J(%
%TCIMACRO{\U{2124} }%
%BeginExpansion
\mathbb{Z}
%EndExpansion
_{p})=\{X\in GU_{g}(\mathfrak{R}_{p};I_{g}):X\Phi =\Phi X\}.
\end{equation*}

\noindent Via the embedding $\iota $, we identify $\varphi _{\alpha }(\func{%
mod}\Pi )$ with $\sqrt{\alpha }(\func{mod}p)\in \mathbb{F}_{p^{2}}$;\ if $%
X\in J(%
%TCIMACRO{\U{2124} }%
%BeginExpansion
\mathbb{Z}
%EndExpansion
_{p})$, then $\pi (X)\in GU_{g}(\mathbb{F}_{p^{2}})$ and the equation $X\Phi
=\Phi X$ for an $(r,s)$-block matrix $X=\left( \QTATOP{A}{C}\QTATOP{B}{D}%
\right) \in GU_{g}(\mathfrak{R}_{p};I_{g})$ reduces to the equation in $%
M_{g}(\mathbb{F}_{p^{2}})$:%
\begin{equation*}
\left( \QTATOP{A}{C}\QTATOP{B}{D}\right) (\func{mod}\Pi )\cdot \left( 
\QTATOP{-\sqrt{\alpha }\cdot I_{r}}{0_{s,r}}\QTATOP{0_{r,s}}{\sqrt{\alpha }%
\cdot I_{s}}\right) (\func{mod}p)=\left( \QTATOP{-\sqrt{\alpha }\cdot I_{r}}{%
0_{s,r}}\QTATOP{0_{r,s}}{\sqrt{\alpha }\cdot I_{s}}\right) (\func{mod}%
p)\cdot \left( \QTATOP{A}{C}\QTATOP{B}{D}\right) (\func{mod}\Pi )
\end{equation*}

\noindent We deduce that $B(\func{mod}\Pi )=0_{r,s}$ and $C(\func{mod}\Pi
)=0_{s,r}$, since $p\neq 2$, so that $\pi (X)\in G(p).$

On the other side, if $Y=\left( \QTATOP{T}{0_{s,r}}\QTATOP{0_{r,s}}{S}%
\right) \in G(p)$, then we can lift $Y$ to a matrix in $G(U_{r}\times U_{s})(%
\mathcal{O}_{k_{\nu }})$, and we can see $G(U_{r}\times U_{s})(\mathcal{O}%
_{k_{\nu }})\subset J(%
%TCIMACRO{\U{2124} }%
%BeginExpansion
\mathbb{Z}
%EndExpansion
_{p})$ via $\iota $, so that $\pi $ is onto.\ $\blacksquare $

\subsection{Unitary Dieudonn\'{e}\ modules and invariant differentials}

We will need the following general fact:

\begin{lemma}
\label{Dieudonne pairing}Let $L$ be a subfield of $\overline{\mathbb{F}}_{p}$
and let $\sigma $ denote the restriction of the Frobenius of $W$ to $W(L)$.
Let $M$ be a Diudonn\'{e} module over $W(L)$ endowed with the $%
%TCIMACRO{\U{2124} }%
%BeginExpansion
\mathbb{Z}
%EndExpansion
_{(p)}$-linear action of a $%
%TCIMACRO{\U{2124} }%
%BeginExpansion
\mathbb{Z}
%EndExpansion
_{(p)}$-algebra with involution $O$; assume $M$ comes with a principal
polarization $e:M\times M\rightarrow W(L)$ that is skew-Hermitian with
respect to the $O$-action. Then the assignment:%
\begin{equation*}
\left\langle ,\right\rangle :\frac{M}{FM}\times \frac{M}{VM}\rightarrow L,\
\ (\overline{x},\overline{y})\mapsto e(x,Fy)(\func{mod}p)
\end{equation*}

\noindent \noindent is a well-defined, perfect pairing which is $L$-linear
in the first variable and $L$-semilinear (with respect to $\sigma $) in the
second variable. Furthermore, $\left\langle ,\right\rangle $ is
skew-Hermitian with respect to the action of $O$.
\end{lemma}

\begin{proposition}
If $F+V=0$ on $M$, then $\left\langle ,\right\rangle $ defines a $\sigma
^{-1}$-alternating pairing $\frac{M}{VM}\times \frac{M}{VM}\rightarrow L$.
\end{proposition}

\textbf{Proof. }For $x,y,m,m^{\prime }\in M$ we have: 
\begin{eqnarray*}
e(x+Fm,F(y+Vm^{\prime })) &=&e(x,Fy)+e(Fm,Fy)+e(x+Fm,FVm^{\prime })= \\
&=&e(x,Fy)+e(m,VFy)^{\sigma }+e(x+Fm,pm^{\prime })\equiv \\
&\equiv &e(x,Fy)(\func{mod}p),
\end{eqnarray*}

\noindent so that $\left\langle ,\right\rangle $ is well defined; it is
clearly $L$-linear in the first variable, and $\sigma $-semilinear in the
second since $F$ is $\sigma $-semilinear. If $b\in O$ and $x,y\in M$ we have:%
\begin{eqnarray*}
\left\langle b\overline{x},\overline{y}\right\rangle &=&\left\langle 
\overline{bx},y\right\rangle =e(bx,Fy)(\func{mod}p)= \\
&=&e(x,b^{\ast }Fy)(\func{mod}p)=e(x,Fb^{\ast }y)(\func{mod}p)= \\
&=&\left\langle \overline{x},b^{\ast }\overline{y}\right\rangle .
\end{eqnarray*}

\bigskip

To show that $\left\langle ,\right\rangle $ is non-degenerate we need to
show that the $L$-linear map:%
\begin{equation*}
\varepsilon :\frac{M}{FM}\longrightarrow \limfunc{Hom}\nolimits_{\sigma 
\text{-}lin}(\frac{M}{VM},L)
\end{equation*}

\noindent induced by $\left\langle ,\right\rangle $ is an isomorphism of $L$%
-vector spaces. Let $x\in M$ such that $\left\langle \overline{x},\overline{y%
}\right\rangle =0$ for all $y\in M$. By assumption we have $e(y,Vx)\in pW(L)$
for all $y\in M$; If we denote by $\mu :M\rightarrow \limfunc{Hom}%
\nolimits_{W}(M,W)$ the isomorphism defined by $\mu (m):=e(\cdot ,m)$, we
have that $\mu (Vx)$ has image contained inside $pW(L)$, so that there is $%
z\in M$ such that $\frac{1}{p}\mu (Vx)=\mu (z)$; we obtain $Vx=pz$, so that $%
px=FVx=pFz$, hence $x=Fz\in FM.$\noindent\ We deduce that $\varepsilon $ is
injective. Similarly one can show that $\left\langle ,\right\rangle $
induces an injective $\sigma $-semilinear map of $L$-space:%
\begin{equation*}
\frac{M}{VM}\hookrightarrow \limfunc{Hom}\nolimits_{L}(\frac{M}{FM},L).
\end{equation*}

\noindent We conclude that $\dim _{L}\frac{M}{VM}=\dim _{L}\frac{M}{FM}%
<\infty $ and $\varepsilon $ is forced to be an isomorphism.

Finally, if $F+V=0$, we have $FM=VM$ and if $x,y\in M$ one computes:%
\begin{eqnarray*}
\left\langle \overline{x},\overline{y}\right\rangle &=&e(x,Fy)(\func{mod}%
p)=-e(x,Vy)(\func{mod}p)= \\
&=&-e(Fx,y)^{\sigma ^{-1}}(\func{mod}p)=e(y,Fx)^{\sigma ^{-1}}(\func{mod}p)=
\\
&=&\left\langle \overline{y},\overline{x}\right\rangle ^{\sigma ^{-1}}\text{.%
}
\end{eqnarray*}

\noindent This says that $\left\langle ,\right\rangle $ is $\sigma ^{-1}$%
-alternating. $\blacksquare $

\bigskip

Fix a triple $(A_{0},i,\overline{\lambda })$ such that $[(A_{0},i,\overline{%
\lambda },\overline{\alpha })]\in Z^{\prime }(\overline{\mathbb{F}}_{p})$
for some level structure $\overline{\alpha }$; we know that $(A_{0},i,%
\overline{\lambda })$ has a canonical $\mathbb{F}_{p^{2}}$-structure, that
we denote by $(A_{0}^{\prime },i,\overline{\lambda }^{\prime })$; by
functoriality, the various object that we associated to $(A_{0},i,\overline{%
\lambda })$ (as $p$-divisible groups, Dieudonn\'{e} modules, polarizations,
actions of algebras) are obtained as base changes to $\overline{\mathbb{F}}%
_{p}$ (or to $W$) of analogous objects defined over $\mathbb{F}_{p^{2}}$ (or
over $W(\mathbb{F}_{p^{2}})$). It is therefore equivalent to work with $%
(A_{0}^{\prime },i,\overline{\lambda }^{\prime })$ or with $(A_{0},i,%
\overline{\lambda })$.

Let $\mathbf{X}^{\prime }:=A_{0}^{\prime }(p)$ be the $p$-divisible group of 
$A_{0}^{\prime }$: it is defined over $\mathbb{F}_{p^{2}}$; its covariant
Dieudonn\'{e} module $\mathbf{M}^{\prime }=M_{\ast }(\mathbf{X}^{\prime })$
is a $W(\mathbb{F}_{p^{2}})$-module, and $\mathbf{M}^{\prime
}=(A_{1/2}^{\prime })^{g}$, (cf. \ref{superspecial-supersingular});
furthermore $\mathbf{M=M}^{\prime }\otimes _{W(\mathbb{F}_{p^{2}})}W$. The
symbols $i_{\mathbf{M}^{\prime }}$ and $e_{\mathbf{M}^{\prime }}$ have the
obvious meanings, relatively to the given action $i$ and polarization $%
\lambda ^{\prime }$.

By \cite{BBM} 3.3.1, there is a positive integer $m$ such that the canonical
map of cotangent spaces (at the origin) 
\begin{equation*}
\mathfrak{t}_{A_{0}^{\prime }(p)}^{\ast }=\mathfrak{t}_{\mathbf{X}^{\prime
}}^{\ast }\longrightarrow \mathfrak{t}_{A_{0}^{\prime }[p^{m}]}^{\ast }
\end{equation*}

\noindent is an isomorphism. Furthermore, the closed immersion $%
A_{0}^{\prime }[p^{m}]\hookrightarrow A_{0}^{\prime }$ of $\mathbb{F}%
_{p^{2}} $-group schemes induces an epimorphism of $\mathbb{F}_{p^{2}}$%
-vector spaces $\mathfrak{t}_{A_{0}^{\prime }}^{\ast }\longrightarrow 
\mathfrak{t}_{A_{0}^{\prime }[p^{m}]}^{\ast }=\mathfrak{t}_{\mathbf{X}%
^{\prime }}^{\ast }.$ \noindent Since $A_{0}^{\prime }$ is superspecial, we
have $\dim \mathbf{X}^{\prime }=g=\dim A_{0}^{\prime }$, so that we obtain
canonical identifications $\mathfrak{t}_{A_{0}^{\prime }}^{\ast }=\mathfrak{t%
}_{\mathbf{X}^{\prime }}^{\ast }$ and $\limfunc{Lie}A_{0}^{\prime }=\limfunc{%
Lie}\mathbf{X}^{\prime }.$

By covariant Dieudonn\'{e} theory, we have a canonical isomorphism of $%
\mathbb{F}_{p^{2}}$-vector spaces:%
\begin{equation*}
\limfunc{Lie}(\mathbf{X}^{\prime })=\frac{\mathbf{M}^{\prime }}{V\mathbf{M}%
^{\prime }}.
\end{equation*}

All the above isomorphisms and identifications respect the $\mathcal{O}%
_{k,(p)}$-structures of the modules considered, and also the polarizations
induced by $\lambda $.

By our assumptions, $F+V=0$ on $\mathbf{M}^{\prime }$, so that by
Proposition \ref{Dieudonne pairing}, the principal polarization $e_{\mathbf{M%
}^{\prime }}:\mathbf{M}^{\prime }\times \mathbf{M}^{\prime }\rightarrow W(%
\mathbb{F}_{p^{2}})$ induces a non-degenerate pairing of $\mathbb{F}_{p^{2}}$%
-spaces:%
\begin{equation*}
\left\langle ,\right\rangle :\frac{\mathbf{M}^{\prime }}{V\mathbf{M}^{\prime
}}\times \frac{\mathbf{M}^{\prime }}{V\mathbf{M}^{\prime }}\rightarrow 
\mathbb{F}_{p^{2}}
\end{equation*}

\noindent which is linear in the first argument, $\sigma $-linear in the
second argument, and $\sigma $-alternating (i.e. $\left\langle
x,y\right\rangle =\left\langle y,x\right\rangle ^{\sigma }$; recall that $%
\sigma ^{2}=1$ here); hence $(\frac{\mathbf{M}^{\prime }}{V\mathbf{M}%
^{\prime }};i_{\mathbf{M}^{\prime }},\left\langle ,\right\rangle )$ is a 
\textit{Hermitian space over }$\mathbb{F}_{p^{2}}$\textit{\ of dimension }$%
g=2n$\textit{, endowed with an action of }$\mathcal{O}_{k,(p)}$ with respect
to which the pairing $\left\langle ,\right\rangle $ is skew-symmetric. Since 
$e_{\mathbf{M}^{\prime }}$ is determined only up to a constant in $%
%TCIMACRO{\U{2124} }%
%BeginExpansion
\mathbb{Z}
%EndExpansion
_{p}^{\times }$, the pairing $\left\langle ,\right\rangle $ is determined up
to a constant in $\mathbb{F}_{p}^{\times }$.

We need to study a bit more the Dieudonn\'{e} module $\mathbf{M}^{\prime }$
(cf. for example \cite{Yu} \S 3, or \cite{BW06} 2.1). Via the fixed
embedding $\nu $ we obtain the decomposition:%
\begin{equation*}
\mathbf{M}^{\prime }=\mathbf{M}_{-}^{\prime }\oplus \mathbf{M}_{+}^{\prime },
\end{equation*}

\noindent where:%
\begin{equation*}
\mathbf{M}_{\pm }^{\prime }:=\{m\in \mathbf{M}^{\prime }:i_{\mathbf{M}%
^{\prime }}(\sqrt{\alpha })m=\pm \sqrt{\alpha }m\},\text{ \ }rk_{W(\mathbb{F}%
_{p^{2}})}\mathbf{M}_{\pm }^{\prime }=g\text{\ .}
\end{equation*}

\noindent One easily sees that $V\mathbf{M}_{\pm }^{\prime }\subseteq 
\mathbf{M}_{\mp }^{\prime }$, $F\mathbf{M}_{\pm }^{\prime }\subseteq \mathbf{%
M}_{\mp }^{\prime }$ and that $e_{\mathbf{M}^{\prime }}(\mathbf{M}%
_{+}^{^{\prime }},\mathbf{M}_{+}^{^{\prime }})=0,e_{\mathbf{M}^{\prime }}(%
\mathbf{M}_{-}^{^{\prime }},\mathbf{M}_{-}^{^{\prime }})=0$ (i.e. $\mathbf{M}%
_{-}^{\prime }$ and $\mathbf{M}_{+}^{\prime }$ are totally isotropic with
respect to the principal polarization). We have:%
\begin{equation*}
\frac{\mathbf{M}^{\prime }}{V\mathbf{M}^{\prime }}=\frac{\mathbf{M}%
_{-}^{\prime }}{V\mathbf{M}_{+}^{\prime }}\oplus \frac{\mathbf{M}%
_{+}^{\prime }}{V\mathbf{M}_{-}^{\prime }}.
\end{equation*}

\noindent Notice that this is a decomposition as $\mathbb{F}_{p^{2}}$-vector
spaces with action of $\mathcal{O}_{k,(p)}$, where $\sqrt{\alpha }$ acts as $%
-\sqrt{\alpha }(\func{mod}p)\in \mathbb{F}_{p^{2}}$ on the first summand,
and as $\sqrt{\alpha }(\func{mod}p)$ on the second. Furthermore:

\begin{equation*}
\dim _{\mathbb{F}_{p^{2}}}\frac{\mathbf{M}_{+}^{\prime }}{V\mathbf{M}%
_{+}^{\prime }}=r\text{, \ }\dim _{\mathbb{F}_{p^{2}}}\frac{\mathbf{M}%
_{+}^{\prime }}{V\mathbf{M}_{-}^{\prime }}=s.
\end{equation*}

\begin{proposition}
\label{determination of aut}Let $(A_{0},i,\overline{\lambda })$ be a triple
over $\overline{\mathbb{F}}_{p}$ such that for some level structure $%
\overline{\alpha }$ we have $[(A_{0},i,\overline{\lambda },\overline{\alpha }%
)]\in Z^{\prime }(\overline{\mathbb{F}}_{p})$; let $(A_{0}^{\prime },i,%
\overline{\lambda }^{\prime })$ be the canonical $\mathbb{F}_{p^{2}}$%
-structure of $(A_{0},i,\overline{\lambda })$. The automorphism group of the
Hermitian space with $\mathcal{O}_{k,(p)}$-action $(\frac{\mathbf{M}^{\prime
}}{V\mathbf{M}^{\prime }};i_{\mathbf{M}^{\prime }},\left\langle
,\right\rangle )$ is isomorphic to the finite group $G(p)=G(U_{r}\times
U_{s})(\mathbb{F}_{p^{2}})$.
\end{proposition}

\textbf{Proof. }Let $\mathfrak{L}_{\pm }:=\frac{\mathbf{M}_{\pm }^{\prime }}{%
V\mathbf{M}_{\mp }^{\prime }}$. Let $\mathcal{B}^{+}$ (resp. $\mathcal{B}%
^{-} $)\ be a fixed ordered basis of $\mathfrak{L}_{+}$ (resp. $\mathfrak{L}%
_{-}$). If $\mathfrak{X}$ is an automorphism of $\frac{\mathbf{M}^{\prime }}{%
V\mathbf{M}^{\prime }}$ which commutes with the action of $\mathcal{O}%
_{k,(p)} $, we have $\mathfrak{X}$ $\mathfrak{L}_{\pm }\subseteq \mathfrak{L}%
_{\pm }$, so that the matrix representing $\mathfrak{X}$ with respect to $%
\mathcal{B}:=\mathcal{B}^{-}\mathcal{\cup B}^{+}$ is of the form:%
\begin{equation*}
X=\left( 
\begin{array}{cc}
X_{-} & 0_{r,s} \\ 
0_{s,r} & X_{+}%
\end{array}%
\right) \in GL_{g}(\mathbb{F}_{p^{2}}).
\end{equation*}

\noindent Viceversa, any such matrix represents - with respect to $\mathcal{B%
}$ - an automorphism of $\frac{\mathbf{M}^{\prime }}{V\mathbf{M}^{\prime }}$
commuting with the action of $\mathcal{O}_{k,(p)}$.

Since $F\mathbf{M}_{+}^{\prime }\subseteq \mathbf{M}_{-}^{\prime }$ and $e_{%
\mathbf{M}^{\prime }}(\mathbf{M}_{-}^{^{\prime }},\mathbf{M}_{-}^{^{\prime
}})=0$, we deduce that $\left\langle \mathfrak{L}_{-},\mathfrak{L}%
_{+}\right\rangle =0$, by definition of the pairing $\left\langle
,\right\rangle $. This implies that $\left\langle ,\right\rangle $ is
represented, with respect to $\mathcal{B}$, by a Hermitian diagonal matrix $%
\left( \QDATOP{U_{-}}{0_{s,r}}\QDATOP{0_{r,s}}{U_{+}}\right) \in GL_{g}(%
\mathbb{F}_{p^{2}})$, so that if $X$ is as above, we have that $X$
represents an automorphism of $(\frac{\mathbf{M}^{\prime }}{V\mathbf{M}%
^{\prime }};i_{\mathbf{M}^{\prime }},\left\langle ,\right\rangle )$ with
respect to $\mathcal{B}$ if and only if:%
\begin{equation*}
X_{\pm }^{\ast }\cdot U_{\pm }\cdot X_{\pm }=cU_{\pm },
\end{equation*}

\noindent where $c\in \mathbb{F}_{p}^{\times }$ is a scalar depending only
on $X$.

We conclude that the automorphism group of $(\frac{\mathbf{M}^{\prime }}{V%
\mathbf{M}^{\prime }};i_{\mathbf{M}^{\prime }},\left\langle ,\right\rangle )$
is isomorphic to the group:%
\begin{equation*}
\mathcal{G=\{(}X_{-},X_{+}\mathcal{)}\in GU_{r}(\mathbb{F}%
_{p^{2}};U_{-})\times GU_{s}(\mathbb{F}%
_{p^{2}};U_{+}):c_{-}(X_{-})=c_{+}(X_{+})\mathcal{\}}\text{,}
\end{equation*}

\noindent where $c_{\pm }$ is the similitude factor homomorphism in $%
GU_{\cdot }(\mathbb{F}_{p^{2}};U_{\pm })$. The unitary spaces $(\mathbb{F}%
_{p^{2}}^{r};U_{-})$ and $(\mathbb{F}_{p^{2}}^{r};I_{r})$ are isomorphic,
hence we can find an isomorphism $GU_{r}(\mathbb{F}_{p^{2}};U_{-})\simeq
GU_{r}(\mathbb{F}_{p^{2}})$ preserving the similitude factor of
corresponding matrices in each group; similarly we can do for $GU_{s}(%
\mathbb{F}_{p^{2}};U_{+})$. Putting things together we obtain an isomorphism 
$\mathcal{G\simeq }G(p)$. $\blacksquare $

\bigskip

We now need to switch to cotangent spaces in our considerations. As usual, $%
\mathfrak{t}_{A_{0}^{\prime }}^{\ast }$ denotes the cotangent space (at the
origin) of $A_{0}^{\prime }$. As vector spaces over $\mathbb{F}_{p^{2}}$, we
have:%
\begin{equation*}
\mathfrak{t}_{A_{0}^{\prime }}^{\ast }=\limfunc{Hom}\nolimits_{\mathbb{F}%
_{p^{2}}}\left( \frac{\mathbf{M}^{\prime }}{V\mathbf{M}^{\prime }},\mathbb{F}%
_{p^{2}}\right) \text{.}
\end{equation*}

\noindent For simplicity, let $\mathfrak{L}:=\frac{\mathbf{M}^{\prime }}{V%
\mathbf{M}^{\prime }}=\mathfrak{L}_{-}\oplus \mathfrak{L}_{+}$, where $%
\mathfrak{L}_{\pm }:=\frac{\mathbf{M}_{\pm }^{\prime }}{V\mathbf{M}_{\mp
}^{\prime }}$, so that $\mathfrak{t}_{A_{0}^{\prime }}^{\ast }=\mathfrak{L}%
^{\ast }$. The action of $\mathcal{O}_{k,(p)}$ on $\mathfrak{L}$ induces by
functoriality an algebra homomorphism:%
\begin{equation*}
i^{\vee }:\mathcal{O}_{k,(p)}^{opp}=\mathcal{O}_{k,(p)}\longrightarrow 
\limfunc{End}\nolimits_{\mathbb{F}_{p^{2}}}(\mathfrak{L}^{\ast })
\end{equation*}

\noindent defined by $i^{\vee }(b)(\eta ):=\eta \circ i(b)$ for all $b\in 
\mathcal{O}_{k,(p)}$ and $\eta \in \mathfrak{t}_{A_{0}^{\prime }}^{\ast }$.
Notice that $\sqrt{\alpha }\in \mathcal{O}_{k,(p)}$ acts on $\mathfrak{L}%
_{\pm }^{\ast }$ as $\pm \sqrt{\alpha }(\func{mod}p)$, via $i^{\vee }$.

\noindent The non-degenerate Hermitian pairing $\left\langle ,\right\rangle $
on $\mathfrak{L}$ induces a $\sigma $-semilinear isomorphism of $\mathbb{F}%
_{p^{2}}$-spaces $\varepsilon :\mathfrak{L\rightarrow L}^{\ast }$ by setting 
$\varepsilon _{v}:w\mapsto \left\langle w,v\right\rangle $ for all $v,w\in 
\mathfrak{L}$. This allows us to define a pairing $\left( ,\right) $ on $%
\mathfrak{L}^{\ast }$ by setting $(\varepsilon _{v_{1}},\varepsilon
_{v_{2}}):=\left\langle v_{1},v_{2}\right\rangle $ for all $v_{1},v_{2}\in 
\mathfrak{L}$. One can check that we have obtained a non-degenerate pairing $%
\left( ,\right) :\mathfrak{L}^{\ast }\times \mathfrak{L}^{\ast
}\longrightarrow \mathbb{F}_{p^{2}},$ which is $\sigma $-semilinear in the
first variable, linear in the second, and such that $\left( \eta _{1},\eta
_{2}\right) =(\eta _{2},\eta _{1})^{\sigma }$ for all $\eta _{1},\eta
_{2}\in \mathfrak{L}^{\ast }$. Furthermore $(i^{\vee }(b)\eta _{1},\eta
_{2})=(\eta _{1},i^{\vee }(\overline{b})\eta _{2})$ for all $b\in \mathcal{O}%
_{k,(p)}$, and $\left( \mathfrak{L}_{-}^{\ast },\mathfrak{L}_{+}^{\ast
}\right) =0$. We have therefore obtained a $\mathbb{F}_{p^{2}}$-Hermitian
space:%
\begin{equation*}
(\mathfrak{t}_{A_{0}^{\prime }}^{\ast }=\mathfrak{L}^{\ast };i^{\vee
},\left( ,\right) )
\end{equation*}

\noindent of dimension $g$, endowed with an action $i^{\vee }$ of $\mathcal{O%
}_{k,(p)}$ with respect to which the pairing $\left( ,\right) $ is
skew-Hermitian.

\begin{lemma}
There is an isomorphism of groups:%
\begin{equation*}
\limfunc{Aut}\nolimits_{\mathbb{F}_{p^{2}}}(\mathfrak{t}_{A_{0}^{\prime
}}^{\ast };i^{\vee },\left( ,\right) )\simeq G(p).
\end{equation*}
\end{lemma}

\textbf{Proof. }The result follows from Proposition \ref{determination of
aut}, since the map $\mathfrak{X\mapsto (X}^{\ast })^{-1}$ defines a
canonical isomorphism of groups $\limfunc{Aut}\nolimits_{\mathbb{F}_{p^{2}}}(%
\mathfrak{t}_{A_{0}^{\prime }}^{\ast };i^{\vee },\left( ,\right)
)\rightarrow \limfunc{Aut}\nolimits_{\mathbb{F}_{p^{2}}}(\limfunc{Lie}%
A_{0}^{\prime };i,\left\langle ,\right\rangle )$. $\blacksquare $

\bigskip

We can now give the following:

\begin{definition}
\label{inv-diff}Let $(A_{0},i,\overline{\lambda })$ be a triple over $%
\overline{\mathbb{F}}_{p}$ such that for some level structure $\overline{%
\alpha }$ we have $[(A_{0},i,\overline{\lambda },\overline{\alpha })]\in
Z^{\prime }(\overline{\mathbb{F}}_{p})$; let $(A_{0}^{\prime },i,\overline{%
\lambda }^{\prime })$ be the canonical $\mathbb{F}_{p^{2}}$-structure of $%
(A_{0},i,\overline{\lambda })$. A \textbf{basis of invariant differentials
of }$(A_{0}^{\prime },i,\overline{\lambda }^{\prime })$\textbf{\ }(over $%
\mathbb{F}_{p^{2}}$) is a choice of an ordered (similitude) Hermitian basis $%
\eta =(\eta _{-},\eta _{+})$ of the Hermitian module $(\mathfrak{t}%
_{A_{0}^{\prime }}^{\ast };i^{\vee },\left( ,\right) )$ such that $\eta
_{\pm }$ is a basis for $(\mathfrak{t}_{A_{0}^{\prime }}^{\ast })_{\pm
}:=\left( \frac{\mathbf{M}_{\pm }^{\prime }}{V\mathbf{M}_{\mp }^{\prime }}%
\right) ^{\ast }$.
\end{definition}

We have:

\begin{lemma}
Let $(A_{0}^{\prime },i,\overline{\lambda }^{\prime })$ be as above:

\begin{description}
\item[(a)] there is a basis of $\mathfrak{t}_{A_{0}^{\prime }}^{\ast }$ with
respect to which the automorphisms of $(\mathfrak{t}_{A_{0}^{\prime }}^{\ast
};i^{\vee },\left( ,\right) )$ are represented by the matrices of $%
G(p)=G(U_{r}\times U_{s})(\mathbb{F}_{p^{2}})$;

\item[(b)] there are basis of invariant differentials for $(A_{0}^{\prime
},i,\overline{\lambda }^{\prime })$;

\item[(c)] let $\mathcal{B}$ be a basis of $\mathfrak{t}_{A_{0}^{\prime
}}^{\ast }$ as in $(a)$ and let $\eta \in \mathbb{F}_{p^{2}}^{g}$ be the
coordinate column vector of a basis of invariant differentials for $%
(A_{0}^{\prime },i,\overline{\lambda }^{\prime })$. Then any other
coordinate vector (with respect to $\mathcal{B}$) of a basis of invariant
differentials for $(A_{0}^{\prime },i,\overline{\lambda }^{\prime })$ is of
the form $M\eta $ for some $M\in G(p)$.
\end{description}
\end{lemma}

\subsection{Superspecial modular forms}

We assume fixed from now on a basis $\eta _{0}$ of invariant differentials
for $(A_{0}^{\prime },i_{0},\overline{\lambda }_{0}^{\prime })$.

\begin{definition}
The symbol $Z_{diff}^{\prime }(\overline{\mathbb{F}}_{p})$ denotes the set
of equivalence classes of tuples $(A,i,\overline{\lambda },\overline{\alpha }%
,\eta )$, where $(A,i,\overline{\lambda },\overline{\alpha })$ is a
representative for an equivalence class $[(A,i,\overline{\lambda },\overline{%
\alpha })]\in Z^{\prime }(\overline{\mathbb{F}}_{p})$, and $\eta $ is a
choice of basis of invariant differentials for the triple $(A^{\prime
},i^{\prime },\overline{\lambda }^{\prime })$ defined over $\mathbb{F}%
_{p^{2}}$. Two tuples $(A,i,\overline{\lambda },\overline{\alpha },\eta )$
and $(A_{1},i_{1},\overline{\lambda }_{1},\overline{\alpha }_{1},\eta _{1})$
are equivalent if there is an isomorphism $f:(A,i,\overline{\lambda },%
\overline{\alpha })\rightarrow (A_{1},i_{1},\overline{\lambda }_{1},%
\overline{\alpha }_{1})$ such that $f^{\ast }(\eta _{1})=\eta $, where $%
f^{\ast }:\mathfrak{t}_{A_{1}}^{\ast }\rightarrow \mathfrak{t}_{A}^{\ast }$
is the cotangent map induced by $f$.
\end{definition}

\begin{remark}
Let $g\in J(%
%TCIMACRO{\U{2124} }%
%BeginExpansion
\mathbb{Z}
%EndExpansion
_{p})\subset GL_{g}(\mathfrak{R}_{p})$ and let $v\in \limfunc{Lie}%
A_{0}^{\prime }\simeq \mathbb{F}_{p^{2}}^{g}$, $\omega \in \mathfrak{t}%
_{A_{0}^{\prime }}^{\ast }=\left( \limfunc{Lie}A_{0}^{\prime }\right) ^{\ast
}$. Then $g$ acts on $v$ and on $\omega $ as follows:%
\begin{eqnarray*}
g\cdot v &:&=g(\func{mod}\Pi )v \\
g\cdot \omega &:&=\omega \circ g(\func{mod}\Pi ),
\end{eqnarray*}

\noindent where $\Pi $ is a uniformizer for $\mathfrak{R}_{p}$.
\end{remark}

\begin{proposition}
\label{alg-ss}For any fixed choice of left transversal $\mathcal{Y}$ (resp. $%
\mathcal{G}$) for $J(%
%TCIMACRO{\U{2124} }%
%BeginExpansion
\mathbb{Z}
%EndExpansion
_{p})$ in $J(%
%TCIMACRO{\U{211a} }%
%BeginExpansion
\mathbb{Q}
%EndExpansion
_{p})$ (resp. $U_{p}$ in $J(%
%TCIMACRO{\U{2124} }%
%BeginExpansion
\mathbb{Z}
%EndExpansion
_{p})$), the uniformization morphism for the superspecial locus $\Theta
^{\prime }(\overline{\mathbb{F}}_{p})$ induces a Hecke-equivariant
isomorphism $\Theta _{diff}^{\prime }(\overline{\mathbb{F}}_{p}):=\Theta
_{diff}^{\prime }(\overline{\mathbb{F}}_{p})_{U(N)}$:%
\begin{equation*}
\Theta _{diff}^{\prime }(\overline{\mathbb{F}}_{p}):I(%
%TCIMACRO{\U{211a} }%
%BeginExpansion
\mathbb{Q}
%EndExpansion
)\backslash \left( J(%
%TCIMACRO{\U{211a} }%
%BeginExpansion
\mathbb{Q}
%EndExpansion
_{p})/U_{p}\mathcal{\times }G(\mathbb{A}_{f}^{p})/U(N)\right) \rightarrow
Z_{diff}^{\prime }(\mathcal{\overline{\mathbb{F}}}_{p}).
\end{equation*}
\end{proposition}

\textbf{Proof. }(In this proof, for $\xi \in I(%
%TCIMACRO{\U{211a} }%
%BeginExpansion
\mathbb{Q}
%EndExpansion
)$ we will sometimes write $\xi $ to denote $\alpha _{0,p}(\xi )\in J(%
%TCIMACRO{\U{211a} }%
%BeginExpansion
\mathbb{Q}
%EndExpansion
_{p})$, or viceversa, if no ambiguity arises, cf. \ref{groups I and J}). Let 
$I(%
%TCIMACRO{\U{211a} }%
%BeginExpansion
\mathbb{Q}
%EndExpansion
)\cdot \left( \rho ^{-1}U_{p}\mathcal{\times }xU(N)\right) $ be a fixed
element in the above left hand side;\ there are uniquely determined $%
y^{-1}\in \mathcal{Y}$ and $g^{-1}\in \mathcal{G}$ such that $\rho
^{-1}U_{p}=y^{-1}g^{-1}U_{p}$. By the definition of $\Theta ^{\prime }(%
\overline{\mathbb{F}}_{p})$, we obtain a well defined tuple $(y_{\ast
}A_{0},y_{\ast }i_{0},y_{\ast }\overline{\lambda }_{0},x^{-1}\cdot y_{\ast }%
\overline{\alpha }_{0})$ representing a class in $Z^{\prime }(\overline{%
\mathbb{F}}_{p})$. Since the $p$-divisible group of $(y_{\ast }A_{0},y_{\ast
}i_{0},y_{\ast }\overline{\lambda }_{0})$ coincides with $(\mathbf{X},i_{%
\mathbf{X}},\overline{\lambda }_{\mathbf{X}})$, $\eta _{0}$ is a basis of
invariant differentials for the $\mathbb{F}_{p^{2}}$-model of $(y_{\ast
}A_{0},y_{\ast }i_{0},y_{\ast }\overline{\lambda }_{0})$, via the canonical
identification: 
\begin{equation*}
\left( \limfunc{Lie}y_{\ast }A_{0}^{\prime }\right) ^{\ast }=\left( \limfunc{%
Lie}\mathbf{X}^{\prime }\right) ^{\ast };
\end{equation*}%
therefore $\eta _{0}g$ is also a basis of invariant differentials for the $%
\mathbb{F}_{p^{2}}$-model of $(y_{\ast }A_{0},y_{\ast }i_{0},y_{\ast }%
\overline{\lambda }_{0})$. We set:%
\begin{equation*}
\Theta _{diff}^{\prime }(\overline{\mathbb{F}}_{p}):I(%
%TCIMACRO{\U{211a} }%
%BeginExpansion
\mathbb{Q}
%EndExpansion
)\cdot \left( y^{-1}g^{-1}U_{p}\mathcal{\times }xU(N)\right) \longmapsto
\lbrack (y_{\ast }A_{0},y_{\ast }i_{0},y_{\ast }\overline{\lambda }%
_{0},x^{-1}\cdot y_{\ast }\overline{\alpha }_{0},\eta _{0}g)].
\end{equation*}

Notice the above assignment is well defined, once we fixed the transversals $%
\mathcal{Y}$ and $\mathcal{G}$: we only need to check that the map
constructed factors through $I(%
%TCIMACRO{\U{211a} }%
%BeginExpansion
\mathbb{Q}
%EndExpansion
)$. Let $\xi \in I(%
%TCIMACRO{\U{211a} }%
%BeginExpansion
\mathbb{Q}
%EndExpansion
)$, so that $\xi _{\ast }A_{0}=A_{0}.$ Then $\xi
y^{-1}g^{-1}U_{p}=y_{1}^{-1}g_{1}^{-1}U_{p}$ for uniquely determined $%
y_{1}^{-1}\in \mathcal{Y}$ and $g_{1}^{-1}\in \mathcal{G}$ (notice that here 
$\xi y^{-1}g^{-1}$ should be more properly be written as $\alpha _{0,p}(\xi
)y^{-1}g^{-1}$). Write $y_{1}=f\cdot y\xi ^{-1}$ with $f=g_{1}^{-1}ug\in J(%
%TCIMACRO{\U{2124} }%
%BeginExpansion
\mathbb{Z}
%EndExpansion
_{p})$, for some $u\in U_{p}$. The isomorphism $f$ induces an isomorphism:%
\begin{equation*}
f_{ab}:(y_{\ast }A_{0},y_{\ast }i_{0},y_{\ast }\overline{\lambda }_{0})%
\overset{\simeq }{\rightarrow }(y_{1,\ast }A_{0},y_{1,\ast }i_{0},y_{1,\ast }%
\overline{\lambda }_{0}).
\end{equation*}

\noindent By definition of the action of $I(%
%TCIMACRO{\U{211a} }%
%BeginExpansion
\mathbb{Q}
%EndExpansion
)$ on $G(\mathbb{A}_{f}^{p})/U(N)$ we have:%
\begin{equation*}
\xi \cdot xU(N)=\alpha _{0}^{p}(\xi )xU(N)=\alpha _{0}\circ H_{1}(\xi )\circ
\alpha _{0}^{-1}\circ xU(N).
\end{equation*}

\noindent The level structure on $(y_{1,\ast }A_{0},y_{1,\ast
}i_{0},y_{1,\ast }\overline{\lambda }_{0})$ associated to the pair $%
y_{1}^{-1}g_{1}^{-1}U_{p}\times \alpha _{0}^{p}(\xi )xU(N)$ is therefore
induced by$:$%
\begin{eqnarray*}
x^{-1}\alpha _{0}^{p}(\xi ^{-1})y_{1,\ast }(\alpha _{0}) &=&x^{-1}\alpha
_{0}H_{1}(\xi ^{-1})\alpha _{0}^{-1}\alpha _{0}H_{1}(y_{1}^{-1})= \\
&=&x^{-1}\alpha _{0}H_{1}(y^{-1}f^{-1})= \\
&=&x^{-1}y_{\ast }(\alpha _{0})\cdot H_{1}(f^{-1}),
\end{eqnarray*}

\noindent so that $f_{ab}$ is an isomorphism 
\begin{equation*}
f_{ab}:(y_{\ast }A_{0},y_{\ast }i_{0},y_{\ast }\overline{\lambda }%
_{0},x^{-1}\cdot y_{\ast }\overline{\alpha }_{0})\overset{\simeq }{%
\rightarrow }(y_{1,\ast }A_{0},y_{1,\ast }i_{0},y_{1,\ast }\overline{\lambda 
}_{0},x^{-1}\alpha _{0}^{p}(\xi ^{-1})y_{1,\ast }(\alpha _{0}))
\end{equation*}

\noindent The cotangent map induced by $f_{ab}$ gives: $f_{ab}^{\ast }(\eta
_{0}g_{1})=\eta _{0}g_{1}f=\eta _{0}ug=\eta _{0}g$, so that $f_{ab}$
preserves the choices of invariant differentials.

The map $\Theta _{diff}^{\prime }(\overline{\mathbb{F}}_{p})$ is surjective:
let $[(A,i,\overline{\lambda },\overline{\alpha },\eta )]\in
Z_{diff}^{\prime }(\overline{\mathbb{F}}_{p})$; we can find $y^{-1}\in 
\mathcal{Y}$, $\widetilde{g}^{-1}\in J(%
%TCIMACRO{\U{2124} }%
%BeginExpansion
\mathbb{Z}
%EndExpansion
_{p})$ and $x\in G(\mathbb{A}_{f}^{p})$ such that $(A,i,\overline{\lambda },%
\overline{\alpha },\eta )$ is isomorphic to a tuple of the form:%
\begin{equation*}
(y_{\ast }A_{0},y_{\ast }i_{0},y_{\ast }\overline{\lambda }_{0},x^{-1}\cdot
y_{\ast }\overline{\alpha }_{0},\eta _{0}\widetilde{g}).
\end{equation*}

\noindent Let $g^{-1}\in \mathcal{G}$ such that $\widetilde{g}%
^{-1}U_{p}=g^{-1}U_{p}$; then $[(A,i,\overline{\lambda },\overline{\alpha }%
,\eta )]=\Theta _{diff}^{\prime }(\overline{\mathbb{F}}_{p})(I(%
%TCIMACRO{\U{211a} }%
%BeginExpansion
\mathbb{Q}
%EndExpansion
)\cdot \left( y^{-1}g^{-1}U_{p}\mathcal{\times }xU(N)\right) )$.

The map $\Theta _{diff}^{\prime }(\overline{\mathbb{F}}_{p})$ is injective:
let us fix $y^{-1},y_{1}^{-1}\in \mathcal{Y}$, $g^{-1},g_{1}^{-1}\in 
\mathcal{G}$ and $x,x_{1}\in G(\mathbb{A}_{f}^{p})$ such that there is an
isomorphism: 
\begin{equation*}
f_{ab}:(y_{\ast }A_{0},y_{\ast }i_{0},y_{\ast }\overline{\lambda }%
_{0},x^{-1}\cdot y_{\ast }\overline{\alpha }_{0},\eta _{0}g)\overset{\simeq }%
{\rightarrow }(y_{1,\ast }A_{0},y_{1,\ast }i_{0},y_{1,\ast }\overline{%
\lambda }_{0},x_{1}^{-1}\cdot y_{1,\ast }\overline{\alpha }_{0},\eta
_{0}g_{1}).
\end{equation*}

\noindent Denote by $f\in J(%
%TCIMACRO{\U{2124} }%
%BeginExpansion
\mathbb{Z}
%EndExpansion
_{p})$ the automorphism induced by $f_{ab}$ on $(\mathbf{X},i_{\mathbf{X}},%
\overline{\lambda }_{\mathbf{X}})$ and let $\xi =(y^{-1}f^{-1}y_{1})_{ab}\in
I(%
%TCIMACRO{\U{211a} }%
%BeginExpansion
\mathbb{Q}
%EndExpansion
)$ be the auto-quasi-isogeny of $(A_{0},i_{0},\overline{\lambda }_{0})$
inducing $y^{-1}f^{-1}y_{1}$ on $(\mathbf{X},i_{\mathbf{X}},\overline{%
\lambda }_{\mathbf{X}})$. Since $f$ is an isomorphism, we have:

\begin{description}
\item[(i)] $x_{1}^{-1}y_{1,\ast }(\overline{\alpha }_{0})\cdot
H_{1}(f)=x^{-1}y_{\ast }(\overline{\alpha }_{0})$, so that $x\equiv \alpha
_{0}H_{1}(y^{-1}f^{-1}y_{1})\alpha _{0}^{-1}x_{1}(\func{mod}U(N))$;

\item[(ii)] $\eta _{0}g_{1}f=\eta _{0}g$, hence $g_{1}f\equiv g(\func{mod}%
U_{p})$.
\end{description}

\noindent Then we have:%
\begin{align*}
I(%
%TCIMACRO{\U{211a} }%
%BeginExpansion
\mathbb{Q}
%EndExpansion
)\cdot \left( y_{1}^{-1}g_{1}^{-1}U_{p}\mathcal{\times }x_{1}U(N)\right) & =
\\
& =I(%
%TCIMACRO{\U{211a} }%
%BeginExpansion
\mathbb{Q}
%EndExpansion
)\cdot \left( \xi y_{1}^{-1}g_{1}^{-1}U_{p}\mathcal{\times \alpha }%
_{0}^{p}(\xi )x_{1}U(N)\right) = \\
& =I(%
%TCIMACRO{\U{211a} }%
%BeginExpansion
\mathbb{Q}
%EndExpansion
)\cdot \left( y^{-1}(g_{1}f)^{-1}U_{p}\mathcal{\times \alpha }%
_{0}H_{1}(y^{-1}f^{-1}y_{1})\alpha _{0}^{-1}x_{1}U(N)\right) \overset{(ii)}{=%
} \\
& =I(%
%TCIMACRO{\U{211a} }%
%BeginExpansion
\mathbb{Q}
%EndExpansion
)\cdot \left( y^{-1}g^{-1}U_{p}\mathcal{\times \alpha }%
_{0}H_{1}(y^{-1}f^{-1}y_{1})\alpha _{0}^{-1}x_{1}U(N)\right) \overset{(i)}{=}
\\
& =I(%
%TCIMACRO{\U{211a} }%
%BeginExpansion
\mathbb{Q}
%EndExpansion
)\cdot \left( y^{-1}g^{-1}U_{p}\mathcal{\times }xU(N)\right) .
\end{align*}

The Hecke-equivariance of $\left\{ \Theta _{diff}^{\prime }(\overline{%
\mathbb{F}}_{p})_{K^{p}}\right\} _{K^{p}}$ (with respect to the projective
systems of domain and codomain obtained by varying the level structures)\ is
an easy consequence of the definition of the Hecke operators in this
context: for $\gamma \in G(\mathbb{A}_{f}^{p})$ let us denote by $H_{\gamma
} $ the corresponding Hecke operator acting on the domain or codomain of $%
\Theta _{diff}^{\prime }(\overline{\mathbb{F}}_{p})$, for a suitable level
subgroup$.$ For $y^{-1}\in \mathcal{Y}$, $g^{-1}\in \mathcal{G}$ and $x\in G(%
\mathbb{A}_{f}^{p})$ we have:%
\begin{eqnarray*}
&&I(%
%TCIMACRO{\U{211a} }%
%BeginExpansion
\mathbb{Q}
%EndExpansion
)\cdot \left( y^{-1}g^{-1}U_{p}\mathcal{\times }xU(N)\right) \overset{%
H_{\gamma }}{\longmapsto }I(%
%TCIMACRO{\U{211a} }%
%BeginExpansion
\mathbb{Q}
%EndExpansion
)\cdot \left( y^{-1}g^{-1}U_{p}\mathcal{\times \gamma }^{-1}x\gamma \cdot
\gamma ^{-1}U(N)\gamma \right) \\
&&\overset{\Theta _{diff}^{\prime }(\overline{\mathbb{F}}_{p})}{\longmapsto }%
[(y_{\ast }A_{0},y_{\ast }i_{0},y_{\ast }\overline{\lambda }_{0},\mathcal{%
\gamma }^{-1}x^{-1}\gamma \cdot y_{\ast }(\gamma ^{-1}\overline{\alpha }%
_{0}),\eta _{0}g); \\
&&I(%
%TCIMACRO{\U{211a} }%
%BeginExpansion
\mathbb{Q}
%EndExpansion
)\cdot \left( y^{-1}g^{-1}U_{p}\mathcal{\times }xU(N)\right) \overset{\Theta
_{diff}^{\prime }(\overline{\mathbb{F}}_{p})}{\longmapsto }[(y_{\ast
}A_{0},y_{\ast }i_{0},y_{\ast }\overline{\lambda }_{0},x^{-1}\cdot y_{\ast }(%
\overline{\alpha }_{0}),\eta _{0}g) \\
&&\overset{H_{\gamma }}{\longmapsto }[(y_{\ast }A_{0},y_{\ast }i_{0},y_{\ast
}\overline{\lambda }_{0},\mathcal{\gamma }^{-1}x^{-1}\cdot y_{\ast }(%
\overline{\alpha }_{0}),\eta _{0}g),
\end{eqnarray*}

\noindent and $\mathcal{\gamma }^{-1}x^{-1}\gamma \cdot y_{\ast }(\gamma
^{-1}\overline{\alpha }_{0})=\mathcal{\gamma }^{-1}x^{-1}\cdot y_{\ast }(%
\overline{\alpha }_{0})$. $\blacksquare $

\bigskip

The isomorphism $\Theta _{diff}^{\prime }(\overline{\mathbb{F}}_{p})$
described above depends upon the choices of transversal $\mathcal{G}$ and $%
\mathcal{Y}$. We assume from now on that such transversal have been fixed.

For an algebraic $\mathbb{F}_{p}$-representation $\rho :GL_{g}\rightarrow
GL_{d}$ of the group $GL_{g}$, we denote by $M_{\rho }(N;\overline{\mathbb{F}%
}_{p}):=M_{\rho }(\mathcal{D}_{(r,s),p}^{U};\overline{\mathbb{F}}_{p})$ the $%
\overline{\mathbb{F}}_{p}$-vector space of unitary modular forms $(\func{mod}%
p)$ of signature $(r,s)$\ for the field $k$, having weight $\rho $ and level 
$N$. Denote by $\iota :Z^{\prime }\hookrightarrow \mathcal{S}_{\mathcal{D}%
,U(N)}\otimes \overline{\mathbb{F}}_{p}$ the closed immersion of $\overline{%
\mathbb{F}}_{p}$-schemes associated to the inclusion of sets $Z^{\prime }(%
\overline{\mathbb{F}}_{p})\subset \mathcal{S}_{\mathcal{D},U(N)}(\overline{%
\mathbb{F}}_{p})$. We want to consider "restrictions" of modular forms to
the superspecial locus, as in \cite{Gh04}, 4. At this purpose, let $\tau
:G(U_{r}\times U_{s})\rightarrow GL_{m}$ be an algebraic $\mathbb{F}_{p}$%
-representation of the $\mathbb{F}_{p}$-group $G(U_{r}\times U_{s})$.

\begin{definition}
\label{def sspecial mod forms}The space $M_{\tau }^{ss}(N;\overline{\mathbb{F%
}}_{p})$ of $(r,s)$-\textbf{unitary superspecial modular forms }$(\func{mod}%
p)$\textbf{\ }for the field\textbf{\ }$k$, having weight $\tau $ and level $%
N $ is the finite dimensional $\mathcal{\overline{\mathbb{F}}}_{p}$-vector
space whose elements $f$ are rules that assign, to any tuple $(A,i,\overline{%
\lambda },\overline{\alpha },\eta )/\overline{\mathbb{F}}_{p}$ such that $%
[(A,i,\overline{\lambda },\overline{\alpha })]$ is an element of $Z^{\prime
}(\overline{\mathbb{F}}_{p})$ and $\eta $ is an ordered basis of invariant
differentials for $(A^{\prime },i,\overline{\lambda }^{\prime })$, an
element $f(A,i,\overline{\lambda },\overline{\alpha },\eta )\in \overline{%
\mathbb{F}}_{p}^{m}$ in such a way that:
\end{definition}

\begin{description}
\item[(a)] $f(A,i,\overline{\lambda },\overline{\alpha },\eta M)=\tau
(M)^{-1}f(A,i,\overline{\lambda },\overline{\alpha },\eta )$ for all $M\in
G(p)\simeq \limfunc{Aut}\nolimits_{\mathbb{F}_{p^{2}}}(\mathfrak{t}%
_{A^{\prime }}^{\ast };i^{\vee },\left( ,\right) );$

\item[(b)] if $(A,i,\overline{\lambda },\overline{\alpha },\eta )\simeq
(A_{1},i_{1},\overline{\lambda }_{1},\overline{\alpha }_{1},\eta _{1})$ then 
$f(A,i,\overline{\lambda },\overline{\alpha },\eta )=f(A_{1},i_{1},\overline{%
\lambda }_{1},\overline{\alpha }_{1},\eta _{1})$.
\end{description}

\noindent We have another description of the above space:

\begin{proposition}
\label{down to earth sspec} For any algebraic $\mathbb{F}_{p}$%
-representation $\rho :GL_{g}\rightarrow GL_{d}$, denote by $\limfunc{Res}%
\rho $ its restriction to $G(U_{r}\times U_{s})$. Then:%
\begin{equation*}
M_{\limfunc{Res}\rho }^{ss}(N;\overline{\mathbb{F}}_{p})=H^{0}(Z^{\prime }(%
\overline{\mathbb{F}}_{p}),\iota ^{\ast }\mathbb{E}_{\rho }).
\end{equation*}
\end{proposition}

\textbf{Proof. }By Proposition \ref{down to earth mod forms}, if $f\in
H^{0}(Z^{\prime }(\overline{\mathbb{F}}_{p}),\iota ^{\ast }\mathbb{E}_{\rho
})$, then $f$ satisfies $(a)$ and $(b)$ of Def. \ref{def sspecial mod forms}%
. Viceversa, let $f\in M_{\limfunc{Res}\rho }^{ss}(N;\overline{\mathbb{F}}%
_{p})$ so that $f$ is an assignment on tuples $(A,i,\overline{\lambda },%
\overline{\alpha },\eta )/\overline{\mathbb{F}}_{p}$ as in Def. \ref{def
sspecial mod forms}; in particular $\eta $ here is a basis of invariant
differentials for $(A^{\prime },i,\overline{\lambda }^{\prime })$. If $%
\omega $ is any ordered basis of $\mathfrak{t}_{A}^{\ast }$, there is a
unique $X_{\omega ,\eta }\in GL_{g}(\overline{\mathbb{F}}_{p})$ such that $%
\omega =\eta X_{\omega ,\eta }$ and we define $f(A,i,\overline{\lambda },%
\overline{\alpha },\omega ):=\rho (X_{\omega ,\eta })^{-1}\cdot f(A,i,%
\overline{\lambda },\overline{\alpha },\eta )$. This assignment is well
defined as $\rho $ is a representation of $GL_{g}$, and allows us to view $f$
as an element of $H^{0}(Z^{\prime }(\overline{\mathbb{F}}_{p}),\iota ^{\ast }%
\mathbb{E}_{\rho })$. $\blacksquare $

\bigskip

The definition of $M_{\tau }^{ss}(N;\overline{\mathbb{F}}_{p})$ depends upon 
$Z^{\prime }(\overline{\mathbb{F}}_{p})$, hence upon the choice of $%
(A_{0},i_{0},\overline{\lambda }_{0},\overline{\alpha }_{0})$ that we have
fixed at the beginning. It is clear that the Hecke operators act upon $%
M_{\tau }^{ss}(N;\overline{\mathbb{F}}_{p})$.

\subsubsection{Algebraic modular forms}

We briefly recall the definition of algebraic modular forms $\func{mod}p$ in
our context (cf. \cite{Gr}). \noindent By Proposition \ref{RZ}, we have
identifications $I(%
%TCIMACRO{\U{211a} }%
%BeginExpansion
\mathbb{Q}
%EndExpansion
_{p})=J(%
%TCIMACRO{\U{211a} }%
%BeginExpansion
\mathbb{Q}
%EndExpansion
_{p})$ and $I(\mathbb{A}_{f}^{p})=G(\mathbb{A}_{f}^{p})$; we set $U:=U_{p}%
\mathcal{\times }U(N)$ and view it as an open compact subgroup of $I(\mathbb{%
A}_{f})$. The group $I(%
%TCIMACRO{\U{211a} }%
%BeginExpansion
\mathbb{Q}
%EndExpansion
)$ is discrete inside $I(\mathbb{A}_{f})$, by \cite{RZ}, 6.23, so that the
double coset space $I(%
%TCIMACRO{\U{211a} }%
%BeginExpansion
\mathbb{Q}
%EndExpansion
)\backslash I(\mathbb{A}_{f})/U$ \noindent is finite, because $I(%
%TCIMACRO{\U{211a} }%
%BeginExpansion
\mathbb{Q}
%EndExpansion
)\backslash I(\mathbb{A}_{f})$ is compact (cf. \cite{Gr}, Prop. 1.4). Assume
fixed an algebraic $\mathbb{F}_{p}$-representation $\tau :G(U_{r}\times
U_{s})\rightarrow GL_{m}$.

\begin{definition}
The space of \textbf{algebraic modular forms }$(\func{mod}p)$\textbf{\ for
the group }$I$, having level $U$ and weight $\tau $ is the $\overline{%
\mathbb{F}}_{p}$-vector space (of finite dimension):%
\begin{eqnarray*}
M_{\tau }^{\limfunc{alg}}(N;\overline{\mathbb{F}}_{p}) &:&=\{f:I(%
%TCIMACRO{\U{211a} }%
%BeginExpansion
\mathbb{Q}
%EndExpansion
)\backslash I(\mathbb{A}_{f})/U\rightarrow \overline{\mathbb{F}}_{p}^{m}: \\
&:&f(\overline{g}M)=\tau (M)^{-1}f(\overline{g})\text{, }M\in G(p),\text{ }%
\overline{g}\in I(%
%TCIMACRO{\U{211a} }%
%BeginExpansion
\mathbb{Q}
%EndExpansion
)\backslash I(\mathbb{A}_{f})/U\},
\end{eqnarray*}

\noindent where the right action of $G(p)$ on $I(%
%TCIMACRO{\U{211a} }%
%BeginExpansion
\mathbb{Q}
%EndExpansion
)\backslash I(\mathbb{A}_{f})/U$ is induced by the identification $G(p)=J(%
%TCIMACRO{\U{2124} }%
%BeginExpansion
\mathbb{Z}
%EndExpansion
_{p})/U_{p}$.
\end{definition}

The space $M_{\tau }^{\limfunc{alg}}(N;\overline{\mathbb{F}}_{p})$ is
endowed with a natural action of the Hecke algebra and our previous
computations give:

\begin{proposition}
\label{alg-superspecial}There is an isomorphism of\ finite dimensional $%
\overline{\mathbb{F}}_{p}$-vector spaces endowed with Hecke action: 
\begin{equation*}
M_{\tau }^{\limfunc{alg}}(N;\overline{\mathbb{F}}_{p})\simeq M_{\tau
}^{ss}(N;\overline{\mathbb{F}}_{p}).
\end{equation*}
\end{proposition}

\textbf{Proof. }By Proposition \ref{alg-ss}, we have an isomorphism of $%
M_{\tau }^{\limfunc{alg}}(N;\overline{\mathbb{F}}_{p})$ with the space of
functions $Z_{diff}^{\prime }(\mathcal{\overline{\mathbb{F}}}%
_{p})\rightarrow \overline{\mathbb{F}}_{p}^{m}$ satisfying condition $(a)$
of Def. \ref{def sspecial mod forms}. (Note that if $g^{-1}\in \mathcal{G}$, 
$mU_{p}\in J(%
%TCIMACRO{\U{2124} }%
%BeginExpansion
\mathbb{Z}
%EndExpansion
_{p})/U_{p}$ then by definition $\eta _{0}g\cdot mU_{p}=\eta _{0}m^{-1}g$). $%
\blacksquare $

\bigskip

\subsection{The Hecke eigensystems correspondence\label{hecke correspondence}%
}

To establish the correspondence between Hecke eigensystems coming from
algebraic modular forms and PEL-modular forms in the unitary settings, we
begins by following \cite{Gh04}, 4.

Let $\mathcal{I}$ be the ideal sheaf associated to the closed immersion of $%
\overline{\mathbb{F}}_{p}$-schemes $\iota :Z^{\prime }\hookrightarrow 
\mathcal{S}_{\mathcal{D},U(N)}\otimes \overline{\mathbb{F}}_{p}$, so that
the following is an exact sequence of sheaves over $\mathcal{S}:=\mathcal{S}%
_{\mathcal{D},U(N)}\otimes \overline{\mathbb{F}}_{p}$:%
\begin{equation*}
0\rightarrow \mathcal{I\rightarrow O}_{\mathcal{S}}\rightarrow \iota _{\ast }%
\mathcal{O}_{Z^{\prime }}\rightarrow 0.
\end{equation*}

\noindent \noindent Since $\mathcal{S}$ is of finite type over $\overline{%
\mathbb{F}}_{p}$, it is noetherian, so that $\mathcal{I}$ is a coherent
ideal sheaf on $\mathcal{S}$. \noindent After tensoring the above sequence
with the locally free sheaf of $\mathcal{O}_{\mathcal{S}}$-modules $\mathbb{E%
}_{\rho }$ and then taking cohomology, we obtain the exact sequence:%
\begin{equation*}
0\rightarrow H^{0}(\mathcal{S},\mathcal{I\otimes }\mathbb{E}_{\rho
})\rightarrow H^{0}(\mathcal{S},\mathbb{E}_{\rho })\rightarrow H^{0}(%
\mathcal{S},\iota _{\ast }\mathcal{O}_{Z^{\prime }}\mathcal{\otimes }\mathbb{%
E}_{\rho }).
\end{equation*}

\noindent By \cite{Ha} III, 2.10 we have $H^{0}(\mathcal{S},\iota _{\ast }%
\mathcal{O}_{Z^{\prime }}\mathcal{\otimes }\mathbb{E}_{\rho
})=H^{0}(Z^{\prime },\iota ^{\ast }\mathbb{E}_{\rho })$, since the
projection formula gives $\iota _{\ast }\iota ^{\ast }\mathbb{E}_{\rho
}=\iota _{\ast }\mathcal{O}_{Z^{\prime }}\mathcal{\otimes }\mathbb{E}_{\rho
} $. The above exact sequence of vector spaces can therefore be written, by
Prop. \ref{down to earth sspec}, as:%
\begin{equation*}
0\rightarrow H^{0}(\mathcal{S},\mathcal{I\otimes }\mathbb{E}_{\rho
})\rightarrow M_{\rho }(N;\overline{\mathbb{F}}_{p})\overset{r}{\rightarrow }%
M_{\limfunc{Res}\rho }^{ss}(N;\overline{\mathbb{F}}_{p}),
\end{equation*}

\noindent where the last map $r$ need not to be surjective. \noindent Recall
that $\limfunc{Res}\rho $ is the restriction of $\rho $ to the algebraic
group $G(U_{r}\times U_{s})$. Furthermore, observe that $r$ is
Hecke-equivariant.

\begin{proposition}
\noindent \label{epi}There exists a positive integer $m_{0}$ such that the
above map $r$ is a surjection $M_{\rho \otimes \det^{m}}(N;\overline{\mathbb{%
F}}_{p})\overset{r}{\twoheadrightarrow }M_{\limfunc{Res}(\rho \otimes
\det^{m})}^{ss}(N;\overline{\mathbb{F}}_{p})$ for any $m>m_{0}$.
\end{proposition}

\textbf{Proof. }The invertible sheaf of $\mathcal{O}_{\mathcal{S}}$-modules $%
\tbigwedge\nolimits^{g}\mathbb{E=E}_{\det }$ is ample over $\mathcal{S}$
(cf. for example \cite{Lan}, Th. 7.24.1). The proposition now follows in the
same way as \cite{G}, Prop. 24. $\blacksquare $

\bigskip

With the notation introduced in this section, we have:

\begin{theorem}
\label{main thm}Let $p$ be an odd prime and $k/%
%TCIMACRO{\U{211a} }%
%BeginExpansion
\mathbb{Q}
%EndExpansion
$ be a quadratic imaginary field extension in which $p$ is inert. Let $n$ be
a positive integer and let $r,s$ be non-negative integers such that $%
r+s=g:=2n.$ Let furthermore $N\geq 3$ be an integer not divisible by $p$.
Denote by $I$ the reductive $%
%TCIMACRO{\U{211a} }%
%BeginExpansion
\mathbb{Q}
%EndExpansion
$-group whose $%
%TCIMACRO{\U{211a} }%
%BeginExpansion
\mathbb{Q}
%EndExpansion
$-rational points are given by $I(%
%TCIMACRO{\U{211a} }%
%BeginExpansion
\mathbb{Q}
%EndExpansion
)=\{X\in GU_{g}(\mathfrak{B};I_{g}):X\Phi =\Phi X\}$, where $\mathfrak{B}$
is the quaternion algebra over $%
%TCIMACRO{\U{211a} }%
%BeginExpansion
\mathbb{Q}
%EndExpansion
$ ramified at $p$ and $\infty $, and $\Phi $ is as in \ref{I, J, G and phi}.

The systems of Hecke eigenvalues coming from $(r,s)$-unitary PEL-modular
forms $(\func{mod}p)$\textbf{\ }for the quadratic imaginary field\textbf{\ }$%
k$, having genus $g$, fixed level $N$ and any possible weight $\rho
:GL_{g}\rightarrow GL_{m(\rho )}$, are the same as the systems of Hecke
eigenvalues coming from $(\func{mod}p)$ algebraic modular forms for the
group $I$ having level $U=U_{p}\times U(N)\subset I(\mathbb{A}_{f})$ and any
possible weight $\rho ^{\prime }:G(U_{r}\times U_{s})\rightarrow
GL_{m^{\prime }(\rho ^{\prime })}$.
\end{theorem}

\textbf{Proof. }To prove the result we first show that \textit{any} system
of Hecke eigenvalues occurring in the spaces $\{M_{\rho }(N;\overline{%
\mathbb{F}}_{p})\}_{\rho }$ for variable $\rho :GL_{g}\rightarrow GL_{m(\rho
)}$ also occurs in the spaces $\{M_{\tau }^{ss}(N;\overline{\mathbb{F}}%
_{p})\}_{\tau }$ for variable $\tau :G(U_{r}\times U_{s})\rightarrow
GL_{m(\tau )}$. Then we follow the proof of Th. 28 in\ \cite{G}, to show
that the converse is also true, and finally we apply Proposition \ref%
{alg-superspecial}. Notice that the first part of this proof is different
from the one given in \cite{G}.

For any integer $i\geq 0$ we have an exact sequence of $\mathcal{O}_{%
\mathcal{S}}$-modules:%
\begin{equation*}
0\rightarrow \mathcal{I}^{i+1}\mathcal{\rightarrow I}^{i}\rightarrow 
\mathcal{I}^{i}/\mathcal{I}^{i+1}\rightarrow 0
\end{equation*}

\noindent giving rise to the exact sequence in cohomology:%
\begin{equation*}
0\rightarrow H^{0}(\mathcal{S},\mathcal{I}^{i+1}\otimes \mathbb{E}_{\rho
})\rightarrow H^{0}(\mathcal{S},\mathcal{I}^{i}\otimes \mathbb{E}_{\rho })%
\overset{r_{i}}{\rightarrow }H^{0}(\mathcal{S},\mathcal{I}^{i}/\mathcal{I}%
^{i+1}\otimes \mathbb{E}_{\rho }),
\end{equation*}

\noindent which defines the homomorphisms $r_{i}$ for any $i\geq 0$ (notice
that $r_{0}=r$).

For any $j\geq 1$ we also have the exact sequence of sheaves of $\mathcal{O}%
_{\mathcal{S}}$-modules:%
\begin{equation*}
\mathcal{I\otimes I}^{j}/\mathcal{I}^{j+1}\mathcal{\rightarrow O}_{\mathcal{S%
}}\otimes \mathcal{I}^{j}/\mathcal{I}^{j+1}\rightarrow \iota _{\ast }%
\mathcal{O}_{Z^{\prime }}\otimes \mathcal{I}^{j}/\mathcal{I}%
^{j+1}\rightarrow 0.
\end{equation*}

\noindent Notice that the image of the first map is zero, so that we obtain
isomorphisms of $\mathcal{O}_{\mathcal{S}}$-modules $\mathcal{I}^{j}/%
\mathcal{I}^{j+1}\simeq \iota _{\ast }\mathcal{O}_{Z^{\prime }}\otimes 
\mathcal{I}^{j}/\mathcal{I}^{j+1}$ for any $j\geq 1$. In cohomology we have
therefore canonical isomorphisms $\xi _{j}$ for any $j\geq 1$:%
\begin{equation*}
\xi _{j}:H^{0}(\mathcal{S},\mathcal{I}^{j}/\mathcal{I}^{j+1}\otimes \mathbb{E%
}_{\rho })\overset{\simeq }{\rightarrow }H^{0}(Z^{\prime },\iota ^{\ast }(%
\mathcal{I}^{j}/\mathcal{I}^{j+1}\otimes \mathbb{E}_{\rho })).
\end{equation*}

Let $0\neq f\in M_{\rho }(N;\overline{\mathbb{F}}_{p})=H^{0}(\mathcal{S},%
\mathbb{E}_{\rho })$ be a Hecke eigenform for some fixed weight $\rho
:GL_{g}\rightarrow GL_{m}$. If $r(f)\neq 0$ then the system of Hecke
eigenvalues associated to $f$ occurs in $M_{\limfunc{Res}\rho }^{ss}(N;%
\overline{\mathbb{F}}_{p})$, since $r$ is Hecke-equivariant. If $r(f)=0$,
then $f\in H^{0}(\mathcal{S},\mathcal{I\otimes }\mathbb{E}_{\rho })$. We
claim that there is a \textit{positive} integer $h$ such that $f\in H^{0}(%
\mathcal{S},\mathcal{I}^{h}\otimes \mathbb{E}_{\rho })$ and $r_{h}(f)\neq 0$%
.\ Assume not, then $r_{1}(f)=0$ and $f\in H^{0}(\mathcal{S},\mathcal{I}%
^{2}\otimes \mathbb{E}_{\rho })$, so that we can compute $r_{2}(f)$ and we
need to have $r_{2}(f)=0$; therefore $f\in H^{0}(\mathcal{S},\mathcal{I}%
^{3}\otimes \mathbb{E}_{\rho })$, etc. Reiterating this procedure we deduce
that $f\in H^{0}(\mathcal{S},\mathcal{I}^{i}\otimes \mathbb{E}_{\rho })$ for
all integers $i>0$ (recall that $...\subseteq H^{0}(\mathcal{S},\mathcal{I}%
^{i+1}\otimes \mathbb{E}_{\rho })\subseteq H^{0}(\mathcal{S},\mathcal{I}%
^{i}\otimes \mathbb{E}_{\rho })\subseteq ...\subseteq H^{0}(\mathcal{S},%
\mathcal{I}\otimes \mathbb{E}_{\rho })$), so that $f=0$, contradicting our
assumption $f\neq 0$. Then there exists a positive integer $h$ such that $%
f\in H^{0}(\mathcal{S},\mathcal{I}^{h}\otimes \mathbb{E}_{\rho })$ and $%
r_{h}(f)\neq 0$. Let:%
\begin{equation*}
f^{ss}:=\xi _{h}(r_{h}(f))\in H^{0}(Z^{\prime },\iota ^{\ast }(\mathcal{I}%
^{h}/\mathcal{I}^{h+1}\otimes \mathbb{E}_{\rho })).
\end{equation*}%
Since $\xi _{h}$ is injective, $f^{ss}$ is non-zero. Observe that $\mathcal{I%
}^{h}/\mathcal{I}^{h+1}=\limfunc{Sym}^{h}(\mathcal{I}/\mathcal{I}^{2})$ and
that $\iota ^{\ast }(\mathcal{I}/\mathcal{I}^{2})=\iota ^{\ast }(\Omega _{%
\mathcal{S}}^{1})$ (cf. \cite{Ha} II, 8.17).

\noindent We now need to use the Kodaira-Spencer isomorphism for our PEL
variety : let $\mathcal{X}_{/\overline{\mathbb{F}}_{p}}$ denotes the
universal abelian scheme over $\mathcal{S}$, endowed with the polarization $%
\lambda _{univ}$ and action $i_{univ}$ of the ring $\mathcal{O}_{k,(p)}$;
let $\widehat{\mathcal{X}}_{/\overline{\mathbb{F}}_{p}}$ denotes the dual
abelian scheme. By Prop. 2.3.4.2. in \cite{Lan}, the Kodaira-Spencer map
induces an isomorphism of $\mathcal{O}_{\mathcal{S}}$-sheaves:%
\begin{equation*}
\underline{KS}:\frac{\mathbb{E}_{\mathcal{X}}\mathbb{\otimes }_{\mathcal{O}_{%
\mathcal{S}}}\mathbb{E}_{\widehat{\mathcal{X}}}}{J^{\prime }}\longrightarrow
\Omega _{_{\mathcal{S}}}^{1},
\end{equation*}

\noindent where:%
\begin{equation*}
J^{\prime }=\left( 
\begin{array}{c}
\lambda _{univ}^{\ast }(y)\otimes z-\lambda _{univ}^{\ast }(z)\otimes y \\ 
i_{univ}(b)^{\ast }(x)\otimes y-x\otimes \widehat{i}_{univ}(b)^{\ast }(y)%
\end{array}%
:x\in \mathbb{E}_{\mathcal{X}};y,z\in \mathbb{E}_{\widehat{\mathcal{X}}%
};b\in \mathcal{O}_{k,(p)}\right) ;
\end{equation*}

\noindent (the map $\lambda _{univ}^{\ast }$ denotes the pull-back
isomorphism $\lambda _{univ}^{\ast }:\mathbb{E}_{\widehat{\mathcal{X}}%
}\rightarrow \mathbb{E}_{\mathcal{X}}$, and $\widehat{i}_{univ}(b)^{\ast }$
denotes the endomorphism of $\mathbb{E}_{\widehat{\mathcal{X}}}$ induced by $%
\widehat{i}_{univ}(b)$). \noindent Precomposing $\underline{KS}$ with $%
id\otimes \lambda ^{\ast }$ we get the isomorphism of sheaves:%
\begin{equation*}
\frac{\limfunc{Sym}^{2}\mathbb{E}}{J}\longrightarrow \Omega _{_{\mathcal{S}%
}}^{1},
\end{equation*}

\noindent where $J=\left( i_{univ}(b)^{\ast }(x)\otimes y-x\otimes
i_{univ}(b)^{\ast }(y):x,y\in \mathbb{E},b\in \mathcal{O}_{k,(p)}\right) $
and $\mathbb{E}:\mathbb{=E}_{\mathcal{X}}$ as usual. Write $\limfunc{Sym}%
_{\flat }^{2}\mathbb{E}:=\left( \limfunc{Sym}^{2}\mathbb{E}\right) /J$ and
notice that $J$ is not preserved by the group $GL_{g}$ (if $rs>0$), but it
has an action of $GL_{r}\times GL_{s}$, because of the determinant condition
imposed in the definition of the moduli problem.

\noindent We have:%
\begin{equation*}
H^{0}(Z^{\prime },\iota ^{\ast }(\mathcal{I}^{h}/\mathcal{I}^{h+1}\otimes 
\mathbb{E}_{\rho }))=H^{0}(Z^{\prime },\iota ^{\ast }(\limfunc{Sym}%
\nolimits^{h}\left( \limfunc{Sym}\nolimits_{\flat }^{2}\mathbb{E}\right) 
\mathbb{\otimes E}_{\rho })).
\end{equation*}

\noindent The group$\ GL_{r}\times GL_{s}$ acts on the sheaf $\iota ^{\ast }(%
\limfunc{Sym}\nolimits^{h}\left( \limfunc{Sym}\nolimits_{\flat }^{2}\mathbb{E%
}\right) \mathbb{\otimes E}_{\rho })$, and this is enough for our purposes,
as the space of superspecial modular forms is defined for representations $%
\tau $ of $G(U_{r}\times U_{s})\subset GL_{r}\times GL_{s}$. We conclude
that: 
\begin{equation*}
f^{ss}\in H^{0}(Z^{\prime },\iota ^{\ast }(\limfunc{Sym}\nolimits^{h}\left( 
\limfunc{Sym}\nolimits_{\flat }^{2}\mathbb{E}\right) \mathbb{\otimes E}%
_{\rho }))=M_{\limfunc{Sym}\nolimits^{h}\limfunc{Sym}\nolimits_{\flat
}^{2}(std)\otimes \limfunc{Res}\rho }^{ss}(N;\overline{\mathbb{F}}_{p}),
\end{equation*}

\noindent where we are viewing $\limfunc{Sym}\nolimits^{h}\limfunc{Sym}%
\nolimits_{\flat }^{2}(std)\otimes \limfunc{Res}\rho $ as a representation
of $G(U_{r}\times U_{s})$ by restriction ($std:GL_{g}\rightarrow GL_{g}$ is
the standard representation of $GL_{g}$).

\noindent The maps $r_{h}$ and $\xi _{h}$ are Hecke equivariant; as N.
Fakhruddin showed in \cite{NF}, the Kodaira-Spencer map is also
Hecke-equivariant, modulo a rescaling on the Hecke operators acting on $%
\limfunc{Sym}_{\flat }^{2}(\mathbb{E})$. We deduce that, after performing
the mentioned rescaling, the non-zero form $f^{ss}$ is an Hecke-eigenform
with the same eigenvalues than our original $f$.

On the other side, let assume we are given a non-zero eigenform $f^{ss}\in
M_{\rho ^{\prime }}^{ss}(N;\overline{\mathbb{F}}_{p})$ for some weight $\rho
^{\prime }:G(U_{r}\times U_{s})\rightarrow GL_{m^{\prime }}$. There is a
rational $\mathbb{F}_{p}$-representation $\widetilde{\rho }%
:GL_{g}\rightarrow GL_{m}$ whose restriction to $G(U_{r}\times U_{s})$
contains $\rho ^{\prime }$. Indeed the algebraically induced representation $%
\rho ^{\prime \prime }:=\limfunc{Ind}_{G(U_{r}\times U_{s})}^{GL_{g}}\rho
^{\prime }$ contains (non-canonically) a finite dimensional $G(U_{r}\times
U_{s})$-invariant subspace $\tau $ that is $G(U_{r}\times U_{s})$-isomorphic
to $\rho ^{\prime }$; by local finiteness there is a finite dimensional $%
GL_{g}$-submodule $\widetilde{\rho }$ of $\rho ^{\prime \prime }$ containing 
$\tau $ as a $G(U_{r}\times U_{s})$-submodule.

\noindent By Proposition \ref{epi}, there is an integer $c>0$ divisible by $%
p^{2}-1$ such that the map:%
\begin{equation*}
r:M_{\widetilde{\rho }\otimes \det^{c}}(N;\overline{\mathbb{F}}%
_{p})\longrightarrow M_{\limfunc{Res}(\widetilde{\rho }\otimes
\det^{c})}^{ss}(N;\overline{\mathbb{F}}_{p})=M_{\limfunc{Res}\widetilde{\rho 
}}^{ss}(N;\overline{\mathbb{F}}_{p})
\end{equation*}%
is surjective; since $M_{\rho ^{\prime }}^{ss}(N;\overline{\mathbb{F}}%
_{p})\subseteq M_{\limfunc{Res}\widetilde{\rho }}^{ss}(N;\overline{\mathbb{F}%
}_{p})$ and since $r$ is Hecke-equivariant, we see that a system of Hecke
eigenvalues occurring in $M_{\rho ^{\prime }}^{ss}(N;\overline{\mathbb{F}}%
_{p})$ also occurs in $M_{\widetilde{\rho }\otimes \det^{c}}(N;\overline{%
\mathbb{F}}_{p})$.

We conclude that the system of Hecke eigenvalues arising from our spaces of
modular forms $M_{\rho }(N;\overline{\mathbb{F}}_{p})$ for varying $\rho
:GL_{g}\rightarrow GL_{m}$, coincide with the systems of Hecke eigenvalues
arising from the spaces $M_{\rho ^{\prime }}^{ss}(N;\overline{\mathbb{F}}%
_{p})$ for varying $\rho ^{\prime }:G(U_{r}\times U_{s})\rightarrow
GL_{m^{\prime }}$. The theorem now follows from Proposition \ref%
{alg-superspecial}. $\blacksquare $

\bigskip

We presented here the construction of the Hecke correspondence in the PEL
unitary case. One obtains the result of Ghitza (for $g>1$) by forgetting
about the action of the algebra with involution that appears in our
computations; observe that for Siegel modular forms, the superspecial locus
has an easier shape, as explained by Remark \ref{gh is ok}.

Under suitable assumptions that guarantee that the superspecial locus is
non-empty, one obtains a result similar to the above Theorem \ref{main thm}
in the context of Hilbert modular forms, i.e. for PEL-data of type C-I with $%
e\geq 1$ (cf. \ref{Examples}). One condition that needs to be satisfied in
this case is the existence of an embedding of the fixed totally real field
into the subset of symmetric matrices of $M_{g}(\mathfrak{B})$, for $%
\mathfrak{B}$ the quaternion algebra over $%
%TCIMACRO{\U{211a} }%
%BeginExpansion
\mathbb{Q}
%EndExpansion
$ ramified at $p$ and $\infty $ (here a matrix $X\in M_{g}(\mathfrak{B})$ is
said to be symmetric if $X=\overline{X}^{t}$).

\subsection{Computing the number of Hecke eigensystems}

We continue to assume fixed the notation introduced in the previous section,
in particular we work with unitary $(\func{mod}p)$ PEL-modular forms of
signature $(r,s)$, with $p>2$. We want to give an estimate of the number of
distinct $(\func{mod}p)$ Hecke eigensystems occurring in the spaces $M_{\rho
}(N;\overline{\mathbb{F}}_{p})$ for $N$ fixed, and varying $\rho $; we
follow the lines of \cite{Gh4}.

Denote by $\mathcal{N}:=\mathcal{N}(p;k,r,s;N)$ the number of distinct Hecke
eigensystems occurring in the totality of spaces $M_{\rho }(N;\overline{%
\mathbb{F}}_{p})$'s for $\rho $ varying over the set of representations of $%
GL_{g}$ over $\overline{\mathbb{F}}_{p}$; by Th. \ref{main thm} and Prop. %
\ref{alg-superspecial}, $\mathcal{N}$ is the number of distinct Hecke
eigensystems occurring in the totality of spaces $M_{\rho }^{ss}(N;\overline{%
\mathbb{F}}_{p})$ where $\rho $ now runs over the finite set $\limfunc{Irr}%
(G(p))$ of isomorphism classes of irreducible finite-dimensional
representations of $G(p):=G(U_{r}\times U_{s})(\mathbb{F}_{p^{2}})$ over $%
\overline{\mathbb{F}}_{p}$. If $\rho :G(p)\rightarrow GL(W_{\rho })$ is any
fixed element representing a class in $\limfunc{Irr}(G(p))$, we have:%
\begin{eqnarray*}
M_{\rho }^{ss}(N;\overline{\mathbb{F}}_{p}) &=&\{f:Z_{diff}^{\prime }(%
\overline{\mathbb{F}}_{p})\rightarrow W_{\rho }:f([(A,i,\overline{\lambda },%
\overline{\alpha },\eta M)])= \\
\rho (M)^{-1}f([(A,i,\overline{\lambda },\overline{\alpha },\eta )])\text{,
all }M &\in &G(p)\text{, }[(A,i,\overline{\lambda },\overline{\alpha },\eta
)]\in Z_{diff}^{\prime }(\overline{\mathbb{F}}_{p})\},
\end{eqnarray*}

\noindent so that, by definition of $Z^{\prime }(\overline{\mathbb{F}}_{p})$
we have $\dim _{\overline{\mathbb{F}}_{p}}M_{\rho }^{ss}(N;\overline{\mathbb{%
F}}_{p})\leq \#Z^{\prime }(\overline{\mathbb{F}}_{p})\cdot \dim _{\overline{%
\mathbb{F}}_{p}}W_{\rho }$, \noindent and:%
\begin{equation}
\mathcal{N}\leq \#Z^{\prime }(\overline{\mathbb{F}}_{p})\cdot \sum\nolimits_{%
\left[ \rho \right] \in \limfunc{Irr}(G(p))}\dim _{\overline{\mathbb{F}}%
_{p}}W_{\rho }.  \label{to estimate}
\end{equation}

\noindent We now give estimates of the two factors appearing in the right
hand side of the last inequality.

\subsubsection{Estimate of the cardinality of the superspecial locus}

In order to compute $\#Z^{\prime }(\overline{\mathbb{F}}_{p})$, one would
like to have an explicit mass formula for superspecial varieties of the
PEL-type considered here;\ lacking such an explicit formula, we can instead
using what is known for Siegel varieties.\ Let us denote by $\mathcal{A}:=%
\mathcal{A}_{g,1,N}$ \ the Siegel moduli scheme over $\mathcal{O}_{k,(p)}$
classifying prime--to-$p$ isogeny classes of tuples $(A,\overline{\lambda },%
\overline{\overline{\alpha }})$, where $A$ is an abelian projective scheme
over some $S\in SCH_{\mathcal{O}_{k,(p)}}$ of relative dimension $g$, $%
\overline{\lambda }$ is a principal homogeneous polarization of $A$, and $%
\overline{\overline{\alpha }}$ is a full level $N$ structure on $(A,%
\overline{\lambda })$; there is a natural mapping $j$ from the moduli $%
\mathcal{O}_{k,(p)}$-scheme $\mathcal{S}:=\mathcal{S}_{\mathcal{D}}$
associated to the PEL-datum $\mathcal{D}$ of type $A$ that we fixed here, to 
$\mathcal{A}$. More precisely, by fixing an isomorphism of $%
%TCIMACRO{\U{211a} }%
%BeginExpansion
\mathbb{Q}
%EndExpansion
$-vector spaces $V:=k^{g}\simeq 
%TCIMACRO{\U{211a} }%
%BeginExpansion
\mathbb{Q}
%EndExpansion
^{2g}$ we obtain a monomorphism of $%
%TCIMACRO{\U{211a} }%
%BeginExpansion
\mathbb{Q}
%EndExpansion
$-groups $GU_{g}(k;r,s)\hookrightarrow GSp_{2g}(J)\simeq GSp_{2g}$, where $J$
is some symplectic form on $%
%TCIMACRO{\U{211a} }%
%BeginExpansion
\mathbb{Q}
%EndExpansion
^{2g}$; then by definition, if $S$ is some locally noetherian $\mathcal{O}%
_{k,(p)}$-scheme, $j$ sends the equivalence class $[(A,i,\overline{\lambda },%
\overline{\alpha })]\in \mathcal{S}(S)$ to the equivalence class $[(A,%
\overline{\lambda },\overline{\overline{\alpha }})]\in $ $\mathcal{A}(S)$,
where $\overline{\overline{\alpha }}$ is the $U^{\prime }(N)$ orbit of the
symplectic isomorphism $\alpha :H_{1}(A,\mathbb{A}_{f}^{p})\rightarrow
V\otimes _{%
%TCIMACRO{\U{211a} }%
%BeginExpansion
\mathbb{Q}
%EndExpansion
}\mathbb{A}_{f}^{p}$, with $U^{\prime }(N):=Ker(GSp_{2g}(\hat{%
%TCIMACRO{\U{2124}}%
%BeginExpansion
\mathbb{Z}%
%EndExpansion
}^{p};J)\rightarrow GSp_{2g}(\hat{%
%TCIMACRO{\U{2124}}%
%BeginExpansion
\mathbb{Z}%
%EndExpansion
}^{p}/N\hat{%
%TCIMACRO{\U{2124}}%
%BeginExpansion
\mathbb{Z}%
%EndExpansion
}^{p};J))$ (notice that $U(N)=U^{\prime }(N)\cap GU_{g}(\mathcal{O}%
_{k}\otimes \hat{%
%TCIMACRO{\U{2124}}%
%BeginExpansion
\mathbb{Z}%
%EndExpansion
}^{p};r,s)$).

By works of\ Kisin and Vasiu, $j$ is a closed immersion. For our purposes,
we content ourselves with the fact that $j$ induces an injection on closed $%
\overline{\mathbb{F}}_{p}$-points, and sends injectively the set of closed
points $Z^{\prime }(\overline{\mathbb{F}}_{p})$ of the superspecial locus of
the unitary PEL-variety $\mathcal{S}\otimes \overline{\mathbb{F}}_{p}$ -
relative to our choice of $(A_{0},i_{0},\overline{\lambda }_{0},\overline{%
\alpha }_{0})$ - into the superspecial locus of $\mathcal{A}(\overline{%
\mathbb{F}}_{p})$. We can use the estimate given in \cite{Gh4}, based on the
explicit mass formula for superspecial principally polarized abelian $%
\overline{\mathbb{F}}_{p}$-varieties due to Ekedahl (\cite{Ek}) and based on
work of Hashimoto-Ibukiyama (\cite{IK}):%
\begin{equation}
\#Z^{\prime }(\overline{\mathbb{F}}_{p})\leq C_{g}\cdot \#GSp_{2g}(%
%TCIMACRO{\U{2124} }%
%BeginExpansion
\mathbb{Z}
%EndExpansion
/N%
%TCIMACRO{\U{2124} }%
%BeginExpansion
\mathbb{Z}
%EndExpansion
)\cdot \prod\nolimits_{i=1}^{g}(p^{i}+(-1)^{i}),  \label{counting Z'}
\end{equation}

\noindent where the constant $C_{g}$ is defined to be: 
\begin{equation*}
C_{g}\mathcal{=}\frac{(-1)^{g(g+1)/2}}{2^{g}}\cdot
\prod\nolimits_{i=1}^{g}\zeta (1-2i)=\frac{1}{2^{2g}g!}\cdot
\prod\nolimits_{i=1}^{g}B_{2i},
\end{equation*}

\noindent (here $\zeta $ is the Riemann zeta function, and $B_{2i}$ denotes
the $2i$-th Bernoulli number).

\subsubsection{Estimates on the size of irreducible representations of $G(p)$%
}

All the representations we consider in this section are finite dimensional
over the appropriate field. The number of pairwise non-isomorphic
irreducible representations of the finite group $G(p)$ over $\overline{%
\mathbb{F}}_{p}$ coincides with the number $k^{p}(G(p))$ of $p$-regular
conjugacy classes of $G(p)$; a matrix element $X$ of $G(p)$ is $p$-regular
if and only if its minimal polynomial has only simple roots over $\overline{%
\mathbb{F}}_{p}$, that is to say if and only if $X$ is semi-simple (over $%
\overline{\mathbb{F}}_{p}$).

The group $G(p)$ is the set of $\mathbb{F}_{p}$-points of the connected
reductive algebraic group $\mathbb{G}:=G(U_{r}\times U_{s})$ defined over $%
\mathbb{F}_{p}$; one can compute the center $Z$ and the derived subgroup $%
\mathbb{G}^{\prime }$ of $\mathbb{G}$ and find:%
\begin{eqnarray*}
Z &=&Z^{0}\simeq \left\{ 
\begin{array}{c}
G(U_{1}\times U_{1}) \\ 
GU_{1}%
\end{array}%
\begin{array}{c}
\text{if }rs\neq 0, \\ 
\text{if }rs=0,%
\end{array}%
\right. \\
\mathbb{G}^{\prime } &=&SU_{r}\times SU_{s}.
\end{eqnarray*}

\noindent Since $\mathbb{G}^{\prime }$ is connected, simply-connected and
semi-simple with rank:%
\begin{equation*}
rk(\mathbb{G}^{\prime })=\left\{ 
\begin{array}{c}
g-2 \\ 
g-1%
\end{array}%
\begin{array}{c}
\text{if }rs\neq 0, \\ 
\text{if }rs=0,%
\end{array}%
\right.
\end{equation*}%
by applying Theorem 3.7.6 of \cite{Car}, we have:%
\begin{equation}
k^{p}(G(p))=\#Z^{0}(\mathbb{F}_{p})\cdot p^{rk(\mathbb{G}^{\prime
})}=\left\{ 
\begin{array}{c}
p^{g-2}\cdot (p-1)(p+1)^{2} \\ 
p^{g-1}\cdot (p-1)(p+1)%
\end{array}%
\begin{array}{c}
\text{if }rs\neq 0, \\ 
\text{if }rs=0.%
\end{array}%
\right.  \label{conj classes}
\end{equation}

If $t\geq 2$, then $SU_{t}(\mathbb{F}_{p^{2}})$ is the set of $\mathbb{F}%
_{p} $-points of a group of a simply connected group of type $^{2}A_{t-1}(p)$%
, and its order is (cf. \cite{Car} 2.9) $\#SU_{t}(\mathbb{F}_{p^{2}})=p^{%
\frac{t(t-1)}{2}}\cdot \prod\nolimits_{i=2}^{t}(p^{i}-(-1)^{i}).$ (We set $%
SU_{0}(\mathbb{F}_{p^{2}})=SU_{1}(\mathbb{F}_{p^{2}}):=\{1\}$). \noindent
Using the exactness of the sequence $1\rightarrow SU_{t}(\mathbb{F}%
_{p^{2}})\rightarrow U_{t}(\mathbb{F}_{p^{2}})\overset{\det }{\rightarrow }%
U_{1}(\mathbb{F}_{p^{2}})\rightarrow 1$ \noindent for $t>0$, one deduces
that $\#U_{t}(\mathbb{F}_{p^{2}})=\#SU_{t}(\mathbb{F}_{p^{2}})\cdot (p+1)$.
We conclude that for any choice of non-negative integers $r$ and $s$ such
that $r+s=g$ we have:%
\begin{eqnarray*}
\#G(U_{r}\times U_{s})(\mathbb{F}_{p^{2}}) &=&\#U_{r}(\mathbb{F}%
_{p^{2}})\cdot \#U_{s}(\mathbb{F}_{p^{2}})\cdot (p-1)= \\
&=&p^{\frac{r(r-1)+s(s-1)}{2}}\cdot
\tprod\nolimits_{i=1}^{r}(p^{i}-(-1)^{i})\cdot
\tprod\nolimits_{i=1}^{s}(p^{i}-(-1)^{i})\cdot (p-1).
\end{eqnarray*}

\noindent In particular, a $p$-Sylow subgroup of $G(U_{r}\times U_{s})(%
\mathbb{F}_{p^{2}})$ has order $p^{\frac{r(r-1)+s(s-1)}{2}}$. Since $G(p)$
is a group with a split $(B,N)$-pair (cf. \cite{Car} 1.18), we deduce that
if $\rho :G(p)\rightarrow GL(W_{\rho })$ is an irreducible representation of 
$G(p)$ over $\overline{\mathbb{F}}_{p}$, then:%
\begin{equation}
\dim _{\overline{\mathbb{F}}_{p}}W_{\rho }\leq p^{\frac{r(r-1)+s(s-1)}{2}}.
\label{dimension}
\end{equation}

\noindent (The proof of this fact is contained in \cite{Cur}; cf. esp.
corollaries 3.5 and 5.11). \noindent Putting together formulae (\ref{conj
classes}) and (\ref{dimension}) we obtain:%
\begin{eqnarray}
\sum\nolimits_{\left[ \rho \right] \in \limfunc{Irr}(G(p))}\dim _{\overline{%
\mathbb{F}}_{p}}W_{\rho } &\leq &  \label{counting irr} \\
&\leq &\left\{ 
\begin{array}{c}
p^{\frac{r(r-1)+s(s-1)}{2}}\cdot p^{g-2}(p-1)(p+1)^{2} \\ 
p^{\frac{g(g-1)}{2}}\cdot p^{g-1}(p-1)(p+1)%
\end{array}%
\begin{array}{c}
\text{if }rs\neq 0, \\ 
\text{if }rs=0.%
\end{array}%
\right. \cdot  \notag
\end{eqnarray}

\subsubsection{Upper bound for the number of Hecke eigensystems}

\noindent Putting formulae (\ref{counting Z'}) and (\ref{counting irr})
together into formula (\ref{to estimate}), we obtain:

\begin{theorem}
\label{counting theorem}Let $p>2,k,r,s,N$ be fixed as above (in particular $%
r,s\geq 0$ and $r+s=g$) and set $C_{g}:=2^{-2g}(g!)^{-1}\prod%
\nolimits_{i=1}^{g}B_{2i}.$ The number $\mathcal{N}:=\mathcal{N}(p;k,r,s;N)$
of distinct $(\func{mod}p)$ Hecke eigensystems occurring in the spaces $%
M_{\rho }(N;\overline{\mathbb{F}}_{p})$ for varying $\rho $ satisfies the
following inequality: 
\begin{eqnarray*}
\mathcal{N} &\leq &C_{g}\cdot \#GSp_{2g}(%
%TCIMACRO{\U{2124} }%
%BeginExpansion
\mathbb{Z}
%EndExpansion
/N%
%TCIMACRO{\U{2124} }%
%BeginExpansion
\mathbb{Z}
%EndExpansion
)\cdot \\
&&\cdot \prod\nolimits_{i=1}^{g}(p^{i}+(-1)^{i})\cdot \left\{ 
\begin{array}{c}
p^{\frac{r(r-1)+s(s-1)}{2}}\cdot p^{g-2}(p-1)(p+1)^{2} \\ 
p^{\frac{g(g-1)}{2}}\cdot p^{g-1}(p-1)(p+1)%
\end{array}%
\begin{array}{c}
\text{if }rs\neq 0, \\ 
\text{if }rs=0.%
\end{array}%
\right.
\end{eqnarray*}%
In particular, if we keep $k,r,s,N$ fixed and let $p>2$ vary:%
\begin{equation*}
\mathcal{N}\in O(p^{g^{2}+g+1-rs}),\text{ \ for }p\rightarrow \infty .
\end{equation*}
\end{theorem}

\bigskip

For an estimate of $\mathcal{N}$ in the\ case of Siegel modular forms, cf. 
\cite{Gh4}; for elliptic modular forms, one can find a conjectural mass
formula for the asymptotic with respect to $p$ of two-dimensional odd and
irreducible Galois representations of $%
%TCIMACRO{\U{211a} }%
%BeginExpansion
\mathbb{Q}
%EndExpansion
$ in \cite{Cen}.

\bigskip

\newpage

\bibliographystyle{amsalpha}
\bibliography{acompat,davide}

\end{document}